\documentclass[11pt]{article} % for LaTex2e

\usepackage{amsfonts}
\usepackage{latexsym}
\font\tinybbfont=msbm6
\font\scriptsizebbfont=msbm7 scaled \magstep 1
 1
\font\footnotesizebbfont=msbm9 scaled \magstep 0
\font\smallbbfont=msbm7 scaled \magstep 2
\font\bbfont=msbm9 scaled \magstep1  % 1.2 pt
\font\largebbfont=msbm10 scaled \magstep 1
 2
 3
 4

\def\tinyBbb#1{\hbox{\tinybbfont #1}}
\def\scriptsizeBbb#1{\hbox{\scriptsizebbfont #1}}

\def\footnotesizeBbb#1{\hbox{\footnotesizebbfont #1}}
\def\smallBbb#1{\hbox{\smallbbfont #1}}
\def\Bbb#1{\hbox{\bbfont #1}}
\def\largeBbb#1{\hbox{\largebbfont #1}}

\newcommand{\Aut}{\mbox{\it Aut}\,}
\newcommand{\Average}{\mbox{\it Average}\,}

\newcommand{\Def}{\mbox{\it Def}\,}
 \newcommand{\scriptsizeDef}{\mbox{\scriptsize\it Def}\,}
 \newcommand{\tinyDef}{\mbox{\tiny\it Def}\,}

\newcommand{\Div}{\mbox{\rm Div}\,}
\newcommand{\End}{\mbox{\it End}\,}

\newcommand{\Ext}{\mbox{\it Ext}\,}
\newcommand{\Extsheaf}{\mbox{\it ${\cal E}$xt}\,}

\newcommand{\Glue}{\mbox{\it Glue}\,}

\newcommand{\Hom}{\mbox{\it Hom}\,}

\newcommand{\Id}{\mbox{\it Id}\,}
\newcommand{\Image}{\mbox{\it Im}\,}
 \newcommand{\footnotesizeImage}{\mbox{\it\footnotesize Im}\,}

\newcommand{\IP}{\mbox{\it IP}\,}
\newcommand{\Isom}{\mbox{\it Isom}\,}
 \newcommand{\footnotesizeIsom}{\mbox{\it\footnotesize Isom}\,}

 \newcommand{\footnotesizeIsombf}{\mbox{\bf\footnotesize Isom}\,}

\newcommand{\Ker}{\mbox{\it Ker}\,}
\newcommand{\Kuranishi}{\,\mbox{\rm\scriptsize Kuranishi}}

\newcommand{\Neck}{\mbox{\it Neck}\,}
\newcommand{\scriptsizeNeck}{\mbox{\it\scriptsize Neck}\,}

\newcommand{\Per}{\mbox{\it Per}\,}
 \newcommand{\smallPer}{\mbox{\it\small Per}\,}

\newcommand{\Real}{\mbox{\it Re}\,}

\newcommand{\Spec}{\mbox{\it Spec}\,}
 \newcommand{\footnotesizeSpec}{\mbox{\it\footnotesize Spec}\,}

\newcommand{\Sym}{\mbox{\it Sym}}

\newcommand{\Trunk}{\mbox{\it Trunk}\,}

\newcommand{\approxi}{\mbox{\rm\scriptsize approx}}
 \newcommand{\tinyapproxi}{\mbox{\rm\tiny approx}}

\newcommand{\aux}{\,\mbox{\rm\scriptsize aux}}
 \newcommand{\tinyaux}{\,\mbox{\rm\tiny aux}}
\newcommand{\bboxtimes}{\Box\hspace{-1.75ex}\raisebox{.15ex}{$\times$}\,}
\newcommand{\blue}{\mbox{\rm blue}}
\newcommand{\bn}{\mbox{\rm\scriptsize b.n.}}

\newcommand{\codimm}{\mbox{\it codim}\,}
\newcommand{\deform}{\mbox{\rm\scriptsize deform}}
\newcommand{\degree}{\mbox{\it deg}\,}
 \newcommand{\footnotesizedeg}{\mbox{\it\footnotesize deg}\,}
 
\newcommand{\dimm}{\mbox{\it dim}\,}

\newcommand{\divisor}{\mbox{\it div}}
\newcommand{\domain}{\mbox{\rm\scriptsize domain}}
\newcommand{\ev}{\mbox{\it ev}\,}
\newcommand{\fiber}{\,\mbox{\rm\scriptsize fiber}\,}

\newcommand{\hand}{\mbox{\rm\scriptsize hand}\,}
\newcommand{\im}{\mbox{\rm Im}\,}

\newcommand{\ind}{\mbox{\it ind}\,}

\newcommand{\intnode}{\mbox{\rm\scriptsize i.n.}}
\newcommand{\jet}{\mbox{\it jet}\,}
\newcommand{\leg}{\mbox{\rm\scriptsize leg}\,}

\newcommand{\loc}{\mbox{\rm\scriptsize loc}}

\newcommand{\layer}{\mbox{\it layer}\,}
\newcommand{\map}{\,\mbox{\scriptsize\rm map}\,}

\newcommand{\node}{\,\mbox{\rm\scriptsize node}}

\newcommand{\nonrigid}{\,\mbox{\rm\scriptsize non-rigid}}
 \newcommand{\tinynonrigid}{\,\mbox{\rm\tiny non-rigid}}
\newcommand{\oin}{\mbox{\rm\scriptsize o.i.n}}

\newcommand{\order}{\mbox{\it ord}\,}
 \newcommand{\footnotesizeord}{\mbox{\it\footnotesize ord}\,}
\newcommand{\pd}{\mbox{\rm\scriptsize pd}}
\newcommand{\pr}{\mbox{\it pr}}
\newcommand{\pt}{\mbox{\it pt}}
\newcommand{\rank}{\mbox{\it rank}\,}
\newcommand{\rb}{\mbox{\it rb}\,}

\newcommand{\red}{\mbox{\rm red}}
\newcommand{\reg}{\mbox{\rm\scriptsize reg}}
\newcommand{\rel}{\mbox{\it\scriptsize rel}}
\newcommand{\relpd}{\mbox{\scriptsize {\it rel}-{\rm pd}}}

\newcommand{\rigid}{\,\mbox{\rm\scriptsize rigid}}
\newcommand{\rigidified}{\mbox{\rm\scriptsize rigidified}}
\newcommand{\rigidifying}{\mbox{\rm\scriptsize rigidifying}}

\newcommand{\sing}{\mbox{\scriptsize\it sing}}

\newcommand{\smooth}{\mbox{\scriptsize\rm smooth}}
\newcommand{\shiftproduct}{\mbox{\it sp}}
\newcommand{\spsccw}{\mbox{\scriptsize\rm spsccw}}
\newcommand{\st}{\mbox{\it st}}
 \newcommand{\scriptsizest}{\mbox{\rm\scriptsize st}}

\newcommand{\target}{\mbox{\rm\scriptsize target}}

\newcommand{\trivial}{\mbox{\rm\scriptsize trivial}\,}
\newcommand{\vdim}{\mbox{\it vdim}\,}
\newcommand{\virt}{\mbox{\scriptsize\it virt}}

\newcommand{\disjointunion}{\raisebox{.3ex}{\scriptsize $\,\coprod\,$}}
\newcommand{\tinybullet}{\mbox{\tiny $\bullet$}}

\begin{document}

\enlargethispage{23cm}

\begin{titlepage}

$ $

\vspace{-1cm} % Re: -1.5cm for PC; -2.5cm for UT-Math-system

\noindent\hspace{-1cm}
\parbox{6cm}{\small September 2006}\
   \hspace{6.5cm}\
   \parbox[t]{5cm}{math.SG/0609483 \newline
                   OGW: $1/4$}

\vspace{1cm}
%\vspace{2cm}

%title
\centerline{\large\bf
 Degeneration and gluing of Kuranishi structures in Gromov-Witten theory}
%\vspace{1ex}
%{\large \bf -----\hspace{1ex} \parbox[t]{28em}{a
%  degeneration axiom and a gluing axiom under a symplectic cut\\[.6ex]
%  for open Gromov-Witten invariants of a symplectic manifold\\[.6ex]
%  with a decorated Lagrangian submanifold}}
\vspace{1ex}
\centerline{\large\bf
 and the degeneration/gluing axioms for open Gromov-Witten invariants}
\vspace{1ex}
\centerline{\large\bf
 under a symplectic cut}
\vspace{1ex}
% \centerline{\large\bf
%   for open Gromov-Witten invariants}
% \vspace{1ex}
% \centerline{\large\bf
%  a degeneration axiom and a gluing axiom under a symplectic cut}
% \vspace{1ex}
% \centerline{\large\bf
%   for open Gromov-Witten invariants}
% \vspace{1ex}
% \centerline{\large\bf
%  of a symplectic manifold with a decorated Lagrangian submanifold}

%\vspace{1cm}
\vspace{3em}
%authors-'n-addresses
\centerline{\large
  Chien-Hao Liu
  \hspace{1ex} and \hspace{1ex}
  Shing-Tung Yau
}

\vspace{2em}
\centerline{\small
 ($\,${\it In memory of Professor Raoul Bott.}$\,$)}
\vspace{2em}

%abstract%
\begin{quotation}
\centerline{\bf Abstract}
\vspace{0.3cm}
\baselineskip 12pt  %13pt for [12pt] style
{\small
  We construct a family Kuranishi structure in the Fukaya-Ono format
   on the moduli space
   $\overline{\cal M}_{\,\bullet\,}(W/B,L\,|\,\bullet\,)$
   of open stable $J$-holomorphic maps to the fibers of
   an almost-complex degeneration family $W/B$ that arises from
   a symplectic cut.
  The degenerate fiber of the family Kuranishi structure defines
   a Kuranishi structure on the moduli space of open stable maps
   to a singular symplectic space of the gluing form
   $Y_1\cup_D Y_2$ from a symplectic cut,
    with a Lagrangian submanifold $L$ contained in the smooth locus.
  The same discussion and construction apply also to
   relative open Gromov-Witten theory for a relative pair $(Z,L;D)$,
   where $D$ is a codimension-$2$ symplectic submanifold of $Z$,
    disjoint from the Lagrangian submanifold $L$.
  We derive then the degeneration-gluing relations of these
   Kuranishi structures.
  The good flat behavior of the family Kuranishi structure
   on $\overline{\cal M}_{\,\bullet\,}(W/B,L\,|\,\bullet\,)$
   motivates both a degeneration axiom and a gluing axiom
   for open Gromov-Witten invariants of a symplectic manifold $X$
   with a decorated Lagrangian submanifold $L^{\alpha}$.
  When a symplectic cut at the boundary of a tubular neighborhood
   of $L$ exists, the construction of open Gromov-Witten invariants
   of $(X,L^{\alpha})$ can then be put in two steps:
   (1) use the degeneration axiom and the gluing axiom to fix
       the ambiguity in the choice of fundamental chain class;
   (2) intersection theory on the specific kind of singular Kuranishi
       space with the induced decoration on the moduli space of
       relative maps to the relative pairs from the degenerate target.
  Step (1) is analytical and is dealt with in this work.
  In the appendix we comment on the equivalence of
   Li-Ruan/Li's degeneration formula and
   Ionel-Parker's degeneration formula in closed Gromov-Witten theory.
} %endsmall
\end{quotation}

%\smallskip
\vspace{2em}

\baselineskip 12pt
{\footnotesize
} %endfootnotesize

\medskip

\noindent {\small
MSC number 2000$\,$:
 53D45; 14N35, 81T30.
} % end-small

\medskip

\baselineskip 11pt
{\footnotesize
\noindent{\bf Acknowledgements.}
 We thank
  Kenji Fukaya, Kefeng Liu, Cumrun Vafa
   for discussions.
 C.-H.L.\ thanks also
  K.L., Chiu-Chu Liu, Xiaowei Wang, Rugang Ye
   for discussions on symplectic Gromov-Witten theory,
   spring 2003;
  Yongbin Ruan
   for discussions on degeneration formula,
   April and November 2005;
  C.V.
   for the semester-long topic course on string theory, spring 2006,
   the numerous thought-provoking explanations and
   answer-to-questions along the lectures,
   and the literature guide;
  these three events help C.-H.L.\ in understanding the various
   themes in the project;
  Samit Dasgupta, Eaman Eftekhary, Dennis Gaitsgory, Joe Harris,
  Eleny-Nicoleta Ionel, Albrecht Klemm, Nikolai Krylov, Yum-Tong Siu,
  Chin-Lung Wang, Ilia Zharkov
   for lectures/discussions;
  A.K., Y.R., Department of Mathematics/Physics of
   U.\ Wisconsin - Madison for hospitality;
 Ling-Miao~Chou for moral support.
 The project is supported by NSF grants DMS-9803347 and DMS-0074329.
} %endfootnotesize

\end{titlepage}

\newpage
\begin{titlepage}

$ $

\vspace{2em}

{\footnotesize

\begin{flushleft}
{\bf
 Professor Raoul Bott and mirror symmetry
 - a reminiscence of a curious mind}, by C.-H.L..
\end{flushleft}
 In the fall 2000, after  the semester-long lectures of
  Prof.~Cumrun Vafa on string theory and stringy duality for
  both mathematicians and physicists in the spring that year,
 Prof.~Raoul~Bott got intrigued in mirror symmetry and
  had a conversation with Prof.~Yau, the second author.
 It ended in a surprise e-mail from Prof.~Yau to me one day
  with an assignment:
  to teach Prof.~Bott what mirror symmetry is about.
 To ``teach" a then-78-year-old legendary mathematician?!
 Unable to turn it away, I thus took the task as a new hand,
  expecting that my ``student" would get bored very soon and
  I could resume the full focus on real projects.
 Amazingly, the outlining lecture for Prof.~Bott was extended
  to weekly meetings for two intensive months in that semester.

 Now, stringy duality, including mirror symmetry,
  is a very broad and technical subject.
 Its true/best explanation
  as yet remains largely physical, rather than mathematical.
 Its foundation lies on quantum field theory (QFT) and
  the rigidity of supersymmetric QFT's,
 together with numerous other mathematical and physical notions,
  objects, structures, and moduli problems that are incorporated into
  superstring theory along its continual fast-paced developments.
 With such an origin from physics, statements from stringy duality are
  unavoidably mysterious, shocking, and awe-inspiring to mathematicians.
 Anyone who ventures to lecture on such a subject before a mathematics
  master like Prof.~Bott should expect to face many
  legitimate-yet-hard-to-give-a-round-off-answer questions,
  making him/her ``hanging on the blackboard" forever.
 Although I focused only on the much limited topic on toric mirror
  symmetry, such embarrassing moments still happened
  no matter how complete-in-a-small-range I thought I had prepared.
 Yet, this is indeed how Prof.~Bott in turn started to ``teach" me
  what mirror symmetry is about!
 He rejected assumptions without sound reasons.
 He liked to see things derived from low/dirty scratches
  rather than from some high end.
 He constantly asked, ``{\it WHY?}".
 With that energetic mind, he even attempted to provide his own pictures
  or explanations after listening to what I had presented.
 While each great mind has his/her own way of functioning,
  which can only inspire and is almost always unlearnable,
 it remains quite an experience to see how a great mind functions
  as he digests raw materials, thinks, polishes, and comments on them.
 His questions become a guide toward a deeper understanding.

 Prof.~Bott impressed me that he is not inclined to read a lot of
  literatures.
 This is very different from those from Yau's school.
 Once I brought for the lecture a pile of related papers
  marked with red under-lines and margin notes,
 he stared at them and asked me:
  ``{\it How much time have you left for thinking?}"
 Actually, one reason string theory is demanding is that
  no matter how many notions/techniques one has finally
   brought to his/her mastery and employs them for fruitful results,
  there are always things that remain to be learned/understood when
   one attempts to reach a fuller/more-comprehensive picture.
 As so many intelligent people are devoted diligently and intensively
   to it,
  the growth and diversity and broadness of stringy literatures,
   including both mathematics and physics, can be terrifying.
 That particular question of his reminds me of the necessary balance
  between reading and independent thinking
  - a lesson I should keep in mind for good.
 He once said in a lecture at U.C.~Berkeley:
 ``{\it Doing mathematics should be like paddling a canoe downstream
      - natural and effortless}."
 Most of us who study his works will never be able to reach
  such a Zen-like level of doing mathematics;
 yet perhaps this is part of what he meant to teach us
  through the insight, beauty, and elegance of his works.
 Among his far-reaching influences in mathematics,
 the orbifold/stack version of his joint work with
  Prof.~Michael Atiyah that gives the
  {\it Atiyah-Bott Localization Formula}
  has been used again and again in the exact computations of
  Gromov-Witten invariants, a topic within mirror symmetry as well.
 The formula can be interpreted as a special
  mathematical version of Feynman's path-integral.
 Its format of localization can be generalized to other equivariant
  (co)homology theories, including equivariant K-theory,
  that can be used for gauge instanton counting for
  $d=4$, $N=2$ super Yang-Mills theory.
 Such theory (i.e.\ Seiberg-Witten theory) can be linked, too,
  to Gromov-Witten theory and mirror symmetry picture,
  e.g.\ with the mirror geometry encoded in the complex geometry
  of a family of Seiberg-Witten curves embeddable in a family
  of Calabi-Yau manifolds!

 The news of Prof.~Bott's passing away came in December 2005
  while this work was being written with full vigor.
 These unforgettable hours with him on mirror symmetry
  are like a gift from him in his later years ---
  completely unexpected, yet marking my mind deep.
 We thus dedicate this work to the memory of Prof.~Bott,
  an inspiring and forever curious/learning mind.

} % end-footnotesize

\end{titlepage}

%paper
\newpage
$ $

\vspace{-4em}  % Re: -4cm for PC; -6cm for UT-Math-system

%short heading
\centerline{\sc
 Degeneration of Kuranishi Structure and Axioms for Open GW-Invariants}

\vspace{2em}

\baselineskip 13pt  %Re: 14pt for [11pt] style
                    %Re: 15pt for [12pt] style.

\begin{flushleft}
{\Large\bf 0. Introduction and outline.}
\end{flushleft}
The moduli space of prestable labelled-bordered Riemann surfaces
 is an Artin stack locally modelled on a quotient of
 manifolds-with-corners.
This leads to the singular real codimension-$1$ boundary
 in the Kuranishi structure ${\cal K}$ for the moduii space
 $\overline{\cal M}_{\,\bullet\,}(X,L\,|\,\bullet\,)$
 of open stable maps to a symplectic manifold $X$ with boundary
 confined in a Lagrangian submanifold $L$.
Such boundary gives rise to an ambiguity in choosing the virtual
 fundamental chain on the Kuranishi structure ${\cal K}$ for defining
 open Gromov-Witten invariants of $(X,L)$.
To fix the ambiguity,
an extra data (i.e.\ a ``decoration") $\alpha$ on $L$
 has to be added to the problem and
the induced effect of the decoration on $L$ to the whole
 $\overline{\cal M}_{\,\bullet\,}(X,L\,|\,\bullet\,)$
 and ${\cal K}$ has to be understood.
Examples of such decoration $\alpha$ are
 a group action on $L$,
 a bundle map on the restriction $T_{\ast}X|_L$, or
 a diffeomorphism on a neighborhood of $L$ in $X$ that leaves $L$
  invariant.
However, unless this decoration is extendable to the whole $X$,
 there is no obvious way to go from ``$\alpha$ on $L$" to
 ``an associated extra structure on
 $\overline{\cal M}_{\,\bullet\,}(X,L\,|\,\bullet\,)$ and
 ${\cal K}$" to help fix the choice of the virtual fundamental chain
 $[\overline{\cal M}_{\,\bullet\,}(X,L^{\alpha})\,|\,\bullet\,]^{\virt}$
 on ${\cal K}$.
The main goal of this work is to propose and explain
 a degeneration axiom and a gluing axiom under a symplectic cut
 for open Gromov-Witten invariants of a symplectic manifold with
 a decorated Lagrangian submanifold $(X,L^{\alpha})$ to take care of
 the above technical issue for an important class of $(X,L^{\alpha})$
 that occurs in the compact version of conifold
 transitions of Calabi-Yau $3$-folds in
 open/closed string duality in string theory ([Va1]).

Technically, we construct a Kuranishi structure for moduli spaces in
 \begin{itemize}
  \item[(1)]
   a {\it family open Gromov-Witten theory}
   for a symplectic/almost-complex degeneration associated
   to a symplectic cut, and

  \item[(2)]
   a {\it relative open Gromov-Witten theory} for
   a symplectic/almost-complex manifold $X$
    with a Lagrangian/totally-real submanifold $L$
   relative to a codimension-$2$ symplectic/almost-complex submanifold
   $D$ that is disjoint from $L$.
 \end{itemize}
Notions, constructions, and techniques developed by our predecessors
 in the various formats/\\ settings/categories are uniformized/merged
 into the present study.
Such structure
 extends [F-O] and [Liu(C)] to
  a degeneration-family open Gromov-Witten theory and
  a relative open Gromov-Witten theory.
In the case that $L$ is empty,
 the study re-writes both the symplecto-analytic [L-R], [I-P1], [I-P2]
  and the algebro-geometric [Li1], [Li2]
  in the symplecto-analytic Fukaya-Ono format.
For the technical step of constructing the transition data
 in the Kuranishi structure, we bring in also [Sie1].
How these Kuranishi structures are relevant to the construction of
 open Gromov-Witten invariants can be summarized by:

{\small
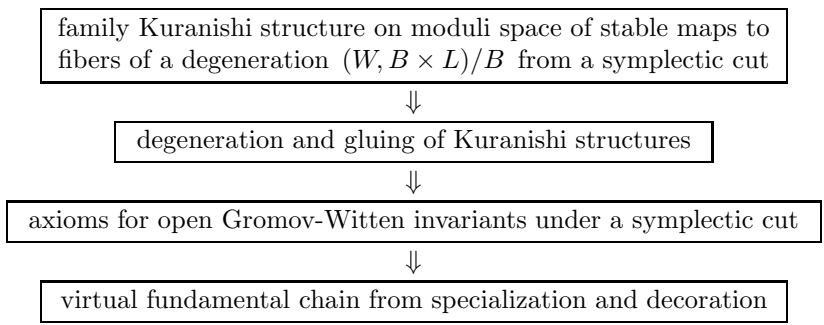
 $$
  \begin{array}{c}
   % \begin{array}{crcccrl}
   %  \hspace{10ex}
   %  && \mbox{[F-O]} & \mbox{[Sie1]} & \hspace{3em}
   %   & \mbox{[I-P1, I-P2], [L-R]}   & \mbox{(symplecto-analytic)}
   %   \\[1.2ex]
   %  & \swarrow   &&&
   %   & \mbox{[Li1,Li2]}  & \mbox{(algebro-geometric)}  \\[1.2ex]
   %  & \mbox{[Liu(C)]}  & & & & \\[1.2ex]
   %  & \searrow & \downarrow & \downarrow
   %   && \downarrow\hspace{3em}\downarrow\hspace{1.6em} \\[1.2ex]
   % \end{array} \\
   \fbox{\parbox{64ex}{\centering
     family Kuranishi structure on moduli space of stable maps to fibers of
     a degeneration $\,(W,B\times L)/B\,$ from a symplectic cut}}\\[2.2ex]
   \Downarrow\\[.2ex]
   \fbox{\parbox{51ex}{\centering
     degeneration and gluing of Kuranishi structures}} \\[.8ex]
   \Downarrow\\[.2ex]
   \fbox{\parbox{70ex}{\centering
    axioms for open Gromov-Witten invariants
                                under a symplectic cut}}\\[.8ex]
   \Downarrow\\[.2ex]
   \fbox{\parbox{64ex}{\centering
    virtual fundamental chain from specialization and decoration}}
  \end{array}
 $$
} % end-small

\bigskip

Compared with closed Gromov-Witten theory,
it may look at first surprising that
{\it in order to understand
 absolute open Gromov-Witten theory one has to understand
  both degeneration and relative open Gromov-Witten theory as well}.
It could be true that this is not the only way.
However, from the viewpoint of algebraic geometry, the route we take
 in this project, of which the current work is a part,
 is an elaborate adoption of the
 {\it deformation-specialization technique} already long in use
 in enumerative (algebraic) geometry;
 (see, e.g.\ [Fu: Sec.~10.4] for an introduction).
Furthermore, the conjectural {\it open/closed string duality} on
 Calabi-Yau three-folds that differ by an extremal transition
 ([Go-V], [O-V1], [O-V2], [Va1])
 almost selects/specifies for us this route uniquely
 among other possible candidate constructions.
Particularly for the motivation and the constant strong drive behind,
 we owe the credits of this work to
  enumerative algebraic geometers and string theorists --
   especially, our teachers Joe Harris and Cumrun Vafa
   and their respective school.
The current work is a step toward a mathematical understanding of
 the compact version of the open/closed string duality for
 Calabi-Yau $3$-folds in [Go-V], [O-V2], and [Va1]
 at the level of moduli spaces/stacks of stable maps,
 cf.\ the diagram in [L-Y2: Introduction].
(See also [D-F] for related discussions.)

\bigskip

\noindent
{\it Convention.}
 Standard notations, terminology, operations, facts in
  (1) symplectic geometry;
  (2) algebraic geometry; (3) Sobolev theory; (4) topology
  can be found respectively in (1) [MD-S2], [G-S], [Woo];
  (2) [Hart], [G-H];
  (3) [MD-S3: Appendix.~B], [Au: Chap.~2 - Chap.~3];
  (4) [Sp].
 \begin{itemize}
  \item[$\cdot$]
   All {\it dimension}, {\it codimension}, {\it rank}, {\it index},
    ..., etc.\ are with respect to ${\Bbb R}$ unless otherwise noted.

  \item[$\cdot$]
   $|\bullet|$ stands
    for the {\it cardinality} of $\,\bullet\,$
      when $\,\bullet\,$ is a finite set or a finite group,
    for the {\it absolute value} or {\it norm} of $\,\bullet\,$
     when $\,\bullet\,$ is a real or complex number or a vector,
    for the {\it sum} of the entries
     when $\,\bullet\,$ is a vector of integers referring to
      some combinatorial quantity
      (like number of marked points or contact order).

  \item[$\cdot$]
   The {\it complex projective space} of complex dimension $n$
    is denoted by ${\Bbb P}^n$.

  \item[$\cdot$]
   In denoting a stable map $f:\Sigma\rightarrow (X,L)$
    to a symplectic space $X$ with a Lagrangian submanifold $L$,
   it is assumed that $f(\partial\Sigma)\subset L$.
   Similarly, for a relative map $f:\Sigma\rightarrow (Z,L;D)$.
   When $L$ is empty, so is $\partial\Sigma$.

  \item[$\cdot$]
   Properties of a map from a nodal (bordered) curve $\Sigma$ to
    another is imposed on its normalization $\tilde{\Sigma}$;
   e.g.\ a $C^{\infty}$ map $f$ from $\Sigma$ to $Y:=Y_1\cup_D Y_2$
    is a continuous map $f:\Sigma\rightarrow Y$  such that
    its lift to $\tilde{\Sigma_0}$  is
    a $C^{\infty}$ map to either $Y_1$ or $Y_2$,
    where $\Sigma_0$ runs over all irreducible components of $\Sigma$.

  \item[$\cdot$]
   Almost-complex (resp.\ complex) structures on different target
    (resp.\ domains) spaces are usually denoted by the same $J$
    (resp.\ $j$) unless the distinction is crucial to the discussion.

  \item[$\cdot$]
   The term ``{\it orbifolds}" and ``{\it sub-orbifolds}"
    are not restricted only to smooth ones.

  \item[$\cdot$]
   Omitted superscripts or subscripts are often denoted by
    $\,^{\cdot}\,$, $\,^{\mbox{\tiny $\bullet$}}\,$ or
    $\,_{\cdot}\,$, $\,_{\mbox{\tiny $\bullet$}}\,$.
 \end{itemize}
 \begin{itemize}
  \item[$\cdot$]
  Commonly used notations for different objects that have no chances
   of confusion:
  \begin{itemize}
   \item
    ${\Bbb C}$ and ${\Bbb R}\,$: as the {\it complex plane} and
      the {\it real line} in differential geometry   vs.\
    as {\it ground fields} in algebraic geometry;

   \item
    {\it curve class} $\beta$ vs.\
    {\it isomorphism} $(\alpha,\beta)$;

   \item
    {\it isomorphism} $\alpha$ of curves or graphs vs.\
    {\it decoration}
          $\alpha$ on a Lagrangian/almost-complex submanifold;

   \item
    {\it universal curve} ${\cal C}$
      over different bases vs.\
    {\it category} ${\cal C}$;

   \item
    {\it genus} $g$ vs.\ {\it map} $g$;

   \item
    {\it index set} $I$ vs.\
    {\it gluing map} $I_{\tinybullet}$ that {\it i}dentifies subsets.
  \end{itemize}
 \end{itemize}

%\bigskip
\newpage
\begin{flushleft}
{\bf Outline.}
\end{flushleft}
{\small
\baselineskip 11pt  %13pt
\begin{itemize}
 \item[1.]
  Symplectic cut and
  the direct system of expanded degenerations
   in the almost-\\ complex category.
  \vspace{-.6ex}
  \begin{itemize}
   \item[1.1]
    Symplectic cut and the associated expanded degenerations.

   \item[1.2]
    Symplectic/almost-complex relative pairs and their expansions.
  \end{itemize}

 \item[2.]
  Prestable labelled-bordered Riemann surfaces.

 \item[3.]
  The moduli space
   $\overline{\cal M}_{(g,h), (n, \vec{m})}
                    (W/B,L\,|\,[\beta],\vec{\gamma},\mu)$
   of stable maps.
  \vspace{-.6ex}
  \begin{itemize}
   \item[3.1]
    Maslov index of a map to a singular space or a relative pair.

   \item[3.2]
    Monodromy effect and the choice of curve class data in $H_2$.

   \item[3.3]
    The moduli space
     $\overline{\cal M}_{(g,h), (n, \vec{m})}
                        (W/B,L\,|\,[\beta],\vec{\gamma},\mu)$
     of stable maps to fibers of
     $(\widehat{W},\widehat{L})/\widehat{B}$.
  \end{itemize}

 \vspace{-1.6ex}
 \item[4.]
  The moduli space
   $\check{\cal W}^{1,p}_{(g,h),(n,\vec{m})}
    ((\widehat{W},\widehat{L})/\widehat{B}\,|\,[\beta],\vec{\gamma},\mu)$
   of stable $\check{W}^{1,p}$-maps.
  \vspace{-1.2ex}
  \begin{itemize}
   \item[4.1]
    The moduli space
     $\check{\cal W}^{1,p}_{(g,h), (n, \vec{m})}
            (W[k],L[k]\,|\,[\beta],\vec{\gamma},\mu)$
     of stable $\check{W}^{1,p}$-maps to $(W[k],L[k])$,\\
    its relative tangent and relative obstruction bundles.

   \item[4.2]
    The moduli space
     $\check{\cal W}^{1,p}_{(g,h),(n,\vec{m})}
      ((\widehat{W},\widehat{L})/\widehat{B}\,|\,[\beta],\vec{\gamma},\mu)$
     of stable $\check{W}^{1,p}$-maps to fibers of
     $(\widehat{W},\widehat{L})/\widehat{B}$,\\
    the relative $\check{W}^{1,p}$-tangent-obstruction fibration complex.
  \end{itemize}

 \item[5.]
  Construction of a Kuranishi structure for
            $\overline{\cal M}_{(g,h),(n,\vec{m})}
                   (W/B,L\,|\,[\beta],\vec{\gamma},\mu)$.
  \vspace{-.6ex}
  \begin{itemize}
   \item[5.1]
    Family Kuranishi structure modelled in the category
     ${\cal C}_{\mbox{\rm\tiny spsccw}}/{\smallBbb C}$.

   \item[5.2]
    Local transversality and locally regular almost-complex structures.

   \item[5.3]
    Construction of family Kuranishi neighborhoods.

   \item[5.4]
    Construction of a family Kuranishi structure.
  %%%%%%%%%%%%%%%%%%%%%%%%%%%%%%%%
  % \item[5.5]
  %  Family orientations of the family Kuranishi structure.
  %%%%%%%%%%%%%%%%%%%%%%%%%%%%%%%%
  \end{itemize}

 \item[6.]
  The moduli space
   $\overline{\cal M}_{(g,h),(n+l(\vec{s}),\vec{m})}
    (Z,L;D\,|\,\beta^{\prime},\vec{\gamma},\mu^{\prime};\vec{s})$
   of relative stable maps  and\\
  its Kuranishi structure.
  \vspace{-.6ex}
  \begin{itemize}
   \item[6.1]
    The moduli space
     $\overline{\cal M}_{(g,h),(n+l(\vec{s}),\vec{m})}
      (Z,L;D\,|\,\beta^{\prime},\vec{\gamma},\mu^{\prime};\vec{s})$
     of relative stable maps.

   \item[6.2]
    A Kuranishi structure for
      $\overline{\cal M}_{(g,h),(n+l(\vec{s}),\vec{m})}
       (Z,L;D\,|\,\beta^{\prime},\vec{\gamma},\mu^{\prime};\vec{s})$.
  \end{itemize}

 \item[7.]
  Degeneration and gluing of Kuranishi structures  and
  axioms of open Gromov-Witten invariants under a symplectic cut.
  \vspace{-.6ex}
  \begin{itemize}
   \item[7.1]
     The degeneration-gluing relations of Kuranishi structures.

   \item[7.2]
    A degeneration axiom and a gluing axiom for
     open Gromov-Witten invariants under\\   a symplectic cut.
  \end{itemize}

  \item[]\hspace{-1.6em}
  Appendix.
   \parbox[t]{32em}{The equivalence of Li-Ruan/Li's degeneration formula
       and Ionel-Parker's \\ degeneration formula.}
\end{itemize}
} %endsmall

\bigskip

\newpage
\baselineskip 13pt  %Re: 14pt for [11pt] style
                    %Re: 15pt for [12pt] style.

\section{Symplectic cut and the direct system of
         expanded degenerations in the almost-complex category.}

A direct system of expanded degenerations of almost-complex spaces
 that merges the symplectic construction (via multi- symplectic cut)
  in [I-P2: Sec.~2 and Sec.~12] and [L-R: Sec.~3]
  with the algebro-geometric construction (via blow-ups and
  blow-downs) in [Li: Sec.~1]
  without using the full language of stacks is given in this section.
The fibers, up to a relative isomorphism, of the families
 in the system will occur as the targets of open stable maps
 in the problem.
The same construction gives also a direct system of
 expanded relative pairs
 (cf.\ [I-P1: Remark~7.7], [L-R: Sec.~4];
       [Gr-V: Sec.~2], [Li1: Sec.~4.1])
 needed for relative open Gromov-Witten theory.

\bigskip

\subsection{Symplectic cut and the associated expanded degenerations.}

\subsubsection{Expanded degenerations from a symplectic cut.}

\begin{flushleft}
{\bf Symplectic cut and a compatible almost-complex degeneration.}
\end{flushleft}
{\it Symplectic cut} was introduced in [Le] and used in [I-P1], [I-P2],
 and [L-R].
We review it here to fix notations.
Given a free Hamiltonian $S^1$ action on a connected open set
 $U$ of a symplectic manifold $X$ that separates $X$.
Fix a Hamiltonian function $h: U\rightarrow (-l, l)$ of the
 $S^1$-action and let $X-h^{-1}(0)=X_+\disjointunion X_-$.
Then the manifold with boundary
 $\overline{X}_+ := X_+\cup h^{-1}(0)$
 (resp.\ $\overline{X}_- := X_-\cup h^{-1}(0)$) gives rise to
 a symplectic manifold $Y_1$ (resp.\ $Y_2$)
  by taking the quotient of the $S^1$-action on the boundary, and
the boundary $h^{-1}(0)$ descends to a codimension-$2$
 symplectic submanifold $D$ in $Y_1$ (resp.\ $Y_2$).
Let $Y$ be the singular symplectic space from gluing $Y_1$ and $Y_2$
 canonically along $D$.
Then, there is a natural map $\xi: X\rightarrow Y$ that
 is modelled on a symplectic reduction (and hence an $S^1$-bundle)
  over the singular locus $D:= Y_1\cap Y_2$ on $Y$ and
 is a symplectomorphism from $X-\xi^{-1}D$ to $Y-D$.

\bigskip

\noindent
{\bf Definition 1.1.1.1 [symplectic cut].} {\rm
With an abuse of language and a different naming than
  the original work [Le] of Lerman,
 we will call both the map ${\xi}:X\rightarrow Y=Y_1\cup_D Y_2$
  and the singular symplectic space $Y$
  a {\it symplectic cut} of $X$.
} % end-definition

\bigskip

Given a symplectic cut $\xi:X\rightarrow Y=Y_1\cup_D Y_2$,
 one can identify a small neighborhood of $D$ in $Y_1$ (resp.\ $Y_2$)
  with a neighborhood of the zero-section of a complex line bundle
  ${\Bbb L}$ (resp.\ the dual complex line bundle ${\Bbb L}^{\ast}$)
  over $D$ and
 construct a complex $1$-parameter family $\pi: W\rightarrow B$ of
  symplectic spaces $W_{\lambda}:=\pi^{-1}(\lambda)$, $\lambda\in B$,
  with a compatible almost-complex structure $J_{W_{\lambda}}$
  such that
   $W_0$ is symplecto-isomorphic to $Y$, with the restriction of
    $J_{W_0}$ to a neighborhood of $D$ almost-complex-isomorphic to
    the gluing of a neighborhood of the zero-section in ${\Bbb L}$ and
     a neighborhood of the zero-section in ${\Bbb L}^{\ast}$ along $D$,
    and
   $W_{\lambda}$, $\lambda\ne 0$, is symplecto-isomorphic to $X$.
Here
 $B$ is a small neighborhood of $0$ in ${\Bbb C}$ and
 the total space of ${\Bbb L}$ and ${\Bbb L}^{\ast}$ are equipped
  with an $U(1)$-invariant almost-complex structure that combines
   an almost-complex structure $J_D$ on $D$ and
   the complex structure on fiber ${\Bbb C}$
   via a $U(1)$-connection on ${\Bbb L}$ and ${\Bbb L}^{\ast}$.
See [I-P2: Sec.~2] (and also [Go] and [MC-W]) for an explicit
 construction.
The total space $W$ is equipped with a symplectic structure
$\omega_W$
 and a compatible almost-complex structure $J_W$ that gives
 $(\omega_{W_{\lambda}},J_{W_{\lambda}})$ when restricted to $W_{\lambda}$.
We will denote the family $\pi:W\rightarrow B$ also by $W/B$ as
  in algebraic geometry and
 call $W/B$ a {\it compatible almost-complex degeneration} associated
  to the symplectic cut $\xi:X\rightarrow Y$.

Fix and denote a local fiber complex coordinate of ${\Bbb L}$
 (resp.\ ${\Bbb L}^{\ast}$) by $w$ (resp.\ $w^{\prime}$) and
 treat both ${\Bbb L}$ and ${\Bbb L}^{\ast}$ as a $U(1)$-bundle.
Let $0<\varepsilon<1$ be sufficiently small.
Then, possibly after shrinking, we may assume that
 $$
  B\;=\; \{ \lambda\in {\Bbb C}\,:\, |\lambda| < \varepsilon^2/2 \}\,.
 $$
The following defines a subset of ${\Bbb L}\oplus{\Bbb L}^{\ast}$
 $$
  ({\Bbb L}\oplus{\Bbb L}^{\ast})_{\le {\varepsilon}}\;
  =\;
    \{(\,\cdot\,,w, w^{\prime})\,:\,
     |w|\le \varepsilon\,,\,|w^{\prime}|\le \varepsilon\,,\,
     |ww^{\prime}|\le \varepsilon^2/2\}\,.
 $$
It admits a fibration
 $({\Bbb L}\oplus{\Bbb L}^{\ast})_{\le {\varepsilon}} \rightarrow B$
 defined by $(\,\cdot\,, w,w^{\prime})\mapsto ww^{\prime}$.
With this fibration and an adjustment in the construction of $W$,
 there is a decomposition
 $$
  W\;=\; (B\times\overline{U_1}) \cup
         ({\Bbb L}\oplus{\Bbb L}^{\ast})_{\le {\varepsilon}} \cup
         (B\times\overline{U_2})
 $$
 over $B$, where
  $B$ is taken to be
   $\{\lambda\in {\Bbb C}: |\lambda|< \varepsilon^2/2\}$,
  $U_1= Y_1 -$($\varepsilon$-neighborhood
              of the zero-section in ${\Bbb L}$),
  $U_2= Y_2 -$($\varepsilon$-neighborhood
              of the zero-section in ${\Bbb L}^{\ast}$), and
  the gluing is along the related boundary circle-bundle over $B\times D$
   in a way that respects the fibration of and the $U(1)$-action on
    these boundaries over $B\times D$.
This decomposition allows us to construct an expanded almost-complex
 degeneration associated to $W/B$, which we explain in the next two
 themes.

\bigskip

\begin{flushleft}
{\bf Local expanded degenerations in the almost-complex category.}
\end{flushleft}
The expanded degeneration around $D$ in the almost-complex category can
 be described by
  a finite collection of almost-complex manifolds together with
  a collection of gluing isomorphisms
   between open dense almost-complex submanifolds therein
 as follows.

Let
 ${\Bbb L}$ be a complex line bundle on the almost-complex
  manifold with the ${\Bbb C}^{\times}$-structure
  reduced to a $U(1)$-structure, and let
 $\alpha$ be a $U(1)$-connection on ${\Bbb L}$.
This induces a unique $U(1)$-structure on the complex dual line bundle
 ${\Bbb L}^{\ast}$ to ${\Bbb L}$ on $D$ with a $U(1)$-connection
 $\alpha^{\ast}$.
Denote the almost-complex structure on $D$ be $J_D$;
 then the pairs $(J_D,\alpha)$ and $(J_D,\alpha^{\ast})$, together with
 the fiberwise complex structures, determines almost-complex structures
 $J_{\scriptsizeBbb L}$, $J_{{\scriptsizeBbb L}^{\ast}}$
 on the total space (still denoted by the same notation)
 of ${\Bbb L}$ and ${\Bbb L}^{\ast}$ respectively.
$J_{\scriptsizeBbb L}$ and $J_{{\scriptsizeBbb L}^{\ast}}$
 together induce an almost-complex structure
 $J_{\,{\scriptsizeBbb L}\oplus{\scriptsizeBbb L}^{\ast}}$
 on (the total space of) ${\Bbb L}\oplus {\Bbb L}^{\ast}$.

Denote
 the zero-section of ${\Bbb L}$ or ${\Bbb L}^{\ast}$ by ${\mathbf 0}$
  and
 the projection map
  ${\Bbb L}\oplus{\Bbb L}^{\ast}\rightarrow {\Bbb L}$
  (resp.\ ${\Bbb L}\oplus{\Bbb L}^{\ast}\rightarrow {\Bbb L}^{\ast}$)
 by $\pr$ (resp.\ $\pr^{\prime}$).
Fix a system of local trivializations and
 $U(1)$-valued transition functions for ${\Bbb L}$.
They induce a system of local trivializations and
 $U(1)$-valued transition functions on ${\Bbb L}^{\ast}$, and
then a system of local trivializations and
 transition functions on ${\Bbb L}\oplus {\Bbb L}^{\ast}$.
With respect to this, the map
 $\pi: {\Bbb L}\oplus{\Bbb L}^{\ast}\rightarrow {\Bbb C}$
  given by $(x; w,w^{\prime})\mapsto ww^{\prime}$ is well-defined
  and compatible with
  $J_{\,{\scriptsizeBbb L}\oplus{\scriptsizeBbb L}^{\ast}}$,
 where $(x,w)$ (resp.\ $(x,w^{\prime})$) are local coordinates
        for ${\Bbb L}$ (resp.\ ${\Bbb L}^{\ast}$)
        in the specified local trivialization.
Let $M_{\lambda}\subset {\Bbb L}\oplus{\Bbb L}^{\ast}$ be
 the preimage $\pi^{-1}(\lambda)$ of $\lambda\in{\Bbb C}$.
For $\lambda\in {\Bbb C}-\{0\}$, $M_{\lambda}$ is isomorphic to
 ${\Bbb L}-{\mathbf 0}\;(\simeq {\Bbb L}^{\ast}-{\mathbf 0})$
 as an almost-complex submanifold.
For $\lambda=0$, $M_0={\Bbb L}\vee {\Bbb L}^{\ast}$,
 the union of ${\Bbb L}$ and ${\Bbb L}^{\ast}$ with
 the zero-sections glued by the canonical isomorphism with $D$.
Thus, the family
 $\pi:{\Bbb L}\oplus{\Bbb L}^{\ast}\rightarrow {\Bbb C}$
 is a smoothing of $M_0$ over $D$ in the almost-complex category.

\bigskip

\noindent
{\it Notation 1.1.1.2
     $[\, M_{\lambda}$, $\lambda\ne 0$, from gluing$\,]$.}
 Associated to the $U(1)$-structure on ${\Bbb L}$ and ${\Bbb L}^{\ast}$
  is a well-defined norm function $|\cdot|$ on fibers of
  ${\Bbb L}$ and ${\Bbb L}^{\ast}$.
 Let
  ${\Bbb L}_{>\delta}=\{\,|w|>\delta\,\}\subset {\Bbb L}$,
  ${\Bbb L}^{\ast}_{>\delta}=\{\,|w^{\prime}|>\delta\,\}
    \subset {\Bbb L}^{\ast}$, and,
 similarly, for
  ${\Bbb L}_{\le\delta}$,
  ${\Bbb L}_{[\delta_1\,,\,\delta_2]}$,
  ${\Bbb L}^{\ast}_{\le \delta}$
  ${\Bbb L}^{\ast}_{[\delta_1\,,\,\delta_2]}$, $\cdots$, etc..
 These bundles over $D$ are equipped with the $U(1)$-connection
  (still denoted by $\alpha$ and $\alpha^{\ast}$) from the restriction
  of that on ${\Bbb L}$ and ${\Bbb L}^{\ast}$ respectively.
 The local fiberwise maps
  $(x,w^{\prime})\mapsto (x,w)=(x,\lambda/w^{\prime})$,
  glue to a bundle isomorphism
  $$
   \varphi_{\lambda}\; :\;
    {\Bbb L}^{\ast}_{[|\lambda|/\delta\,,\,\delta]}\;
      \stackrel{\sim}{\longrightarrow}\;
      {\Bbb L}_{[|\lambda|/\delta\,,\,\delta]}\,,\hspace{2em}
   (x,w^{\prime})\mapsto (x,w)=(x,\lambda/w^{\prime})
  $$
  for $0\le |\lambda|<\delta^2$
  such that $\varphi_{\lambda}^{\ast}\,\alpha=-\,\alpha^{\ast}$.
 Thus, $\varphi_{\lambda}$ is an isomorphism in the category of
  almost-complex manifolds as well.
 In terms of this, $M_{\lambda}$, $\lambda\ne 0$,
  is the almost-complex manifold obtained from gluing
  ${\Bbb L}_{>|\lambda|/\delta}$ and
   ${\Bbb L}^{\ast}_{>|\lambda|/\delta}$, with $|\lambda|<\delta^2$,
  by $\varphi_{\lambda}$.
 The maps
  $$
   \begin{array}{ccccccccccc}
    \theta_{\lambda} & :
     & {\Bbb L}_{>0} & \longrightarrow & M_{\lambda}
     & \hspace{2em}\mbox{and}\hspace{2em}
    & \theta^{\prime}_{\lambda}
     & : & {\Bbb L}^{\ast}_{>0} & \longrightarrow & M_{\lambda}\\[.6ex]
    && (\,\cdot\,,w) & \longmapsto & (\,\cdot\,,w,\frac{\lambda}{w})
     &&&& (\,\cdot\,,w^{\prime})   & \longmapsto
        & (\,\cdot\,,\frac{\lambda}{w^{\prime}}, w^{\prime})
   \end{array}
  $$
  are almost-complex isomorphisms.
 We will denote the restriction of $\theta_{\lambda}$
  (resp.\ $\theta^{\prime}_{\lambda}$,
          $\theta_{\lambda}\cup\theta^{\prime}_{\lambda}$) to
  the subsets
   ${\Bbb L}_{[\delta_1,\delta_2]}$, ..., etc.\ of ${\Bbb L}_{>0}$
  (resp.\ ${\Bbb L}^{\ast}_{[\delta_1^{\prime},\delta_2^{\prime}]}$ ...
           of ${\Bbb L}^{\ast}_{>0}$,
          ${\Bbb L}_{[\delta_1,\delta_2]}
            \cup {\Bbb L}^{\ast}_{[\delta_1^{\prime},\delta_2^{\prime}]}$
           of ${\Bbb L}_{>0}\cup{\Bbb L}^{\ast}_{>0}$)
  by $\theta_{\lambda\,;\,[\delta_1,\delta_2]}$
  (resp.\ $\theta^{\prime}
            _{\lambda\,;\,[\delta_1^{\prime},\delta_2^{\prime}]}$,
          $\theta_{\lambda\,;\,[\delta_1,\delta_2]} \cup
           \theta^{\prime}
             _{\lambda\,;\,[\delta_1^{\prime},\delta_2^{\prime}]}$).

\bigskip

For $k\in {\Bbb Z}_{\ge 0}$, let
 $B[k]= {\Bbb C^{k+1}}$,
   with coordinates $(\lambda_0,\,\ldots,\,\lambda_k)$, and
 $\pr_i: B[k]\rightarrow {\Bbb C}$ be the $i$-th coordinate projection map.
Let $({\Bbb L}\oplus {\Bbb L}^{\ast})_i
       =\pr_i^{\ast}({\Bbb L}\oplus{\Bbb L}^{\ast}) $,
 $i=0,\,\ldots,\,k$, be the pulled-back of
 $\pi:{\Bbb L}\oplus{\Bbb L}^{\ast}\rightarrow {\Bbb C}$
 to $B[k]$ via $\pr_i$.
The local coordinates of
 $({\Bbb L}\oplus {\Bbb L}^{\ast})_i$ will be denoted by
 $(\lambda_0,\,\ldots,\,\lambda_i,\,\ldots,\,\lambda_k;
   x, w_i, w_i^{\prime})$ with $\lambda_i=w_iw_i^{\prime}$.
Let
 $({\Bbb L}\oplus{\Bbb L}^{\ast})_i^0
   := ({\Bbb L}\oplus {\Bbb L}^{\ast})_i-\pr_i^{\ast}{\Bbb L}^{\ast}$,
   which is $\{w_i\ne 0\}$ in local coordinates,
  and
 $({\Bbb L}\oplus{\Bbb L}^{\ast})_i^{\infty}
   := ({\Bbb L}\oplus {\Bbb L}^{\ast})_i-\pr_i^{\ast}{\Bbb L}$,
  which is $\{w_i^{\prime}\ne 0\}$ in local coordinates.
We will use coordinates
 $(\lambda_0,\,\ldots,\,\lambda_i,\,\ldots,\,\lambda_k;\,
   x,\,w_i,\,\lambda_i/w_i)$ and
 $(\lambda_0,\,\ldots,\,\lambda_i,\,\ldots,\,\lambda_k;\,
   x,\,\lambda_i/w_i^{\prime},\,w_i^{\prime})$ respectively
 for these two open dense almost-complex submanifolds of
 $({\Bbb L}\oplus{\Bbb L}^{\ast})_i$.
In terms of these, the following map

 \vspace{-2ex}
 {\small
 $$
  \begin{array}{ccc}
   ({\smallBbb L}\oplus {\smallBbb L}^{\ast})_{i-1}^{\infty}
    & \stackrel{\varphi_{i-1,i}\raisebox{-1.2ex}{}}{\longrightarrow}
    & ({\smallBbb L}\oplus {\smallBbb L}^{\ast})_i^0          \\[.6ex]
   (\lambda_0,\,\ldots,\,\lambda_{i-1},\,\lambda_i,\,
       \ldots,\,\lambda_k;
       x, \frac{\lambda_{i-1}}{w_{i-1}^{\prime}}, w_{i-1}^{\prime})
    & \longmapsto
    & (\lambda_0,\,\ldots,\,\lambda_{i-1},\,\lambda_i,\,
          \ldots,\,\lambda_k;
        x,\frac{1}{w_{i-1}^{\prime}}, \lambda_i w_{i-1}^{\prime})
  \end{array}
 $$
 } % end-small

 \vspace{-2ex}
 \noindent
 is an isomorphism in the almost-complex category for $i=1,\,\ldots,\,k$.
The system
 $$
  \left( \{({\Bbb L}\oplus{\Bbb L}^{\ast})_i\}_{i=0}^k\,,\,
   \{\varphi_{\,i-1,\,i}\}_{i=1}^k \right)
 $$
 of almost-complex manifolds and gluing data
 determines an almost-complex manifold
 $({\Bbb L}\oplus{\Bbb L}^{\ast})[k]$ that fibers over $B[k]$.
 %%%%%%%%%%%%%%%%%%%%%%%
 % ({\sc Figure} ???.)
 %
 % \marginpar{\raggedright\tiny $\bullet$
 % {\sc Figure ???} after the definition.}
 %%%%%%%%%%%%%%%%%%%%%%%

\bigskip

\noindent
{\bf Definition 1.1.1.3
     [expanded degeneration of
                $({\Bbb L}\oplus{\Bbb L}^{\ast})/{\Bbb C}$].} {\rm
 We will call the family of almost-complex spaces as constructed above,
  $\pi[k]: ({\Bbb L}\oplus{\Bbb L}^{\ast})[k] \rightarrow B[k]$
  (in short hand: $({\Bbb L}\oplus{\Bbb L}^{\ast})[k]/B[k]$),
  the {\it $k$-th expanded degeneration} of
  the degeneration $\pi:{\Bbb L}\oplus{\Bbb L}^{\ast}\rightarrow {\Bbb C}$.
} % end-definition

%\begin{figure}[htbp]
% \setcaption{\small
% {\sc Figure} ?-?-?.
%  \baselineskip 12pt
%  How the expanded degeneration
%   $({\smallBbb L}\oplus{\Bbb L}^{\ast})[k]$ over $B[k]$ is constructed
%   via gluing is indicated.
%  By construction,
%   each $({\smallBbb L}\oplus{\smallBbb L}^{\ast})_i/B[k]$
%   is open dense in $({\smallBbb L}\oplus{\smallBbb L}^{\ast})[k]/B[k]$.
% } % end-setcaption
%\centerline{\psfig{figure=loc-exp.eps,width=13cm,caption=}}
%\end{figure}

\bigskip

We will use the above gluing construction of
 $({\Bbb L}\oplus{\Bbb L}^{\ast})[k]/B[k]$ as the foundation
 for the rest of the discussion on expanded degenerations.

The natural maps from pull-backs can be re-scaled to give maps
 $\tilde{\mathbf p}[k]_i :
  ({\Bbb L}\oplus{\Bbb L}^{\ast})_i
  \rightarrow {\Bbb L}\oplus{\Bbb L}^{\ast}$ defined by
 \begin{eqnarray*}
  \lefteqn{
   (\lambda_0,\,\ldots,\,\lambda_{i-1},\,\lambda_i,\,\lambda_{i+1},\,
               \ldots,\,\lambda_k; x,\,w_i,\,w_i^{\prime})  }\\[.6ex]
  &&
   \longmapsto\;
   (\lambda_0\,\cdots\,\lambda_k; x,\,
     (\lambda_0\,\cdots\,\lambda_{i-1})\,w_i,\,
     (\lambda_{i+1}\,\cdots\,\lambda_k)\,w_i^{\prime})
 \end{eqnarray*}
 with $\,w_i w_i^{\prime}=\lambda_i\,$, for $i=0,\,\ldots,\,k$.
These maps glue to a map
 $\tilde{\mathbf p}[k]: ({\Bbb L}\oplus{\Bbb L}^{\ast})[k]
                      \rightarrow {\Bbb L}\oplus{\Bbb L}^{\ast}$
 over ${\mathbf p}[k]:B[k]\rightarrow {\Bbb C}$
 in the almost-complex category.

All the gluing isomorphisms $\varphi_{\,i-1,\,i}$, $i=1,\,\ldots,\,k$,
 are maps over $B[k]$.
Thus, the fibers of
 $\pi[k]:({\Bbb L}\oplus{\Bbb L}^{\ast})[k]\rightarrow B[k]$
 can be described by the corresponding gluing over a fixed values of
 $\vec{\lambda}=(\lambda_0,\,\ldots,\,\lambda_k)$, as follows.
First, note that the $\varphi_{\lambda}$,
 $\lambda\in {\Bbb C}^{\times}$, defined in Notation 1.1.1.2
   % Notation [$M_{\lambda}$, $\lambda\ne 0$, from gluing]
 gives as well an isomorphism
 $\varphi_{\lambda} : {\Bbb L}^{\ast}_{>0}\;
                          \rightarrow\; {\Bbb L}_{>0}$
 in the almost-complex category.
The gluing of ${\Bbb L}$ and ${\Bbb L}^{\ast}$ by $\varphi_{\lambda}$
 gives a ruled (i.e.\ ${\Bbb P}^1$-fibered) manifold $\Delta$ over $D$
 with a well-defined almost-complex
 structure that contains ${\Bbb L}$ and ${\Bbb L}^{\ast}$ as
 open almost-complex submanifolds.
Different choices of $\lambda$ give rise to isomorphic almost-complex
 manifolds over $D$ with a such isomorphism provided by the identity
 map on ${\Bbb L}$ (and hence on $\Delta$).
We will take $\Delta$ as from the gluing $\varphi_1$.
Denote the zero-section ${\mathbf 0}$ of ${\Bbb L}$
 (resp.\ ${\Bbb L}^{\ast}$) by $D_0$ (resp.\ $D_{\infty}$) in $\Delta$.
Let
 $({\Bbb L}\vee{\Bbb L}^{\ast})_i := {\Bbb L}_i\vee {\Bbb L}^{\ast}_i$,
  $i=0,\,\ldots,\,k^{\prime}$, be identical copies of
  ${\Bbb L}\vee {\Bbb L}^{\ast}$ and
 $({\Bbb L}\vee {\Bbb L}^{\ast})_{[k^{\prime}]}$ be the gluing of
  $({\Bbb L}\vee {\Bbb L}^{\ast})_i$, $i=0,\,\ldots,\,k^{\prime}$, by
  $\varphi_{i;1}:
    {\Bbb L}^{\ast}_{i-1}\rightarrow {\Bbb L}_{i}$,
  $(x, w_{i-1}^{\prime})\mapsto (x, w_i)=(x,1/w_{i-1}^{\prime})$.
Then, as an almost-complex space,
 $$
  ({\Bbb L}\vee {\Bbb L}^{\ast})_{[k^{\prime}]}\;
  =\;
    {\Bbb L}\cup_{\,{\mathbf 0}=D_{1,\infty}}\Delta_1
       \cup_{D_{1,0} = D_{2,\infty}}\,\cdots\,
       \cup_{D_{k^{\prime}-1,0}=D_{k^{\prime},\infty}}\Delta_{k^{\prime}}
       \cup_{D_{k^{\prime},0}={\mathbf 0}} {\Bbb L}^{\ast}\,,
 $$
 where $(\Delta_i;D_{i,0},D_{i,\infty})=(\Delta; D_0,D_{\infty})$.
There is a natural map
 $({\Bbb L}\vee{\Bbb L}^{\ast})_{[k^{\prime}]}
                         \rightarrow {\Bbb L}\vee{\Bbb L}^{\ast}$
 that restricts to the identity map on ${\Bbb L}$ and ${\Bbb L}^{\ast}$,
      and collapses all $\Delta_i$ to $D$.
The natural ${\Bbb G}_m := {\Bbb C}^{\times}$-action on ${\Bbb L}$
 extends to a ${\Bbb G}_m$-action on $\Delta$ as a group of
 automorphisms of $\Delta$ over $D$ in the almost-complex category.
For $\sigma\in{\Bbb G}_m$, the induced action coincides with
 the composition $\varphi_{\sigma}\circ {\varphi_1}^{-1}$
 on $\Delta- D_0\cup D_{\infty}$.
It leaves $D_0\cup D_{\infty}$ fixed.
The relative automorphism group
 $\Aut( ({\Bbb L}\vee{\Bbb L}^{\ast})_{[k^{\prime}]}/
                                    {\Bbb L}\vee{\Bbb L}^{\ast} )$
 in the almost-complex category is the product
 $\prod_{\,i=1}^{\,k^{\prime}}\Aut(\Delta_i/D)
                          =({\Bbb C}^{\times})^{\,k^{\prime}}$.
Now let
 $I=\{i_0,\,\ldots,\,i_{k^{\prime}}\}$ be a subset of
  $\{0,\,\ldots,\,k\}$ and
 $\dot{H}_I$ be the locally closed submanifold of $B[k]$,
  whose points have coordinates $\lambda_i=0$ exactly when $i\in I$.
$B[k]$ is the disjoint union of $\dot{H}_I$, where $I$ runs over
 all the subsets of $\{0,\,\ldots,\,k\}$.
Let $\vec{\lambda}=(\lambda_0,\,\ldots,\,\lambda_k)\in \dot{H}_I$.
Then,
 $\pi[k]^{-1}(\vec{\lambda})$ is the almost-complex space from
 the system
 $$
  \left( \{M_{\lambda_i}\}_{i=0}^k\,,\,
   \{\varphi_{\,i-1,\,i\,;\,\vec{\lambda}}\}_{i=1}^k \right)\,,
 $$
 where
  $M_{\lambda_i}=\{w_iw_i^{\prime}=\lambda_i\}
     \subset {\Bbb L}_i\oplus{\Bbb L}_i^{\ast}$ and
  $\varphi_{\,i-1,\,i\,;\,\vec{\lambda}}:
     M_{\lambda_{i-1}}-{\Bbb L}_{i-1}
               \rightarrow M_{\lambda_i}-{\Bbb L}_i^{\ast}$
   is given by
   $(\lambda_{i-1}/w_{i-1}, w_{i-1}^{\prime})
     \mapsto
     (w_i,\lambda_i/w_i)=(1/w_{i-1}^{\prime}, \lambda_i w_{i-1}^{\prime})$.
This system can be reduced to the following system
 $$
  \left( \{M_{\lambda_{\,i_j}}\}_{j=-1}^{k^{\prime}+1}\,,\,
         \{ \tilde{\varphi}_{\,j-1,\,j\,;\,\vec{\lambda}} \}
                                       _{j=0}^{k^{\prime}+1} \right)\,,
 $$
 where
  $M_{\lambda_{\,i_{-1}}}={\Bbb L}_0-\{\mathbf 0\}$ and
  $M_{\lambda_{\,i_{k^{\prime}+1}}}={\Bbb L}_k^{\ast}-\{\mathbf 0\}$
   by convention, and
  $\tilde{\varphi}_{\,j-1,\,j\,;\,\vec{\lambda}}:
     M_{\lambda_{\,i_{j-1}}}-{\Bbb L}_{i_{j-1}}
        \rightarrow M_{\lambda_{\,i_j}}-{\Bbb L}_{i_j}^{\ast}$
   is the composition
   $\varphi_{i_j-1\,,\,i_j}\, \circ\,
     \cdots\,\circ\, \varphi_{i_{j-1}+1\,,\,i_{j-1}+2}\,
     \circ\, \varphi_{i_{j-1}\,,\,i_{j-1}+1}$
   with $\varphi_{i_{-1}\,,0}$ and $\varphi_{k\,,\,k+1}$
    being identity maps by convention.

In summary,

\bigskip

\noindent
{\bf Lemma 1.1.1.4 [natural map and its fibers].} {\it
 Let ${\mathbf p}[k]:B[k]\rightarrow {\Bbb C}$ be the product map
  defined by
 $(\lambda_0,\,\ldots,\,\lambda_k)\mapsto \lambda_0\,\cdots\,\lambda_k$.
 Then there is a natural map
  $\tilde{\mathbf p}[k]:
   ({\Bbb L}\oplus{\Bbb L}^{\ast})[k]
    \rightarrow ({\Bbb L}\oplus{\Bbb L}^{\ast})$
  in the almost-complex category that covers ${\mathbf p}[k]$.
 The fiber of $\pi[k]$ at $\vec{\lambda}\in \dot{H}_I$ is isomorphic
  to $M_{\lambda_0\,\cdots\,\lambda_k}$, if $I=\emptyset$, and
  to $({\Bbb L}\vee{\Bbb L}^{\ast})_{[k^{\prime}]}$,
     if $I=\{i_0,\,\ldots,\,i_{k^{\prime}}\}$ is non-empty.
 In particular, $\tilde{\mathbf p}[k]$ is an isomorphism when restricted
  to the fiber over a point in the complement of complex codimension-$2$
  coordinate subspaces.
} % end-lemma

\bigskip

\begin{flushleft}
{\bf Expanded almost-complex degenerations associated to $W/B$.}
\end{flushleft}
Given a fibered almost-complex space $W/B$ from a symplectic cut
 as constructed above,
recall the decomposition over $B$
 $$
  W\;=\; (B\times\overline{U_1}) \cup
         ({\Bbb L}\oplus{\Bbb L}^{\ast})_{\le {\varepsilon}} \cup
         (B\times\overline{U_2})\,.
 $$

Let
 $$
  \begin{array}{lllll}
   \overline{U_1}[k]
    & :=      & \pr_0^{\ast}((B\times \overline{U_1})/B)
    & \simeq  & B[k]\times \overline{U_1}\,, \\[.6ex]
   \overline{U_2}[k]
    & :=      & \pr_k^{\ast}((B\times \overline{U_2})/B)
    & \simeq  & B[k]\times \overline{U_2}\,,
  \end{array}
 $$
 and
define the {\it $\varepsilon$-truncation}
 $({\Bbb L}\oplus{\Bbb L}^{\ast})[k]_{\le\varepsilon}$  of
 $({\Bbb L}\oplus{\Bbb L}^{\ast})[k]$ as follows.
\begin{itemize}
 \item[$\cdot$]
 First consider the preliminary truncation of
  $({\Bbb L}\oplus{\Bbb L}^{\ast})_0$ and
  $({\Bbb L}\oplus{\Bbb L}^{\ast})_k$
  defined respectively by
  $$
   {\cal U}_0^{\prime}\;:=\; \{|w_0|\le \varepsilon\}
     \subset ({\Bbb L}\oplus{\Bbb L}^{\ast})_0
    \hspace{1em}\mbox{and}\hspace{1em}
   {\cal U}_k^{\prime}\;:=\; \{|w_k^{\prime}|\le \varepsilon\}
     \subset ({\Bbb L}\oplus{\Bbb L}^{\ast})_k\,.
  $$

 \item[$\cdot$]
 Then the pair of gluing maps
  $$
   (\, \varphi_{i-1,i}\circ\,\cdots\,\circ\varphi_{0,1}\,,\,
       \varphi_{i,i+1}\circ\circ\,\cdots\,\circ\varphi_{k-1,k}^{-1} \,)\,,\;
  $$
  from $({\cal U}_0^{\prime},{\cal U}_k^{\prime})$ to
  ${\cal U}_0^{\prime}$,
  $({\Bbb L}\oplus{\Bbb L}^{\ast})_i$, $i=1,\,\ldots\,,k-1$, and
  ${\cal U}_k^{\prime}$ respectively induces a truncation thereon
  defined by
  $$
   \begin{array}{lcl}
    {\cal U}_0 & =
     & \{\,  (\lambda_0,\,\ldots,\,\lambda_k; x,\,w_0,\,w_0^{\prime})
         \in ({\Bbb L}\oplus{\Bbb L}^{\ast})_0\;:\;
         |w_0| \le \varepsilon\,,\;
         |w_0^{\prime}| \le
          \frac{\varepsilon}{|\lambda_1\,\cdots\,\lambda_k|}
       \, \}\,,    \\[1.6ex]
    {\cal U}_i & =
     & \{\,  (\lambda_0,\,\ldots,\,\lambda_k; x,\,w_i,\,w_i^{\prime})
         \in ({\Bbb L}\oplus{\Bbb L}^{\ast})_i\;:\;
         |w_i| \le
          \frac{\varepsilon}{|\lambda_0\,\cdots\,\lambda_{i-1}|}
            \,,\;
         |w_i^{\prime}| \le
          \frac{\varepsilon}{|\lambda_{i+1}\,\cdots\,\lambda_k|}
       \, \}\,,    \\[1ex]
    && \hspace{24em}i\;=\;1\,,\,\ldots\,,\,k-1\,,\\[.6ex]
    {\cal U}_k & =
     & \{\,  (\lambda_0,\,\ldots,\,\lambda_k; x,\,w_k,\,w_k^{\prime})
         \in ({\Bbb L}\oplus{\Bbb L}^{\ast})_k\;:\;
         |w_k| \le
          \frac{\varepsilon}{|\lambda_0\,\cdots\,\lambda_{k-1}|}
            \,,\;
         |w_k^{\prime}| \le \varepsilon
       \, \}\,.
   \end{array}
  $$

 \item[$\cdot$]
 The gluing map $\varphi_{i-1,i}$ sends
  ${\cal U}_{i-1}\cap
   ({\Bbb L}\oplus{\Bbb L}^{\ast})_i^{\infty}$  to
  ${\cal U}_i\cap ({\Bbb L}\oplus{\Bbb L}^{\ast})_i^0$.
 Thus the collection $\{\varphi_{i-1,i}\}_{i=1}^k$ of gluing maps
  glue the collection $\{{\cal U}_i\}_{i=0}^k$ of almost-complex
  manifolds to an almost- complex manifold over $B[k]$,
  which will be denoted
   $({\Bbb L}\oplus{\Bbb L}^{\ast})[k]_{\le \varepsilon}$  and
   called the {\it $\varepsilon$-truncation} of
   $({\Bbb L}\oplus{\Bbb L}^{\ast})[k]$.
\end{itemize}

\noindent
Then the gluing
 $W\;=\; (B\times\overline{U_1}) \cup
          ({\Bbb L}\oplus{\Bbb L}^{\ast})_{\le {\varepsilon}} \cup
          (B\times\overline{U_2})$
 induces via $\pr_1^{\ast}$ and $\pr_k^{\ast}$
 a canonical gluing of
 $$
  \overline{U_1}[k] \cup
  ({\Bbb L}\oplus{\Bbb L}^{\ast})[k]_{\le\varepsilon} \cup
  \overline{U_2}[k]\; =:\; W[k]
 $$
 over $B[k]$, that goes with a map
 $\pi[k]:W[k]\rightarrow B[k]$.
Here we shrink and re-define $B[k]$ to be
 $$
  B[k]\;:=\;
  \{(\lambda_0,\,\ldots,\,\lambda_k)\;:\;
   |\lambda_i| < \varepsilon^2/2\,,\;
    i\,=\,0\,,\,\ldots\,,\,k \}\,.
 $$
The fiber of
  $\overline{U_1}[k]$,
  $({\Bbb L}\oplus{\Bbb L}^{\ast})[k]_{\le {\varepsilon}}$, and
  $\overline{U_2}[k]$ over the same point in $B[k]$
 glue to an almost-complex space.

\bigskip

\noindent
{\bf Definition 1.1.1.5 [expanded degeneration of $W/B$].} {\rm
 The family $\pi[k]: W[k]\rightarrow B[k]$, in short $W[k]/B[k]$,
 of almost-complex spaces is called an
 {\it expanded almost-complex degeneration} of $W/B$.
} % end-definition

\bigskip

Let $I\subset \{1,\,\ldots,\,n\}$ be non-empty.
Then it follows from the construction that,
 for $\vec{\lambda}\in B[k]\cap \dot{H}_I$,
 the fiber almost-complex space
  $W[k]_{\vec{\lambda}}:=\pi[k]^{-1}(\vec{\lambda})$ is
  almost-complex-isomorphic to
  $$
   Y_{[k^{\prime}]}\; :=\;
    Y_1\cup_{D=D_{1,\infty}}\Delta_1
       \cup_{D_{1,0} = D_{2,\infty}}\,\cdots\,
       \cup_{D_{k^{\prime}-1,0}=D_{k^{\prime},\infty}}\Delta_{k^{\prime}}
       \cup_{D_{k^{\prime},0}=D} Y_2
  $$
  with $k^{\prime}=|I|$.
By construction, there is an almost-complex morphism
 $$
  \tilde{\mathbf p}[k]\;:\; W[k]/B[k]\;\longrightarrow\; W/B\,,
 $$
 cf.\ Lemma 1.1.1.4.
      % Lemma [natural map and its fibers]

\bigskip

\bigskip

\begin{flushleft}
{\bf Neck-trunk decompositions of $W[k]/B[k]$ and re-forgings.}
\end{flushleft}
We introduce here
  neck-trunk decompositions of $W[k]/B[k]$ and re-forging morphisms
 for the discussion of rigidification in Sec.~4.2 and
     the gluing construction of a Kuranishi neighborhood in Sec.~5.3.

Recall $0<\varepsilon<1$ and consider the decomposition over $B\,$:
 $$
  ({\Bbb L}\oplus{\Bbb L}^{\ast})_{\le 1/\varepsilon}\;
  =\; ({\Bbb L}\oplus{\Bbb L}^{\ast})
                       _{[\varepsilon, 1/\varepsilon]\,;\, 1}
      \cup ({\Bbb L}\oplus{\Bbb L}^{\ast})_{\le \varepsilon}
      \cup ({\Bbb L}\oplus{\Bbb L}^{\ast})
                       _{[\varepsilon, 1/\varepsilon]\,;\, 2}\,,
 $$
 where
 $$
  \begin{array}{lll}
   ({\Bbb L}\oplus{\Bbb L}^{\ast})
                     _{[\varepsilon, 1/\varepsilon]\,;\, 1}
    & =
    & \{(\,\cdot\,,w, w^{\prime})\,:\,
      \varepsilon \le |w|\,\le 1/\varepsilon\,,\,
      |ww^{\prime}|\le \varepsilon^2/2\}\,, \\[.6ex]
   ({\Bbb L}\oplus{\Bbb L}^{\ast})
                      _{[\varepsilon, 1/\varepsilon]\,;\,2}
    & =
    & \{(\,\cdot\,,w, w^{\prime})\,:\,
      \varepsilon \le |w^{\prime}|\le 1/\varepsilon\,,\,
      |ww^{\prime}|\le \varepsilon^2/2\}\,,
  \end{array}
 $$
 all three components fiber over $B$ via
  $(\,\cdot\,,w,w^{\prime})\mapsto ww^{\prime}$, and
 the gluing is along their horizonal boundary over $B$.
Let
 $$
  \begin{array}{llll}
  \Neck[k]_i  & =
   & \mbox{the image of
      $\pr_i^{\ast} (({\Bbb L}\oplus{\Bbb L}^{\ast})_{\le \varepsilon})$
      in $W[k]$}\,,
   & i\;=\;0\,,\,\ldots\,,\,k\,,                         \\[.6ex]
  \Trunk[k]_{i;1}  & =
   & \mbox{the image of
      $\pr_i^{\ast} (({\Bbb L}\oplus{\Bbb L}^{\ast})
        _{[\varepsilon, 1/\varepsilon]\,;\,1})$ in $W[k]$}\,,
   & i\;=\; 1\,,\,\ldots\,,\,k\,,                        \\[.6ex]
  \Trunk[k]_{i;2}  & =
   & \mbox{the image of
      $\pr_i^{\ast} (({\Bbb L}\oplus{\Bbb L}^{\ast})
        _{[\varepsilon, 1/\varepsilon]\,;\,2})$ in $W[k]$}\,,
   & i\;=\; 0\,,\,\ldots\,,\,k-1\,.
  \end{array}
 $$
Then all these spaces fiber over $B[k]$.
Furthermore,
since $0<\varepsilon<1$ and
 $|\lambda_i|<\varepsilon^2$ for all
 $\vec{\lambda}=(\lambda_0,\,\ldots,\,\lambda_k)\in B[k]$,
one has
 $$
  \Trunk[k]_{i-1;2}\;
   =\; \Trunk[k]_{i;1}\; =:\; \Trunk[k]_i
   \hspace{10em}\mbox{and}\hspace{5.4em}
 $$

 \vspace{-4ex}
 {\small
 \begin{eqnarray*}
  \lefteqn{W[k]/B[k]\; =} \\[.6ex]
  &&
   (\, \overline{U_1}[k]  \,\cup\,
       \Neck[k]_0  \,\cup\,  \Trunk[k]_1  \,\cup\,  \Neck[k]_1
         \,\cup\,\cdots\, \cup\,
       \Trunk[k]_k  \,\cup\,  \Neck[k]_k  \,\cup\,
       \overline{U_2}[k] \,)/B[k]\,,
 \end{eqnarray*}
 } % end-small

 \vspace{-2ex}
 \noindent
 where the gluings are along the horizontal boundary of each component
   over $B[k]$ and
  are all induced by the gluing maps $\varphi_{i-1,i}$'s.
 We shall call this a ($\varepsilon$-){\it neck-trunk decomposition},
  $\Neck[k]_i$ a ($\varepsilon$-){\it neck region}, and
  $\Trunk[k]_i$ a ($\varepsilon$-){\it trunk region}
  of $W[k]/B[k]$.
When in need of expressing $\varepsilon$ explicitly, we will denote
 a neck (resp.\ trunk) by $\Neck_{\varepsilon}[k]_i$
 (resp.\ $\Trunk_{\varepsilon}[k]_i$).

Denote the fiber of
 $\Neck[k]_i$, $\Trunk[k]_j$,
  $i=0\,,\,\ldots\,,\,k$, $j=1\,,\,\ldots\,,\,k$,
 over $\vec{\lambda}\in B[k]$ by
 $\Neck[k]_{i\,,\,\vec{\lambda}}$, $\Trunk[k]_{j\,,\,\vec{\lambda}}$
 respectively. Then $W[k]_{\vec{\lambda}}$ is divided to
  a gluing-along-boundary:

 \vspace{-2ex}
 {\small
 \begin{eqnarray*}
  \lefteqn{W[k]_{\vec{\lambda}} =} \\[.6ex]
  && \overline{U_1}  \,\cup\,
       \Neck[k]_{0\,,\,\vec{\lambda}}  \,\cup\,
       \Trunk[k]_{1\,,\,\lambda}  \,\cup\,
       \Neck[k]_{1\,,\,\lambda}   \,\cup\,\cdots\, \cup\,
       \Trunk[k]_{k\,,\,\lambda}  \,\cup\,
       \Neck[k]_{k\,,\,\lambda}  \,\cup\,
       \overline{U_2}\,.
 \end{eqnarray*}
 } % end-small

\vspace{-2ex}
\noindent
Recall Notation 1.1.1.2
         % Notation [$M_{\lambda}$, $\lambda\ne 0$, from gluing]
 and that $W[k]_{\vec{0}}=Y_{[k]}$ and denote
  $D_i=\Delta_i\cap \Delta_{i+1}$, $i=0,\,\ldots,\,k$.
Let
 $\vec{\lambda}\in B[k]$  and
 $0\le i_0 <\,\cdots\, < i_{k^{\prime}} \le k$ be the associated
  indices so that $\lambda_{i_j}=0$.
Then there are canonical almost-complex morphisms built-in to
 the construction:

 \vspace{-1ex}
 {\small
 $$
  \begin{array}{cl}
   \Neck[k]_{i\,,\,\vec{0}}\;  =\; \Neck[k]_{i\,,\,\vec{\lambda}}\,,
    \hspace{1em}
   \Trunk[k]_{j\,,\,\vec{0}}\; =\; \Trunk[k]_{j\,,\,\vec{\lambda}}\,,
     & \hspace{-1.4em}
       i\,,\,j\in \{i_0,\,\ldots,\,i_{k^{\prime}}\}\,;        \\[1.6ex]
   \pr_i^{\ast}(\theta
                  _{\lambda_i\,;\,[\sqrt{|\lambda_i|}, \varepsilon]} \cup
                \theta^{\prime}
                  _{\lambda_i\,,\,[\sqrt{|\lambda_i|},\varepsilon]})
     :\Neck[k]_{i\,,\,\vec{0}}-N_{\sqrt{|\lambda_i|}}(D_i)
      \rightarrow \Trunk[k]_{i\,,\,\vec{\lambda}}\,,          \\[1.6ex]
   \pr_i^{\ast}(\theta
                  _{\lambda_i\,;\,[|\lambda_i|/\varepsilon, \varepsilon]}
                \cup
                \theta^{\prime}
                  _{\lambda_i\,,\,[|\lambda_i|/\varepsilon, \varepsilon]})
     :\Neck[k]_{i\,,\,\vec{0}}-N_{|\lambda_i|/\varepsilon}(D_i)
      \rightarrow \Trunk[k]_{i\,,\,\vec{\lambda}}\,,
     & i\notin \{i_0,\,\ldots,\,i_{k^{\prime}}\}\,;           \\[1.6ex]
   \pr_j^{\ast} \theta_{\lambda_j\,;\,[\varepsilon,1/\varepsilon]}
    =\pr_{j-1}^{\ast}
       \theta^{\prime}_{\lambda_{j-1}\,;\,[\varepsilon,1/\varepsilon]}
    : \Trunk[k]_{j\,,\,\vec{0}}\; \stackrel{\sim}{\rightarrow}
         \Trunk[k]_{j\,,\,\vec{\lambda}}\,,
     & j\notin \{i_0,\,\ldots,\,i_{k^{\prime}}\}\,.
  \end{array}
 $$
 } % end-small

\vspace{-1ex}
\noindent
Here
 $N_{\,\cdot\,}(D_i)$ is the (open) tubular neighborhood of $D_i$
 in $Y_{[k]}=W[k]_{\vec{0}}$ of the specified radious from
 the norm on ${\Bbb L}$ and ${\Bbb L}^{\ast}$.
The collection of these morphisms glue/descend to
 two almost-complex morphisms
 $$
  \begin{array}{lllll}
   I_{\vec{\lambda}}     & :
    &  Y_{[k]}-\cup_{\,i=0}^{\,k}\,N_{\sqrt{|\lambda_i|}}(D_i)
    & \longrightarrow    & W[k]_{\vec{\lambda}}   \;,\\[1ex]
   I_{\vec{\lambda},\varepsilon}      & :
    &  Y_{[k]}-\cup_{\,i=0}^{\,k}\,N_{|\lambda_i|/\varepsilon}(D_i)
    & \longrightarrow  & W[k]_{\vec{\lambda}}     \;,
  \end{array}
 $$
 both of which shall be called a {\it re-forging morphism} from
 $W[k]_{\vec{0}}$ to $W[k]_{\vec{\lambda}}$.
Note that
 $I_{\vec{\lambda}}$ glues along the paired boundary of
  the connected components of
  $Y_{[k]}-\cup_{\,i=0}^{\,k}\,N_{\sqrt{|\lambda_i|}}(D_i)$
while
 $I_{\vec{\lambda},\varepsilon}$
  glues along the paired boundary of
   the connected components of
   $Y_{[k]}-\cup_{\,i=0}^{\,k}\,N_{|\lambda_i|/\varepsilon}(D_i)$
  but along a collar of non-paired boundary associated to
   $i\notin\{i_0\,,\,\ldots\,,\,i_{k^{\prime}}\}$.

\bigskip

\noindent
{\it Remark 1.1.1.6 $[$trunk region$]$.}
 The discussion implies that
  $\Trunk[k]_j\simeq B[k]\times \Trunk[k]_{j\,,\,\vec{0}}$
  canonically for $j=1\,,\,\ldots\,,\,k$.

\bigskip

\noindent
{\it Remark 1.1.1.7 $[$gluing map$]$.}
With Notation 1.1.1.2,
  % Notation [$M_{\lambda}$, $\lambda\ne 0$, from gluing]
 % Express $N_{\delta}(D_i)=N_{\delta}(D_i)_1\cup_{D_i}N_{\delta}(D_i)_2$
 % with $N_{\delta}(D_i)_1\subset \Delta_i$ and
 %      $N_{\delta}(D_i)_2\subset \Delta_{i+1}$.
 $$
  I_{\vec{\lambda},\varepsilon} \circ \pr_i^{\ast}
    \left(\rule{0ex}{2ex}\right.
      \varphi_{\lambda_i}\::\:
      {\Bbb L}^{\ast}
        _{[|\lambda_i|/\varepsilon\,,\,\varepsilon]}\;
      \longrightarrow\;
      {\Bbb L}_{[|\lambda_i|/\varepsilon\,,\,\varepsilon]}
     \left.\rule{0ex}{2ex}\right)\;
   =\; \Id_{\scriptsizeNeck_{\varepsilon}[k]_{i\,,\,\vec{\lambda}}}\,,
 $$
 where both ${\Bbb L}$ and ${\Bbb L}^{\ast}$ are regarded as
  canonically embedded in ${\Bbb L}\oplus{\Bbb L}^{\ast}$.

\bigskip

\noindent
{\it Remark 1.1.1.8
     $[$neck-trunk decomposition of $W[k]_{\vec{\lambda}}$$]$.}
Let $0\le i_0 <\,\cdots\, < i_{k^{\prime}} \le k$ be the associated
 indices to a $\vec{\lambda}\in B[k]$ so that $\lambda_{i_j}=0$.
Then,

 \vspace{-1.6ex}
 {\footnotesize
 $$
 \begin{array}{lll}
  &&
   \left(
    \overline{U_1}\,\cup\,
       \Neck[k]_{0\,,\,\vec{\lambda}}  \,\cup\,
       \Trunk[k]_{1\,,\,\vec{\lambda}}  \,\cup\,
       \Neck[k]_{1\,,\,\vec{\lambda}}
         \,\cup\,\cdots\, \cup\,
       \Trunk[k]_{i_0\,,\,\vec{\lambda}}
   \right)   \\[1.6ex]
  \bigcup\; \Neck[k]_{i_0\,,\,\vec{\lambda}}
   & \bigcup
   & \left(
       \Trunk[k]_{i_0+1\,,\,\vec{\lambda}}  \,\cup\,
       \Neck[k]_{i_0+1\,,\,\vec{\lambda}}
         \,\cup\,\cdots\, \cup\,
       \Trunk[k]_{i_1\,,\,\vec{\lambda}}
     \right)   \\[1.6ex]
  \bigcup\; \Neck[k]_{i_1\,,\,\vec{\lambda}}\;
   & \bigcup
   & \cdots\hspace{1em}
     \bigcup\hspace{1em}
     \left(
       \Trunk[k]_{i_{k^{\prime}-1}+1\,,\,\vec{\lambda}}  \,\cup\,
       \Neck[k]_{i_{k^{\prime}-1}+1\,,\,\vec{\lambda}}
         \,\cup\,\cdots\, \cup\,
       \Trunk[k]_{i_{k^{\prime}}\,,\,\vec{\lambda}}
      \right)   \\[1.6ex]
  \bigcup\; \Neck[k]_{i_{k^{\prime}}\,,\,\vec{\lambda}}
   & \bigcup
   & \left(
      \Trunk[k]_{i_{k^{\prime}}+1\,,\,\vec{\lambda}}  \,\cup\,
      \Neck[k]_{i_{k^{\prime}}+1\,,\,\vec{\lambda}}
        \,\cup\,\cdots\, \cup\,
      \Trunk[k]_{k\,,\,\vec{\lambda}} \,\cup\,
      \Neck[k]_{k\,,\,\vec{\lambda}} \,\cup\,
      \overline{U_2}
      \right)
 \end{array}
 $$
 } % end-footnotesize

 \vspace{-1ex}
 \noindent
 defines a {\it neck-trunk decomposition} of
 $W[k]_{\vec{\lambda}}\simeq Y_{[k^{\prime}]}$.

\bigskip

\subsubsection{The pseudo-${\Bbb G}_m[k]$-action on $W[k]/B[k]$
    in almost-complex category.}

Let ${\Bbb G}_m[k]
      := {\Bbb C}^{\times}\times\,\cdots\,\times{\Bbb C}^{\times}$
  ($k$ times) with coordinates $(\sigma_1,\,\ldots,\,\sigma_k)$.
It pseudo-acts\footnote{For non-algebraic-geometers:
  here ${\footnotesizeBbb G}_m$ means the {\it multiplicative group}
   of the ground field
   (e.g., ${\footnotesizeBbb C}^{\times}$ in ${\footnotesizeBbb C}$
          in our case) and
   is a standard notation from algebraic geometry.
  Also, given a group $G$ with the identity $e$,
   a {\it pseudo-group action} of $G$ on a space $M$ is a map from a
   neighborhood of $e\times M$ in $G\times M$ to $M$ that satisfies all
   the group-action axioms whenever items in the axioms are defined.}
 on $B[k]$ by
 $$
  (\lambda_0,\,\ldots,\,\lambda_i,\,\ldots,\, \lambda_k)
   \longmapsto\;
  (\sigma_0\sigma_1^{-1}\lambda_0,\,\ldots,\,
   \sigma_i\sigma_{i+1}^{-1}\lambda_i,\,\ldots,\,
   \sigma_k\sigma_{k+1}^{-1}\lambda_k)\,,
 $$
 where $\sigma_0=\sigma_{k+1} = 1$ by convention.
It admits a lifting to a pseudo-action on $W[k]/B[k]$ as follows.

Consider first the lifting of this pseudo-action to
 $({\Bbb L}\oplus{\Bbb L}^{\ast})_i$, $i=0,\,\ldots,\,k$, over $B[k]$
 by
 \begin{eqnarray*}
  \lefteqn{
   (\lambda_0,\,\ldots,\,\lambda_i,\,\ldots,\, \lambda_k;
     x,w_i,w_i^{\prime}) }                               \\[.6ex]
   && \longmapsto\;
      (\sigma_0\sigma_1^{-1}\lambda_0,\,\ldots,\,
       \sigma_i\sigma_{i+1}^{-1}\lambda_i,\,\ldots,\,
       \sigma_k\sigma_{k+1}^{-1}\lambda_k;
       x,\,\sigma_iw_i,\,\sigma_{i+1}^{-1}w_i^{\prime})\,.
 \end{eqnarray*}
This is well-defined since
 $(\sigma_i w_i)(\sigma_{i+1}^{-1} w_i^{\prime})
                         = \sigma_i\sigma_{i+1}^{-1}\lambda_i$.
This pseudo-action leaves both $({\Bbb L}\oplus{\Bbb L}^{\ast})_i^0$
 and $({\Bbb L}\oplus{\Bbb L}^{\ast})_i^{\infty}$ invariant, and
it follows from the explicit expression in Sec.~1.1.1
 that the gluing map
 $\varphi_{i-1,i}:({\Bbb L}\oplus{\Bbb L}^{\ast})_{i-1}^{\infty}
    \rightarrow ({\Bbb L}\oplus{\Bbb L}^{\ast})_i^0$
 is ${\Bbb G}_m[k]$-equivariant, for $i=1,\,\ldots,\,k$.
Consequently, the pseudo-${\Bbb G}_m[k]$-actions on
 $({\Bbb L}\oplus{\Bbb L}^{\ast})_i$, $i=0,\,\ldots,\,k$, glue
 to a pseudo-${\Bbb G}_m[k]$-action on
 $({\Bbb L}\oplus{\Bbb L}^{\ast})[k]$
 that lifts the pseudo-${\Bbb G}_m[k]$-action on $B[k]$.
This pseudo-action embeds ${\Bbb G}_m[k]$ into
 $\Aut(({\Bbb L}\oplus{\Bbb L}^{\ast})[k])$ in the almost
 complex category;
the isotropy group of $\vec{\lambda}\in B[k]$ under this
 pseudo-action coincides with
 $\Aut(\pi[k]^{-1}(\vec{\lambda})/{\Bbb L}\vee{\Bbb L}^{\ast})$.

By construction, the pseudo-${\Bbb G}_m[k]$-action on
 $({\Bbb L}\oplus{\Bbb L}^{\ast})[k]/B[k]$ descends to the trivial
 action on $({\Bbb L}\oplus{\Bbb L}^{\ast})/{\Bbb C}$ under
 $(\tilde{\mathbf p}[k], {\mathbf p}[k])$.
It follows that ${\Bbb G}_m[k]$ leaves
 $({\Bbb L}\oplus{\Bbb L}^{\ast})[k]_{\le \varepsilon}$ invariant and
its restriction to the horizontal boundary
 $\partial_{/B[k]}({\Bbb L}\oplus{\Bbb L}^{\ast})[k]_{\le \varepsilon}
   =B[k]\times (\partial\overline{U_1}\amalg\partial\overline{U_2})$
 of $({\Bbb L}\oplus{\Bbb L}^{\ast})[k]_{\le \varepsilon}$
 over $B[k]$ acts purely on the $B[k]$-factor.
This together with the gluing form
 $W[k]/B[k]
  = ( \overline{U_1}[k]  \cup
      ({\Bbb L}\oplus{\Bbb L}^{\ast})[k]_{\le \varepsilon}  \cup
      \overline{U_2}[k] )/B[k]$
 of $W[k]/B[k]$
implies that the pseudo-${\Bbb G}_m[k]$-action on
 $({\Bbb L}\oplus{\Bbb L}^{\ast})[k]_{\le \varepsilon}$
 extends to a pseudo-${\Bbb G}_m[k]$-action
 on $W[k]/B[k]$ such that its restriction on
 $\overline{U_1}[k]=B[k]\times\overline{U_1}$ and
 $\overline{U_2}[k]=B[k]\times\overline{U_2}$
 acts only on the $B[k]$-factor.

The following lemma follows immediately from the gluing construction
 of $W[k]/B[k]$ in Sec.~1.1.1.

\bigskip

\noindent
{\bf Lemma 1.1.2.1 [${\Bbb T}^k$-action on $W[k]/B[k]$].} {\it
 The restriction of the pseudo-${\Bbb G}_m[k]$-action on $W[k]/B[k]$
  to its maximal compact subgroup ${\Bbb T}^k:=U(1)^k$ gives
  an honest ${\Bbb T}^k$-action on $W[k]/B[k]$.
 This ${\Bbb T}^k$-action leaves the neck-trunk decomposition of
  $W[k]/B[k]$ invariant;
 the two re-forging morphisms
   $I_{\vec{\lambda}}$ and $I_{\vec{\lambda},\varepsilon}$
  are equivariant with respect to the stabilizer of the fiber
  $W[k]_{\vec{\lambda}}$ under the ${\Bbb T}^k$-action on $W[k]$.
} % end-lemma

\bigskip

\subsubsection{The topological quotient space $\widehat{W}/\widehat{B}$
               associated to $W/B$.}

We now construct a topological space $\widehat{W}/\widehat{B}$
 with charts that accommodates all the fibers
  $\{W_{\lambda}\}_{\lambda\in B}
    \cup \{Y_{[k]}\}_{k\in{\scriptsizeBbb Z}_{>0}}$
  that occur in an expanded degeneration of $W/B$.
For notation, given fibered spaces $W^{\prime}$ over $B^{\prime}$
 and $W^{\prime\prime}$ over $B^{\prime\prime}$,
a map $\varphi:W^{\prime}/B^{\prime}
           \rightarrow W^{\prime\prime}/B^{\prime\prime}$
 means a map $\varphi:W^{\prime}\rightarrow W^{\prime\prime}$
 that is descendable to a map
 $\underline{\varphi}:B^{\prime}\rightarrow B^{\prime\prime}$
 on the base.
Similarly, for a {\it pseudo-map}\footnote{Here a
        {\it pseudo-map} $f:A_1\rightarrow A_2$ means a map $f$ from a
         subset of $A_1$ to $A_2$. Similarly, a {\it pseudo-embedding}
         $f:A_1\rightarrow A_2$ means a pseudo-map  $f:A_1\rightarrow A_2$
         that is an embedding on where $f$ is defined.
        For non-algebraic-geometers: the reason for introducing such notion
         here is as follows.
        In the full construction of a moduli stack via
         the $\Isom$-functor, for two families of geometric objects
         in question
          (e.g.\ all the almost-complex isomorphism classes
           of fibers that occur in expanded degenerations of $W/B$)
         $\pi_1:W_1/B_1$ and $\pi_2:W_2\rightarrow B_2$, one
         constucts/defines a universal "overlapping" family
         $\pi:W\rightarrow \footnotesizeIsombf(\pi_1,\pi_2)$.
        Encoded into the construction of the family $\pi$ are
         natural morphisms
          $p_1:\footnotesizeIsombf(\pi_1,\pi_2)\rightarrow B_1$ and
          $p_2:\footnotesizeIsombf(\pi_1,\pi_2)\rightarrow B_2$, and
         tautological isomorphisms
          $p_1^{\ast}W_1\simeq W \simeq p_2^{\ast} W_2$ over
          $\footnotesizeIsombf(\pi_1,\pi_2)$.
        In Grothendieck's picture, illuminated by Mumford, each of
         $\pi_1$ and $\pi_2$ gives a local chart of the ``moduli space"
         behind, and the data from the $\footnotesizeIsom$-construction
         gives the Grothendieck's generalized notion of
         ``gluing" local charts $B_1$ and $B_2$ of the ``moduli space".
        As we mean to avoid the distraction of such formality,
        in our case it happens that one may relate $B_1$ and $B_2$
         instead by directly choosing (non-canonically and non-uniquely)
           a section to $p_1$,
          which is only defined on $\footnotesizeImage p_1\subset B_1$,
          and then post-compose it with $p_2$. This gives then
          a substitute ``transition" map
          $\underline{\varphi}:\footnotesizeImage p_1\rightarrow B_2$.
        Furthermore, as long as the
         ``quotient topology on the moduli space" is concerned,
         all that matters is that
          the domain $\footnotesizeImage p_1$ of $\underline{\varphi}$
           contains an {\it open} neighborhood of the point in $B_1$
           over which the central fiber in question sits;
        the precise tracking of $\footnotesizeImage p_1$ is irrelevant.
        Thus we directly re-denote $\underline{\varphi}$ as
         a {\it pseudo-map} $\underline{\varphi}:B_1\rightarrow B_2$.
        Via the canonical isomorphism
         $p_1^{\ast}W_1\simeq W \simeq p_2^{\ast} W_2$, accompanying
         the construction of $\underline{\phi}$ is also the pseudo-map
         $\varphi:W_1/B_1\rightarrow W_2/B_2$ that covers
         $\underline{\varphi}$.
        See [L-MB] for details on stacks and
            [L-L-Y: Sec.~1] for a literature guide.
      Similar use of ``{\it pseudo-}" applies to terms:
       {\it pseudo-embedding}, {\it pseudo-isomorphisms}, ..., etc.,
       and their compositions.}
 $W^{\prime}/B^{\prime}\rightarrow W^{\prime\prime}/B^{\prime\prime}$.

Recall the base $B[k]$ with coordinates
 $(\lambda_0,\,\ldots,\,\lambda_k)$ from the product ${\Bbb C}^{k+1}$.
To make the discussion more specific/concrete,
for a subset $I=\{\,i_0,\,\ldots,\,i_{k^{\prime}}\,\}$ of
 $\{\,0,\,\ldots, k\,\}$
 let $B[k]_I^{\varepsilon^2/4}$ be the affine coordinate subspace
  of $B[k]$,
   whose points have coordinates $\lambda_i=\varepsilon^2/4$
   for $i\notin I$ and
 denote $\pi[k]^{-1}(B[k]_I^{\varepsilon^2/4})$ by
  $W[k]_{B[k]_I^{\varepsilon^2/4}}$.
Then one has a pseudo-embedding of almost-complex spaces
 via the composition
 $$
  \varphi_{k^{\prime},k; I}\;:\; W[k^{\prime}]/B[k^{\prime}]\;
   \stackrel{\sim}{\longrightarrow}\;
  W[k]_{B[k]_I^{\varepsilon^2/4}}/B[k]_I^{\varepsilon^2/4}\;
  \hookrightarrow W[k]/B[k]\,,
 $$
 where
  $W[k^{\prime}]/B[k^{\prime}]
    \stackrel{\sim}{\longrightarrow}
     W[k]_{B[k]_I^{\varepsilon^2/4}}/B[k]_I^{\varepsilon^2/4}$
  is the almost-complex pseudo-isomorphism that lifts
  the pseudo-isomorphism
  $B[k^{\prime}]\rightarrow B[k]_I^{\varepsilon^2/4}$
  defined by
  $(\lambda^{\prime}_0\,,\,\ldots\,,\,\lambda^{\prime}_{k^{\prime}})
   \mapsto (\lambda_0\,,\,\ldots\,,\,\lambda_k)$ with
   $\lambda_i
     =(\frac{4}{\varepsilon^2})^{k-k^{\prime}}\lambda^{\prime}_j$,
    for $i=i_j\in I$, and $=\varepsilon^2/4$, for $i\notin I$.
The defining domain of $\varphi_{k^{\prime},k; I}$ contains
 an open neighborhood of the central fiber $\simeq Y_{[k^{\prime}]}$
 of $W[k^{\prime}]/B[k^{\prime}]$.
$\varphi_{k^{\prime},k; I}$ is equivariant with respect to
 ${\Bbb G}_m[k^{\prime}]\hookrightarrow {\Bbb G}_m[k]$ with
 $$
  (\sigma^{\prime}_1,\,\cdots,\, \sigma^{\prime}_{k^{\prime}})\;
  \longmapsto\;
  ( \underbrace{\sigma^{\prime}_0,\,\cdots,\sigma^{\prime}_0}
      _{i_0}\,,\,
    \underbrace{\sigma^{\prime}_1,\,\cdots,\,\sigma^{\prime}_1}
      _{i_1-i_0}\,,\,
    \cdots\,,
    \underbrace{\sigma^{\prime}_{k^{\prime}},\,
                  \cdots,\,\sigma^{\prime}_{k^{\prime}}}
      _{i_{k^{\prime}}-i_{k^{\prime}-1}}\,,\,
    \underbrace{\sigma^{\prime}_{k^{\prime}+1},\,
                  \cdots,\,\sigma^{\prime}_{k^{\prime}+1}}
      _{k-i_{k^{\prime}}} )\,,
 $$
 where
  $\sigma^{\prime}_0=\sigma^{\prime}_{k^{\prime}+1}=1$
   by convention and
  the multiplicity of each repeated entry is indicated.
Let $W_{(k)}/B_{(k)}$ be the quotient space of
 $W[k]/B[k]$ by ${\Bbb G}_m[k]$ with the quotient topology.
Then
 $(W-W_0)/(B-\{0\})$ embeds in $W_{(k)}/B_{(k)}$ canonically for all
  $k\in{\Bbb Z}_{\ge 0}$ and
 $\varphi_{k^{\prime},k;I}$ induces an embedding
  $$
   \varphi_{(k^{\prime},k;I)}\;:\;
    W_{(k^{\prime})}/B_{(k^{\prime})}\;
    \hookrightarrow\; W_{(k)}/B_{(k)}
  $$
  over $W/B$, for all $k^{\prime}<k$, that restricts to
  the identity map on $(W-W_0)/(B-\{0\})$.

Let $\widehat{B}=B\cup {\Bbb Z}_{>0}$ with the topology
 generated by the open subsets of $B$ and the subsets of
 $\widehat{B}$ of the form $U\cup\{1,\,\ldots,\,k\}$, where $U$
 is an open neighborhood of $0\in B$ and $k\in {\Bbb Z}_{>0}$.
Define the set
 $$
  \widehat{W}/\widehat{B}\;
   :=\; \left.\left(\amalg_{k\in{\scriptsizeBbb Z}_{\ge 0}}\,
          W_{(k)}/B_{(k)}\right)\right/\raisebox{-.6ex}{$\sim$}\,,
 $$
 where
  $p\in W_{(k)}$ and $p^{\prime}\in W_{(k^{\prime})}$
    with $k>k^{\prime}$ are defined to be equivalent
    (in notation, $p\sim p^{\prime}$)
   if $p$ is the image of $p^{\prime}$ under some
    $\varphi_{(k^{\prime},k;I)}$
  (this defines $\widehat{W}$)
    and
  $p\in B[k]$ and $p^{\prime}\in B[k^{\prime}]$ are equivalent
   if $p$ is the image of $p^{\prime}$ under some
   $\underline{\varphi_{(k^{\prime},k;I)}}$
  (this reproduces $\widehat{B}$).
As indicated, the fibrations $W_{(k)}/B_{(k)}$, $k\in {\Bbb Z}_{\ge 0}$,
 induce a fibration of $\widehat{W}$ over $\widehat{B}$.
By construction, there are natural embeddings (of sets)
 $$
  \varphi_{(k)}\;:\;
   W_{(k)}/B_{(k)}\;\hookrightarrow\; \widehat{W}/\widehat{B}\,,
    \hspace{1em} k\;\in\; {\Bbb Z}_{\ge 0}\,.
 $$
Equip $\widehat{W}$ with the topology that specifies
 a subset $\widehat{U}$ of $\widehat{W}$ to be open
 if and only if $\widehat{U}=\cup_{\alpha}\widehat{U}_{\alpha}$
 such that for each $\alpha$
  there exists $k_{\alpha}\in {\Bbb Z}_{\ge 0}$ so that
  $\widehat{U}_{\alpha}=\varphi_{(k_{\alpha})}(U_{\alpha})$
  for some open subset $U_{\alpha}$ of $W_{(k_{\alpha})}$.
We will call this topology the {\it quotient topology}
 on $\widehat{W}$.
Note that this topology involves all $\varphi_{(k^{\prime},k;I)}$
 so that the information of how one $Y_{[k^{\prime}]}$ or $W_{\lambda}$
 degenerates to another $Y_{[k]}$ with $k>k^{\prime}$ is all kept.
By construction,
 both the natural map $\widehat{W}\rightarrow \widehat{B}$ and
 the defining maps
 $$
  \varphi[k]\; :\; W[k]/B[k]\;\longrightarrow \;\widehat{W}/\widehat{B}
 $$
 from the composition
  $W[k]/B[k]\rightarrow W_{(k)}/B_{(k)}
                   \rightarrow \widehat{W}/\widehat{B}$
 are continuous.
$(W[k]/B[k],\varphi[k])$ is named a {\it standard local chart}
 on $\widehat{W}/\widehat{B}$ and
the collection $\{(W[k]/B[k],\varphi[k]): k\in {\Bbb Z}_{\ge 0}\}$
 the {\it standard atlas} for $\widehat{W}/\widehat{B}$.

Finally, note that the collection of maps
 $\{ \tilde{\mathbf p}[k]:W[k]/B[k]\rightarrow W/B \}
                             _{k\in{\scriptsizeBbb Z}_{\ge 0}}$
 descends to a (continuous) tautological map
 $$
  \widehat{\mathbf p}\;:\;
    \widehat{W}/\widehat{B}\; \longrightarrow\; W/B\,.
 $$

\bigskip

\noindent
{\it Remark 1.1.3.1 $[$quotient topology versus stack$\,]$.}
 To identify consistently isomorphic fibers (as almost-complex spaces)
  in the collection $\{W[k]/B[k]\}_{k\in{\scriptsizeBbb Z}_{\ge 0}}$
  and make the final family universal,
 one has to employ Grothendieck's generalized notion in algebraic geometry
  of ``{\it gluing}" via the Isom-functor construction,
  of a ``{\it space}" as a collection of local charts together
      with a gluing data in the generalized sense, and
  of a ``{\it global structure}" as a descent datum.
 Following this, the collection
  $\{W[k]/B[k]\}_{k\in {\scriptsizeBbb Z}_{\ge 0}}$
  would be glued to an Artin stack ${\cal B}$,
  together with a universal expanded degeneration ${\cal W}$
  over ${\cal B}$.
 The set of geometric points of ${\cal B}$ would be $B\cup{\Bbb Z}_{>0}$
  with the corresponding set of isomorphism class
   of fibers of ${\cal W}/{\cal B}$ being
   $\{W_{\lambda}\}_{\lambda\in B}
    \cup \{Y_{[k]}\}_{k\in{\scriptsizeBbb Z}_{>0}}$.
 (Cf.\ [Li1: Sec.~1]; see [L-L-Y: Sec.~1] for a brief tour on stacks).
 Since it is the stable maps,
  i.e.\ triples $(\Sigma, W[k]_{\lambda},
                 f:\Sigma\rightarrow W[k]_{\vec{\lambda}})$,
  that we want to study in this work,
 it turns out that what we finally need most essentially
  is a structure that describes the ``nearness" between
  a $W_{\lambda}$ or $Y_{[k]}$ and
  another $W_{\lambda^{\prime}}$ or $Y_{[k^{\prime}]}$.
 For this reason, the space $\widehat{W}/\widehat{B}$
  with the quotient topology and the standard atlas
  as constructed above that accommodates all
  $\{W_{\lambda}\}_{\lambda\in B}
   \cup \{Y_{[k]}\}_{k\in{\scriptsizeBbb Z}_{>0}}$
  suffices.

\bigskip

\subsection{Symplectic/almost-complex relative pairs and their expansions.}

A {\it symplectic} (resp.\ {\it almost-complex}) {\it relative pair}
 $(Z;D)$ is a symplectic (resp.\ almost-complex) manifold $Z$
 together with a real codimension-$2$ symplectic (resp.\ almost-complex)
 submanifold $D$.
Given a symplectic relative pair $(Z;D)$ with a Hamiltonian
 $U(1)$-action on a (open) tubular neighborhood $N(D)$ of $D$ in $Z$
 that fixes $D$,
define $Z[1]$ to be the total space of a compatible almost-complex
 degeneration of a symplectic cut on $Z$ associated to the given
 local $U(1)$-action around $D$.
By construction, $Z[1]$ fibers over
 $A[1]:=B=\{\lambda\in{\Bbb C}\,:\,|\lambda|<\varepsilon^2/2\}$,
 $0<\varepsilon<1$,
 with the singular fiber $Z[1]_0=Z\cup_{D=D_{1,\infty}}\Delta_1$.
Since the pinched locus of the symplectic cut is disjoint from $D$ and
 it separates $D$ with $Z-N_{\varepsilon}(D)$,
 $D[1]:= A[1]\times D$ embeds canonically in $Z[1]$ over $A[1]$
 with $D[1]_0:= \{0\}\times D$ identical to $D_{1,0}$ in $\Delta_1$.

The construction in Sec.~1.1.1 applied to
 the almost-complex degeneration $Z[1]/A[1]$ then gives rise to
 an {\it almost-complex expanded relative pair} $(Z[k];D[k])/A[k]$
 with $A[k]=B[k-1]$ and $D[k]=A[k]\times D$, for $k\in{\Bbb Z}_{\ge 0}$.
Its fiber, e.g., at $\vec{0}\in A[k]$ is the almost-complex
 relative pair
 $$
  \begin{array}{cll}
   (Z[k];D[k])_{\vec{0}}
    & =   & (Z\cup_{D=D_{1,\infty}}\Delta_1
            \cup_{D_{1,0} = D_{2,\infty}}\,\cdots\,
            \cup_{D_{k-1,0}=D_{k,\infty}}\Delta_k\,;\, D_{k,0}) \\[.6ex]
    & =:  & (Z_{[k]};D_{[k]})\,.
  \end{array}
 $$
There is also the almost-complex morphism
 $$
  \tilde{\mathbf p}[k]\;:\; (Z[k];D[k])/A[k]\;\longrightarrow\; (Z;D)/\pt
 $$
 from the construction.

Let $\overline{U}=Z-N_{\varepsilon}(D)$,
 where $N_{\varepsilon}(D)$ is the open $\varepsilon$-neighborhood
       of $D$ in $Z$ with respect to the norm on ${\Bbb L}$.
Then $(Z[k];D[k])/A[k]$ admits a {\it neck-trunk decomposition}:

 \vspace{-1.6ex}
 {\small
 \begin{eqnarray*}
  \lefteqn{Z[k]/A[k]\; =} \\[.6ex]
  &&
   (\, \overline{U}[k]  \,\cup\,
       \Neck[k]_0  \,\cup\,  \Trunk[k]_1  \,\cup\,  \Neck[k]_1
         \,\cup\,\cdots\, \cup\,
       \Trunk[k]_k  \,\cup\,  N_{\varepsilon}(D)[k])/A[k]\,,
 \end{eqnarray*}
 } % end-small

\vspace{-1ex}
\noindent
where
 $\overline{U}[k]$, $\Neck[k]_i$, $i=0\,,\,\ldots\,,\,k-1$,
  $\Trunk[k]_j$, $j=1\,,\,\ldots\,,\,k$,
  here are similar to their counterpart:
   $\overline{U_1}[k-1]$, $\Neck[k-1]_i$, and
   $\Trunk[k-1]_j$, in Sec.~1.1.1  and
 $N_{\varepsilon}(D)[k]=A[k]\times N_{\varepsilon}(D)$,
  which contains $D[k]$.
This induces a neck-trunk decomposition to the fiber
 $(Z[k];D[k])_{\vec{\lambda}}$ of $(Z[k];D[k])$
 at $\vec{\lambda}\in A[k]$,
 cf.\ Remark 1.1.1.8.
      % Remark [neck-trunk decomposition of $W[k]_{\vec{\lambda}}$]
There are {\it re-forging morphisms} from $Z[k]_{\vec{0}}=Z_{[k]}$
 to $Z[k]_{\vec{\lambda}}$ constructed in the same way as earlier:
 $$
  \begin{array}{llllll}
   I_{\vec{\lambda}}     & :
    &  Z_{[k]}-\cup_{\,i=0}^{\,k-1}\,N_{\sqrt{|\lambda_i|}}(D_i)
    & \longrightarrow    & Z[k]_{\vec{\lambda}}   \;,\\[1ex]
   I_{\vec{\lambda},\varepsilon}      & :
    &  Z_{[k]}-\cup_{\,i=0}^{\,k-1}\,N_{|\lambda_i|/\varepsilon}(D_i)
    & \longrightarrow  & Z[k]_{\vec{\lambda}}     \;,
    & \hspace{1em}\vec{\lambda}\;\in\; A[k]\,.
  \end{array}
 $$

The group ${\Bbb G}_m[k]$ now pseudo-acts on $A[k]$ by
 $$
  (\lambda_0,\,\ldots,\,\lambda_i,\,\ldots,\, \lambda_{k-1})
   \longmapsto\;
  (\sigma_0\sigma_1^{-1}\lambda_0,\,\ldots,\,
   \sigma_i\sigma_{i+1}^{-1}\lambda_i,\,\ldots,\,
   \sigma_{k-1}\sigma_k^{-1}\lambda_{k-1})\,,
 $$
 where $\sigma_0=1$ by convention.
Similar to Sec.~1.1.2,
 it lifts to a {\it pseudo-${\Bbb G}_m[k]$-action} on $Z[k]$
  that leaves $D[k]$ invariant in such a way that
   the pseudo-action on $D[k]=A[k]\times D$ acts only on
   the $A[k]$-factor.
As a parallel to Lemma 1.1.2.1,
  % Lemma [${\Bbb T}^k$-action on $W[k]/B[k]$]
 the restriction of the pseudo-${\Bbb G}_m[k]$-action on
  $(Z[k];D[k])/A[k]$ to its maximal compact subgroup
  ${\Bbb T}^k$ gives an honest ${\Bbb T}^k$-action on
  $(Z[k];D[k])/A[k]$.
This ${\Bbb T}^k$-action leaves the neck-trunk decomposition of
  $(Z[k];D[k])/A[k]$ invariant and
 the two re-forging morphisms
   $I_{\vec{\lambda}}$ and $I_{\vec{\lambda},\varepsilon}$
  are equivariant with respect to the stabilizer of
  $Z[k]_{\vec{\lambda}}$ under the ${\Bbb T}^k$-action on $Z[k]$.

To connect the various expanded relative pairs,
each
 $I^{\prime}=\{\,i_0,\,\ldots,\,i_{k^{\prime}-1}\,\}
  \subset \{\,0,\,\ldots, k-1\,\}$
is associated to a pseudo-embedding of almost-complex spaces
 $$
  \varphi^{\prime}_{k^{\prime},k; I^{\prime}}\;:\;
   (Z[k^{\prime}];D[k^{\prime}])/A[k^{\prime}]\;
  \hookrightarrow (Z[k];D[k])/A[k]\,,
 $$
 which covers the pseudo-embedding
  $A[k^{\prime}]\rightarrow A[k]$,
  defined by
  $(\lambda^{\prime}_0\,,\,\ldots\,,\,\lambda^{\prime}_{k^{\prime}-1})
   \mapsto (\lambda_0\,,\,\ldots\,,\,\lambda_{k-1})$
  with
   $\lambda_i
     =(\frac{4}{\varepsilon^2})^{k-k^{\prime}}\lambda^{\prime}_j$,
     for $i=i_j\in I^{\prime}$, and
    $=\varepsilon^2/4$, for $i\notin I^{\prime}$.
 and
 is equivariant with respect to the group homomorphism
  ${\Bbb G}_m[k^{\prime}]\hookrightarrow {\Bbb G}_m[k]$
  defined by
  $$
   (\sigma^{\prime}_1,\,\cdots,\, \sigma^{\prime}_{k^{\prime}})\;
   \longmapsto\;
   ( \underbrace{1,\,\cdots\,,1}
       _{i_0}\,,\,
     \underbrace{\sigma^{\prime}_1,\,\cdots,\,\sigma^{\prime}_1}
       _{i_1-i_0}\,,\,
    \cdots\,,
     \underbrace{\sigma^{\prime}_{k^{\prime}-1},\,
                   \cdots,\,\sigma^{\prime}_{k^{\prime}-1}}
       _{i_{k^{\prime}-1}-i_{k^{\prime}-2}}\,,\,
     \underbrace{\sigma^{\prime}_{k^{\prime}},\,
                   \cdots,\,\sigma^{\prime}_{k^{\prime}}}
       _{k-i_{k^{\prime}-1}} )\,.
  $$
Let $(Z_{(k)};D_{(k)})/A_{(k)}$ be the quotient space of
 $(Z[k];D[k])/A[k]$ by ${\Bbb G}_m[k]$ with the quotient topology.
Then
 $(Z;D)$ embeds in $(Z_{(k)};D_{(k)})/A_{(k)}$ canonically for all
  $k\in{\Bbb Z}_{\ge 0}$ and
 $\varphi^{\prime}_{k^{\prime},k;I^{\prime}}$ induces an embedding
  $$
   \varphi^{\prime}_{(k^{\prime},k;I)}\;:\;
    (Z_{(k^{\prime})};D_{(k^{\prime})})/A_{(k^{\prime})}\;
    \hookrightarrow\; (Z_{(k)};D_{(k)})/A_{(k)}\,,
  $$
  for all $k^{\prime}<k$, that restricts to
  the identity map on $(Z;D)$.

Let $\widehat{A}={\Bbb Z}_{\ge 0}$ with the topology generated
 by the defining open subsets
  $\{i\in {\Bbb Z}_{\ge 0}\,:\, 0\le i\le n\}$, $n\in {\Bbb Z}_{\ge 0}$.
Then, the construction in Sec.~1.1.3 applied to
 $\{(Z[k];D[k])/A[k]\}_{k\in{\scriptsizeBbb Z}_{\ge 0}}$,
  where $(Z[0];D[0])/A[0]=(Z;D)$ by convention,
 gives rise to a {\it topological relative pair}
 $(\widehat{Z};\widehat{D})$ over $\widehat{A}$
  with the {\it quotient topology},
 the natural embeddings
 $$
  \varphi_{(k)}\;:\;
   (Z_{(k)};D_{(k)})/A_{(k)}\;
    \hookrightarrow\; (\widehat{Z};\widehat{D})/\widehat{A}\,,
    \hspace{1em} k\;\in\; {\Bbb Z}_{\ge 0}\,,
 $$
 the {\it standard local charts}
 $$
  \varphi[k]\; :\; (Z[k];D[k])/A[k]\;
    \longrightarrow\; (\widehat{Z};\widehat{D})/\widehat{A}\,,
   \hspace{1em} k\;\in\; {\Bbb Z}_{\ge 0}\,,
 $$
 and a (continuous) tautological map
 $$
  \widehat{\mathbf p}\;:\;
   (\widehat{Z};\widehat{D})/\widehat{A}\; \longrightarrow\; (Z;D)\,.
 $$
The topological relative pair $(\widehat{Z};\widehat{D})/\widehat{A}$
 equipped with the standard local charts substitutes the stack
 of expanded relative pairs obtained by gluing $(Z[k];D[k])/A[k]$'s
 via the Isom-functor construction.

Readers are referred also to
 [I-P1: Sec.~3 and Sec.~6], [L-R: Sec.~3], and [Li1: Sec.~4]
 for related discussions.

\bigskip

\section{Prestable labelled-bordered Riemann surfaces.}

In this section we review/rephrase/modify definitions/facts of
 labelled-bordered Riemann surfaces with marked points to
 introduce and fix terminologies and notations that we will use.
This is a classical topic with long history.
Readers are referred to
 [Sie1: Sec.~2], [F-O: Sec.~9 and pp.~988~-~991], and
 [Liu(C): Sec.~2 - Sec.~4]
 for related discussions and guide to literatures.
See also
 [Ab], [A-G], [D-M], [H-M], [I-S2], [Kn], [Ma], [Se], [Sil], and [Wol].

\bigskip

\begin{flushleft}
{\bf Prestable labelled-bordered Riemann surfaces with marked points.}
\end{flushleft}
%
% \noindent
{\bf Definition 2.1 [prestable labelled-bordered Riemann surface].}
{\rm
 A {\it prestable labelled-bordered Riemann surface of
  $($combinatorial$)$ type $((g,h),(n,\vec{m}))$
  $($with labelled boundary and marked points$)$}\footnote{The
                definition here is based on [Liu(C): Definition 3.9].
                We phrase it to make it manifest that
                 an interior marked point on a nodal bordered Riemann
                  surface is allowed to move and land on the boundary
                  to become a {\it double} boundary point.
                This freedom is required to obtain a compact moduli space
                 of stable bordered Riemann surfaces with marked points.
                We avoid the term {\it marked bordered Riemann surface}
                 to reserve its more traditional meaning in
                 the Teichm\"{u}ller theory of Riemann surfaces.},
  where $\vec{m}=(m_1\,,\,\ldots\,,\,m_h)$,
 consists of the following data:
 \begin{itemize}
  \item[$\cdot$]
   a {\it compact connected nodal bordered Riemann surface} $\Sigma$,
    whose points are locally modelled at $0$ or $(0,0)$
     in the following holomorphic models:
    \begin{itemize}
     \item[(i)] {\it interior point}$\,$:
      \begin{itemize}
       \item[(i1)]
        $\{z\in {\Bbb C}\,:\, |z|<1\}$
         for a {\it smooth interior point},

       \item[(i2)]
        $\{(z_1,z_2)\in{\Bbb C}^2\,:\,
                  |z_1|<1\,,\, |z_2|<1\,,\,z_1z_2=0\}$
         for an {\it interior node};
      \end{itemize}

     \item[(b)] {\it boundary point}$\,$:
      \begin{itemize}
       \item[(b1)]
        $\{z\in {\Bbb C}\,:\,|z|<1\,,\im(z)\ge 0 \}$
         for a {\it smooth boundary point},

       \item[(b2)]
        $\{(z_1,z_2)\in{\Bbb C}^2\,:\,
                    |z_1|<1\,,\,|z_2|<1\,,\, z_1z_2=0\}
             /(z_1,z_2)\sim (\overline{z_2},\overline{z_1})$\\
         for a {\it boundary node of type E},

       \item[(b3)]
        $\{(z_1,z_2)\in{\Bbb C}^2\,:\,
                    |z_1|<1\,,\,|z_2|<1\,,\, z_1z_2=0\}
            /(z_1,z_2)\sim (\overline{z_1},\overline{z_2})$\\
         for a {\it boundary node of type H};
      \end{itemize}

     \item[]
     \hspace{-1.8em}
     the number of interior (resp.\ boundary) node will be denoted
       $n_{\intnode}$ (resp.\ $n_{\bn}$).
    \end{itemize}
 \end{itemize}
 \begin{itemize}
  \item[$\cdot$]
   {\it labelled boundary} and $h\,$:
   a {\it boundary component} of $\Sigma$ is either
    the image of an {\it embedding} of $S^1$ in $\partial\Sigma$  or
    a boundary node of type E;
   $\Sigma$ has $h$-many boundary components and
   they are labelled from $1$ to $h$;
   the labelled boundary of $\Sigma$ will be denoted
    by $\dot{\partial}\Sigma$ (or simply $\partial\Sigma$
    when the labelling is understood);
   note that different boundary components of $\Sigma$
    may intersect at a boundary node of type H.

  \item[$\cdot$]
   {\it genus $g$}$\,$:
   each boundary component of $\Sigma$ can be capped by a $2$-disc;
   let $\hat{\Sigma}$ be the nodal Riemann surface without boundary
    obtained by capping all the boundary components of $\Sigma$
    by discs,
   then $\hat{\Sigma}$ has {\it arithmetic genus} $g$.

  \item[$\cdot$]
   {\it free marked points}$\,$:
   an $n$-tuple $\vec{p}=(p_1\,,\,\cdots\,,\,p_n)$
    of {\it smooth interior points} or
       {\it double boundary points}\footnote{For non-algebraic-geometers:
                     in the affine ${\footnotesizeBbb R}$-scheme model
                      a {\it smooth interior point}
                        (resp.\ {\it smooth boundary point})
                       on $\Sigma$ is modelled on
                      a {\it complex closed point},
                       $(x^2-(c+\bar{c})x+c\bar{c})$,
                        $c\in{\footnotesizeBbb C}-{\footnotesizeBbb R}$,
                       (resp.\ {\it real closed point} $(x-a)$,
                                    $a\in{\footnotesizeBbb R}$)  in
                       $\footnotesizeSpec {\footnotesizeBbb R}[x]$.
                     A complex closed point in
                      $\footnotesizeSpec {\footnotesizeBbb R}[x]$
                      can be deformed to a {\it double real point},
                       described by an ideal $(x-a)^2$
                       for some $a\in{\footnotesizeBbb R}$, in
                      ${\footnotesizeBbb R}[x]$.
                     While a real double point as above can be deformed
                       to a complex closed point,
                      a closed real point can only be deformed
                       to another closed real point.
                     In other words, an interior marked point in $\Sigma$
                      can be deformed to a double boundary point
                      on $\partial\Sigma$ and vice versa.
                     Together, we name them {\it free} marked points on
                      $\Sigma$.
                     Thus, a free marked point,
                      whether in the interior or on boundary, always have
                      real $2$-dimensional family of deformations.
                     In particular, fixing a complex point always
                      contributes {\it two} real constraints whether
                      that point is in the interior or on the boundary.
                     In contrast, a boundary marked point on $\Sigma$
                      can move around only in the boundary $\partial\Sigma$
                      and contributes only one real condition.},
             on $\Sigma$; the support of the latter free points is
             required to be smooth boundary points.
           The notation $n\doteq n^{\prime}+n^{\prime\prime}$ means that
            there are $n^{\prime}$-many interior marked points and
            $n^{\prime\prime}$-many free marked points supported
             in $\partial\Sigma$, when the distinction is needed.

  \item[$\cdot$]
   {\it boundary marked points}$\,$:
   an $m_i$-tuple of smooth boundary points
    $\vec{p}_i=(p_{i1}\,,\,\cdots\,,\,p_{i\,m_i})$
    on the boundary component of $\Sigma$ labelled by $i$
    for $i=1\,,\,\ldots\,,\,h$;
   we require that the set of boundary marked points is disjoint
    from the support of free marked points that land on the boundary.
 \end{itemize}

 By definition, the set of nodes and the set of marked points
  on $\Sigma$ are disjoint from each other.
 Any point in the union of the two is called
  a {\it special point} on $\Sigma$.

 A {\it regular} or {\it smooth point} on $\Sigma$ is either
  a smooth interior point or a smooth boundary point on $\Sigma$.
 The set of regular points on $\Sigma$ with the induces topology
  and holomorphic/complex structure is denoted by $\Sigma_{\reg}$
  and called the {\it regular} or {\it smooth locus} of $\Sigma$.

 From the local model of points on $\Sigma$, one can define
  the {\it normalization} $\tilde{\Sigma}$ of $\Sigma$
  as in algebraic geometry.
 Topological, $\tilde{\Sigma}$ is obtained by  first
  removing all the (interior as well as boundary) nodes on $\Sigma$
   and then
  filling all the resulting (interior as well as boundary) punctures
   by distinct points.
 $\tilde{\Sigma}$ is a possibly disconnected bordered Riemann surface
 (with neither interior nor boundary nodes).
 Let
  $\nu: \tilde{\Sigma}\rightarrow \Sigma$ be the normalization
   of $\Sigma$ and
  $\tilde{\Sigma}=\amalg_{i}\tilde{\Sigma}_i$
   be the disjoint union of connected components;
 then each $\nu(\tilde{\Sigma}_i)$ in $\Sigma$ is called
  an {\it irreducible component} of $\Sigma$.

 Let $\overline{\Sigma}$ be the nodal bordered Riemann surface
  with the same topology as $\Sigma$ but
  with the complex-conjugated holomorphic structure
   from that of $\Sigma$.
 Then
  $\Sigma_{\scriptsizeBbb C}
   := \Sigma
      \cup_{\partial\Sigma=\partial\overline{\Sigma}}
      \overline{\Sigma}$
  has a canonically induced nodal Riemann surface structure
  without boundary.
 It is called the {\it Schottky/complex double} of $\Sigma$.
 By construction, there is an {\it involution} $\tau$ that acts on
  $\Sigma_{\scriptsizeBbb C}$ by complex conjugation.

 An {\it isomorphism}
  $h: (\Sigma, \dot{\partial}\Sigma,
       \vec{p}, \vec{p}_1,\,\cdots\,,\,\vec{p}_h)  \rightarrow
      (\Sigma^{\prime}, \dot{\partial}\Sigma^{\prime},  \vec{p}^{\prime},
      \vec{p}^{\prime}_1,\,\cdots\,,\,\vec{p}^{\prime}_h)$
  from a labelled-border-\\  ed Riemann surface to
  another of the same type is a bi-holomorphic map
   $h:(\Sigma,\partial\Sigma)
      \rightarrow (\Sigma^{\prime},\partial\Sigma^{\prime})$
  that
   preserves the label of the boundary components  and
   sends $p_i$ to $p_i^{\prime}$, $q_{ij}$ to $q_{ij}^{\prime}$.
 An {\it automorphism} of
  $(\Sigma, \dot{\partial}\Sigma,
     \vec{p}, \vec{p}_1,\,\cdots\,,\,\vec{p}_h)$
  is an isomorphism from
  $(\Sigma, \dot{\partial}\Sigma,
       \vec{p}, \vec{p}_1,\,\cdots\,,\,\vec{p}_h)$
  to itself.
 $(\Sigma, \dot{\partial}\Sigma,
   \vec{p}, \vec{p}_1,\,\cdots\,,\,\vec{p}_h)$  is called {\it stable}
  if its group
      $\Aut(\Sigma, \dot{\partial}\Sigma,
        \vec{p}, \vec{p}_1,\,\cdots\,,\,\vec{p}_h)$
     of automorphisms is finite.

 We will denote the data
  $(\Sigma, \dot{\partial}\Sigma,
     \vec{p}, \vec{p}_1,\,\cdots\,,\,\vec{p}_h)$
  also by $(\Sigma,\partial\Sigma)$ or $\Sigma$ in short.
 The isomorphism class of labelled-bordered Riemann surfaces isomorphic
  to $\Sigma$ will be denoted $[\Sigma]$.
 When there is no chance of confusion, we will call
  $\Sigma$ also a {\it curve} and denote it by $C$,
  as a $1$-dimensional scheme over $\Spec{\Bbb C}$ or $\Spec{\Bbb R}$
   in algebraic geometry
  with labelled irreducible components of ${\Bbb R}$-locus
   and marked points.
 The {\it moduli space} of isomorphism classes of {\it stable}
  (resp.\ {\it prestable}) labelled-bordered Riemann surfaces of
  type $((g,h),(n,\vec{m}))$ will be denoted
   $\overline{\cal M}_{(g,h), (n,\vec{m})}$
   (resp.\ $\widetilde{\cal M}_{(g,h), (n,\vec{m})})$.
} % end-definition

\bigskip

\noindent
{\bf Theorem 2.2 [$\overline{\cal M}_{(g,h), (n,\vec{m})}$].} {\it
 The moduli space $\overline{\cal M}_{(g,h), (n,\vec{m})}$
  of stable labelled-bordered Riemann surfaces with marked points
    of type $((g,h),(n,\vec{m}))$,
   with its topology defined via the dilatation of
    quasi-conformal maps and their composition with
    circle/arc-with-ends-in-boundary pinching maps  or
   via the local Fenchel-Nielsen coordinates
    associated to pants-decompositions,
  is a compact, Hausdorff, orientable orbifold-with-corners.
} % end-theorem

\bigskip

\noindent
See [Liu(C): Theorem 4.9, Theorem 4.14] and the quoted references there.

The universal deformation ${\cal C}/Def(\Sigma)$ of $\Sigma$,
  canonically acted upon by $\Aut(\Sigma)$,
 provides a local orbifold-chart
 $\psi_{[\Sigma]}:
  \Def(\Sigma)\rightarrow\overline{\cal M}_{(g,h), (n,\vec{m})}$
 around $[\Sigma]$ in $\overline{\cal M}_{(g,h), (n,\vec{m})}$.
Topologically this is a quotient of a neighborhood of the origin
 in the manifold-with-corners
 $$
  \begin{array}{c}
   \Ext^1_{\Sigma_{\tinyBbb C}}
    \left(
      \Omega_{\Sigma_{\tinyBbb C}}
      ( \sum_{i=1}^n(p_i+\overline{p_i})\,  +\,
        \sum_{j=1}^h\sum_{k=1}^{m_j}p_{jk} )\,,\,
      {\cal O}_{\Sigma_{\tinyBbb C}}
     \right)^{\tau}\;
      \\[1.6ex]
   \simeq\;
    {\Bbb C}^{3g-3+h+n^{\prime}}  \times
    \overline{\Bbb H}^{n^{\prime\prime}}  \times
    {\Bbb R}^{h-n_{bn}+m_1+\,\cdots\,+m_h}  \times
    ({\Bbb R}_{\ge 0})^{n_{bn}}
  \end{array}
 $$
 by $\Aut(\Sigma)$,
 where
  $\,\mbox{\tiny $\bullet$}\,^{\tau}$ is the fixed-point locus of the
   induced action of $\tau$ on $\,\mbox{\tiny $\bullet$}\,$,
  $\overline{\Bbb H} =$
   the closed upper half-plane $\{z\in {\Bbb C}: \im(z)\ge 0\}$,  and
  $n\doteq n^{\prime}+n^{\prime\prime}$.
As an orbifold, $\overline{\cal M}_{(g,h),(n,\vec{m})}$
 goes with a {\it universal family}, denotes also by
 ${\cal C}/\overline{\cal M}_{(g,h), (n,\vec{m})}$.
We will call this ${\cal C}$ the {\it universal curve}
 of type $((g,h),(n,\vec{m}))$.
$\overline{\cal M}_{(g,h), (n,\vec{m})}$ is {\it naturally stratified}
 by a finite collection of locally closed sub-orbifolds-with-corners.
The stratification is governed by
  the {\it topological type} (i.e.\ equivalence up to homeomorphisms
    of the underlying topology of punctured bordered Riemann surfaces)
    and
  the {\it degeneration patterns} of
 a labelled-bordered Riemann surface with marked points.
See, e.g., [Liu(C): Figures 1, 2, 3, 9, 10, 11] for illustrations of
 such stratifications.

\bigskip

\begin{flushleft}
{\bf Local chart on $\widetilde{\cal M}_{(g,h),(n,\vec{m})}$.}
\end{flushleft}
There are $11$ types of
 {\it unstable {\rm (}irreducible{\rm )} components}
 that can happen for a prestable labelled-bordered Riemann surface:
 (1) {\it closed} component$\,$:
     ($g=0$) ${\Bbb P}^1$ with $0$, $1$, or $2$ special points;
     ($g=1$) torus without special points or
             nodal torus with one node and without marked points;
 (2) {\it bordered} component$\,$: (all with $g=0$)
     $2$-disc $D^2$ with $0$ or $1$ free marked point; or
      $D^2$ with $1$ or $2$ boundary marked points;
     annulus without special point or
     nodal annulus with one node and without marked points.
These components contribute positive-dimensional subgroups
 to $\Aut(\Sigma)$.
The following discussion is an immediate generalization of
 [F-O: pp.~989~-~990] and [Sie1: Sec.~2.2] to
 the case of labelled-bordered Riemann surfaces.
The moduli space $\widetilde{M}_{(g,h),(n,\vec{m})}$
  of isomorphism classes of {\it prestable} labelled-bordered
   Riemann surfaces of type $((g,h),(n,\vec{m}))$
  can be associated to an Artin stack.
The discussion below gives a substitute quotient topology structure.

A {\it semi-universal deformation} ${\cal C}/Def(\Sigma)$
 of $\Sigma$, together with a specification of
  an {\it approximate pseudo-$\Aut(\Sigma)$-action}
  on ${\cal C}/Def(\Sigma)$, defines a local chart
  $$
   \psi_{[\Sigma]}\;:\; \Def(\Sigma)\;
      \longrightarrow\; \widetilde{\cal M}_{(g,h),(n,\vec{m})}
  $$
 of $[\Sigma]\in\widetilde{M}_{(g,h),(n,\vec{m})}$.
Such a pair of data can be constructed as follows:
 \begin{itemize}
  \item[$(1)$]
  {\it the defining family ${\cal C}/\Def(\Sigma)$ of the chart}$\,$:
   Let
    $\Sigma^{\prime}
    = (\Sigma, (p^{\prime}_{\,\cdot\,})_{\,\cdot\,})$,
     where $(p^{\prime}_{\,\cdot\,})_{\,\cdot\,}$ is a minimal
      tuple of rigidifying additional marked points on $\Sigma$
      that are disjoint from all the existing special points
      of $\Sigma$.
   Take ${\cal C}/\Def(\Sigma)$ to be the universal deformation
    ${\cal C}^{\prime}/Def(\Sigma^{\prime})$ of $\Sigma^{\prime}$
    with the sections $s^{\prime}_{\,\cdot\,}$
    associated to $p^{\prime}_{\,\cdot\,}$'s removed.

  \item[(2)]
  {\it the approximate pseudo-$\Aut(\Sigma)$-action}$\,$:
   Let $e$ be the identity element of $\Aut(\Sigma)$
    and recall that the central fiber ${\cal C}^{\prime}_0$
     of ${\cal C}^{\prime}/\Def(\Sigma^{\prime})$ is $\Sigma^{\prime}$.
   Consider the product family
    $(\Aut(\Sigma)\times {\cal C}^{\prime})/
                 (\Aut(\Sigma)\times \Def(\Sigma^{\prime}))$.
   First, extend the section $s_{\,\cdot\,}$ over
    $\{e\}\times \Def(\Sigma^{\prime})$ to over
    $\Aut(\Sigma)\times \{0\}$ by setting
   $s^{\prime}_{\,\cdot\,}(\sigma,0)
         =\sigma\cdot p^{\prime}_{\,\cdot\,}$.
   Then, further extend them to a collection of sections
    $s^{\prime}_{\,\cdot\,}$ over a neighborhood
     (still denoted by $\Aut(\Sigma)\times\Def(\Sigma^{\prime})$,
        though in general it may not be a product)  of
     $\Aut(\Sigma)\times\{0\}
           \subset \Aut(\Sigma)\times\Def(\Sigma^{\prime})$
    whose image in a fiber are disjoint from each other and
     from the special points and the image of the existing sections
     associated $\Sigma$ on that fiber.
   This can always be done but is non-canonical/non-unique.
   Denote
    the resulting family by
     $( (\Aut(\Sigma)\times{\cal C})/(\Aut(\Sigma)\Def(\Sigma))\,,\,
        (s^{\prime}_{\,\cdot\,})_{\,\cdot\,} )$  and
    the restriction of $s^{\prime}_{\,\cdot\,}$ to over
     $\{\sigma\}\times \Def(\Sigma)$ by $s^{\prime}_{\,\cdot\,,\sigma}$.

  \item[$\cdot$]
   From the universal property of the family
     ${\cal C}^{\prime}/\Def(\Sigma^{\prime})$
   the unique isomorphism
     from the central fiber
      $(\Sigma, (\sigma\cdot p^{\prime}_{\,\cdot\,})_{\,\cdot\,})$
      of the family
      $( (\{\sigma\}\times{\cal C})/(\{\sigma\}\times\Def(\Sigma))\,,\,
         (s_{\,\cdot\,,\,\sigma})_{\,\cdot\,} )$
     to the central fiber $\Sigma^{\prime}$ of
      ${\cal C}^{\prime}/\Def(\Sigma^{\prime})$
    extends to a unique isomorphism
    $$
     \Phi^{\prime}_{\sigma}\; :\;
      (\{\sigma\}\times{\cal C})/(\{\sigma\}\times\Def(\Sigma))\;
       \longrightarrow\; {\cal C}^{\prime}/Def(\Sigma^{\prime})\,,
    $$
    assuming that the neighborhood of
     $\Aut(\Sigma)\times\{0\}$ in $\Aut(\Sigma)\times\Def(\Sigma)$
     we chose is small enough.

  \item[$\cdot$]
   Let
    $$
     F_{\sigma}\;:\;
      ( (\{\sigma\}\times{\cal C})/(\{\sigma\}\times\Def(\Sigma))\,,\,
          (s^{\prime}_{\,\cdot\,,\sigma})_{\,\cdot\,}  ) \;
      \longrightarrow\; {\cal C}/Def(\Sigma)\,,
    $$
    be the forgetful isomorphism that forgets the tuple
     $(s^{\prime}_{\,\cdot\,,\sigma})_{\,\cdot\,}$
     of rigidifying section.
   The morphism
    $$
     \begin{array}{ccccc}
      \Phi_{[\Sigma]}  & :
       & (\Aut(\Sigma)\times {\cal C})/(\Aut(\Sigma)\times\Def(\Sigma))
       & \longrightarrow  & {\cal C}/Def(\Sigma) \\[.6ex]
      && (\sigma, x)  & \longmapsto
       & \sigma\cdot x\; :=\;
         (F_{\sigma}\circ \Phi_{\sigma}^{\prime\,-1}\circ F_e^{-1}) (x)
     \end{array}
    $$
    defines then an {\it approximate}\footnote{Here,
       the term ``{\it approximate}" is referring to the fact that
       the composition law
       $\Phi_{[\Sigma]}(\sigma_1,\Phi_{[\Sigma]}(\sigma_2,x))
          =\Phi_{[\Sigma]}(\sigma_1\sigma_2,x)$
       may not hold but, for $\Def(\Sigma)$ small enough,
       $\Phi_{[\Sigma]}(\sigma_1,\Phi_{[\Sigma]}(\sigma_2,x))$
       is always in a small neighborhood of
       $\Phi_{[\Sigma]}(\sigma_1\sigma_2,x)$.}
     pseudo-$\Aut(\Sigma)$-action on
    ${\cal C}/Def(\Sigma)$.

  \item[(3)]
   {\it The coordinate map $\psi_{[\Sigma]}$}$\,$:
    The family ${\cal C}/\Def(\Sigma)$ specifies a map
     $\psi_{[\Sigma]}:
       \Def(\Sigma) \rightarrow \widetilde{M}_{(g,h),(n,\vec{m})}$
    by sending $b\in\Def(\Sigma)$ to the isomorphism class
     $[C_b]\in \widetilde{\cal M}_{(g,h),(n,\vec{m})}$
     of the fiber $C_b$ of ${\cal C}$ over $b$.
 \end{itemize}
Topologically, $\psi_{[\Sigma]}$ is a quotient of
 a neighborhood of the origin of the manifold-with-corners
 $$
  \begin{array}{c}
   \Ext^1_{\Sigma_{\tinyBbb C}}
    \left(
      \Omega_{\Sigma_{\tinyBbb C}}
      ( \sum_{i=1}^n(p_i+\overline{p_i})\,  +\,
        \sum_{j=1}^h\sum_{k=1}^{m_j}p_{jk} + D_{\rigidifying})\,,\,
      {\cal O}_{\Sigma_{\tinyBbb C}}
     \right)^{\tau}\;
      \\[1.6ex]
   \simeq\;
    {\Bbb C}^{3g-3+h+n^{\prime}+d_c}  \times
    \overline{\Bbb H}^{n^{\prime\prime}}  \times
    {\Bbb R}^{h-n_{bn}+m_1+\,\cdots\,+m_h+d_b}  \times
    ({\Bbb R}_{\ge 0})^{n_{bn}}
  \end{array}
 $$
 by the induced $\Aut(\Sigma)$-action,
 where
  $D_{\rigidifying}$ is a minimal $\tau$-invariant rigidifying
   divisor on $\Sigma_{\scriptsizeBbb C}$ whose support is disjoint
   from the existing special points on $\Sigma_{\scriptsize C}$,
  $n=\doteq n^{\prime}+n^{\prime\prime}$,  and
  $d_c$ (resp.\ $d_b$) is the complex (resp.\ real) dimension of
   the product of the automorphism group of the closed
   (resp.\ bordered) unstable components of $\Sigma$.
The stacky (real) dimension of these charts, i.e.\
 $\dimm\Def(\Sigma)-\dimm Aut(\Sigma)$, remains
 $6g-6+3h+2n+m_1+\,,\cdots\,+m_h$.

\bigskip

\noindent
{\bf Definition 2.3
  [standard local chart of $\widetilde{M}_{(g,h),(n,\vec{m})}$].}
{\rm
 We will call the tuple
   $(\Def(\Sigma),\Phi_{[\Sigma]},\psi_{[\Sigma]})$,
   in short $\Def(\Sigma)$, a {\it standard local chart} of
   $[\Sigma]\in\widetilde{\cal M}_{(g,h),(n,\vec{m})}$ and
  the ${\cal C}$ that accompanies $\Def(\Sigma)$ in the construction
    and is equipped with the approximate pseudo-$\Aut(\Sigma)$-action
  the {\it universal curve} over the chart $\Def(\Sigma)$.
} % end-definition

\bigskip

\begin{flushleft}
{\bf Resemblance of the approximate pseudo-action with a pseudo-action.}
\end{flushleft}
$\Phi_{[\Sigma]}$ defines a relation $\sim$ on $\Def(\Sigma)$
 generated by
  $b_1\sim b_2$ if there exists a
  $\sigma\in \Aut(\Sigma)$ such that $b_2=\sigma\cdot b_1$.
As the major step of the construction is a morphism to
 the universal deformation space of $\Sigma$ with added rigidifying
 marked points,
it remains true that two fibers $C_{b_1}$ and $C_{b_2}$ of
 ${\cal C}/\Def(\Sigma)$ are isomorphic if and only if
 $b_1\sim b_2$;
and, in this case,
 an isomorphism $C_{b_2}\simeq C_{b_1}$ can be given by
 the composition $\sigma_1\cdot\ldots\cdot\sigma_k\cdot$ for some
 $\sigma_1\,,\,\ldots\,,\,\sigma_k\in\Aut(\Sigma)$.
Furthermore,
 as long as $\Def(\Sigma)$ in the construction is small enough,
 the map $\sigma:{\cal C}/\Def(\Sigma)\rightarrow {\cal C}/\Def(\Sigma)$
  is bijective on the domain it is defined.
These two properties make the approximate pseudo-$\Aut(\Sigma)$-action
 on ${\cal C}/\Def(\Sigma)$ equally good as a genuine one.

\bigskip

\noindent
{\bf Definition 2.4 [$\Aut(\Sigma)$-orbit].} {\rm
 An equivalence class of $\sim$ in $\Def(\Sigma)$
  is called an {\it $\Aut(\Sigma)$-orbit} on $\Def(\Sigma)$.
 Similarly for the approximate pseudo-$\Aut(\Sigma)$-action
  on ${\cal C}$.
} % end-definition

\bigskip

$\Def(\Sigma)$ admits a stratification by locally closed subsets
 such that points in the same stratum have the corresponding fibers
 in ${\cal C}$ of the same topological type.
It follows that
 the approximate pseudo-$\Aut(\Sigma)$-action leaves each stratum
  invariant  and
 points of $\Def(\Sigma)$ in the same fiber have their
  $\Aut(\Sigma)$-orbits of the same dimension.
When not of the finitely many exceptional types,
a general point $b\in\Def(\Sigma)$ has the $\Aut(\Sigma)$-orbit
 $\Aut(\Sigma)\cdot b$ of the same dimension as $\Aut(\Sigma)$,
while $0\in Def(\Sigma)$, which corresponds to the fiber $\Sigma$,
 is always a fixed point of $\Aut(\Sigma)$.

\bigskip

\noindent
{\it Remark 2.5 $[$abelian $\Aut(\Sigma)$$]$.}
When $\Aut(\Sigma)$ is abelian, a similar construction as in Sec.~1.1.2
 shows that $\Aut(\Sigma)$ does pseudo-acts on ${\cal C}/\Def(\Sigma)$
 in this case.

\bigskip

\noindent
{\bf Lemma 2.6 [pseudo-$\Gamma\!\cdot\!\Aut_e(\Sigma)^{\circ}$-action].}
{\it
 Let
  $\Gamma$ be a finite subgroup of $\Aut(\Sigma)$,
  $\Aut_e(\Sigma)^{\circ}$ be a small enough neighborhood of
   the identity element $e$ of $\Aut(\Sigma)$, and\\
  $\Gamma\!\cdot\!\Aut_e(\Sigma)
   =\cup_{\sigma\in\Gamma}\,\sigma\cdot\Aut_e(\Sigma)^{\circ}$.
 Then, possibly after shrinking $\Def(\Sigma)$,
  the defining $\Gamma\!\cdot\!\Aut_e(\Sigma)^{\circ}$-action on the
  center fiber $\Sigma$ of ${\cal C}/\Def(\Sigma)$ extends to a
  pseudo-action on ${\cal C}/\Def(\Sigma)$ by isomorphisms.
 This pseudo-$\Gamma\!\cdot\!\Aut_e(\Sigma)^{\circ}$-action extends
  to an approximate pseudo-$\Aut(\Sigma)$-action on
  ${\cal C}/\Def(\Sigma)$ by isomorphisms.
} % end-lemma

\bigskip

\noindent {\it Proof.}
 Fix a rigidifying devisor $\sum_{\,\cdot\,}p^{\prime}_{\,\cdot\,}$
   on $\Sigma$ away from the nodes and
 let $\Sigma=(\cup_{q_i}N_i)\cup(\cup_jV_j)$
  be a neck-trunk decomposition of $\Sigma$
   (cf.\ the thick-thin decomposition of $\Sigma$ when $\Sigma$
         is of hyperbolic type),
   where
    $N_i$ is a neck on $\Sigma$ in a small neighborhood of node $q_i$
     with $q_i$ running over the set of nodes of $\Sigma$,  and
    $V_j$ be a connected component of $\Sigma-\cup_{q_i}N_i$,
  such that $\Gamma\!\cdot\!\Aut_e(\Sigma)^{\circ}(\cup_i\partial N_i)$
   remains in a tubular neighborhood of $\cup_i\partial N_i$ in $\Sigma$
  and the $\Gamma\!\cdot\!\Aut_e(\Sigma)^{\circ}$-orbits of
   all marked points, including the added regidifying ones
   $p^{\prime}_{\,\cdot\,}$, are away from this tubular neighborhood.
 As $\Gamma$ sends nodes to nodes, this can be realized as long as
  $\Aut_e(\Sigma)^{\circ}$ is small enough.
 Extend this neck-trunk decomposition of $\Sigma$ to a neck-trunk
  decomposition
  $$
   {\cal C}/\Def(\Sigma)\;=\;
    \left(\rule{0ex}{2ex}\cup_{q_i}\Neck(q_i)  \right)\bigcup
    \left(\rule{0ex}{2ex} \cup_j Trunk_j\right)
  $$
  of ${\cal C}/\Def(\Sigma)$,
  where
   $\{q_i\}_i$ is the set of nodes of $\Sigma$;
   $\Neck(q_i)$ is a neck region in ${\cal C}$ associated to $q_i$; and
   $\{\Trunk_j/\Def(\Sigma)\}_j$ is the set of connected components
    of ${\cal C}/\Def(\Sigma)-\Neck(q_i)$,
    equipped with a fixed product decomposition
    $\Trunk_j= \Def(\Sigma)\times V_j$.
 This can be realized as long as $\Def(\Sigma)$ is small enough.
 Denote the section of ${\cal C}/\Def(\Sigma)$ associated to
  $p^{\prime}_{\,\cdot\,}$ by $s^{\prime}_{\,\cdot\,}$.
 %%%%%%%%%%%%%%%%%%%%
 % To fix a Dehn twist ambiguity that will occur, let us digress here
 %  to take a look at the pseudo-group
 %  $\Gamma\!\cdot\!\Aut_e(\Sigma)^{\circ}$.
 % Let $\Gamma_0:=\Gamma\cap\Aut_e(\Sigma)$.
 % It is an abelian subgroup of $\Gamma$ and it naturally acts
 %  on ${\cal C}/\Def(\Sigma)$.
 % The quotient group  $\Gamma/(\Gamma\cap\Aut)$ parameterizes
 %  the set of homotopy equivalence classes of automorphisms
 %  of $\Sigma$ parameterized by $\Gamma\!\cdot\!\Aut_e(\Sigma)^{\circ}$.
 % The exact sequence
 %  $$
 %   1\rightarrow\Gamma_0\rightarrow \Gamma\rightarrow \Gamma/\Gamma_0
 %      \rightarrow 1
 %  $$
 %  realizes $\Gamma$ as a semi-direct product
 %   $\Gamma=\Gamma_0\rsemiproduct (\Gamma/\Gamma_0)$,
 %   where $\Gamma/\Gamma_0$ acts on $\Gamma_0$ by conjugation.
 %%%%%%%%%%%%%%%%%%%%
 The specification of a neck-trunk decomposition of
  ${\cal C}/\Def(\Sigma)$ specifies simultaneously how each fiber
  of ${\cal C}/\Def(\Sigma)$ is obtained from a cut-and-paste of
  $\Sigma$, (cf.\ the re-forging morphisms in Sec.~1.1.1).
 This then induces a pseudo-$\Gamma\!\cdot\!\Aut_e(\Sigma)^{\circ}$-action
  $$
   \Phi_{[\Sigma]}^{\circ}\;:\;
    (\Gamma\!\cdot\!\Aut_e(\Sigma)^{\circ})\times
     ({\cal C}/\Def(\Sigma))\;\longrightarrow\;
     {\cal C}/\Def(\Sigma)
  $$
  on ${\cal C}/\Def(\Sigma)$ as the cut-and-paste region remain near
  the neck region of $\Sigma$ under the smallness assumption of
  $\Aut_e(\Sigma)^{\circ}$.
 This proves the first statement of the lemma.

 To extend this to an approximate pseudo-$\Aut(\Sigma)$-action on
  ${\cal C}/\Def(\Sigma)$,
 consider the product family
  $(\Aut(\Sigma)\times{\cal C})/(\Aut(\Sigma)\times\Def(\Sigma))$.
 Recall $s^{\prime}_{\,\cdot\,}$ the sections of ${\cal C}/\Def(\Sigma)$
  that correspond to the added rigidifying points
  $p^{\prime}_{\,\cdot\,}$ on $\Sigma$.
 Their image lies in the trunk region of ${\cal C}/\Def(\Sigma)$.
 Extend these sections first to over
  $\Gamma\!\cdot\!\Aut_e(\Sigma)^{\circ}\times\Def(\Sigma)$
  by setting
  $s^{\prime}_{\,\cdot\,,\,\sigma}
         = \sigma\cdot s^{\prime}_{\,\sigma\,}$
  over $\{\sigma\}\times\Def(\Sigma)$,
  where $(\{\sigma\}\times{\cal C})/(\{\sigma\}\times\Def(\Sigma))$
   is canonically identified with ${\cal C}/\Def(\Sigma)$.
 These sections again have their image in the trunk region of
  $(\{\sigma\}\times {\cal C})/(\{\sigma\}\times\Def(\Sigma))$.
 Extend these sections next to over $\Aut(\Sigma)\times\{0\}$ as well
  by the $\Aut(\Sigma)$-action on $\Sigma$.
 Finally extend the resulting sections to over $\Aut(\Sigma)\times \Def(\Sigma)$.
 This then defines an approximate pseudo-$\Aut(\Sigma)$-action
  on ${\cal C}/\Def(\Sigma)$ by isomorphisms
  that extends the pseudo-$\Gamma\!\cdot\!\Aut_e(\Sigma)^{\circ}$-action
  constructed.
 This concludes the proof.

\noindent\hspace{15cm}$\Box$

\bigskip

The same argument gives also:

\bigskip

\noindent
{\bf Lemma 2.7 [finite group].} {\it
 Any finite group action on $\Sigma$ by automorphisms extends to
  an action on ${\cal C}/\Def(\Sigma)$ by isomorphisms.
 This action extends to
  a pseudo- $\Gamma\!\cdot\!\Aut_e(\Sigma)^{\circ}$-action on
  ${\cal C}/\Def(\Sigma)$ and then to an approximate
  pseudo-$\Aut(\Sigma)$ on ${\cal C}/\Def(\Sigma)$,
  both by isomorphisms.
} % end-corollary

%%%%%%%%%%%%%%%%%%%%%%%%%%%%%%%%%%%%%%%%%%%%%%%%%%%%%%%%%%%%%%%%%
% ??????????.
% The maximal abelian subgroup $\Aut_0^{\circ,\ab}$ of
%  $\Aut_0^{\circ}(\Sigma)$ adapted to a choice of disc neighborhoods
%  at nodes of unstable components of $\Sigma$.
% The restriction of the approximate pseudo-
%  $\Aut_0^{\circ}(\Sigma)$-action on ${\cal C}/Def(\Sigma)$
%  to the subgroup $\Aut_0^{\circ,\ab}$ is indeed a group action and
%  can be described explicitly from the gluing construction of
%  deformations of labelled-bordered Riemann surfaces.
% ????????????.
%%%%%%%%%%%%%%%%%%%%%%%%%%%%%%%%%%%%%%%%%%%%%%%%%%%%%%%%%%%%%%%%%

\bigskip

\begin{flushleft}
{\bf The quotient topology on $\widetilde{\cal M}_{(g,h),(n,\vec{m})}$
      and the stabilization morphism.}
\end{flushleft}
The {\it quotient topology} on $\widetilde{\cal M}_{(g,h),(n,\vec{m})}$
 is defined by setting
a subset $U\subset\widetilde{\cal M}_{(g,h),(n,\vec{m})}$
  to be {\it open}
 if $U=\cup_{\alpha}U_{\alpha}$ such that
  there exist a collection of standard local charts
   $(V_{\alpha},\Phi_{\alpha},\psi_{\alpha})$ of
   $\widetilde{M}_{(g,h),(n,\vec{m})}$
  such that
   $U_{\alpha}\subset \psi_{\alpha}(V_{\alpha})$  and that
   $\psi_{\alpha}^{-1}(U_{\alpha})$ is open in $V_{\alpha}$.
This is similar to the construction in Sec.~1.1.3 and Sec.~1.2
 for the quotient topology on $\widehat{B}$ and $\widehat{A}$.

For $((g,h),(n,\vec{m}))$ with $2(2g+h+n)+m_1+\,\cdots\,+m_h\ge 5$,
 {\it stabilization} of prestable labelled-bordered Riemann surfaces
  by contracting the unstable components
 gives rise to
  a {\it flat local complete intersection morphism}\footnote{See
    [Fu] for a general definition of
    {\it local complete intersection morphism}.
    Such a morphism has a well-defined Gysin map,
    and hence push-pull, on cycles.}
   $\st:{\cal C}/\Def(\Sigma)  \rightarrow
        {\cal C}_{\scriptsizest}/\Def(\Sigma_{\scriptsizest})$,
  together with a group homomorphism
   $\Aut(\Sigma)\rightarrow \Aut(\Sigma_{\scriptsizest})$
   that makes $\st$ equivariant,
 for each $[\Sigma]\in\widetilde{\cal M}_{(g,h),(n,\vec{m})}$.
The collection of these pairs of morphisms on
  local charts-with-structure-group
 descend to the {\it stabilization morphism}
 $\widetilde{\st}: \widetilde{\cal M}_{(g,h),(n,\vec{m})}
        \rightarrow \overline{M}_{(g,h),(n,\vec{m})}$.
We say that  {\it
 $\widetilde{\st}$ is a local complete intersection morphism
 in the stacky sense}.
It is continuous with respect to the quotient topology
 on $\widetilde{M}_{(g,h),(n,\vec{m})}$.
The inclusion
 $\overline{M}_{(g,h),(n,\vec{m})}
   \hookrightarrow \widetilde{M}_{(g,h,(n,\vec{m})}$
 is a section to $\widetilde{\st}$ with open-dense image.

\bigskip

\noindent
{\it Remark 2.8 $[$local factorization of $\st\,$$]$.}
Assume that $2(2g+h+n)+m_1+\,\cdots\,+m_h\ge 5$.
Let $\Sigma=\Sigma^s\cup\Sigma^u$,
 where
  the subcurve $\Sigma^u$ consists
   of all the unstable irreducible components of $\Sigma$  and
  $\Sigma^s$ is the union of the remaining irreducible components.
Then a connected component of $\Sigma^u$ may intersect
 $\Sigma^s$ at either $1$ or $2$ nodes of $\Sigma$;
it is called a {\it tree} in the formal case and
 a {\it chain} in the latter case,
 in which it can only be either
  a chain of ${\Bbb P}^1$ of the form
   ${\Bbb P}^1_{(1)}\cup\,\cdots\,\cup{\Bbb P}^1_{(k)}$
   with $0$ of ${\Bbb P}^1_{(i)}$ glued to
        $\infty$ of ${\Bbb P}^1_{(i+1)}$,
  or a chain of discs $D^2=\{z\in{\Bbb C}\,:\, |z|\le 1\}$
   of the form $D^2_{(1)}\cup\,\,\cdots\,\cup D^2_{(k)}$
   with $-\sqrt{-1}$ of $D^2_{(i)}$ glued to
        $\sqrt{-1}$ of $D^2_{(i+1)}$.
These are reflected to the stabilization map:
 locally $\st$ can be factorized to a composition of
  a projection map of a product space,
   for no collapsing or collapsing a tree of unstable components;
  a map of the form ${\mathbf \pi}[k]:B[k]\rightarrow B$
   in Lemma 1.1.1.4,  % Lemma [natural map and its fibers]
   for collapsing a chain of unstable ${\Bbb P}^1$ components;  and
  a map of the form
   $$
    ({\Bbb R}_{\ge 0})^{k+1}\;\longrightarrow\; {\Bbb R}_{\ge 0}\,,
      \hspace{2em}
      (t_0,\,\ldots\,,t_k)\;\longmapsto\; t_0\,\cdots\,t_k\,,
   $$
   for collapsing a chain of unstable discs.
Cf.\ [Sie1: end of Sec.~2.2].

\bigskip

\section{The moduli space
         $\overline{\cal M}_{(g,h),(n,\vec{m})}
                  (W/B,L\,|\,[\beta],\vec{\gamma},\mu)$
         of stable maps.}

In the previous two sections, we discuss respectively
 the targets and the domains of the maps we want to study.
However, as a lesson from the various standard moduli problems
 in algebraic geometry, which can almost always be traced back to
 the complicated problem of Hilbert-schemes,
to render a reasonable moduli space of maps from
  bordered Riemann surfaces to fibers of $\widehat{W}/\widehat{B}$,
 we need to fix some combinatorial quantities of such maps
  that are constant for a continuous/flat family.
The closed Gromov-Witten theory indicates a partial set of such data:
 the {\it combinatorial type of domain curves},
 the image curve class $\beta\in H_2(X,L;{\Bbb Z})$,  and
 {\it boundary loop class $\vec{\gamma}$} from $H_1(L;{\Bbb Z})$.
The study of [Liu(C)] implies that for open Gromov-Witten theory
 the boundary effect is reflected also in the
 {\it Maslov index} $\mu\in {\Bbb Z}$,
 which is not fixed by $\beta$ in general.
This quantity thus has to be generalized to our case and
 be included in the combinatorial data.
This is done in Sec.~3.1  and
 the generalized Maslov index does enter the operator index
  in Sec.~5.3.1.
However, this addition of data is not enough.
While it turns out that the
 datum $\vec{\gamma}$ from $H_1(L;{\Bbb Z})$ is not influenced,
the datum $\beta\in H_2(X,L;{\Bbb Z})$
 is not the correct choice of the image curve class datum in our case
 since in general it is not well-defined to all fibers
  in the family $W[k]/B[k]$, which contains $X$ as a fiber,
 due to the monodromy effect.
It thus has to be enlarged to and replace by
 the {\it minimal common monodromy-invariant curve-class subset}
  $[\beta]\subset H_2(X,L;{\Bbb Z})$,
 generated by $\beta$ under the monodromy of $W[k]/B[k]$,
  for all $k\in {\Bbb Z}_{\ge 0}$.
This is done in Sec.~3.2.
Once these combinatorial data are identified,
 one can then define the related moduli space
 $\overline{\cal M}_{(g,h),(n, \vec{m})}
          (W/B,L\,|\,[\beta],\vec{\gamma},\mu)$  of maps accordingly.
This is done in Sec.~3.3.

\bigskip

\subsection{Maslov index of a map to a singular space or
            a relative pair.}

A generalization of the notion of Maslov index to a map from
 a bordered Riemann surface to a relative pair or
 a singular space from a symplectic cut
 is given in this subsection.
This quantity is needed to select a reasonable (union of)
 component(s) of the moduli space of stable maps in question.

%%%%%%%%%%%%%%%%%%%%%%%%%%%%%%%%%%%%%%%%%%%%%%%%%%%%%%%
% The Maslov index of a $C^{\infty}$
%  relative map $f:\Sigma\rightarrow (Y,L;D)$ from a prestable
%   bordered Riemann surface $\Sigma$ to $(Y,L;D)$
%   with $f({\partial\Sigma})\subset L$
%  is meant to encode
%   how the Lagrangian submanifold $L$ rotates around from the point
%   of view of the domain surface $\Sigma$.
% This quantity will influence the deformation properties of $f$,
%  e.g.~the virtual dimension around $f$ in the moduli space of open
%       relative stable maps,
%  and is not determined solely by the image class
%   $f_{\ast}[\Sigma] \in H_2(X,L;{\Bbb Z})$.
% For our study, we have to define the Maslov index for $f$ with
%  target $Y$ being the singular space from symplectic cuts as well.
%%%%%%%%%%%%%%%%%%%%%%%%%%%%%%%%%%%%%%%%%%%%%%%%%%%%%%%
%
Given a $C^{\infty}$ map $f:(\Sigma,\partial\Sigma)\rightarrow (X,L)$
 from a prestable bordered Riemann surface $\Sigma$ to a smooth
 symplectic manifold $X$.
Endow $X$ with a compatible almost-complex structure $J$ that renders
 $T_{\ast}X$ a complex vector bundle with
 $T_{\ast}L\hookrightarrow (T_{\ast}X)|_L$ as a totally real subbundle.
Then $E:=f^{\ast}(\det(T_{\ast}X))$ is a complex line bundle on
 $\Sigma$ whose restriction to $\partial\Sigma$ contains a real
 line subbundle $E_{\scriptsizeBbb R}(L)$ associated to
 $f^{\ast}(T_{\ast}L)$.
The Maslov index of $f$ in this case
 (cf.\ [K-L: Definition 3.7.2]) is defined by:

\bigskip

\noindent
{\bf Definition 3.1.1 [Maslov index - smooth target].} {\rm
 The {\it Maslov index} $\mu(f)$ of the $C^{\infty}$ map $f$ above
   is twice the index of a general extension of
   $E_{\scriptsizeBbb R}(L)\subset E|_{\partial\Sigma}$
  to a real-line subbundle with isolated singularities in $E$,
   still denoted by $E_{\scriptsizeBbb R}(L)$,
  over the whole $\Sigma$.
 For convenience, we set
  $\mu(f)=2\,\degree (f^{\ast}\det(T_{\ast}X))$ if
  either $L$ or $\partial\Sigma$ is empty.
} % end-definition

\bigskip

\noindent
Note that this definition is more in the almost-complex category than
 in the symplectic category.
However, $\mu(f)$ thus defined is independent of
 the choice of $\omega$-compatible $J$ on $X$ and
 the general extension $E_{\scriptsizeBbb R}$ on $\Sigma$.
To turn the real line field language to the more convenient real vector
 field language, one considers the complex line bundle
 $E^{\,\otimes 2}$ and rephrases $\mu(f)$ as the index
 of a general global section $s$ of $E^{\,\otimes 2}$ that extends
 the section $s_L$ in $E^{\,\otimes 2}|_{\partial\Sigma}$ determined by
 $f^{\ast}(T_{\ast}L)$.

In the complex K\"{a}hler category, the key object in the above
 description of $\mu(f)$, namely the (complex) determinant line bundle
  $K:= \det \Omega_X = (\det T_{\ast}X)^{-1}$,
 can be defined for a {\it singular} $Y$ from a symplectic cut.
Once having this, the Maslov index of a $C^{\infty}$ map
 $f:\Sigma \rightarrow (Y,L)$,
  with $L$ disjoint from the singular locus $Y_{\sing}$ of $Y$,
 can be defined in exactly the same way as above:
 the index of a global section $s$ in
 $f^{\ast}(K^{\,\otimes(-2)})$ that extends a global section
 $s_{L}$ in $(f^{\ast}(K^{\,\otimes(-2)}))|_{\partial \Sigma}$
 determined by $f^{\ast}(T_{\ast}L)$.
Taking $\det$ of a coherent sheaf in algebraic geometry brings in
 a twisting effect from a divisor whose support is contained in
 the non-locally-free locus of the coherent sheaf (cf.~[Kn-M]).
For $\Omega_Y$ in K\"{a}hler category,
 such locus coincides with the singular locus of $Y$.
One can compute such effect explicitly and compare them with
 the contribution to $\mu(f)$ from each individual smooth
 irreducible component of $Y$.
The result can be stated in both the symplectic and
 the almost-complex category.
This gives rise to the following definitions.

\bigskip

\noindent
{\bf Definition 3.1.2
     [Maslov index - relative pair and symplectic gluing].}
{\rm
 The Maslov index of a $C^{\infty}$ map from
  a bordered Riemann surface $\Sigma$
  to a symplectic pair or a symplectic space from a symplectic cut
  is defined as follows:
 \begin{itemize}
  \item[(1)]
   Let
    $(Z,L;D)$ be a smooth symplectic pair $(Z;D)$ with
     a Lagrangian submanifold $L$ disjoint from $D$  and
    $f:(\Sigma,\partial\Sigma)\rightarrow (Z,L)$ be a $C^{\infty}$ map.
   Then, define the {\it Maslov index} of $f$ {\it relative to $D$}
    to be
    $$
     \mu^{\rel}(f) \; =\; \mu(f)\, -\, 2\, f_{\ast}[\Sigma]\cdot D\,,
    $$
    where $\mu(f)$ is the usual Maslov index of $f$
     as defined in Definition 3.1.1.
      % Definition [Maslov index - smooth target]
   (If $L$ is empty, then set
     $\mu(f)=\deg f^{\ast}(K_Z^{\otimes(-2)})
            =-2\, f_{\ast}[\Sigma]\cdot K_Z$.
    Note that both $L$ and $D$ in the definition can be disconnected.)

  \item[(2)]
   Let
    $(Y,L)
     =(Y_1,L_1)\cup_{D_1\simeq D_2} (Y_2,L_2)$
     be the singular symplectic space from gluing of
     two Lagrangian-decorated relative pairs $(Y_1,L_1;D_1)$
     and $(Y_2,L_2;D_2)$ and
    $f=f_1\sqcup f_2:\Sigma:= \Sigma_1\cup\Sigma_2
         \rightarrow (Y_1,L_1)\cup_D (Y_2,L_2)$
     be a $C^{\infty}$ map to $(Y,L)$.
    Then, define the {\it Maslov index} of $f$ to be
     $$
      \mu(f)\; =\; \mu^{\rel}(f_1) + \mu^{\rel}(f_2)\;
         =\; (\mu(f_1)-2f_{1\,\ast}[\Sigma_1]\cdot D_1)\,
             +\, (\mu(f_2)-2f_{2\,\ast}[\Sigma_2]\cdot D_2)\,.
     $$

  \item[(3)]
   For a $C^{\infty}$ map $f$ to a symplectic space from gluing
    a finite collection of Lagrangian-decorated symplectic pairs,
   apply Item (1) and Item (2) above inductively to define the
    Maslov index $\mu(f)$ or $\mu^{\rel}(f)$.
 \end{itemize}

 The same definitions hold in the almost-complex category with
  $L$ replaced by a totally real submanifold  and
  $D$ replaced by a real-codimension-$2$ almost-complex submanifold.
} % end-definition

\bigskip

\noindent
{\bf Example 3.1.3 [relative Maslov index].} (Cf.\ Sec.~1.2.)
 Given $(Z,L;D)$, let $(Z_{[k]},L_{[k]};D_{[k]})$
  be the central fiber of its $k$-th expanded relative pairs.
 For an open relative stable map
  $f:\Sigma\rightarrow (Z_{[k]},L_{[k]};D_{[k]})$
  with the corresponding decomposition
  $f=f_0\sqcup f_1 \sqcup\,\cdots\,\sqcup f_k$,
  where $f_0:\Sigma_0\rightarrow Y$ and
   $f_i:\Sigma_i\rightarrow \Delta_i$, $i=1,\,\ldots,\,k$,
 the Maslov index of $f$ as a relative map is then
 $$
  \mu^{\rel}(f)\;
   = \; \left(\mu(f_0) - 2f_{0\,\ast}[\Sigma_0]\cdot D_0\right)\,
        -\, 2\, \sum_{i=1}^k\, f_{i\,\ast}[\Sigma_i]
                 \cdot (K_{\Delta_i} + D_{i,0}+D_{i,\infty})\,,
 $$
 where $\mu(f_0)$ is defined as in Definition 3.1.1 for smooth target.
  % Definition [Maslov index - smooth target]

\bigskip

We list as lemmas the basic invariance properties of the
 Maslov index of $C^{\infty}$ maps, as defined above,
 that are part of the foundations of later discussions.
The proof of these lemmas are straightforward and hence omitted.

\bigskip

\noindent
{\bf Lemma 3.1.4 [invariance under homotopy and deformation].} {\it
 $(1)$
  Let $Z$ be a smooth manifold of even dimension,
   $L$ be a smooth submanifold of $Z$ of the middle dimension, and
   $D$ be a smooth codimension-$2$ submanifold of $Z$ disjoint from $L$.
  Let
   $f_t:\Sigma\rightarrow (Z,\omega_t)$, $t\in [0,1]$,
    be a homotopy class of $C^{\infty}$ maps
    from a prestable bordered Riemann surface $\Sigma$ to $(Z;D)$
    with $f_t(\partial\Sigma)\subset L$ and
   $\omega_t$ is a $1$-parameter family of symplectic structures
    $($say, of class $C^2$$)$ on $Z$ keeping $L$ a Lagrangian
     submanifold and $D$ a symplectic submanifold.
  Then $\mu^{\rel}(f_0)=\mu^{\rel}(f_1)\,$.
 $(2)$
  Let
   $Y=Y_1\cup_D Y_2$ be a space from gluing smooth even-dimensional
    $(\,$manifold, codimension-$2$ submanifold$\,)$-pairs and
   $L$ be a smooth submanifold of $Y$ of the middle dimension disjoint
    from $D$.
  Let
   $f_t:\Sigma\rightarrow (Y,\omega_t)$, $t\in [0,1]$,
    be a homotopy class of $C^{\infty}$ maps
    from a prestable bordered Riemann surface $\Sigma$ to $Y$
    with $f_t(\partial\Sigma)\subset L$ and
   $\omega_t$ is a $1$-parameter family of symplectic structures
    $($say, of class $C^2$$)$ on $Y$ keeping $L$ a Lagrangian
     submanifold and $D$ a symplectic submanifold.
  Then $\mu(f_0)=\mu(f_1)\,$.
}

% \bigskip
%
% \noindent
% {\it Proof.}
% This follows from the topological constant nature of the index
%  of general sections of a topological complex line bundle on
%  a prestable bordered Riemann surface.
%
% \noindent\hspace{14cm}$\Box$

\bigskip

\noindent
{\bf Lemma 3.1.5 [invariance under domain degeneration].}
{\it
 Let $(X,L)$ be either a smooth symplectic manifold or a singular
  symplectic space from symplectic cut,
  with a Lagrangian submanifold $L$ disjoint from $X_{\sing}$.
 Let $p: \Sigma\rightarrow \underline{\Sigma}$ be a pinching map
  that arise from a degeneration of $\Sigma$ that pinches
  a finite disjoint union of simple loops on $\Sigma$.
 Given a $C^{\infty}$ map $f:\Sigma\rightarrow (X,L)$ and
  a family of deformations of $f$
  to a $g:\underline{\Sigma}\rightarrow (X,L)$,
 Then $\mu(f)=\mu(g)$.
 Similarly for $C^{\infty}$ maps into $(Z,L;D)$.
} % end-lemma

% \bigskip
%
% \noindent
% {\it Proof.}
%  Fix a compatible almost-complex structure on $X$.
%  Let $\det T_{\ast}X$ be the associated complex determinant line bundle.
%  (One can consider only the smooth locus of $X$ if $X$ is singular
%   for the purpose of comparing Maslov indices since all the maps are
%   required to be nondegenerate.)
% For $f:\Sigma\rightarrow X$ deformable
%  (nondegenerately for $X$ singular)
%  to $g:\underline{\Sigma}\rightarrow X$,
% $(f^{\ast}\det T_{\ast}X)^{\otimes 2}$ on $\Sigma$ must be
%  topologically trivial on the annuli around the simple loops pinched
%  by $\Sigma\rightarrow \underline{\Sigma}$.
% The lemma thus follows.
%
% \noindent\hspace{14cm}$\Box$

\bigskip

\noindent
{\bf Lemma 3.1.6 [invariance under symplectic cut on target].} {\it
  Let
   $\xi: (X,L)\rightarrow Y:= (Y_1,L_1)\cup_D (Y_2,L_2)$
   be a symplectic cut with $L_1$ and $L_2$ disjoint form $D$.
  $(1)$
   Let
    $f:\Sigma\rightarrow (X,L)$ be a $C^{\infty}$ map
     that intersects $\xi^{-1}(D)$ at a finite union of $S^1$-orbits
     and
    $g:\underline{\Sigma}\rightarrow Y$
     be the $C^{\infty}$ map descended from $f$,
     where $\underline{\Sigma}$ is obtained from $\Sigma$
      by pinching each connected component of $f^{-1}(\xi^{-1}(D))$
      to a nodal point.
   Then $\mu(g)=\mu(f)$.
  $(2)$
   Conversely, let
    $g:\underline{\Sigma}\rightarrow \underline{X}$
     be a pre-deformable $C^{\infty}$ map
        {\rm (}cf.\ Definition 3.3.1{\rm )}
        % Definition [stable map to fibers of $(W[k],L[k])/B[k]$]
     and
    $f:\Sigma\rightarrow X$ be a lifting of $g$,
     where $\Sigma$ is a deformation of
      $\underline{\Sigma}$ that smoothes exactly the nodes
       $g^{-1}(D)$ in $\underline{\Sigma}$.
   Then $\mu(f)=\mu(g)$.
} % end-lemma

\bigskip

\noindent
We remark that, if one associates the symplectic cut $\xi$ to
 a symplectic deformation family as constructed in [Go], [MC-W],
 and [I-P2],
then Lemma 3.1.6 is a corollary of [I-P2: Lemma 2.2].
  % Lemma [invariance under symplectic cut on target]
The same statements of these lemmas,  with
  $L$ replaced by a totally real submanifold and
  $D$ replaced by a real-codimension-$2$ almost-complex submanifold,
 in the almost-complex category hold as well.

\bigskip

\noindent
{\it Remark 3.1.7 $[$homotopy vs.\ homology$]$.}
 As in the absolute case in [Liu(C)],
 the Maslov index of an open relative stable map
  $f:\Sigma\rightarrow (Z,L;D)$ or the singular $(Y,L)$
  influences the deformation properties of $f$.
 Though a homotopy invariant, it is not determined by the image class
  $f_{\ast}[\Sigma]$ of $f$ in $H_2(Z,L;{\Bbb Z})$ or $H_2(Y,L;{\Bbb Z})$,
  cf.\ [K-L: Remark 4.2.2].

\bigskip

\subsection{Monodromy effect and the choice of curve class data in $H_2$.}

Recall
 the symplectic cut $\xi:X\rightarrow Y=Y_1\cup_D Y_2$  and
 the associated almost-complex degeneration $W/B$.
Let $L$ be an Lagrangian submanifold disjoint from the cutting locus
 $\xi^{-1}(D)$ then it gives rise to $(W,B\times L)/B$,
   where $L$ is totally real in each fiber of $W/B$;  and
the construction in Sec.~1.1 extends immediately to give
  expanded degenerations $(W[k], L[k])/B[k]$
   with the equivariant pseudo-${\Bbb G}_m[k]$-action,
  the topological space $(\widehat{W},\widehat{B}\times L)/\widehat{B}$,
  the standard local charts
   $\varphi[k]:(W[k],L[k])/B[k]
     \rightarrow (\widehat{W},\widehat{L})/\widehat{B}$
   of $(\widehat{W},\widehat{L})/\widehat{B}$
   with the product-induced map
    $\tilde{\mathbf p}[k]:(W[k],L[k])/B[k]\rightarrow (W,L)/B$.
Note that $\widehat{L}=\widehat{B}\times L$.
We remark that $L[k]\simeq B[k]\times L$ is a coisotropic submanifold
 in $W[k]$ and is fiberwise Lagrangian/totally-real over $B[k]$.
We can assume that $L[k]$ is contained in the trunk region
 $\overline{U_1}[k]\cup\overline{U_2}[k]$ of $W[k]/B[k]$.
${\mathbf p}[k]$ sends  the
 discriminant locus $\{\lambda_0\,\cdots\,\lambda_k=0\}\subset B[k]$
   of $W[k]/B[k]$
  to the discriminant locus $\{0\}\subset B$ of $W/B$  and the
 complement $B[k]_{\reg}:= B[k]-\{\lambda_0\,\cdots\,\lambda_k=0\}$
  to the complement $B_{\reg} :=B-\{0\}$.
Note that $\pi_1(B[k]_{\reg})\simeq {\Bbb Z}^{\oplus(k+1)}$  is
 generated by the canonically-oriented meridian $S^1$ of  the
 $(k+1)$-many coordinate hyperplanes of $B[k]$.
Fix topological trivializations
 $$
  W[k]_{{\footnotesizeBbb R}_{\ge 0} \cdot
        (\varepsilon^2/4\,,\,\cdots\,,\,\varepsilon^2/4)}\;
   \simeq\; ({\Bbb R}_{\ge 0} \cdot
             (\varepsilon^2/4\,,\,\cdots\,,\,\varepsilon^2/4)) \times X
 $$
 along the diagonal ray of $B[k]$'s.
This fixes an isomorphism
 $$
  H_2(W[k]_{\mbox{\tiny $\bullet$}},
       L[k]_{\mbox{\tiny $\bullet$}};{\Bbb Z})\;
  \simeq\; H_2(X,L;{\Bbb Z})\,,
  \hspace{2em}\mbox{for $\mbox{\tiny $\bullet$}
     \in {\footnotesizeBbb R}_{\ge 0} \cdot
          (\varepsilon^2/4\,,\,\cdots\,,\,\varepsilon^2/4)$}\,.
 $$
Via these identifications, $\pi_1(B[k]_{\reg})$ acts on
 $H_2(X,L;{\Bbb Z})$ by monodromy.
Furthermore, since $L$ is contained in the truck region of $X$,
 one has:

\bigskip

\noindent
{\bf Lemma 3.2.1 [trivial monodromy on $H_1(L;{\Bbb Z})$].} {\it
 As a fiber of $(W[k],L[k])/B[k]$, the monodromy
  $\pi_1(B[k]_{\reg})$-action on $H_1(L;{\Bbb Z})$
  is well-defined and is trivial;  and
 the connecting homomorphism
  $\partial: H_2(X,L;{\Bbb Z})\rightarrow H_1(L;{\Bbb Z})$
  is equivariant with respect to the $\pi_1(B[k]_{\reg})$-action.
} % end-lemma

\bigskip

\noindent
{\bf Lemma/Definition 3.2.2
[$(\widehat{W},\widehat{L})/\widehat{B}$-monodromy orbit].} {\rm
 {\it For each $\beta\in H_2(X,L;{\Bbb Z})$,
  all the monodromy-orbits $\pi_1(B[k]_{\reg})\cdot\beta$,
   $k\in {\Bbb Z}_{\ge 0}$, coincide.}
 We will name it the
  {\it $(\widehat{W},\widehat{L})/\widehat{B}$-monodromy orbit}
  of $\beta$ and denote it by $[\beta]$.
} % end-corollary

\bigskip

\noindent
{\it Proof.}
 Observe that the following diagram commutes
  $$
   \begin{array}{ccc}
    H_2(X,L;{\Bbb Z})
     & \stackrel{\mbox{\scriptsize $\tilde{\mathbf p}[k]_{\ast}$}}
                {\longrightarrow}
     & H_2(X,L;{\Bbb Z})\\
   \mbox{\scriptsize $a$}\downarrow\hspace{1ex}
    && \hspace{7ex}
       \downarrow\mbox{\scriptsize ${\mathbf p}[k]_{\ast}(a)$}\\
    H_2(X,L;{\Bbb Z})
     & \stackrel{\mbox{\scriptsize $\tilde{\mathbf p}[k]_{\ast}$}}
                {\longrightarrow}
     & H_2(X,L;{\Bbb Z})
   \end{array}
  $$
  for all $\,a\in \pi_1(B[k]_{\reg})$; i.e.\
 $\tilde{\mathbf p}[k]_{\ast}$ is equivariant with respect to
  the monodromy actions.
 As $\tilde{\mathbf p}[k]_{\ast}$ is the identity map under our
    identification  and
   ${\mathbf p}[k]_{\ast}:\pi_1(B[k]_{\reg})\rightarrow \pi_1(B_{\reg})$
    is surjective,
  the lemma follows immediately.

\noindent\hspace{15cm}$\Box$

\bigskip

Since the difference of two elements in a same $[\beta]$
 lies in the kernel of the map
 $$
  \xi_{\ast}\; :\; H_2(X,L;{\Bbb Z})\;
              \longrightarrow\;  H_2(Y,L;{\Bbb Z})\,,
 $$
 each $[\beta]$ determines a class, denoted by $\xi_{\ast}[\beta]$,
 in $H_2(Y,L;{\Bbb Z})$.
For simplicity of notation, we will denote $\xi_{\ast}[\beta]$
 also by $[\beta]$.

\bigskip

\noindent
{\it Comparison 3.2.3 $[$Li-Ruan and Ionel-Parker$]$.}
 Though in different format, it should be noted that
  $(\widehat{W},\widehat{L})/\widehat{B}$-monodromy orbits
  in $H_2(X,L;{\Bbb Z})$
  coincides with $\xi_{\ast}^{-1}(0)$-cosets,
   where\\  $\xi_{\ast}:H_2(X,L;{\Bbb Z})\rightarrow H_2(Y,L;{\Bbb Z})$
    for the moment.
 Thus, the curve class considered here is of the same kind as
  [L-R: Sec.~5] when $L$ is empty.
 Furthermore,
 $\xi_{\ast}^{-1}(0)$ is generated precisely by
   the ``{\it rim tori}" of [I-P1: Sec.~5]
  since the monodromy of all $W[k]/B[k]$ are generated exactly
   by uniform simultaneous Dehn twists over $D$.
 As remarked in ibidem it is with respect to such a collection in
  $H_2(X;{\Bbb Z})$ that one expects to have a
  degeneration-formula/gluing-theorem of Gromov-Witten invariants.
 Thus the combinatorial data we use to restrict the moduli problem
  of maps from bordered Riemann surfaces to fibers of
  $(\widehat{W},\widehat{L})/\widehat{B}$ is the same,
   when $L$ is empty, as those in [L-R], [I-P1], and [I-P2].
 See Appendix for a further comparison of [L-R] versus [I-P1], [I-P2].

\bigskip

\noindent
{\it Comparison 3.2.4 $[$refinement of {\rm [Li1]} and {\rm [Li2]}$]$.}
 In the algebro-geometric setting ([Li1], [Li2]) without $L$,
 one assumes the existence of a relative ample line bundle $H$
   on $W/B$  and
  considers a fixed $H$-degree curve class,
   which in general corresponds to a collection of curve classes
   in $H_2(X;{\Bbb Z})$ (or $A_1(X)$).
 Note that,
  since $H|_{W_b}$, $b\in B$, form a flat family of line bundles
   with base $B$,
  the first Chern class of $H|_X$, and hence the fixed $H$-degree
   class, must be monodromy invariant.
 As the moduli space of maps to fibers of $\widehat{W}/\widehat{B}$
  associated to different monodromy orbits must be disjoint from
  each other,
 Jun Li's degeneration formula in [Li1] and [Li2] indeed always
  splits into a disjoint/independent collection\footnote{The
       moduli stacks involved for different monodromy orbits
        are disjoint from each other.
       They are substacks,
         consisting of disjoint collections of connected components,
        of the moduli stack constructed in [Li1] and are equipped with
        the tangent-obstruction complex  and the virtual fundamental
        class from the restriction of those constructed in [Li2] to
        related connected components.
       See [L-Y1] for an explicit example and discussion.}
  of degeneration formulas,
  one for each monodromy orbit in the fixed $H$-degree curve class.
 Since the discussion in this subsection produces the same
  monodromy on $H_2(X;{\Bbb Z})$ (or $A_1(X)$) as the one
  associated to the Artin stack ${\mathfrak W}/{\mathfrak B}$  of
  expanded degenerations associated to $W/B$, constructed in [Li1],
 the $\widehat{W}/\widehat{B}$-monodromy-orbit refinement of [Li2]
  is the finest refinement of Jun Li's formula
  (and is indeed implicitly already in [Li2],
   had a discussion of monodromy at the level of the stack
   ${\mathfrak W}/{\mathfrak B}$ been made.
  Further, it has to be so
   for [Li1], [Li2] to be consistent with [L-R], [I-P1], [I-P2].
  So this is also a consistency check statement.
  See Comparison 3.2.3 above and Appendix).
 The examples studied in [L-Y1] are [L-Y2]
   are both special cases of such refinement:
  there the ${\mathfrak W}/{\mathfrak B}$-monodromy on
   $H_2(X;{\Bbb Z})$ (or $A_1(X)$) is trivial  and hence
  the degeneration formula of Jun Li refines to one associated
  to each fixed curve class in $H_2(X;{\Bbb Z})$ (or $A_1(X)$).

\bigskip

\subsection{The moduli space
    $\overline{\cal M}_{(g,h), (n, \vec{m})}
                       (W/B,L\,|\,[\beta],\vec{\gamma},\mu)$
    of stable maps to fibers of
    $(\widehat{W},\widehat{L})/\widehat{B}$.}

We now define the moduli space
 $\overline{\cal M}_{(g,h), (n, \vec{m})}
            (W/B,L\,|\,[\beta],\vec{\gamma},\mu)/B$
  of stable maps to fibers of
   $(\widehat{W},\widehat{L})/\widehat{B}$
 and highlight its basic properties.
% This moduli space over $B$ plays a central role in [L-Y4]
%  for a construction of open Gromov-Witten invariants.

\bigskip

\begin{flushleft}
{\bf Moduli space of stable maps to fibers of
      $(\widehat{W},\widehat{L})/\widehat{B}\;$: its topology.}
\end{flushleft}
\noindent
{\bf Definition 3.3.1 [stable map to fibers of $(W[k],L[k])/B[k]$].}
% (cf.\ [Liu(C): Definitions 5.1 - 5.4];
%       [I-P1: Definition 4.1], [I-P2: Lemma 3.3],
%       [L-R: Definition 3.14]; [Li1: Definition 3.1].)
{\rm
 Let
  $[\beta]$ be the $(\widehat{W},\widehat{L})/\widehat{B}$-monodromy
   orbit of $\beta\in H_2(X,L;{\Bbb Z})$,
  $\vec{\gamma}
    =(\gamma_1\,,\,\ldots\,,\gamma_h)\in H_1(L;{\Bbb Z})^{\oplus h}$
  such that $\partial\beta=\gamma_1+\,\cdots\,+\gamma_h$,  and
  $\mu\in{\Bbb Z}$.
 A map $f:(\Sigma,\partial\Sigma)/\pt\rightarrow (W[k], L[k])/B[k]$
  from a bordered Riemann surface $\Sigma$ to a fiber\footnote{When
         the fiber in question is almost-complex isomorphic to a
          $W_{\lambda}$ with $\lambda\ne 0$, the existing definitions
          from Gromov-Witten theory for smooth targets apply.
         Thus, all our focus here is on maps with singular targets.
         Such focus of discussions to singular targets
          prevails the whole manuscript.}
   of $(W[k],L[k])/B[k]$ is called {\it prestable} of
   ({\it combinatorial}) {\it type}
   $((g,h),(n,\vec{m})\,|\,[\beta],\vec{\gamma},\mu)$
  if the following conditions are satisfied:
   \begin{itemize}
    \item[$\cdot$]
     $\Sigma$ is a {\it prestable labelled-bordered Riemann surface}
      of type $((g,h),(n,\vec{m}))$;

    \item[$\cdot$]
     $f$ is {\it continuous}  and
     $\tilde{f}:= \nu\circ f$ is {\it $J$-holomorphic}:
       $J\circ d\tilde{f}= d\tilde{f}\circ j$,
      where $\nu:\tilde{\Sigma}\rightarrow \Sigma$
       is the normalization of $\Sigma$;

    \item[$\cdot$]
     $\tilde{\mathbf p}[k]_{\ast}(f_{\ast}[\Sigma,\partial\Sigma])\;
      \in\; [\beta]\,$;\hspace{1ex}
     $\tilde{\mathbf p}[k]_{\ast}(f_{\ast}[\dot{\partial}\Sigma])\;
        =\; \vec{\gamma}\,$;\hspace{1ex}
     $\mu(f)\;=\;\mu$;

    \item[$\cdot$]
     the automorphism group $\Aut^{\rigid}(f)$ of $f$ as a map to
      (the rigid) $W[k]$ is finite.
   \end{itemize}

 An {\it isomorphism} between two prestable maps
    $f_1:\Sigma_1/\pt\rightarrow W[k]/B[k]$,
    $f_2:\Sigma_2/\pt\rightarrow W[k]/B[k]$
   of the same type
  is a pair $(\alpha,\beta)$\footnote{The use of notation
         $(\alpha,\beta)$ here is so compelling. There should be
         no confusion of this $\beta$ with the curve class $\beta$.
         Similarly, for the occasional use of a {\it map} $g$,
         versus the genus $g$.},
  where
   $\alpha:\Sigma_1\rightarrow \Sigma_2$
     is an isomorphism of prestable labelled-bordered Riemann surfaces
      with marked points and
   $\beta\in{\Bbb G}_m[k]$ acts on $W[k]/B[k]$ as in Sec.~1.1.2
  such that $f_1\circ \beta = f_2\circ\alpha$.
 The isomorphism class associated to a prestable map $f$ will be
  denoted by $[f]$.
 The {\it group of automorphisms} $\Aut(f)$ of a prestable
  $f:\Sigma/\pt\rightarrow W[k]/B[k]$ consists then of
  elements $(\alpha,\beta)\in\Aut(\Sigma)\times {\Bbb G}_m[k]$
  such that $\beta\circ f=f\circ\alpha$.

 A prestable map
    $f:\Sigma/\pt \rightarrow (W[k],L[k])/B[k]$,
    with image in fiber $(W[k]_{\vec{\lambda}},L[k]_{\vec{\lambda}})$,
   is called {\it non-degenerate}
  if no irreducible components of $\Sigma$ are mapped into the
   singular locus $W[k]_{\vec{\lambda}\,,\,\sing}$
   of $W[k]_{\vec{\lambda}}$.
 For $f$ non-degenerate,
  $\Lambda := f^{-1}(W[k]_{\vec{\lambda}\,,\,\sing})$
  consists of interior nodes on $\Sigma$.
 A node $q\in\Lambda$ is called a
  {\it distinguished node} on $\Sigma$ {\it under} $f$.

 Assume that the target fiber
  $W[k]_{\vec{\lambda}}\simeq Y_{[k^{\prime}]}$ for some $k^{\prime}$.
 Decompose a non-degenerate prestable $f$ by
  $$
   f\,=\,\cup_{i=0}^{k^{\prime}}\, f_i\,:\,
    \Sigma\,=\,\cup_{i=0}^{k^{\prime}}\Sigma_{(i)}\;
     \longrightarrow\;
    Y_{[k^{\prime}]}\,=\,\cup_{i=0}^{k^{\prime}}\Delta_{i}
  $$  with
  $f_{(i)}=f_{\Sigma_{(i)}}:
     \Sigma_{(i)}\rightarrow \Delta_i$.
 Recall $D_i := \Delta_i\cap \Delta_{i+1}$.
 Let $\Lambda_i:= f^{-1}(D_i)$ and called it the
  {\it $i$-th subset of distinguished nodes}.
 Associated to $q_{\,\cdot\,}\in \Lambda_i$ are unique
  $q_{\,\cdot\,,1}$ on $\Sigma_{(i)}$
  and $q_{\,\cdot\,,2}$ on $\Sigma_{(i+1)}$.
 From the normal form of $J$-holomorphic map at a point
  ([Ye: Theorem 3.1] and [I-P1: Lemma 3.4]),
  $f_i^{-1}(D_i)$ is a divisor of the form
   $\sum_{q_{ij}\in\Lambda_i} s_{ij,1}\,q_{ij,1}$ on $\Sigma_{(i)}$ and
  $f_{i+1}^{-1}(D_i)$ is a divisor of the form
   $\sum_{q_{ij}\in\Lambda_i} s_{ij,2}\,q_{ij,2}$ on $\Sigma_{(i+1)}$.
 A prestable $f$ is called {\it pre-deformable}
  if it is non-degenerate  and
     $s_{ij,1}=s_{ij,2}$ ($=: s_{ij}$) for all $q_{ij}\in\Lambda_i$,
     $i=0\,,\,\ldots\,,\,k$.
 We call $s_{ij}$ the {\it contact order} of $f$ at $q_{ij}$ along $D_i$.
 Both the non-degeneracy condition and the pre-deformability condition
  are preserved under isomorphisms between prestable maps.

 Finally, a prestable $f:\Sigma/\pt \rightarrow (W[k].L[k])/B[k]$
   is called {\it stable}
  if $f$ is pre-deformable and
     its group $\Aut(f)$ of automorphisms is finite.
 The moduli space of isomorphism classes of stable maps to fibers of
  $(W[k],L[k])/B[k]$  of type
  $((g,h),(n,\vec{m})\,|\,[\beta],\vec{\gamma},\mu)$  is denoted by
  ${\cal M}_{(g,h), (n, \vec{m})}^{\nonrigid}
        ((W[k],L[k])/B[k]\,|\,[\beta],\vec{\gamma},\mu)$.
} % end-definition

\bigskip

We have assumed that the almost-complex structure on $W[k]$ is
 $C^{\infty}$;  thus,
all maps parameterized by
 ${\cal M}_{(g,h), (n, \vec{m})}^{\nonrigid}
        ((W[k],L[k])/B[k]\,|\,[\beta],\vec{\gamma},\mu)$
 are $C^{\infty}$ as well when restricted/lifted to the connected
 components of the normalization of the domains.

For
 $[f:\Sigma/\pt\rightarrow W[k]/B[k]]
   \in {\cal M}_{(g,h), (n, \vec{m})}^{\nonrigid}
         ((W[k],L[k])/B[k]\,|\,[\beta],\vec{\gamma},\mu)$,
 fix a Hermitian metric\footnote{A
          {\it Hermitian metric} on an almost-complex space is a
           (Riemannian) metric so that the almost-complex structure
           is an isometry.
          Different
           choices of such auxiliary metrics on domains and targets
           define the same topology.
          Here for $W[k]$ which is equipped with a compatible pair
           $(J,\omega)$ the metric is chosen to be the one associated
           to the pair $(J, \omega)$.}
  on ${\cal C}/\Def(\Sigma)$ and on $W[k]$.
Define the {\it energy}\footnote{Note that
          % $\footnotesizearea(f)\le E(f)^{1/2}/\sqrt{2}$,
          % where
          %  $\footnotesizearea(f):= \int_{\Sigma}\left(
          %    \|\partial f/\partial x\|^2 \|\partial f/\partial y\|^2
          %   - \langle \partial f/\partial x\,,\, \partial f/\partial y
          %      \rangle^2 \right)^{1/2} d\mu$
          %   is the {\it area} of $f$.
          % For $f$ $J$-holomorphic,
          %  $\footnotesizearea(f)= E(f)^{1/2}/\sqrt{2}$;
          % however, they are ...
          $E(f)$ is conformally invariant with respect to the metric
           on $\Sigma$.
          For $J$ $\omega$-tame and $f$ $J$-holomorphic,
           $E(f)$ coincides with the symplectic area, and is
           determined by
           $[\beta]\subset H_2(X,L;{\footnotesizeBbb Z})$.}
 of $f:\Sigma/\pt \rightarrow W[k]/B[k]$ to be
 $$
  E(f)\;=\; \frac{1}{2}\,\int_{\Sigma}|df|^2\,d\mu\,,
 $$
 where
  $|df|^2$ is the norm-squared of $df$ with respect to the metric
   on $W[k]$ and on $\Sigma$, and
  $d\mu$ is the area-form on $\Sigma$ with respect to the metric
   on $\Sigma$.
Then one can define a topology on
 ${\cal M}_{(g,h), (n, \vec{m})}^{\nonrigid}
         ((W[k],L[k])/B[k]\,|\,[\beta],\vec{\gamma},\mu)$
 similar to [Pa: Sec. 2.1] and [Ye: Definition 0.2];
 see also [Gr2], [P-W], [R-T1], [Sie1]; [Liu(C)]; [I-P1], [L-R].
A point $[f^{\prime}]$ in ${\cal M}_{(g,h), (n, \vec{m})}^{\nonrigid}
              ((W[k],L[k])/B[k]\\  |\,[\beta],\vec{\gamma},\mu)$
 is said to in the $(\varepsilon_1,\varepsilon_2)$-neighborhood
 $U_{\varepsilon_1,\varepsilon_2}([f])$ of $[f]$
 if they have representatives
   $f:\Sigma/\pt \rightarrow W[k]/B[k]$ and
   $f^{\prime}:\Sigma^{\prime}/\pt \rightarrow W[k]/B[k]$
  so that
  \begin{itemize}
   \item[(1)]
   there exists a surjective {\it collapsing/pinching map}
    $c:\Sigma^{\prime}\rightarrow\Sigma$
    that is a diffeomorphism  from
     the complement of a collection of simple loops
      and simple arc with ends on $\partial\Sigma^{\prime}$ on
      $\Sigma^{\prime}$  to
     the complement of the set of nodes on $\Sigma^{\prime}$,  and
     collapses/pinches each simple loop (resp.\ arc) in the collection
      to an interior (resp.\ boundary) node of $\Sigma$
    such that

    \item[$\cdot$]
     \parbox{31ex}{(nearness of domain)}\
      \parbox[t]{22em}{
       $\Sigma^{\prime}$ is isomorphic to a fiber of
         ${\cal C}/\Def(\Sigma)$ with\\
        $\|j-c_{\ast}j \|_{C^{\infty}}<\varepsilon_2$
         on $\Sigma-U_{\varepsilon_1}$  and
        $c(p^{\prime}_{\cdot})$ in the $\varepsilon_2$-neighborhood
         of $p_{\cdot}$,
         where $p_{\cdot}$, $p^{\prime}_{\cdot}$ are
          marked points on $\Sigma$, $\Sigma^{\prime}$ that are paired
          by their label;}

    \item[$\cdot$]
     \parbox{31ex}{(nearness of target and map)}
      $\|f-f^{\prime}\circ c^{-1}\|_{C^{\infty}}<\varepsilon_2$
      on $\Sigma-U_{\varepsilon_1}$,
      as maps to $W[k]$;

   \item[(2)]
   (nearness of energy)\footnote{This
                 condition is redundant here as $E(f)=E(f^{\prime})$
                 currently. We reserve it here to stress its importance
                 to Compactness Theorem.}\hspace{9ex}
    $|E(f)-E(f^{\prime})|<\varepsilon_2\,$.
  \end{itemize}
 Here,
  $U_{\varepsilon_1}$ is the $\varepsilon_1$-neighborhood of
   the set of nodes of $\Sigma$
   that is small enough so that it contains no marked points.
The system  $\{\,U_{\varepsilon_1,\varepsilon_2}([f])\,\}
                       _{f;\,\varepsilon_1,\,\varepsilon_2}$
 of subsets generates the {\it $C^{\infty}$-topology}\footnote{Let
                 ${\cal M}$ be the moduli space of pre-deformable stable
                  maps to fibers of
                  $((W[k],L[k])/B[k])^{\,\mbox{\tiny\rm rigid}}$
                  with the $C^{\infty}$-topology from [Ye: Definition 0.2].
                 Then ${\footnotesizeBbb G}_m[k]$ acts on ${\cal M}$,  and
                 our moduli space
                  ${\cal M}_{\,\cdot\,}^{\,\mbox{\tiny\rm non-rigid}}
                    ((W[k],L[k])/B[k]\,|\,\cdots\,)$
                  is contained in ${\cal M}/{\footnotesizeBbb G}_m[k]$
                  with the quotient topology.
                 The induced subset-topology on
                  ${\cal M}_{\,\cdot\,}^{\,\mbox{\tiny\rm non-rigid}}
                    ((W[k],L[k])/B[k]\,|\,\cdots\,)$
                  coincides with its $C^{\infty}$-topology.}
 on
 ${\cal M}_{(g,h), (n, \vec{m})}^{\nonrigid}
           ((W[k],L[k])/B[k]\,|\,[\beta],\vec{\gamma},\mu)$.

The pseudo-embedding
 $\varphi_{k^{\prime},k; I}:
  (W[k^{\prime}],L[k])/B[k^{\prime}] \hookrightarrow (W[k],L[k])/B[k]$,
  $k^{\prime}<k$, $I\subset\{0\,,\,\ldots\,,\,k\}$,
 from Sec.~1.1.3 induces a {\it pseudo-embedding}
 $$
  \begin{array}{r}
   \varphi_{k^{\prime},k; I}\;  :\;
    {\cal M}_{(g,h), (n, \vec{m})}^{\nonrigid}
             ((W[k^{\prime}],L[k^{\prime}])/B[k^{\prime}]\,
                                |\,[\beta],\vec{\gamma},\mu)
                                   \hspace{6em}\\[1ex]
   \hookrightarrow\;
      {\cal M}_{(g,h), (n, \vec{m})}^{\nonrigid}
            ((W[k],L[k])/B[k]\,|\,[\beta],\vec{\gamma},\mu)\,.
  \end{array}
 $$
Define the set of isomorphism classes of stable maps to fibers of
 $(\widehat{W},\widehat{L})/\widehat{B}\,$:
 $$
  \overline{\cal M}_{(g,h), (n, \vec{m})}
              (W/B,L\,|\,[\beta],\vec{\gamma},\mu)\;
  :=\;
   \left.
     \left( \amalg_{k=0}^{\infty}
      {\cal M}_{(g,h), (n, \vec{m})}^{\nonrigid}
         ((W[k],L[k])/B[k]\,|\,[\beta],\vec{\gamma},\mu) \right)
   \right/ \raisebox{-.6ex}{$\sim$}\,,
 $$
 where the equivalence relation $\sim$ is generated by
  $[f]\sim \varphi_{k^{\prime},k;I}([f^{\prime}])$
  for
   $[f]  \in
    {\cal M}_{(g,h), (n, \vec{m})}^{\nonrigid}
     (\\ (W[k],L[k])/B[k]\,|\,[\beta],\vec{\gamma},\mu)$  and
   $[f^{\prime}]  \in$
    the defining domain of $\varphi_{k^{\prime},k;I}$  on
     ${\cal M}_{(g,h), (n, \vec{m})}^{\nonrigid}
     ((W[k^{\prime}],\\  L[k^{\prime}])/B[k^{\prime}]\,
                            |\,[\beta],\vec{\gamma},\mu)$.
By construction, there are embeddings of sets
 $$
  \varphi_{(k)}\,:\,
      {\cal M}_{(g,h), (n, \vec{m})}^{\nonrigid}
         ((W[k],L[k])/B[k]\,|\,[\beta],\vec{\gamma},\mu)\,
     \hookrightarrow\,
      \overline{\cal M}_{(g,h), (n, \vec{m})}
               (W/B,L\,|\,[\beta],\vec{\gamma},\mu)\,, \hspace{1em}
  k\in{\Bbb Z}_{\ge 0}\,.
 $$
A subset $U$ of
  $\overline{\cal M}_{(g,h), (n, \vec{m})}
             (W/B,L\,|\,[\beta],\vec{\gamma},\mu)$
  is said to be {\it open}
 if $U=\cup_{\alpha}U_{\alpha}$ such that
     $U_{\alpha}$ is contained in the image of some $\varphi_{(k)}$
      and
     $\varphi_{(k)}^{-1}(U_{\alpha})$ is open  in
     ${\cal M}_{(g,h), (n, \vec{m})}^{\nonrigid}
         ((W[k],L[k])/B[k]\,|\,[\beta],\\  \vec{\gamma},\mu)$.
This defines the {\it $C^{\infty}$-topology} on the
 moduli space\footnote{We
             could have used the notation
             ${\cal M}_{(h,h),(n,\vec{m})}(
                   (\widehat{W},\widehat{L})/\widehat{B})\,
                                  |\,[\beta],\vec{\gamma},\mu)$
             for the moduli space
             $\overline{\cal M}_{(g,h), (n, \vec{m})}
                         (W/B,L\,|\,[\beta],\vec{\gamma},\mu)$.
             Our choice of the latter reflects the intention to keep
             in mind that maps to singular fibers are meant to be limited
             to those that are approachable from maps to smooth fibers,
             (reflected, e.g.\ by the pre-deformability condition).
             As we will show that this is indeed so at the level of
             Kuranishi/virtual neighborhoods on the moduli space.}
  $\overline{\cal M}_{(g,h), (n, \vec{m})}
          (W/B,L\,|\,[\beta],\vec{\gamma},\mu)$
 of stable maps to fibers of $\widehat{W}/\widehat{B}$.
By construction,
 $\overline{\cal M}_{(g,h), (n, \vec{m})}
         (W/B,L\,|\,[\beta],\vec{\gamma},\mu)$.
 fibers naturally over $B$;  in notation
  $\overline{\cal M}_{(g,h), (n, \vec{m})}
            (W/B,L\,|\,[\beta],\vec{\gamma},\mu)/B$.

\bigskip

\noindent
{\bf Definition 3.3.2 [tautological cover].} {\rm
 By construction,
  $$
   \left\{\, {\cal M}_{(g,h), (n, \vec{m})}^{\nonrigid}
                ((W[k],L[k])/B[k]\,|\,[\beta],\vec{\gamma},\mu)\,
    \right\}_{k\in{\scriptsizeBbb Z}_{\ge 0}}
  $$
  is an open cover of
  $\overline{\cal M}_{(g,h), (n, \vec{m})}
           (W/B,L\,|\,[\beta],\vec{\gamma},\mu)$.
 We will call it the {\it tautological cover} of
  $\overline{\cal M}_{(g,h), (n, \vec{m})}
                     (W/B,L\,|\,[\beta],\vec{\gamma},\mu)$.
} % end-definition

\bigskip

Indeed, there exists $k_0$
             depending $(W/B,L)$ and
                       $((g,h),(n,\vec{m})|[\beta],\vec{\gamma},\mu)$
 such that

 \vspace{-1ex}
 {\scriptsize
 $$
  \begin{array}{l}
  {\cal M}_{(g,h), (n, \vec{m})}^{\tinynonrigid}
     ((W[k_0],L[k_0])/B[k_0]\,|\,[\beta],\vec{\gamma},\mu)\;
   \supset\;
    {\cal M}_{(g,h), (n, \vec{m})}^{\tinynonrigid}
    ((W[k_0+1],L[k_0+1])/B[k_0+1]\,|\,[\beta],\vec{\gamma},\mu)
                                                        \\[1.6ex]
  \hspace{2em} \supset\;
    {\cal M}_{(g,h), (n, \vec{m})}^{\tinynonrigid}
    ((W[k_0+2],L[k_0+2])/B[k_0+2]\,|\,[\beta],\vec{\gamma},\mu)\;
   \supset\;\cdots\,.
  \end{array}
 $$ } % end-scriptsize

\vspace{-1ex}
\noindent
Thus, the tautological cover of
  $\overline{\cal M}_{(g,h), (n, \vec{m})}
                     (W/B,L\,|\,[\beta],\vec{\gamma},\mu)$
 is finite in effect, cf.~Theorem 3.3.8. % Theorem [compactness$/B$]
The universal maps on the universal curve over each
 ${\cal M}_{(g,h), (n, \vec{m})}^{\nonrigid}
   ((W[k],L[k])/B[k]\,|\,[\beta],\vec{\gamma},\mu)$ are glued
 to give the universal map (between spaces with charts)
  $$
   F\; :\; {\cal C}
    /\overline{\cal M}_{(g,h), (n, \vec{m})}
                        (W/B,L\,|\,[\beta],\vec{\gamma},\mu)\;
   \longrightarrow\; (\widehat{W},\widehat{L})/\widehat{B}\,.
  $$

\bigskip

\noindent
{\it Remark 3.3.3 $[$on Definition 3.3.1$]$.}
For the meaning/reason of the various conditions in Definition 3.3.1:
     % Definition [stable map to fibers of $(W[k],L[k])/B[k]$]
[Liu (C): Lemma 6.13],
  which is generalized to Lemma 5.3.1.1 in Sec.~5.3.1,
    % Lemma [index of $D_f\bar{\partial}_J$ for rigid target]
 explains the role of {\it Maslov index} $\mu$
 on infinitesimal deformations of an open stable map;
[L-R: Lemma 3.11 (3)], [I-P1: Lemma 3.3], and [Li1: Proposition 2.2]
  give the reason to the important {\it pre-deformability condition},
  as we want to single out maps that contribute to
  the degeneration formula;
[I-P1: Sec.~6, Step 3] explains why morphisms of maps in question
 are defined so that the singular targets become {\it non-rigid}
 on the ruled-manifold-components from expansion,
  as it has to so that the choice of complex-scaling
  renormalizations in ``stretching/pulling out" a degenerate component
  that falls into $W[k]_{\lambda\,,\,\sing}$ becomes irrelevant.
Furthermore, we will see in Sec.~5.3.5 that
 it is the {\it combination of all three} that renders
 the moduli space
 $\overline{\cal M}_{(g,h), (n, \vec{m})}
              (W/B,L\,|\,[\beta],\vec{\gamma},\mu)/B$
  ``{\it virtually flat}" over $B$.
Only so can one hope for a degeneration formula.

\bigskip

We now highlight three basic properties of
 $\overline{\cal M}_{(g,h), (n, \vec{m})}
              (W/B,L\,|\,[\beta],\vec{\gamma},\mu)/B$
 in parallel to [L-R: Sec.3.3] and [I-P1: Theorem 7.4]
(and to the existing literature quoted earlier on non-family case
 as well).

\bigskip

\begin{flushleft}
{\bf Hausdorffness.}
\end{flushleft}
Let
 ${\cal M}_{(g,h), (n, \vec{m})}^{\rigid}
         ((W[k],L[k])/B[k]\,|\,[\beta],\vec{\gamma},\mu)$
 be the moduli space of isomorphism classes of stable maps to
  fibers of {\it rigid} $(W[k],L[k])/B[k]$  of type
  $((g,h),(n,\vec{m})\,|\,[\beta],\vec{\gamma},\mu)$.
This is defined the same as in Definition 3.3.1
    % Definition [stable map to fibers of $(W[k],L[k])/B[k]$]
 except that
  a morphism between $f_1:\Sigma_1/\pt \rightarrow W[k]/B[k]$  and
   $f_2:\Sigma_2/\pt \rightarrow W[k]/B[k]$ is taken to be
   an isomorphism $\alpha:\Sigma_1\rightarrow \Sigma_2$ such that
   $f_1=f_2 \circ \alpha$,  and
  the stability condition for $f$ to the rigid $W[k]/B[k]$
   is that $\Aut^{\rigid}(f)$ is finite.
Then [Sie1: proof of Proposition 3.8] (see also [F-O: Lemma 10.4])
 can be applied to show that
 ${\cal M}_{(g,h), (n, \vec{m})}^{\rigid}
         ((W[k],L[k])/B[k]\,|\,[\beta],\vec{\gamma},\mu)$
 is Hausdorff.
This space is indeed a singular subspace of a manifold and hence is
 metrizable.
The moduli space
 ${\cal M}_{(g,h), (n, \vec{m})}^{\nonrigid}
       ((W[k],L[k])/B[k]\,|\,[\beta],\vec{\gamma},\mu)$
 is the quotient space of
 ${\cal M}_{(g,h), (n, \vec{m})}^{\rigid}
        ((W[k],L[k])/B[k]\,|\,[\beta],\vec{\gamma},\mu)$
 by the ${\Bbb G}_m[k]$-action.
Due to the stability condition, all the ${\Bbb G}_m[k]$-orbits on
 ${\cal M}_{(g,h), (n, \vec{m})}^{\rigid}
         ((W[k],L[k])/B[k]\,|\,[\beta],\vec{\gamma},\mu)$
 have the same (real) dimension $2k$.
This implies that
 ${\cal M}_{(g,h), (n, \vec{m})}^{\nonrigid}
        ((W[k],L[k])/B[k]\,|\,[\beta],\vec{\gamma},\mu)$
 is also Hausdorff, for $k\in {\Bbb Z}_{\ge 0}$.

Given now
 $[f]$, $[f^{\prime}]\in \overline{\cal M}_{(g,h), (n, \vec{m})}
           (W/B,L\,|\,[\beta],\vec{\gamma},\mu)$,
 assume, without loss of generality, that
  the image fiber of $f$ (resp.\ $f^{\prime}$) is
   $Y_{[k]}$ or $W_{\lambda}$, $\lambda\ne 0$
  (resp.\ $Y_{[k^{\prime}]}$ or $W_{\lambda^{\prime}}$)
  with $k\ge k^{\prime}$.
Then
 $[f]$, $[f^{\prime}]\in
  {\cal M}_{(g,h), (n, \vec{m})}^{\nonrigid}
        ((W[k],L[k])/B[k]\,|\,[\beta],\vec{\gamma},\mu)$.
As
 ${\cal M}_{(g,h), (n, \vec{m})}^{\nonrigid}
        ((W[k],L[k])/B[k]\,|\,[\beta],\vec{\gamma},\mu)$
 embeds in
 $\overline{\cal M}_{(g,h), (n, \vec{m})}
               (W/B,L\,|\,[\beta],\vec{\gamma},\mu)/B$,
this implies,
  by the way we define
   the $C^{\infty}$-topology on
   $\overline{\cal M}_{(g,h), (n, \vec{m})}
                 (W/B,L\,|\,[\beta],\vec{\gamma},\mu)/B$,
 that there are disjoint open subsets in
 $\overline{\cal M}_{(g,h), (n, \vec{m})}
               (W/B,L\,|\,[\beta],\\  \vec{\gamma},\mu)/B$
 that separate $[f]$ and $[f^{\prime}]$.
It follows that:

\bigskip

\noindent
{\bf Proposition 3.3.4 [Hausdorffness].} {\it
 $\overline{\cal M}_{(g,h), (n, \vec{m})}
              (W/B,L\,|\,[\beta],\vec{\gamma},\mu)/B$
 with the $C^{\infty}$-topology is Hausdorff.
} % end-proposition

\bigskip

This proposition can be regarded as a corollary of the
 stability condition on maps, in much the same reason as in
 geometric-invariant-theory quotients in algebraic geometry.

\bigskip

\begin{flushleft}
{\bf Finite stratification.}
\end{flushleft}
%
%%%%%%%%%%%%%%%%%%%%%%%%%%%%%%%%%%%%%%%%%%%%%%%%%%%%%%%%%
% \noindent
% {\bf Definition/Lemma ??? [area and energy].}
%  ([MD-S: Sec.\ 4.1]; see also [L-R: Sec.~3] and [Liu(C): Sec.~5.2];
%   [Ye: Sec.~4] has a slightly different convention on the $1/2$ factor.)
% {\rm
%  Let $f:\Sigma\rightarrow (Y,L)$ be a $C^{\infty}$ map.
%  Define the {\it area} and the {\it energy} of $f$ to be
%   %
%   \marginpar{\raggedright\tiny $\bullet$
%    Find the precise definition of $|df|_J$.}
%   %
%   $$
%    a(f)\;=\; \int_{\Sigma}f^{\ast}\omega \;
%         (=\; \omega(f_{\ast}[\Sigma]))
%     \hspace{2em}\mbox{and}\hspace{2em}
%    E(f)\; =\; \frac{1}{2}\int_{\Sigma}|df|_J^2\,dA
%   $$
%   respectively.
%  {\it
%   The norm $|df|_J$ and the area form $dA$ depend on a choice of
%    a metric on $\Sigma$ but
%     $|df|_J^2\,dA$ doesn't.
%    When $J$ is $\omega$-tame and $f$ is $J$-holomorphic, $a(f)=E(f)$.}
% } % end-definition/lemma
%%%%%%%%%%%%%%%%%%%%%%%%%%%%%%%%%%%%%%%%%%%%%%%%%%%%%
%
% \bigskip
%
We first generalize a simplified version of the constructions/operations
  of [B-M: Sec.~1]
 to incorporate both the boundary of bordered Riemann surfaces and
 the consideration in [I-P1: Sec.~7].
This defines a category ${\mathfrak G}$ of graphs\footnote{In
                  this simplified presentation, we directly identify
                   a (plain) graph with its geometric realization, i.e.\
                    a simplicial $1$-complex consisting of
                     a finite collection of points (i.e.\ {\it vertices});
                     a finite collection of (un-oriented) line segments
                      (i.e.\ {\it edges}),
                       with both ends attached to vertices;
                     a finite collection of (un-oriented) line segments
                      (i.e.\ {\it legs}, {\it hands} or {\it roots})
                      with only one end attached to vertices;
                     a finite collection of (un-oriented) line segments
                      (i.e.\ {\it bridges}) with both ends attached
                       to free ends of hands, and
                     a finite collection of (un-oriented) line segments
                      (i.e.\ {\it fingers}) with only one end attached
                      to free end of hands.
                  We will denote a graph by $\tau$ (not to be confused
                   with the involution $\tau$ in Definition 2.1).
                        % Definition
                        % [prestable labelled-bordered Riemann surface]
                  The full formal language in [B-M: Sec.~1] can be
                   recovered whenever needed.}
 whose objects label the topological types of stable maps to fibers of
 $(\widehat{W},\widehat{L})/\widehat{B}$.

\bigskip

\noindent
{\bf Definition 3.3.5 [weighted layered $(A_2\rightarrow A_1)$-graph].}
{\rm
 Let $A_2\rightarrow A_1$ be a pair of abelian groups with a morphism.
 A {\it weighted layered $(A_2\rightarrow A_1)$-graph} $\tau$ consists
  of the following data:
  \begin{itemize}
   \item[(1)] \parbox{48ex}{({\it graph with
                                  hands, bridges, legs, and fingers})}
    a {\it graph} $\tau$,
     whose set of vertices, edges, legs, hands, bridges, and fingers
      are denoted by
      $V(\tau)$, $E(\tau)$, $L(\tau)$, $H(\tau)$, $B(\tau)$, and $F(\tau)$
      respectively;
    among them the sets $H(\tau)$, $F(\tau)$, $L(\tau)$ are ordered,
      with $L(\tau)$ also bi-colored by $(\blue,\red)$;

   \item[$\cdot$]
    the gluing of hands to vertices
      (resp.\ edges to vertices, bridges to hands,
              legs to vertices, fingers to hands)
     defines the {\it attaching map} $H(\tau)\rightarrow V(\tau)$
      (resp.\ $E(\tau)\rightarrow \Sym^2(V(\tau))$,
              $B(\tau)\rightarrow \Sym^2(H(\tau))$,
              $L(\tau)\rightarrow V(\tau)$,
              $F(\tau)\rightarrow H(\tau)$,
       where $\Sym^2(\,\cdot\,)$ is the symmetric product of $\,\cdot\,$);
    the attaching map $F(\tau)\rightarrow H(\tau)$,
     together with the ordering on the sets $H(\tau)$ and $F(\tau)$,
     groups elements of $F(\tau)$ into a tuple of tuples;

   \item[(2)] \parbox{19ex}{({\it layer structure})}
    a map
     $\layer:V(\tau)\rightarrow \{\,0\,,\,\cdots\,,\,k+1\,\}$,
     for a $k\in{\Bbb Z}_{\ge 0}$, such that
      $\Image(\layer)$ is either $\{0\}$, $\{k+1\}$, or
       the whole $\{\,0\,,\,\cdots\,,\,k+1\,\}$  and that,
      if $v_1$, $v_2\in V(\tau)$ is connected by an edge
      $e\in E(\tau)$,
     then
       either
       $\layer(v_1)=\layer(v_2)$,
        in which case we call $e$ an {\it ordinary edge}, or
       $|\layer(v_1)-\layer(v_2)|=1$, in which case we call $e$ a
        {\it distinguished edge};
    the set of ordinary (resp.\ distinguished) edges
     is denoted by $E^o(\tau)$ (resp.\ $E^{\dag}(\tau))$;
    by definition, $E(\tau)=E^o(\tau)\amalg E^{\dag}(\tau)$;

   \item[$\cdot$]
    we require that a hand can be attached only to a vertex
      in $\layer^{-1}(\{0\,,\,k+1\})$  and
     a red leg can be attached only to a vertex
      to which there is a hand attached;

   \item[(3)] \parbox{19ex}{({\it weight functions})}
    $$
     \begin{array}{ll}
      g\;:\; V(\tau)\; \longrightarrow\; {\Bbb Z}_{\ge 0}\,;    \\[.6ex]
      b : V(\tau)\rightarrow A_2\,,\;
      \gamma : H(\tau)\rightarrow A_1
       & \mbox{such that the morphism $A_2\rightarrow A_1$}\\
      &  \mbox{takes $\sum_{v\in V(\tau)}b(v)$ to
                     $\sum_{h\in H(\tau)}\gamma(h)$}\,; \\[.6ex]
      \order: E^{\dag}(\tau)\rightarrow {\Bbb Z}_{\ge 1}\,;
    \end{array}
    $$

   \item[(4)]
    an assignment $\tau\mapsto {\mu}(\tau)\in{\Bbb Z}$,
     called the {\it index} of $\tau$.
  \end{itemize}

  An {\it isomorphism} $\alpha: \tau_1\rightarrow \tau_2$
   between two weighted layered $(A_2\rightarrow A_1)$-graphs
   is an isotopy class of isomorphisms $\tau_1\rightarrow \tau_2$
     as a simplicial complex
    that induces  isomorphisms of sets, ordered sets, or
         bi-colored ordered sets whichever applicable:
       $V(\tau_1)\stackrel{\sim}{\rightarrow}V(\tau_2)$,
       $H(\tau_1)\stackrel{\sim}{\rightarrow}H(\tau_2)$,
       $E(\tau_1)\stackrel{\sim}{\rightarrow}E(\tau_2)$,
       $B(\tau_1)\stackrel{\sim}{\rightarrow}B(\tau_2)$,
       $L(\tau_1)\stackrel{\sim}{\rightarrow}L(\tau_2)$,
       $F(\tau_1)\stackrel{\sim}{\rightarrow}F(\tau_2)$   and
    that preserves
     the layer $\layer(\,\cdot\,)$,
     weights $g(\,\cdot\,)$,
             $b(\,\cdot\,)$, $\gamma(\,\cdot\,)$, $\order(\,\cdot\,)$,
      and
     the index $\mu(\,\cdot\,)$.

  Denote by ${\mathfrak G}(A_2\rightarrow A_1)$
   (or simply ${\mathfrak G}$ when $A_2\rightarrow A_1$ is understood)
   the category
    whose objects are weighted layered $(A_2\rightarrow A_1)$-graphs  and
    whose morphisms are given by isomorphisms.
} % end-definition

\bigskip

Define the {\it core} $\tau^0$ of
 a weighted layered $(A_2\rightarrow A_1)$-graph to be the
 (weighted layered) sub-graph of $\tau$
 by removing the hands, bridges, legs, and fingers from $\tau$.
For a connected weighted layered $(A_2\rightarrow A_1)$-graph $\tau$,
define the
 {\it genus} of $\tau$ to be
  $$
   g(\tau)\;=\; 1-\chi(\tau^0) + \sum_{v\in V(\tau)}\,g(v)
  $$
   and the
 {\it $b$-weight} $b(\tau)$ of $\tau$ to be
  $$
   b(\tau)\;=\; \sum_{v\in V(\tau)}\,b(v)\,.
  $$
For general $\tau$, define its genus and $b$-weight
 by summing genus and $b$-weight over its connected components.
Let $\tau_1$, $\tau_2$ be connected weighted layered
 $(A_2\rightarrow A_1)$-graphs of the same index.
A {\it contraction} from $\tau_1$ to $\tau_2$ is a homotopy class
 of surjective simplicial pseudo-maps $c:\tau_1\rightarrow \tau_2$
 such that
 \begin{itemize}
  \item[$\cdot$]
   the defining domain of $c$ contains $\tau_1-B(\tau_1)$;

  \item[$\cdot$]
   let $\layer:V(\tau_1)\rightarrow \{0\,,\,\ldots\,,\,k_1+1\}$ and
    $\layer:V(\tau_2)\rightarrow \{0\,,\,\ldots\,,\,k_2+1\}$
    be the layer structure of $\tau_1$ and $\tau_2$ respectively;
   then
    $k_1\ge k_2$  and
    there exists a non-decreasing map
     $I:\{0\,,\,\ldots\,,\,k_1+1\}\rightarrow \{0\,,\,\ldots\,,\,k_2+1\}$
     such that $I \circ \layer(v)=\layer (c(v))$ for all $v\in V(\tau_1)$;

  \item[$\cdot$]
   $c$ is a deformation retract on its defining domain;
   the induced maps from $H(\tau_1)$ to $H(\tau_2)$,
    $L(\tau_1)$ to $L(\tau_2)$, and $F(\tau_1)$ to $F(\tau_2)$
    are bijective;

  \item[$\cdot$]
   let $v\in V(\tau_2)$,
   then $c^{-1}(v)$ is connected and
    $g(v)=g(c^{-1}(v))$, $b(v)=b(c^{-1}(v))$;

  \item [$\cdot$]
   if $e\in E^{\dag}(\tau_1)$ is not mapped to a vertex of $\tau_2$
    then $c(e)\in E^{\dag}(\tau_2)$ and $\order(e)=\order(c(e))$.
 \end{itemize}
A ({\it red-to-blue}) {\it color change}
 $\rb:\tau_1\rightarrow \tau_2$  is a change of the color of
 some red legs to blue, leaving everything else the same.
Both contractions and color-changes preserve $g$ and $b$-weight
 of weighted layered $(A_2\rightarrow A_1)$-graphs.

\bigskip

Associated to a point
 $[f:\Sigma/\pt \rightarrow (\widehat{W},\widehat{L})/\widehat{B}]
  \in   \overline{\cal M}_{(g,h), (n, \vec{m})}
                 (W/B,L\,|\,[\beta],\vec{\gamma},\mu)$,
  with target isomorphic to $(Y_{[k]},L)$,
 is a weighted layered graph $\tau_{[f]}$
 via the following correspondence

{\footnotesize
\noindent\hspace{1em}
\begin{tabular}{lcl}
  & & \\
 \hspace{2em}
 $f:\Sigma \rightarrow (Y_{[k]}, L)$
  && $(H_2(Y,L;{\footnotesizeBbb Z})\stackrel{\partial}{\rightarrow}
       H_1(L;{\footnotesizeBbb Z}))$-graph $\tau$ \\[.6ex]
 \hline  \\[-1.6ex]
 {\it irreducible component} $\Sigma_v$ of $\Sigma$
   && {\it vertex} $v \in V(\tau)$ \\
 {\it labelled boundary component} $(\partial\Sigma_v)_h$ of $\Sigma_v$
   &&  {\it hand} $h\in H(\tau)$ attached to $v$ \\
 \hspace{1em}
 (including {\it boundary node} $q_h$ of {\it type E})\\
 {\it ordinary interior node} $q$ connecting $\Sigma_{v_1}$, $\Sigma_{v_2}$
   &&  {\it ordinary edge} $e_q$ with ends attached to $(v_1, v_2)$ \\
 {\it distinguished node} $q$ connecting $\Sigma_{v_1}$, $\Sigma_{v_2}$
   && {\it distinguished edge} $e_q$ with ends attached to $(v_1,v_2)$ \\
 {\it boundary node} $q$ of {\it type H} connecting
    $(\partial\Sigma)_{h_1}$  % and $(\partial\Sigma)_{h_2}$
   &&  {\it bridge} $b_q$ attached to the free ends of $(h_1,h_2)$ \\
 \hspace{1em}and $(\partial\Sigma)_{h_2}$ \\
 {\it free marked point} $p$ on $\Sigma_v$
   &&  {\it leg} $l_p$ attached to vertex $v$     \\
  \hspace{1em}{\it interior} marked point
   && \hspace{1em}{\it blue} leg\\
  \hspace{1em}{\it boundary} free marked point
   && \hspace{1em}{\it red} leg\\
 {\it boundary marked point} $p\in (\partial\Sigma)_h$
   &&  {\it finger} $f_p$ attached to the free end of hand $h$ \\
 \hspace{2em}------------------------------------
   && \hspace{2em}------------------------------------------------   \\
 $\Sigma_v$ such that $f(\Sigma_v)\subset \Delta_i$
   && $\layer(v)=i$, $v\in V(\tau)$   \\
 \hspace{2em}------------------------------------
   && \hspace{2em}------------------------------------------------   \\
 $g(\Sigma_v)$
   &&  $g(v)$, $v\in V(\tau)$   \\
 $f_{\ast}[\Sigma_v]$
   &&  $b(v)$, $v\in V(\tau)$   \\
 $f_{\ast}[(\partial\Sigma_v)_h]$ if $\partial\Sigma_v\ne \emptyset$
   &&  $\gamma(h)$, $h\in H(\tau)$  \\
 Maslov index $\mu(f)$
   &&  $\mu(\tau)$ \\
 distinguished node $q$ of {\it contact order} $s$
   && $\order(e_q)=s$, $e_q\in E^{\dag}(\tau)$ \\
 && \\
\end{tabular}
} % end-footnotesize

 \noindent
 where it is understood that, when $\Sigma_{v_1}=\Sigma_{v_2}$, $v_1=v_2$.
It is clear that $\tau$ is defined to the isomorphism class $[f]$ of $f$.
We call $\tau_{[f]}$ the {\it dual} (weighted layered) {\it graph}
 of $f$ or $[f]$.
Two stable maps
  $f_1:\Sigma_1/\pt \rightarrow W[k_1]/B[k_1]$,
  $f_2:\Sigma_2/\pt \rightarrow W[k_2]/B[k_2]$
  are said to be {\it of the same topological type}
 if $\tau_{[f_1]}$ is isomorphic to $\tau_{[f_2]}$ in
  the category ${\mathfrak G}$.
Degenerations of stable maps to fibers of
 $(\widehat{W},\widehat{L})/\widehat{B}$ are reflected contravariantly
 by compositions of contractions and color-changes of their dual graphs.
 %
 % %
 % \begin{itemize}
 %  \item[$\cdot$]
 %   ????????????.
 % \end{itemize}
 % %
 % This defines a {\it partial order} on the set of $???$-graphs
 %  by setting
 %  $\tau^{\prime} \preceq \tau$ if $??????$.

\bigskip

The following fundamental lemma on $J$-holomorphic maps to fibers of
 $(W[k],L[k])/B[k]$ is a consequence of
  [Ye: Lemma 4.1, Lemma 4.3, Lemma 4.5] and
  [I-P2: the explicit construction in Sec.~2],
 (see also
  [Gr2]; [MD-S1: Lemma 4.5.2], [F-O: Lemma 8.1], [Pa: Proposition 3.1.3],
  [P-W]; and [I-P1: Lemma 1.5], [L-R: Lemma 3.8 and Lemma 3.9]):

\bigskip

\noindent
{\bf Lemma 3.3.6 [energy lower bound].} {\it
 One can fix Hermitian metrics on $W[k]$, $k\in{\Bbb Z}_{\ge 0}$,
 so that there exists a $\delta_0 >0$ that depends only on
  $(X,J,\omega)$ such that, for all $k\in {\Bbb Z}_{\ge 0}$,
  \begin{itemize}
   \item[$\cdot$]
    any non-constant $J$-holomorphic map
     $f:\Sigma/\pt \rightarrow (W[k],L[k])/B[k]$
     has $\,E(f)\;\ge\;\delta_0\,;$

   \item[$\cdot$]
    for any sequence $f_i:\Sigma/\pt \rightarrow (W[k],L[k])/B[k]$
     of $J$-holomorphic maps on $\Sigma$ and any
     blow-up point\footnote{This
                     is a point on $\Sigma$ to which a positive energy
                      of $f_i$ condenses/accumulates in the limit.
                     This is where a bubbling occurs.
                     See, e.g.\ [MD-S1: Lemma 4.5.5], [P-W], [Ye: Sec.~4],
                                and [L-R: Sec.~3.2].}
     $z\in \Sigma$,
     $$
      \lim_{r\rightarrow 0}\, \limsup_{i\rightarrow\infty}\,
       E(f|_{B_r(z)})\; \ge\; \delta_0\,.
     $$
  \end{itemize}
} % end-lemma

The following lemma is parallel to [L-R: Lemma 3.15].
It follows from Lemma 3.3.6.
  % Lemma [energy lower bound]

\bigskip

\noindent
{\bf Lemma 3.3.7 [finite stratification].} {\it
 The classification of stable maps by their topological types
  gives rise to a finite stratification of
  $\overline{\cal M}_{(g,h), (n, \vec{m})}
                      (W/B,L\,|\,[\beta],\vec{\gamma},\mu)/B$,
  with each stratum $S_{\tau}$ labelled by
  a weighted layered $(H_2(Y,L;{\Bbb Z}), H_1(L;{\Bbb Z}))$-graph
  $\tau\in{\mathfrak G}$.
 % The closure $\overline{S_{\tau_1}}$ intersects $S_{\tau_2}$
 %  if and only if $??????$.
} % end-lemma

\bigskip

\begin{flushleft}
{\bf Compactness.}
\end{flushleft}
The following fundamental compactness result of Gromov-Witten
 theory
 in the current contents is closely related to
 [L-R: Theorem 3.16, Corollary 3.17, Theorem 3.20, Theorem 3.21],
 and [I-P1: Theorem 7.4].
It follows from
  Lemma 3.3.6, Lemma 3.3.7,
   % Lemma [energy lower bound],  Lemma [finite stratification]
  the compactness technique/results in [Ye],
  the compactness techniques/results in [L-R: Sec. 3.2],  and
  the compactness technique/results in
   [I-P1, particularly Sec.~6, Step 3] and [I-P2],
 as the effect around the boundary of domains that is mapped to $L$
     is taken care of in [Ye],
    the effect for degeneration of domains due around the degeneration
    of the neck regions of targets is taken care of in [I-P1],
    and these two regions are disjoint from each other in our situation.
See also [Gr2], [F-O], [I-S1], [Pa], [P-W], and [R-T1]
 for the non-family case.

%%%%%%%%%%%%%%%%%%%%%%%%%%%%%%%%%%%%%%%%%%%%%%%%%%%%%%%%%%%%%%%%%%%
%
% \marginpar{\raggedright\tiny $\bullet$
%  State a lemma that degeneracy for maps to target $\Delta_i$
%   can be arranged so that
%   all the degeneracy occur in $D_{i,0}=D_i$ and
%   none in $D_{i,\infty}=D_{i-1}$.
%  This is needed to justify the way we define the $C^{\infty}$-topology
%   via the tail-collapsing maps $\varphi_I$.
%  If this is not true, one needs to introduce more kinds of $\varphi_I$
%   that distinguished the degeneracy types. }
%%%%%%%%%%%%%%%%%%%%%%%%%%%%%%%%%%%%%%%%%%%%%%%%%%%%%%%%%%%%%%%%%%%

\bigskip

\noindent
{\bf Theorem 3.3.8 [compactness$/B$].} {\it
 The moduli space
  $\overline{\cal M}_{(g,h), (n, \vec{m})}
      (W/B,L\,|\,[\beta],\vec{\gamma},\mu)/B$
  of stable maps to fibers of
   $(\widehat{W},\widehat{L})/\widehat{B}$
  of the specified type, with the $C^{\infty}$-topology,
 is compact over a compact subset of $B$.
} % end-theorem

\bigskip

\noindent
{\it Remark 3.3.9 $[$finiteness of curve classes in $[\beta]$$]$.}
It should be noted that,
 while $[\beta]$ corresponds to an element in $H_2(Y,L:{\Bbb Z})$
  under the map
   $\xi_{\ast}: H_2(X,L;{\Bbb Z})\rightarrow H_2(Y,L;{\Bbb Z})$
  from the symplectic cut $\xi:X\rightarrow Y$,
 the $(W,B\times L)/B$-monodromy orbit $[\beta]$ could be
  an infinite subset in $H_2(X,L;{\Bbb Z})$.
However, only finitely many
  $\overline{\cal M}_{(g,h), (n, \vec{m})}
        (X,L\,|\,\beta^{\prime},\vec{\gamma},\mu)$,
  $\beta^{\prime}\in [\beta]$, can be non-empty
 since all the related stable maps to $X$ of the specified type
  have the same energy and one has the compactness result of [Ye]
  in this case.

\bigskip

\bigskip

Before leaving this section, we should mention that
 the moduli problem of stable maps to fibers of
 $(\widehat{W}, \widehat{L})/\widehat{B}$ has non-trivial obstructions.
The space $\overline{\cal M}_{(g,h), (n, \vec{m})}
              (W/B,L\,|\,[\beta],\vec{\gamma},\mu)/B$
 is very singular in general.
The construction of a family Kuranishi structure on
 $\overline{\cal M}_{(g,h), (n, \vec{m})}
            (W/B,\\  L\,|\,[\beta],\vec{\gamma},\mu)/B$,
 to be done in Sec.~5.3 and Sec.~5.4,
 is meant to accommodate such singularities due to obstructions.

\bigskip

\section{The moduli space
    $\check{\cal W}^{1,p}_{(g,h),(n,\vec{m})}
     ((\widehat{W},\widehat{L})/\widehat{B}\,|\,[\beta],\vec{\gamma},\mu)$
    of stable $\check{W}^{1,p}$-maps.}

In this section we introduce the moduli space
 $\check{\cal W}^{1,p}_{(g,h),(n,\vec{m})}
   ((\widehat{W},\widehat{L})/\widehat{B}\,|\,[\beta],\vec{\gamma},\mu)$
 of $\check{W}^{1,p}$-maps to fibers of
 $(\widehat{W},\widehat{L})/\widehat{B}$.
This space fibers over $B$  and
 is locally embeddable into a Banach orbifold-with-corners;
it contains
 $\overline{\cal M}_{(g,h), (n, \vec{m})}
              (W/B,L\,|\,[\beta],\vec{\gamma},\mu)$
 as a finite dimensional, compact-over-$B$,
      singular sub-orbifold-with-corners.
Members of the system of Kuranishi neighborhoods for
  $\overline{\cal M}_{(g,h), (n, \vec{m})}
              (W/B,L\,|\,[\beta],\vec{\gamma},\mu)/B$,
  to be constructed in Sec.~5.3,
 are embedded in the local singular-orbifold-charts of
  $\check{\cal W}^{1,p}_{(g,h),(n,\vec{m})}
   ((\widehat{W},\widehat{L})/\widehat{B}\,|\,[\beta],\vec{\gamma},\mu)/B$
  as finite dimensional, locally closed, algebraic-type subsets
     that are flat over $B$.
The related Banach relative tangent-space fibration
   $T^1_{\check{\cal W}^{1,p}_{\tinybullet}
         ((\widehat{W},\widehat{L})/\widehat{B}\,|\,\tinybullet)
         /(\widetilde{\cal M}_{\tinybullet}\times\widehat{B})}$  and
 the related Banach relative obstruction-space fibration
   $T^2_{\check{\cal W}^{1,p}_{\tinybullet}
         ((\widehat{W},\widehat{L})/\widehat{B}\,|\,\tinybullet)
          /(\widetilde{\cal M}_{\tinybullet}\times\widehat{B})}$  on
 $\check{\cal W}^{1,p}_{(g,h),(n,\vec{m})}
  ((\widehat{W},\widehat{L})/\widehat{B}\,|\,[\beta],\vec{\gamma},\mu)$
 over $\widetilde{\cal M}_{(g,h),(n,\vec{m})}\times \widehat{B}$
 and  their flattening stratification are also given.

The foundation (Sec.~4.1) of the construction of these moduli spaces
 and fibrations are in [Sie1: Sec.~4 - Sec.~6]; see also [Ru] and [L-R].
These spaces will be used to study the gluability and the gluing of
 family Kuranishi neighborhoods in Sec.~5.4.

We assume throughout the work that $2<p<\infty$ to ensure
 the continuity of $W^{1,p}$- and $\check{W}^{1,p}$-maps
 (to be defined below) on a bordered Riemann surface.

\bigskip

\subsection{The moduli space
    $\check{\cal W}^{1,p}_{(g,h),(n,\vec{m})}
                (W[k],L[k]\,|\,[\beta],\vec{\gamma},\mu)$  of
     stable $\check{W}^{1,p}$-maps to\\  $(W[k],L[k])$,
    its relative tangent and relative obstruction bundles.}
%
%%%%%%%%%%%%%%%%%%%%%%%5
% The techniques and many results of Siebert on the
%  $\bar{\partial}$-operators in families of prestable curves and
%  the Banach orbifold structure on the moduli space of closed
%  $\check{W}^{1,p}$ maps in [Sie1: Sec.~4 - Sec.~6] can be generalized
%  readily to open $\check{W}^{1,p}$ maps from prestable bordered Riemann
%  surfaces to an almost-complex manifold $X$ with the boundary of domains
%  being mapped to a submanifold in $X$.
% Only mild additions to the details in ibidem are required to take
%  care of the effect from the boundary and the boundary nodes on the
%  domains and from the submanifold in the target $X$.
% As a preparation to Sec.~???, we state the results for the case
%  $X=W[k]$ with the coisotropic submanifold $L[k]$.
%%%%%%%%%%%%%%%%%%%%%%%%
%
For a labelled-bordered Riemann surface $\Sigma$ with marked points
 with a K\"{a}hler metric
 so that all the boundary components are geodesics,
let $U_{\varepsilon}$ be the $\varepsilon$-neighborhood of
 the set of nodes on $\Sigma$ with respect to the metric.
Consider a measure $\mu$ on $\Sigma$ defined as follows.
 \begin{itemize}
  \item[$\cdot$]
   on $\Sigma-U_{\varepsilon}$, $\mu$ coincides with the area-form
    associated to the metric;

  \item[$\cdot$]
   around the $\varepsilon/2$-neighborhood of a node
    with local polar coordinates $(r,\theta)$,
     $\mu=drd\theta$, where
    $r$ is the distance function to the node and
    $\theta$ parameterizes the angular direction;
    $\theta$ runs
     over $[0,2\pi]$
      for an interior node or a boundary node of type $E$  and
     over the disjoint union of two finite closed intervals
      for a boundary node of type $H$;

  \item[$\cdot$]
   on $U_{\varepsilon}-U_{\varepsilon/2}$,
   $\mu$ is realized as a non-degenerate $2$-form that
    interpolates smoothly the above two $2$-forms.
 \end{itemize}
Define the $\check{L}^p$- (resp.\ $W^{k,p}$)-norm for a function
 $f$ on $\Sigma$ to be the $L^p$- (resp.\ $W^{k,p}$-)norm of $f$
 with respect to $\mu$:
 $$
  \|f\|_{\check{L}}\;
   =\; \left(\int_{\Sigma}\,|f|^p\,\mu\right)^{1/p}\,,
     \hspace{2em}
  \|f\|_{\check{W}^{k,p}}\;
   =\; \left(
        \int_{\Sigma}\,\sum_{|\mu|\le k}\, |\partial^{\nu}f|^p\, \mu
       \right)^{1/p}\,.
 $$
The completion of $C^{\infty}(\Sigma)$ with respect to the norm
 $\|\,\cdot\,\|_{\check{L}^p}$ and $\|\,\cdot\,\|_{\check{W}^{k,p}}$
 is denoted by $\check{L}(\Sigma)$ and $\check{W}^{k,p}(\Sigma)$
 respectively.
For $2<p<\infty$,
 $$
  C^{\infty}(\Sigma)  \subset  \check{W}^{1,p}(\Sigma)  \subset
  W^{1,p}(\Sigma)     \subset  C^0(\Sigma)\,.
 $$
The notion of bounded $\check{L}$- or
 $\check{W}^{k,p}$-norm depends only on the complex structure on
 $\Sigma$, not the K\"{a}hler metric, $\varepsilon$,  or
  the smooth interpolation on $U_{\varepsilon}-U_{\varepsilon/2}$.
In particular, though such measure $\mu$ on $\Sigma$ is not
 invariant under $\Aut(\Sigma)$ in general,
the notion of functions of bounded $\check{L}^p$-  or
 $\check{W}^{1,p}$-norm on $\Sigma$ is invariant $\Aut(\Sigma)$.
The notion generalizes to maps to manifolds or sections of a bundle.
The choice of such Sobolev sections makes the local trivialization
 over a base $S$ of the space of Sobolev sections for vector bundles
 of a family ${\cal C}_S/S$ of prestable labelled-bordered Riemann
 surfaces with marked points over $S$ that occurs in our problem
 possible;
see [Sie1: Sec.~4] for the technical details,
 which can be generalized to our case.

The symplectic cut $\xi:X\rightarrow Y$ extends to  a
 strong deformation retract $r: W/B\rightarrow Y/\{0\}$
 such that the restriction $r_{\lambda}:W_{\lambda}\rightarrow Y$
 is also a symplectic cut.
The post-composition of $r$ with
 $\tilde{\mathbf p}[k]:W[k]/B[k]\rightarrow W/B$
 defines a map ${\mathbf r}[k]: W[k]/B[k]\rightarrow Y/\pt$.

% \noindent $\bullet$
% $\check{W}^{1,p}$ map,
%  its Maslov index (for the coisotropic boundary case?),
%  and its combinatorial type.

\bigskip

\noindent
{\bf Definition 4.1.1 [stable $\check{W}^{1,p}$-map to $(W[k],L[k])$].}
{\rm
 A $\check{W}^{1,p}$-map
  $h:(\Sigma,\partial\Sigma)\rightarrow (W[k],L[k])$ is said to be
  of ({\it combinatorial}) {\it type}
  $((g,h),(n,\vec{m})\,|\,[\beta],\vec{\gamma},\mu)$
 if
  $\Sigma$ is a labelled-bordered Riemann surface of type
   $((g,h),(n,\vec{m}))$,
  ${\mathbf r}[k]_{\ast}h_{\ast}([\Sigma])=\xi_{\ast}([\beta])$,
  ${\mathbf r}[k]_{\ast}(h_{\ast}[\dot{\partial}\Sigma])\;
   =\; \vec{\gamma}$, and
  the relative homotopy class of ${\mathbf r}[k]\circ h$ contains
   a map of Maslov index $\mu$.
 $h$ is called {\it stable} if the restriction of $h$ to each unstable
  component of $\Sigma$ is non-constant.

 An {\it isomorphism} from $h_1:\Sigma_1\rightarrow W[k]$  to
  $h_2:\Sigma_2\rightarrow W[k]$ is an isomorphism
  $\alpha:\Sigma_1\rightarrow \Sigma_2$ such that $h_1=h_2\circ\alpha$.
 The isomorphism class of $h$ is denoted by $[h]$.
 The moduli space of isomorphisms classes of
  stable $\check{W}^{1,p}$-maps to $(W[k],L[k])$  of
  type $((g,h),(n,\vec{m})\,|\,[\beta],\vec{\gamma},\mu)$ is denoted by
 $\check{\cal W}^{1,p}_{(g,h),(n,\vec{m})}
                (W[k],L[k]\,|\,[\beta],\vec{\gamma},\mu)$.
} % end-definition

\bigskip

\noindent
Note that the stability condition implies that $\Aut(h)$ is finite
 for a stable $\check{W}^{1,p}$-map $h$.

As $2<p<\infty$, a $\check{W}^{1,p}$-map is continuous and one
 can define the {\it $C^0$-topology} on
 $\check{\cal W}^{1,p}_{(g,h),(n,\vec{m})}
                (\\  W[k],L[k]\,|\,[\beta],\vec{\gamma},\mu)$  by
 defining the $(\varepsilon_1,\varepsilon_2)$-neighborhood
 $U_{\varepsilon_1,\varepsilon_2}([h])$ of $[h]$ to consist of
 all $[h^{\prime}: (\Sigma^{\prime}, \partial\Sigma^{\prime})
                   \rightarrow (W[k],L[k])]
      \in \check{\cal W}^{1,p}_{(g,h),(n,\vec{m})}
                              (W[k],L[k]\,|\,[\beta],\vec{\gamma},\mu)$
 such that
  there exists a surjective {\it collapsing/pinching map}
   $c:\Sigma^{\prime}\rightarrow\Sigma$
   that is a diffeomorphism  from
    the complement of a collection of simple loops
     and simple arc with ends on $\partial\Sigma^{\prime}$ on
     $\Sigma^{\prime}$  to
    the complement of the set of nodes on $\Sigma^{\prime}$,  and
    collapses/pinches each simple loop (resp.\ arc) in the collection
     to an interior (resp.\ boundary) node of $\Sigma$
   so that
   \begin{itemize}
    \item[$\cdot$]
     \parbox{31ex}{(nearness of domain)}\
      \parbox[t]{22em}{
       $\Sigma^{\prime}$ is isomorphic to a fiber of
         ${\cal C}/\Def(\Sigma)$ with\\
        $\|j-c_{\ast}j \|_{C^{\infty}}<\varepsilon_2$
         on $\Sigma-U_{\varepsilon_1}$  and
        $c(p^{\prime}_{\cdot})$ in the $\varepsilon_2$-neighborhood
         of $p_{\cdot}$,
         where $p_{\cdot}$, $p^{\prime}_{\cdot}$ are
          marked points on $\Sigma$, $\Sigma^{\prime}$ that are paired
          by their label;}

    \item[$\cdot$]
     \parbox{31ex}{(nearness of map)}
      $\|h-h^{\prime}\circ c^{-1}\|_{C^0}<\varepsilon_2$
      on $\Sigma_{\reg}$.
  \end{itemize}
 Here,
  $U_{\varepsilon_1}$ is the $\varepsilon_1$-neighborhood of
   the set of nodes of $\Sigma$
   that is small enough so that it contains no marked points.
 This topology is equivalent to the $L^{\infty}$-topology,
 [Sie1: Proposition 5.3].

A {\it Banach space-with-corners} is the direct product of
 a Banach space and a polyhedral cone at the origin in
  a finite-dimensional (real) vector space.
A {\it Banach orbifold-with-corners} is an orbifold locally modelled
 on a finite quotient of a neighborhood of the origin of  a
 Banach space-with-corners.
The same techniques for the proof of
 [Sie1: Proposition 3.8 and Theorem 5.1]
 can be applied to prove the following theorem:

\bigskip

\noindent
{\bf Theorem 4.1.2 [$\check{\cal W}^{1,p}_{(g,h),(n,\vec{m})}
                    (W[k],L[k]\,|\,[\beta],\vec{\gamma},\mu)$].}
{\it
 The $C^0$-topology on the moduli space\\
  $\check{\cal W}^{1,p}_{(g,h),(n,\vec{m})}
              (W[k],L[k]\,|\,[\beta],\vec{\gamma},\mu)$  is Hausdorff.
 There exists a refinement of the $C^0$-topology on
  $\check{\cal W}^{1,p}_{(g,h),(n,\vec{m})}
                (W[k],L[k]\,|\,[\beta],\vec{\gamma},\mu)$
  so that it becomes a $($Hausdorff$)$ Banach orbifold-with-corners.
} % end-theorem

\bigskip

We shall call the refined topology in the above theorem  the
 {\it $\check{W}^{1,p}$-topology} on
 $\check{\cal W}^{1,p}_{(g,h),(n,\vec{m})}
                 (\\  W[k],L[k]\,|\,[\beta],\vec{\gamma},\mu)$.
With this topology, a local orbifold-chart of $[h]$ is modelled
 on the quotient of the Banach space-with-corners from  a
 rigidifying slice $V^{\prime\prime}_{[h]}$
 to the approximate pseudo-$\Aut(\Sigma)$-action on
 $\Def(\Sigma)  \times
  W^{1,p}(\Sigma,\partial\Sigma;
          h^{\ast}T_{\ast}W[k],(h|_{\partial})^{\ast}T_{\ast}L[k])$
 by $\Aut(h)$.
Here
 $W^{1,p}(\Sigma,\partial\Sigma;
    h^{\ast}T_{\ast}W[k],\\   (h|_{\partial\Sigma})^{\ast}T_{\ast}L[k])$
 is the Banach space of $\check{W}^{1,p}$-sections $s$ of
 the vector bundle $h^{\ast}T_{\ast}W[k]$ on $\Sigma$
 with $s|_{\partial\Sigma}$ taking values in
 $(h|_{\partial\Sigma})^{\ast}T_{\ast}L[k]$.
We call $(V^{\prime\prime}_{[h]},\Gamma_{[h]}:=\Aut([h]))$,  a
 {\it Banach orbifold-with-corners chart} of $[h]$  in
 $\check{\cal W}^{1,p}_{(g,h),(n,\vec{m})}
                  (W[k],L[k]\,|\,[\beta],\vec{\gamma},\mu)$.
(We will call $[h]$ the center of the chart for convenience.)

% \bigskip
%
% \noindent $\bullet$
% Remark on the construction of local charts from the product
%  deformation space of a domain curve $\Sigma$ and
%  the $\check{W}^{1,p}$ deformation space of a map
%  $h:(\Sigma,\partial\Sigma) \rightarrow (W[k],L[k])$
%  with $\Sigma$ fixed.
% We will call this $h$ the {\it center} of a such chart.

The system of the equivariant relative tangent bundle of the local
 Banach orbifold-with-corners charts over the deformation space
 of the domain curve in the center glue to a Banach orbi-bundle
 $T^1_{\check{\cal W}^{1,p}_{\tinybullet}(W[k],L[k]\,|\,\tinybullet)
           /\widetilde{\cal M}_{\tinybullet}}$  on
 $\check{\cal W}^{1,p}_{(g,h),(n,\vec{m})}
               (W[k],L[k]\,|\,[\beta],\vec{\gamma},\mu)$,
 whose fiber at $[h:(\Sigma,\partial\Sigma)\rightarrow (W[k],L[k])]$
  is given by the Banach $\Aut(h)$-space
  $\check{W}^{1,p}(\Sigma,\partial\Sigma;
        h^{\ast}T_{\ast}W[k],
        (h|_{\partial\Sigma})^{\ast}T_{\ast}(L[k]))$.
We call this orbi-bundle the {\it relative tangent bundle}  of
 $\check{\cal W}^{1,p}_{(g,h),(n,\vec{m})}
               (W[k],L[k]\,|\,[\beta],\vec{\gamma},\mu)$  over
 $\widetilde{M}_{(g,h),(n,\vec{m})}$.

The same construction in [Sie1: Sec.~6.1] gives a Banach orbi-bundle
 $T^2_{\check{\cal W}^{1,p}_{\tinybullet}(W[k],L[k]\,|\,\tinybullet)
       /\widetilde{\cal M}_{\tinybullet}}$  on
 $\check{\cal W}^{1,p}_{(g,h),(n,\vec{m})}
               (W[k],L[k]\,|\,[\beta],\vec{\gamma},\mu)$,
 whose fiber at $[h:(\Sigma,\partial\Sigma)\rightarrow (W[k],L[k])]$
  is given by the Banach $\Aut(h)$-space
  $\check{L}^p(\Sigma;
               \Lambda^{0,1}\Sigma\otimes_J h^{\ast}T_{\ast}W[k])$  of
  $\check{L}^p$-sections of
  $\Lambda^{0,1}\Sigma\otimes_J h^{\ast}T_{\ast}W[k]$.
We call this orbi-bundle the {\it relative obstruction bundle}  of
 $\check{\cal W}^{1,p}_{(g,h),(n,\vec{m})}
               (W[k],L[k]\,|\,[\beta],\vec{\gamma},\mu)$  over
 $\widetilde{M}_{(g,h),(n,\vec{m})}$.

The nonlinear Cauchy-Riemann operator
 $h\mapsto \bar{\partial}_Jh := \frac{1}{2}(dh + J\circ dh\circ j)$
 defines a section (in the sense of orbi-bundle)
 $$
  s_{\bar{\partial}_J}\;:\;
   \check{\cal W}^{1,p}_{(g,h),(n,\vec{m})}
               (W[k],L[k]\,|\,[\beta],\vec{\gamma},\mu)\;
   \longrightarrow\;
   T^2_{\check{\cal W}^{1,p}_{\tinybullet}(W[k],L[k]\,|\,\tinybullet)
          /\widetilde{\cal M}_{\tinybullet}}
 $$
 of the relative obstruction bundle
 $T^2_{\check{\cal W}^{1,p}_{\tinybullet}(W[k],L[k]\,|\,\tinybullet)
         /\widetilde{\cal M}_{\tinybullet}}$.

A connection $\nabla$ on $T_{\ast}W[k]$ induces  an
 {\it partial connection} on the orbi-bundle\\
  $T^2_{\check{\cal W}^{1,p}_{\tinybullet}(W[k],L[k]\,|\,\tinybullet)
           /\widetilde{\cal M}_{\tinybullet}}$,
 using the parallel transport on $T_{\ast}W[k]$ associated to $\nabla$.
Denote this $\nabla$-induced partial connection on
 $T^2_{\check{\cal W}^{1,p}_{\tinybullet}(W[k],L[k]\,|\,\tinybullet)
          /\widetilde{\cal M}_{\tinybullet}}$
 also by $\nabla$;
then its associated horizontal distribution $H^{\nabla}$ at
  a point $([h],\eta)$ over $[h]$ projects isomorphically to  the
  relative tangent space
  $T^1_{\check{\cal W}^{1,p}_{\tinybullet}(W[k],L[k]\,|\,\tinybullet)
            /\widetilde{\cal M}_{\tinybullet}\,,\, [h]}$  of
  $\check{\cal W}^{1,p}_{(g,h),(n,\vec{m})}
                 (W[k],L[k]\,|\,[\beta],\vec{\gamma},\mu)$  over
  $\widetilde{M}_{(g,h),(n,\vec{m})}$ at $[h]$.
One thus has a well-defined {\it vertical projection} $\pi^v$ to the
 tangent space
 $T_{([h],\eta)}\,
   T^2_{\check{\cal W}^{1,p}_{\tinybullet}(W[k],L[k]\,|\,\tinybullet)
            /\widetilde{\cal M}_{\tinybullet}\,,\,[h]}$
 for a tangent vector at $([h], \eta)$ that projects into
 $T^1_{\check{\cal W}^{1,p}_{\tinybullet}(W[k],L[k]\,|\,\tinybullet)
             /\widetilde{\cal M}_{\tinybullet}\,,\, [h]}$.
Together with the vector space translations on fibers of
  $T^2_{\check{\cal W}^{1,p}_{\tinybullet}(W[k],L[k]\,|\,\tinybullet)
           /\widetilde{\cal M}_{\tinybullet}}$,
 the composition $\pi^v\circ d\,s_{\bar{\partial}_J}$ defines
 an orbi-bundle homomorphism
 $$
  D\bar{\partial}_J\;:\;
   T^1_{\check{\cal W}^{1,p}_{\tinybullet}(W[k],L[k]\,|\,\tinybullet)
             /\widetilde{\cal M}_{\tinybullet}}
   \longrightarrow\;
   T^2_{\check{\cal W}^{1,p}_{\tinybullet}(W[k],L[k]\,|\,\tinybullet)
            /\widetilde{\cal M}_{\tinybullet}}\,.
 $$
We shall call $D\bar{\partial}_J$  the
 {\it $\nabla$-induced linearization} of  the
 nonlinear Cauchy-Riemann operator $\bar{\partial}_J$.
The expression for $D\bar{\partial}_J$ can be computed explicitly.
% For $\nabla$ the Levi-Civita connection of the metric associated
%  to $(\omega,J)$, the formula goes as
%  $$
%   (D_h\bar{\partial}_J)(\xi)\;
%    =\;  \frac{1}{2}\,(
%          \nabla\xi\circ dh  +
%          J\circ \nabla\xi\circ dh\circ j  +
%          \nabla_{\xi}J\circ dh\circ j)
%  $$
%  when $h$ is not constant on any open subset of $\Sigma$.
See, e.g.,
 [MD-S1: Eq.~(3.2), Remark~3.3.1], [Liu(C): Proposition 6.12], and
 [Sie1: Sec.~6.3].
Note that, by definition, the $J$-holomorphy locus in
 $\check{\cal W}^{1,p}_{(g,h),(n,\vec{m})}
                  (W[k],L[k]\,|\,[\beta],\vec{\gamma},\mu)$
 is sent by $s_{\bar{\partial}_J}$ to the image of the zero-section of
 $T^2_{\check{\cal W}^{1,p}_{\tinybullet}(W[k],L[k]\,|\,\tinybullet)
           /\widetilde{\cal M}_{\tinybullet}}$;
the linearization $D_h\bar{\partial}_J$ for $h$ $J$-holomorphic
 is thus independent of $\nabla$.

Finally, we remark that,
as $L[k]$ is a coisotropic submanifold that contains properly
  a symplectic submanifold (e.g.\ $\Real(B[k])\times L$) in $W[k]$,
the restriction of the orbi-bundle homomorphism $D\bar{\partial}_J$
 to each fiber is not Fredholm.
Instead, $D\bar{\partial}_J$, when restricted to fibers, has a
 finite-dimensional cokernel  but
 an infinite-dimensional kernel in general.

\bigskip

\subsection{The moduli space
    $\check{\cal W}^{1,p}_{(g,h),(n,\vec{m})}
     ((\widehat{W},\widehat{L})/\widehat{B}\,|\,[\beta],\vec{\gamma},\mu)$
     of stable $\check{W}^{1,p}$-maps to fibers of
     $(\widehat{W},\widehat{L})/\widehat{B}$,
    the relative $\check{W}^{1,p}$-tangent-obstruction fibration complex.}

{\bf Definition 4.2.1
     [stable $\check{W}^{1,p}$-map to $(W[k],L[k])/B[k]$].}
{\rm
 A $\check{W}^{1,p}$-map
  $h:(\Sigma,\partial\Sigma)/\pt\rightarrow (W[k], L[k])/B[k]$
  from a bordered Riemann surface $\Sigma$ to a fiber of
  $(W[k],L[k])/B[k]$ is called {\it prestable} of
   ({\it combinatorial}) {\it type}
   $((g,h),(n,\vec{m})\,|\,[\beta],\vec{\gamma},\mu)$
  if $h$ is a stable $\check{W}^{1,p}$-map from
   $(\Sigma, \partial\Sigma)$ to $(W[k],L[k])$ of type
    $((g,h),(n,\vec{m})\,|\,[\beta],\vec{\gamma},\mu)$
    such that the image of $h$ lies in a fiber of $(W[k],L[k])/B[k]$.
 An {\it isomorphism} between two prestable $\check{W}^{1,p}$-maps
    $h_1:\Sigma_1/\pt\rightarrow W[k]/B[k]$,
    $h_2:\Sigma_2/\pt\rightarrow W[k]/B[k]$
   of the same type is a pair $(\alpha,\beta)$
  where
   $\alpha:\Sigma_1\rightarrow \Sigma_2$
    is an isomorphism of prestable labelled-bordered Riemann surfaces
      with marked points and
   $\beta\in{\Bbb G}_m[k]$ acts on $W[k]/B[k]$ as in Sec.~1.1.3
  such that $f_1\circ \beta = f_2\circ\alpha$.
 The isomorphism class associated to a prestable $\check{W}^{1,p}$-map
  $h$ will be denoted by $[h]$.
 The {\it group of automorphisms} $\Aut(h)$ of a prestable
  $h:\Sigma/\pt\rightarrow W[k]/B[k]$ consists of elements
  $(\alpha,\beta)\in\Aut(\Sigma)\times {\Bbb G}_m[k]$
  such that $\beta\circ h=h\circ\alpha$.

 A prestable $\check{W}^{1,p}$-map
  $h:\Sigma/\pt \rightarrow (W[k].L[k])/B[k]$ is called {\it stable}
  if $\Aut(h)$ is finite.
 The moduli space of isomorphism classes of stable $\check{W}^{1,p}$-maps
  to fibers of $(W[k],L[k])/B[k]$  of type
  $((g,h),(n,\vec{m})\,|\,[\beta],\vec{\gamma},\mu)$  is denoted by
  ${\cal W}^{1,p}_{(g,h), (n, \vec{m})}
        ((W[k],L[k])/B[k]\,|\,[\beta],\vec{\gamma},\mu)$.
} % end-definition

%%%%%%%%%%%%%%%%%%%%%%%%%%%%%%%%%%%%%%%%%%%
% A $\check{W}^{1,p}$ map $h$ from a bordered Riemann surface $\Sigma$
%  to a {\it non-rigid} $Y_{[k]}$,
%  in notation $h:\Sigma\rightarrow Y_{[k]}^{\nonrigid}$,
% is an ordinary $\check{W}^{1,p}$ map $h:\Sigma\rightarrow Y_{[k]}$
%  but with the set $\Mor(h_1,h_2)$ of {\it morphisms} from
%   $h_1:\Sigma_1\rightarrow Y_{[k]}$ to $h_2:\Sigma_2\rightarrow Y_{[k]}$
%  being the set of pairs
%   $(\alpha,\beta)\in\Mor(\Sigma_1,\Sigma_2)\times \Aut^J(Y_{[k]}/Y)$
%   such that the diagram
%   $$
%    \begin{array}{ccc}
%      \Sigma_1 & \stackrel{h_1}{\longrightarrow} & Y_{[k]} \\[.6ex]
%      \mbox{\scriptsize $\alpha$}\downarrow\hspace{2ex} &
%               & \hspace{1ex}\downarrow\mbox{\scriptsize $\beta$} \\
%      \Sigma_2 & \stackrel{h_2}{\longrightarrow} & Y_{[k]}
%    \end{array}
%   $$
%   commutes.
% We say that a such $h$ is {\it stable} if in addition
%  $\Aut(h):=\Mor(h,h)$ defined from above is finite.
% This is a locally closed condition.
%
% The same definitions apply to (open) {\it stable $\check{W}^{1,p}$
%  maps to a non-rigid $(W[k],L[k])$}, with $\Aut(Y_{[k]}/Y)$
%  replaced by ${\Bbb G}_m[k]$ that acts equivariantly on
%  $(W[k],L[k])/B[k]$.
% We will apply such notion particularly to maps whose image lies in
%  a fiber of $W[k]/B[k]$, $k\in {\Bbb Z}_{\ge 0}$.
%%%%%%%%%%%%%%%%%%%%%%%%%%%%%%%%%%%%%

\bigskip

Once having the notion of stable $\check{W}^{1,p}$-maps to the fibers
 of $(W[k],L[k])/B[k]$, one can apply the same procedure/routine of
 gluings as in Sec.~3.3 to define/obtain the moduli space
 $\check{\cal W}^{1,p}_{(g,h),(n,\vec{m})}
    ((\widehat{W},\widehat{L})/\widehat{B}\,
                           |\,[\beta],\vec{\gamma},\mu)$
 of (isomorphism classes of) stable $\check{W}^{1,p}$-maps to
 the fibers of $(\widehat{W},\widehat{L})/\widehat{B}$.

\bigskip

\begin{flushleft}
{\bf The $\check{W}^{1,p}$-topology and the singular
     orbifold-with-corners structure on
      $\check{\cal W}^{1,p}_{(g,h),(n,\vec{m})}
        ((\widehat{W},\widehat{L})/\widehat{B}\,
                               |\,[\beta],\vec{\gamma},\mu)$.}
\end{flushleft}
Let
 $\check{\cal W}^{1,p}_{(g,h),(n,\vec{m})}
        (W[k],L[k]\,|\,[\beta],\vec{\gamma},\mu)^{W[k]/B[k]}$
 be the singular (constructible) sub-orbifold-with-corners of
 the Banach orbifold-with-corners
  $\check{\cal W}^{1,p}_{(g,h),(n,\vec{m})}
                (W[k],L[k]\,|\,[\beta],\vec{\gamma},\mu)$
 whose system of local singular orbifold-with-corners charts
 consists of
 $$
  \left\{\,
    (V^{\prime},\Gamma_{V^{\prime}})\;\left|\;\
    \parbox{30em}{There exists a Banach orbifold-with-corners
     local chart $(V^{\prime\prime},\Gamma_{V^{\prime\prime}})$ of\\
     $\check{\cal W}^{1,p}_{(g,h),(n,\vec{m})}
                     (W[k],L[k]\,|\,[\beta],\vec{\gamma},\mu)$
     such that
     \begin{itemize}
      \item[$\cdot$]
       $V^{\prime}$ is the subset of $V^{\prime\prime}$ parameterizing
        all those $\check{W}^{1,p}$-maps to $(W[k],L[k])$
         parameterized by $V^{\prime\prime}$
        whose image lies completely in a fiber of $(W[k],L[k])/B[k]$
         and
        which are stable in the sense of Definition 4.2.1;
            % Definition
            %    [stable $\check{W}^{1,p}$-map to $(W[k],L[k])/B[k]$]

      \item[$\cdot$]
      $\Gamma_{V^{\prime}}=\Gamma_{V^{\prime\prime}}\,$.
     \end{itemize} } \right.\,
  \right\}\,.
 $$
Note that $V^{\prime}$ is a locally closed subset of
 the corresponding $V^{\prime\prime}$.
The gluing of the system of local charts
 $\{(V^{\prime},\Gamma_{V^{\prime}})\}_{\bullet}$
  for $\check{\cal W}^{1,p}_{(g,h),(n,\vec{m})}
         (W[k],L[k]\,|\,[\beta],\vec{\gamma},\mu)^{W[k]/B[k]}$
 follows from the restriction of the gluing of the subsystem
  $\{(V^{\prime\prime},\Gamma_{V^{\prime\prime}})\}_{\bullet}$
  of charts for $\check{\cal W}^{1,p}_{(g,h),(n,\vec{m})}
                     (W[k],L[k]\,|\,[\beta],\vec{\gamma},\mu)$.
The natural map from
 $\check{\cal W}^{1,p}_{(g,h),(n,\vec{m})}
      (W[k],L[k]\,|\,[\beta],\vec{\gamma},\mu)^{W[k]/B[k]}$
 to $B[k]$ defines the notation
 $\check{\cal W}^{1,p}_{(g,h),(n,\vec{m})}
   (W[k],L[k]\,|\,[\beta],\vec{\gamma},\mu)^{W[k]/B[k]}/B[k]$.
A singular orbifold-with-corners structure  on
 $\check{\cal W}^{1,p}_{(g,h),(n,\vec{m})}
   ((\widehat{W},\widehat{L})/\widehat{B}\,|\,[\beta],\vec{\gamma},\mu)$
 can be obtained by gluing a system of singular charts from further
  orbifolding appropriate subsets of the singular local charts of
  $\check{\cal W}^{1,p}_{(g,h),(n,\vec{m})}
      (W[k], L[k]\,|\,[\beta],\vec{\gamma},\\  \mu)^{W[k]/B[k]}$,
  $k\in {\Bbb Z}_{\ge 0}$, as follows.

Let
 $\rho  \in
  \check{\cal W}^{1,p}_{(g,h),(n,\vec{m})}
   ((\widehat{W},\widehat{L})/\widehat{B}\,|\,[\beta],\vec{\gamma},\mu)$
 be represented by
 $h:(\Sigma,\partial\Sigma)/\pt
       \rightarrow (Y_{[k]},L_{[k]})/\{0\}\subset (W[k],L[k])/B[k]$.
(The case the target is a smooth $W_{\lambda}$, $\lambda\ne 0$,
  is immediate and will be omitted.)
For our final purpose of studying
 $\overline{\cal M}_{(g,h), (n, \vec{m})}
              (W/B,L\,|\,[\beta],\vec{\gamma},\mu)$,
 we will assume that the image of $h$ has non-empty intersection
 with each irreducible component of $Y_{[k]}$.
The following discussion can be adapted to the situation when
 this is not the case as well.
As an element in
 $\check{\cal W}^{1,p}_{(g,h),(n,\vec{m})}
         (W[k],L[k]\,|\,[\beta],\vec{\gamma},\mu)^{W[k]/B[k]}$,
let $(V^{\prime},\Gamma_{V^{\prime}})$ be in the form of a singular local
 chart-with-corners $(V^{\prime}_h,\Aut(h)^{\rigid})$ centered at $h$.
The equivariant pseudo-${\Bbb G}_m[k]$-action on $W[k]/B[k]$
 induces an equivariant pseudo-${\Bbb G}_m[k]$-action on
 $\check{\cal W}^{1,p}_{(g,h),(n,\vec{m})}
   (W[k],L[k]\,|\,[\beta],\vec{\gamma},\mu)^{W[k]/B[k]}/B[k]$
 via post-composition with maps.
Locally this is a (pseudo) ${\Bbb G}_m[k]$-action on the singular
 local chart-with-corners $V^{\prime}$ that commutes with
 the $\Gamma_{V^{\prime}}$-action on $V^{\prime}$.
A $\Gamma_{V^{\prime}}$-invariant slice $V_h$ through $h$ in
 $V^{\prime}$ to rigidify this ${\Bbb G}_m[k]$-action can be
 constructed as follows.

Let ${\cal C}^{\prime}/V^{\prime}$ be the universal bordered Riemann
 surface with marked points over $V^{\prime}$ and
 $F^{\prime}:{\cal C}^{\prime}/V^{\prime}\rightarrow W[k]/B[k]$
 be the universal map.
Both are built-in from the construction of Siebert.
Let $U_{\Sigma}\subset \Sigma$ be an $\Aut(h)^{\rigid}$-invariant
 sub-surface of $\Sigma$ by removing an appropriate small neighborhood
 of all the nodes of $\Sigma$.
Then, for $V^{\prime}$ small enough, the $\Aut(h)^{\rigid}$-equivariant
 embedding $U_{\Sigma}\hookrightarrow \Sigma$ extends to an
 $\Aut(h)^{\rigid}$-equivariant embedding
 $V^{\prime}\times U_{\Sigma} \hookrightarrow {\cal C}^{\prime}$
 over $V^{\prime}$.
This implies that there exist global sections $s_i$, $i=1,\,\ldots,\,k$,
 of ${\cal C}^{\prime}\rightarrow V^{\prime}$ such that
 \begin{itemize}
  \item[$\cdot$]
   $\alpha^{\ast}s_i$, $\alpha\in\Aut(h)^{\rigid}$, takes values
   in the image of $V^{\prime}\times U_{\Sigma}$;

  \item[$\cdot$]
   the image of $F^{\prime}\circ \alpha^{\ast}s_i$,
    $\alpha\in\Aut(h)^{\rigid}$, lies in a neighborhood $U_i[k]$
    of $h\circ s_i(0)$ in
    $\Trunk[k]_i\simeq
     B[k]\times (\Delta_i-N_{\varepsilon}(D_{i-1}\cup D_i))$ of $W[k]$,
    cf.\ Remark 1.1.1.6;

  \item[$\cdot$]
   the finite set
    $\left\{\, \pi_{2,i}\circ F\circ \alpha^{\ast}s_i(0):
       i=1,\,\ldots,\,k;\,\alpha\in\Aut(h)^{\rigid} \,\right\}$
   lies in ${\Bbb C}-{\Bbb R}_{\le 0}$,
    where $\pi_{2,i}:U_i[k]\rightarrow {\Bbb C}-\{0\}$
     is the projection map to the fiber of ${\Bbb L}$ from a local
     trivialization of ${\Bbb L}$, as an embedded submanifold in
     $\Delta_i$, $i=1,\,\ldots\,,\, k$, cf.\ Sec.~1.1.1.
 \end{itemize}

\noindent
Define the average function $\Average$ for a finite subset
 $S$ in ${\Bbb C}-{\Bbb R}_{\le 0}$ by
 $$
  \Average(S)\;
   =\; e^{\, \frac{1}{|S|} \sum_{w\in S}
                \left( \log(|w|) + \sqrt{-1}\arg(w) \right)
         }\,,
 $$
 where $\arg(w)\in (-\pi,\pi)$.
Let
 $$
  \bar{s}_i\;
   :=\;  \Average\left( \pi_{2,i}\circ F\circ \alpha^{\ast}s_i\,
           :\,\alpha\in\Aut(h)^{\rigid} \right)\,.
 $$
For $V^{\prime}$ small enough, this is a well-defined
 $\Aut(h)^{\rigid}$-invariant function on $V^{\prime}$
 for $i=1,\,\ldots,\,k$, with values in ${\Bbb C}-{\Bbb R}_{\le 0}$.
The $k$-tuple
 $$
  R \,:=\, (\bar{s}_1,\,\ldots,\,\bar{s}_k)\;
   :\; V^{\prime}\; \longrightarrow\; ({\Bbb C}-\{0\})^k
 $$
 defines thus an $\Aut(h)^{\rigid}$-invariant ${\Bbb G}_m[k]$-equivariant
 map, where ${\Bbb G}_m[k]$ acts on $({\Bbb C}-\{0\})^k$ by
 $(w_1,,\ldots,\,w_k) \mapsto (\sigma_1w_1,\,\ldots,\,\sigma_kw_k)$,
 $(\sigma_1,\,\ldots,\,\sigma_k)\in {\Bbb G}_m[k]=({\Bbb C}^{\times})^k$.
%%%%%%%%%%%%%%%%%%%%%%%%%%%%%%%%%%%%%5
% (Cf.~{\sc Figure} ???.)
%
% \marginpar{\raggedright\tiny $\bullet$ {\sc Figure} ???.}
%
Let $V_h= R^{-1}(R(0))$.
Then $V_h\subset V^{\prime}$ is a rigidifying slice through $h$ to
 the ${\Bbb G}_m[k]$-action on $V^{\prime}$ and is invariant under
 the $\Gamma_{V^{\prime}}$-action.

By construction, the residual discrete subgroup $\Gamma_{V_h}$ of
 $\Aut(\Sigma)\times {\Bbb G}_m[k]$ that pseudo-acts on $V_h$ is
 an extension of $\Aut(h)^{\rigid}$ by a discrete subgroup of
 ${\Bbb G}_m[k]$ whose elements fix $[h]$
 when they descend to pseudo-act on the quotient space
 $V_h/\Aut(h)^{\rigid}$.
In other words, $(\alpha,\beta)\in\Gamma_{V_h}$
 if and only if $\beta \circ h = h\circ \alpha$.
By shrinking $V_h$ if necessary, one can render the pseudo
 $\Gamma_{V_h}$-action to an honest group action.
This shows that indeed $\Gamma_{V_h}=\Aut(h)$.
Stability of $h$ says that $\Gamma_{V_h}$ is finite.
Thus, $(V_h,\Gamma_{V_h})$ defines a singular orbifold local
 chart-with-corners at
 $\rho=[h]  \in
  \check{\cal W}^{1,p}_{(g,h),(n,\vec{m})}
   ((\widehat{W},\widehat{L})/\widehat{B}\,|\,[\beta],\vec{\gamma},\mu)$.
Re-write $h$ above as $h_{\rho}$ to manifest its representing $\rho$
 and
denote the map
 $V_{h_{\rho}}  \rightarrow
   \check{\cal W}^{1,p}_{(g,h),(n,\vec{m})}
    ((\widehat{W},\widehat{L})/\widehat{B}\,|\,[\beta],\vec{\gamma},\mu)$
 that identifies $V_{h_{\rho}}/\Aut(h_{\rho})$ with
  a neighborhood of $\rho$ in
   $\check{\cal W}^{1,p}_{(g,h),(n,\vec{m})}
    ((\widehat{W},\widehat{L})/\widehat{B}\,|\,[\beta],\vec{\gamma},\mu)$
  by $\psi_{\rho}$,
then a system of singular local charts-with-corners on
 $\check{\cal W}^{1,p}_{(g,h),(n,\vec{m})}
  ((\widehat{W},\widehat{L})/\widehat{B}\,|\,[\beta],\vec{\gamma},\mu)$
 is given by $\{(V_{h_{\rho}},\Aut(h_{\rho}),\psi_{\rho})\}_{\rho}$.
We will identify each $V_{h_{\rho}}/\Aut(h_{\rho})$ directly
 as a subset in
 $\check{\cal W}^{1,p}_{(g,h),(n,\vec{m})}
  ((\widehat{W},\widehat{L})/\widehat{B}\,|\,[\beta],\vec{\gamma},\mu)$.

We next construct the transition data for the local charts.
Given a pair $(p,q)$ with
 $p  \in
  \check{\cal W}^{1,p}_{(g,h),(n,\vec{m})}
  ((\widehat{W},\widehat{L})/\widehat{B}\,|\,[\beta],\vec{\gamma},\mu)$
 and
 $q  \in V_{h_p}/\Gamma_{V_{h_p}}  \subset
  \check{\cal W}^{1,p}_{(g,h),(n,\vec{m})}
  ((\widehat{W},\widehat{L})/\widehat{B}\,|\,[\beta],\vec{\gamma},\mu)$,
there is a $\Gamma_{V_{h_q}}$-invariant neighborhood $V_{qp}$ of $h_q$
 in $V_{h_q}$ such that $\psi_q(V_{qp})\subset V_{h_p}/\Gamma_{h_p}$.
The set of embeddings $\{h_q\}\hookrightarrow V_p$ is parameterized
 by a $\Gamma_{h_p}$-orbit in $V_{h_p}$.
Fixing a such embedding determines an embedding
 $h_{qp}:V_{qp}\rightarrow V_{h_p}$ up to a pre-composition with
 the $\Gamma_{h_q}$-action on $V_{qp}$.
The map $h_{qp}$ determines then an embedding
 $\phi_{qp}:\Gamma_{V_{h_q}}\rightarrow \Gamma_{V_{h_p}}$.
The orbifold cocycle condition (cf.\ Definition 5.1.2 (2))
  % Definition [Kuranishi structure-in-${\cal C}$]
 for a triple $(p,q,r)$ with $(p,q)$ as above and $r\in V_q/\Gamma_{h_q}$
 follows immediately.
Thus, the system $\{(V_{qp},h_{qp},\phi_{qp})\}_{(p,q)}$
 gives a required orbifold transition data.

The two systems $\{(V_{h_{\rho}},\Aut(h_{\rho}),\psi_{\rho})\}_{\rho}$
 and $\{(V_{qp},h_{qp},\phi_{qp})\}_{(p,q)}$ together
 give a singular orbifold-with-corners structure on
 $\check{\cal W}^{1,p}_{(g,h),(n,\vec{m})}
  ((\widehat{W},\widehat{L})/\widehat{B}\,|\,[\beta],\vec{\gamma},\mu)$.
The induces topology on
 $\check{\cal W}^{1,p}_{(g,h),(n,\vec{m})}
  ((\widehat{W},\widehat{L})/\widehat{B}\,|\,[\beta],\vec{\gamma},\mu)$
 from these charts will be called the
 {\it $\check{W}^{1,p}$-topology}  on
 $\check{\cal W}^{1,p}_{(g,h),(n,\vec{m})}
  (\\(\widehat{W},\widehat{L})/\widehat{B}\,|\,[\beta],\vec{\gamma},\mu)$.
Theorem 4.1.2
 % Theorem [$\check{\cal W}^{1,p}_{(g,h),(n,\vec{m})}
 %                  (W[k],L[k]\,|\,[\beta],\vec{\gamma},\mu)$]
 together with the detail above implies:

\bigskip

\noindent
{\bf Proposition 4.2.2 [Hausdorffness].} {\it
 $\check{\cal W}^{1,p}_{(g,h),(n,\vec{m})}
  ((\widehat{W},\widehat{L})/\widehat{B}\,|\,[\beta],\vec{\gamma},\mu)$
 with the $\check{W}^{1,p}$-topology is Hausdorff.
} % end-proposition

\bigskip

Note that there is a natural morphism
 (as topological spaces with a system of local charts and gluing data)
 from
 $\check{\cal W}^{1,p}_{(g,h),(n,\vec{m})}
  ((\widehat{W},\widehat{L})/\widehat{B}\,|\,[\beta],\vec{\gamma},\mu)$
 to $\tilde{\cal M}_{(g,h),(n,\vec{m})}\times\widehat{B}$ that forgets
 the map, keeping only the domain and the target in a stable-map data.
It is with respect to this morphism that we denote
 $\check{\cal W}^{1,p}_{(g,h),(n,\vec{m})}
  ((\widehat{W},\widehat{L})/\widehat{B}\,|\,[\beta],\vec{\gamma},\mu)/
  (\tilde{\cal M}_{(g,h),(n,\vec{m})}\times\widehat{B})$.

\bigskip

\begin{flushleft}
{\bf The relative $\check{W}^{1,p}$-tangent-obstruction fibration
     complex on  $\overline{\cal M}_{(g,h), (n, \vec{m})}
                        (W/B,L\,|\,[\beta],\vec{\gamma},\mu)$.}
\end{flushleft}
The total space of the Banach orbi-bundle
 $T^1_{\check{\cal W}^{1,p}_{\bullet}(W[k],L[k]\,|\,\tinybullet)
                              /\widetilde{\cal M}_{\tinybullet}}$  on
 $\check{\cal W}^{1,p}_{(g,h),(n,\vec{m})}
               (W[k],L[k]\,|\,[\beta],\vec{\gamma},\mu)$
 is itself a Banach orbifold-with-corners.
The system of local trivializations of
 $T^1_{\check{\cal W}^{1,p}_{\bullet}(W[k],L[k]\,|\,\tinybullet)
                               /\widetilde{\cal M}_{\tinybullet}}$
 over the system\footnote{The notation $\{\,\cdots\,\}_{\tinybullet}$
                          to indicate a system of objects of the form
                          $\cdots$ will be used in many places of the
                          work.} % end-footnote
 $\{(V^{\prime\prime},\Gamma^{\prime\prime})\}_{\tinybullet}$
 of local charts on
 $\check{\cal W}^{1,p}_{(g,h),(n,\vec{m})}
               (W[k],L[k]\,|\,[\beta],\vec{\gamma},\mu)$
 provides the Banach orbifold-with-corners charts for
 $T^1_{\check{\cal W}^{1,p}_{\bullet}(W[k],L[k]\,|\,\tinybullet)
                               /\widetilde{\cal M}_{\tinybullet}}$
 with the system of gluing data.
After a refinement if necessary, we may assume that all
 $V^{\prime\prime}$ are small enough, so that the collection
 $$
  \left\{
    \left(\,
      T^1_{V^{\prime\prime}}  :=
       T^1_{\check{\cal W}^{1,p}_{\bullet}(W[k],L[k]\,|\,\tinybullet)
                              /\widetilde{\cal M}_{\tinybullet}}
          \left.\rule{0ex}{2ex}\right|_{V^{\prime\prime}}\,,\,
      \Gamma_{V^{\prime\prime}}
    \right)
  \right\}_{\tinybullet}
 $$
 gives the Banach-orbifold-with-corners local charts for
 $T^1_{\check{\cal W}^{1,p}_{\bullet}(W[k],L[k]\,|\,\tinybullet)
                               /\widetilde{\cal M}_{\tinybullet}}$.

Let $\{(V,\Gamma_{V})\}_{\tinybullet}$ be a (fine enough) system of
 local charts on
 $\check{\cal W}^{1,p}_{(g,h),(n,\vec{m})}
     ((\widehat{W},\widehat{L})/\widehat{B}\,
                            |\,[\beta],\vec{\gamma},\mu)$
 as constructed in the previous theme, with each $V$ admitting
 $$
  V\; \subset\;  V^{\prime}\;  \subset\;  V^{\prime\prime}\,,
 $$
 where, recall that,
  $(V^{\prime},\Gamma_{V^{\prime}})$ is a local chart on
   $\check{\cal W}^{1,p}_{(g,h),(n,\vec{m})}
          (W[k],L[k]\,|\,[\beta],\vec{\gamma},\mu)^{W[k]/B[k]}$,  and
  $(V^{\prime\prime},\Gamma_{V^{\prime\prime}})$
   is a local chart on
   $\check{\cal W}^{1,p}_{(g,h),(n,\vec{m})}
                  (W[k],L[k]\,|\,[\beta],\vec{\gamma},\mu)$,
  for some $k\in{\Bbb Z}_{\ge 0}$ depending on $V$.
Consider the fiberwise-closed singular (constructible) subset
 $T^1_V$  of $T^1_{V^{\prime\prime}}$ defined by
 $$
  T^1_V\; :=\;
  \left\{
    \begin{array}{l}
     ([h:(\Sigma,\partial\Sigma)\rightarrow W[k]/B[k]],\xi)\,
      \in\, T^1_{V^{\prime\prime}}|_V \\[.6ex]
     \hspace{4em}
      :\; \xi  \in
          \check{W}^{1,p}(\Sigma,\partial\Sigma; h^{\ast}T_{W[k]/B[k]},
                               (h|_{\partial\Sigma})^{\ast}T_{\ast}L)\,
    \end{array}
  \right\}\,,
 $$
 where
  $\check{W}^{1,p}(\Sigma,\partial\Sigma;
              h^{\ast}T_{W[k]/B[k]},
              (h|_{\partial\Sigma})^{\ast}T_{\ast}L)$
   is the closed Banach subspace of\\
   $\check{W}^{1,p}(\Sigma,\partial\Sigma;
         h^{\ast}T_{\ast}W[k],
         (h|_{\partial\Sigma})^{\ast}T_{\ast}L[k])$
   that consists of
   $\check{W}^{1,p}$-sections of
   $(h^{\ast}T_{\ast}W[k],
     (h|_{\partial\Sigma})^{\ast}T_{\ast}L[k])$
   that are projected to $0$ under
   $\pi[k]_{\ast}:T_{\ast}W[k]\rightarrow T_{\ast}B[k]$.
Then the $\Gamma_V$-action on $V$ canonically lifts to an action
 on $T^1_V$.
The gluing data of the system
 $\{T^1_{V^{\prime\prime}}\}_{\tinybullet}$ extends to the lifting
 the gluing data on the system $\{(V,\Gamma_V)\}_{\tinybullet}$
 to on $\{(T^1_V,\Gamma_V)\}_{\tinybullet}$.
This gives rise to a {\it singular orbifold-with-corners}
 $T^1_{\check{\cal W}^{1,p}_{\tinybullet}
        ((\widehat{W},\widehat{L})/\widehat{B}|\tinybullet)
        /\widetilde{M}_{\tinybullet}}$.
The system of maps $\{(T^1_V\rightarrow V)\}_{\tinybullet}$
 descends to a morphism of orbifolds
 $$
  T^1_{\check{\cal W}^{1,p}_{\tinybullet}
       ((\widehat{W},\widehat{L})/\widehat{B}\,|\,\tinybullet)
                               /\widetilde{\cal M}_{\tinybullet}}\;
  \longrightarrow\;
   \check{\cal W}^{1,p}_{(g,h),(n,\vec{m})}
    ((\widehat{W},\widehat{L})/\widehat{B}\,|\,[\beta],\vec{\gamma},\mu)\,,
 $$
 whose fiber at $\rho$, represented by
 $h:\Sigma\rightarrow Y_{[k]}$,
 is given by the Banach $\Aut(\rho)$-space\\
  $\check{W}^{1,p}(\Sigma,\partial\Sigma;
       h^{\ast}T_{\ast}Y_{[k]}, (h|_{\partial\Sigma})^{\ast}T_{\ast}L)$
 ($:=$ is the closed Banach subspace of\\
    $\check{W}^{1,p}(\Sigma,\partial\Sigma;
          h^{\ast}T_{\ast}W[k],
          (h|_{\partial\Sigma})^{\ast}T_{\ast}L[k])$
    that consists of
    $\check{W}^{1,p}$-sections of
    $(h^{\ast}T_{\ast}W[k],
      (h|_{\partial\Sigma})^{\ast}T_{\ast}L[k])$
    that are projected to $0$ under
    $\pi[k]_{\ast}:T_{\ast}W[k]\rightarrow T_{\ast}B[k]$).

The same restrict-and-descend construction applied to
 the collection of orbi-bundles:
 $$
  \left\{
   \mbox{
    $T^2_{\check{\cal W}^{1,p}_{\bullet}(W[k],L[k]\,|\,\tinybullet)
           /\widetilde{\cal M}_{\tinybullet}}$
    over
    $\check{\cal W}^{1,p}_{(g,h),(n,\vec{m})}
               (W[k],L[k]\,|\,[\beta],\vec{\gamma},\mu)$ }
  \right\}_{k\in{\scriptsizeBbb Z}_{\ge 0}}
 $$
 gives rise to
the singular orbifold-with-corners
 $T^2_{\check{\cal W}^{1,p}_{\tinybullet}
       ((\widehat{W},\widehat{L})/\widehat{B}\,|\,\tinybullet)
        \widetilde{\cal M}_{\bullet}}$
 with a built-in orbifold morphism
 $$
  T^2_{\check{\cal W}^{1,p}_{\tinybullet}
       ((\widehat{W},\widehat{L})/\widehat{B}\,|\,\tinybullet)
         /\widetilde{\cal M}_{\tinybullet}}\;
   \longrightarrow\;
  \check{\cal W}^{1,p}_{(g,h),(n,\vec{m})}
   ((\widehat{W},\widehat{L})/\widehat{B}\,|\,[\beta],\vec{\gamma},\mu)\,,
 $$
 whose fiber at $\rho$, represented by
   $h:(\Sigma,\partial\Sigma)\rightarrow (Y_{[k]},L)$,
  is given by the Banach $\Aut(\rho)$-space
  $\check{L}^p(\Sigma,\partial\Sigma;
             \Lambda^{0,1}\Sigma\otimes_J h^{\ast}T_{\ast}Y_{[k]})$.

%%%%%%%%%%%%%%%%%%%%%%%%%%%%%%%%%%%%%%%%%%%%%%
% \bigskip
%
% \noindent $\bullet$
% {\bf Q.}\ {\it Do we need to consider their associated sheaves?}
% %
% \marginpar{\raggedright\tiny $\bullet$ Let the work guide itself.}
%
% \bigskip
%%%%%%%%%%%%%%%%%%%%%%%%%%%%%%%%%%%%%%%%%%%%%%

Let $\widehat{B}=(B-\{0\})\amalg{\Bbb Z}_{\ge 0}$
 be the stratification of $\widehat{B}$ by the homeomorphism type
 of the fibers of $\widehat{W}/\widehat{B}$
 (with the stratum $B-\{0\}$ labelled by $-1$).
It induces a stratification
 $\{{\cal S}_k\}_{k\in{\scriptsizeBbb Z}_{\ge -1}}$  on
 $\check{\cal W}^{1,p}_{(g,h),(n,\vec{m})}
                  (W[k],L[k]\,|\,[\beta],\vec{\gamma},\mu)$
 by taking the preimage under the forgetful morphism
 $$
  \check{\cal W}^{1,p}_{(g,h),(n,\vec{m})}
                 (W[k],L[k]\,|\,[\beta],\vec{\gamma},\mu)\;
    \longrightarrow\;  \widehat{B}
 $$
 of each stratum in $\widehat{B}$.
The restriction of
 $T^1_{\check{\cal W}^{1,p}_{\bullet}(W[k],L[k]\,|\,\tinybullet)
        /\widetilde{\cal M}_{\tinybullet}}$ and
 $T^2_{\check{\cal W}^{1,p}_{\bullet}(W[k],L[k]\,|\,\tinybullet)
        /\widetilde{\cal M}_{\tinybullet}}$
 to over ${\cal S}_i$ are orbi-bundles on ${\cal S}_i$.
We say that $\{{\cal S}_i\}_{i\in{\scriptsizeBbb Z}_{\ge -1}}$ gives
 a common {\it flattening stratification} for both
 $T^1_{\check{\cal W}^{1,p}_{\bullet}(W[k],L[k]\,|\,\tinybullet)
        /\widetilde{\cal M}_{\tinybullet}}$ and
 $T^2_{\check{\cal W}^{1,p}_{\bullet}(W[k],L[k]\,|\,\tinybullet)
        /\widetilde{\cal M}_{\tinybullet}}$,
 as for a coherent sheaf in algebraic geometry\footnote{In the
                          algebro-geometric setting, the parallel to
                          the various $T^1_{\tinybullet}$ and
                          $T^2_{\tinybullet}$ here will be constructed
                          as a coherent sheaf
                           on a Deligne-Mumford moduli stack
                          from the deformation-obstruction theory of
                           the moduli problem in question.
                          See, e.g.\
                           [L-T1: Sec.~1] and [Li2: Sec.~1.2, Sec.~1.3].
                          In the analytic category that involves Banach
                           orbifolds, it is simpler to construct
                           directly the associated total space,
                           which are themselves (singular) orbifolds,
                           of the would-be sheaves
                           rather than to construct these sheaves.}.

The system of sections
 $$
  \left\{
    s_{\bar{\partial}_J}\;:\;
     \check{\cal W}^{1,p}_{(g,h),(n,\vec{m})}
                 (W[k],L[k]\,|\,[\beta],\vec{\gamma},\mu)\;
     \longrightarrow\;
     T^2_{\check{\cal W}^{1,p}_{\tinybullet}(W[k],L[k]\,|\,\tinybullet)
            /\widetilde{\cal M}_{\tinybullet}}
  \right\}_{k\in{\scriptsizeBbb Z}_{\ge 0}}
 $$
 restricts and descends to a section (as a morphism of orbifolds)
 $$
  s_{\bar{\partial}_J}\,:
  \check{\cal W}^{1,p}_{(g,h),(n,\vec{m})}
   ((\widehat{W},\widehat{L})/\widehat{B}\,
     |\,[\beta],\vec{\gamma},\mu)\;
  \longrightarrow\;
  T^2_{\check{\cal W}^{1,p}_{\tinybullet}
       ((\widehat{W},\widehat{L})/\widehat{B}\,|\,\tinybullet)
         /\widetilde{\cal M}_{\tinybullet}}
 $$
 of  $T^2_{\check{\cal W}^{1,p}_{\tinybullet}
        ((\widehat{W},\widehat{L})/\widehat{B}\,|\,\tinybullet)
          /\widetilde{\cal M}_{\tinybullet}}$.
In contrast, as the connection $\nabla$ on $W[k]$
 that defines the orbi-bundle homomorphism
 $$
  D\bar{\partial}_J\;:\;
   T^1_{\check{\cal W}^{1,p}_{\tinybullet}(W[k],L[k]\,|\,\tinybullet)
             /\widetilde{\cal M}_{\tinybullet}}
   \longrightarrow\;
   T^2_{\check{\cal W}^{1,p}_{\tinybullet}(W[k],L[k]\,|\,\tinybullet)
            /\widetilde{\cal M}_{\tinybullet}}\,.
 $$
 is not ${\Bbb G}_m[k]$-invariant,
the system of these linearizations does not restrict and descend
 to a linearization of $\bar{\partial}_J$ from
 $T^1_{\check{\cal W}^{1,p}_{\tinybullet}
      ((\widehat{W},\widehat{L})/\widehat{B}\,|\,\tinybullet)
        /\widetilde{\cal M}_{\tinybullet}}$
 to
 $T^2_{\check{\cal W}^{1,p}_{\tinybullet}
       ((\widehat{W},\widehat{L})/\widehat{B}\,|\,\tinybullet)
         /\widetilde{\cal M}_{\tinybullet}}$.
However, as the restriction of the linearization $D\bar{\partial}_J$
 over the $J$-holomorphy locus
 $$
  \overline{\cal M}_{(g,h), (n, \vec{m})}
                        (W/B,L\,|\,[\beta],\vec{\gamma},\mu)\;
  \subset\;
   \check{\cal W}^{1,p}_{(g,h), (n, \vec{m})}
   ((\widehat{W},\widehat{L})/\widehat{B}\,|\,[\beta],\vec{\gamma},\mu)
 $$
 is independent of $\nabla$, one does have a morphism as
 an orbifold map between fibered orbifolds:
 $$
  T^1_{\check{\cal W}^{1,p}_{\tinybullet}
          ((\widehat{W},\widehat{L})/\widehat{B}\,|\,\tinybullet)
            /\widetilde{\cal M}_{\tinybullet}}
    \left.\rule{0ex}{2ex}\right|\,
     _{\overline{\cal M}_{\tinybullet}(W/B,L\,|\,\tinybullet)}\;
  \stackrel{D\bar{\partial}_J}{\longrightarrow}\;
   T^2_{\check{\cal W}^{1,p}_{\tinybullet}
          ((\widehat{W},\widehat{L})/\widehat{B}\,|\,\tinybullet)
            /\widetilde{\cal M}_{\tinybullet}}
    \left.\rule{0ex}{2ex}\right|\,
     _{\overline{\cal M}_{\tinybullet}(W/B,L\,|\,\tinybullet)}\,.
 $$
We will call this the
 {\it relative $\check{W}^{1,p}$-tangent-obstruction fibration complex}
 on  $\overline{\cal M}_{(g,h), (n, \vec{m})}
                          (W/B,L\,|\\[0ex]  \,[\beta],\vec{\gamma},\mu)$.

%%%%%%%%%%%%%%%%%%%%%%%%%%%%%%%%%%
% \bigskip
%
% \noindent $\bullet$
% The family $\bar{\partial}_J$-operator and transversality.
%
% \bigskip
%
% \noindent $\bullet$
% The complex linear operator
%  $\bar{\partial} :
%    {\cal T}_1^{\check{W}^{1,p}} \rightarrow {\cal T}_2^{\check{L}^p}$.
%
% \bigskip
%
% \noindent $\bullet$
% The global section $\bar{\partial}_J$ of ${\cal T}_2^{\check{L}^p}$.
%
% \bigskip
%
% \noindent $\bullet$
% The difference of
%  the linearization $D\bar{\partial}_J$ of $\bar{\partial}_J$ and
%  $\bar{\partial}$.
%
% \bigskip
%
% \noindent $\bullet$
%  ????????????.
%%%%%%%%%%%%%%%%%%%%%%%%%%%%%%%%%%%%%%%%%%%%%%%%

\bigskip

This concludes our discussion for these auxiliary
 $\infty$-dimensional Banach-type orbifolds.
To give an orientation for next, we remark that
 to go from these spaces to a finite-dimensional object
  that serves as local charts for
  $\overline{\cal M}_{(g,h), (n, \vec{m})}
                      (W/B,L\,|\,[\beta],\vec{\gamma},\mu)$
  in a generalized sense and is flat over $B$,
 there are three transversality issues one has to deal with:
 \begin{itemize}
  \item[$\cdot$]
   transversality of the operator $\bar{\partial}_J$,
    or equivalently the section $s_{\bar{\partial}_J}$;

  \item[$\cdot$]
   transversality of matching conditions at distinguished nodes;

  \item[$\cdot$]
   transversality$/S$ of the pre-deformability condition
   at distinguished nodes.
 \end{itemize}
The construction of a Kuranishi structure on
 $\overline{\cal M}_{(g,h), (n, \vec{m})}
               (W/B,L\,|\,[\beta],\vec{\gamma},\mu)/B$
 is guided by the attempt to achieve all three transversality
  conditions simultaneously and in a way that is flat over $B$.
This can be realized by
 a modification of $\bar{\partial}_J$  and
 the take of a system of finite-dimensional
  orbifold-structure-group-invariant subsets
  in the local charts of these auxiliary orbifolds.
%%%%%%%%%%%%%%%%%%%%%
% An essence in the construction of a Kuranishi neighborhood of
%  $\rho\in $ is that all the three transversality conditions can be
%  simultaneously achieved
% if one considers not the bundle
%  $???$ with the section $\bar{\partial}_J$
%  but rather a mild quotient bundle $???$ of $???$
%   with the induced section $\pi_{???}\circ\bar{\partial}_J$,
%   where $\pi_{???}: ???\rightarrow ???$ is the quotient map.
% It follows then from the Implicit Function Theorem that
%  the subset $???$ of the zero-locus of $\pi_{???}\circ\bar{\partial}_J$
%   described by the matching and the pre-deformability condition
%  has manifold fibers over $S$ and can lead to a good virtual
%   (read Kuranishi) neighborhood of $\rho$ in $???$
%  if these fibers pile to form a good family over $B[k]$ (not $S$) for
%   an appropriate choice of $S$ that is in a sense maximal at $\rho$.
%%%%%%%%%%%%%%%%%%%%%

\bigskip

\section{Construction of a family Kuranishi structure on \\
            $\overline{\cal M}_{(g,h),(n,\vec{m})}
             (W/B,L\,|\,[\beta],\vec{\gamma},\mu)$ over $B$.}

In this rather long section,
 we construct a family Kuranishi structure over $B$
 for the moduli space
 $\overline{\cal M}_{(g,h),(n,\vec{m})}
   (W/B,L\,|\,[\beta],\vec{\gamma},\mu)/B$
 of stable maps to fibers of $(\widehat{W},\widehat{L})/\widehat{B}$.
This answers incidentally the simplest case,
  namely the degeneration from a symplectic cut,  of
 the question posed in [F-O: p.~962] on the family version of
  Kuranishi structure for a degeneration.
The detail merges [F-O], [Liu(C)] with [I-P2], [L-R], [Li1]; and
the result gives an almost-complex/analytic/symplectic parallel to
 the algebraic [Li1] and [Li2] when curves are closed.

\bigskip

\subsection{Family Kuranishi structure modelled in the category
    ${\cal C}_{\mbox{\rm\footnotesize spsccw}}/{\largeBbb C}$.}

We extend the notion of Kuranishi structure in [F-O: Sec.~5]
 (and also [F-O-O-O: Appendix 2] and [Liu(C): Sec.~6.1]) and
 define a Kuranishi structure modelled in a specific category of
  topology/geometry that appears in our problem;
see also [Sat: Sec.~1] and [Th: Chapter~13] for related discussions
 on orbifolds.

\bigskip

\begin{flushleft}
{\bf Kuranishi structure modelled in a category of topology/geometry.}
\end{flushleft}
Let ${\cal C}$ be a category of topology/geometry --
 e.g.~smooth manifolds with corners,
      complex spaces of specified type of singularities, or
      fibrations over a fixed topological space  --
 in which the notion of morphisms, embeddings, isomorphisms, bundles,
  and groups actions make sense.
Then, the notion of
 {\it orbifolds}, {\it orbi-bundles} (see also [Th]),
 {\it Kuranishi neighborhoods},
  {\it equivalence of Kuranishi neighborhoods}, and
  {\it Kuranishi structures} in [F-O]
 can be generalized by replacing the model topology/geometry
  in a local chart from domains in ${\Bbb R}^n$
 to objects in ${\cal C}$, with diffeomorphisms
  (resp.\ embeddings; bundle isomorphisms, bundle embeddings)
 that appear in the data of gluing replaced by isomorphisms between
  (resp.\ embeddings of, isomorphisms of bundles over,
          embeddings of bundles over) objects in ${\cal C}$.

\bigskip

\noindent
{\bf Definition 5.1.1 [Kuranishi neighborhood-in-${\cal C}$].}
{\rm
 Let $M$ be a Hausdorff topological space and
  ${\cal C}$ be a category of topology/geometry.
 A {\it Kuranishi neighborhood-in-${\cal C}$} of $p\in M$
  is a $5$-tuple $(V_p,\Gamma_{V_p}, E_{V_p}\,; s_p, \psi_p)$
   (collectively denoted also by $V_p$ for simplicity of notation)
   such that
  \begin{itemize}
   \item[(1)] \parbox{11em}{[$\,${\it neighborhood model}$\:$]}
    \parbox[t]{27em}{$V_p$ is an object in ${\cal C}$,
    $\;\Gamma_{V_p}$ is a finite group that acts on $V_p$
     (as isomorphisms in ${\cal C}$) effectively;
    $\Gamma_{V_p}$ is called the {\it structure group}
     of the Kuranishi neighborhood;}

   \item[(2)] \parbox{11em}{[$\,${\it obstruction bundle}$\,$]}
    $E_{V_p}$ is a $\Gamma_{V_p}$-equivariant vector bundle over $V_p$;

   \item[(3)] \parbox{11em}{[$\,${\it Kuranishi map}$\,$]}
    \parbox[t]{27em}{$s_p: V_p\rightarrow E_{V_p}$
    is a $\Gamma_{V_p}$-equivariant
    continuous section of $E_{V_p}\,$;}

   \item[(4)] \parbox{11em}{[$\,${\it local coordinate map}$\,$]}
    \parbox[t]{27em}{$\psi_p: s_p^{-1}(0)\rightarrow M$
    is a continuous map which induces a homeomorphism from
    $s_p^{-1}(0)/\Gamma_{V_p}$ to a neighborhood of $p$ in $M$.}
  \end{itemize}

 Two Kuranishi neighborhoods-in-${\cal C}$
   $(V_{1,p},\Gamma_{V_{1,p}}, E_{V_{1,p}}\,; s_{1,p}, \psi_{1,p})$,
   $(V_{2,p},\Gamma_{V_{2,p}}, E_{V_{2,p}}\,; s_{2,p}, \psi_{2,p})$
  of $p\in M$ are said to be {\it equivalent},
  in notation\\
   $(V_{1,p},\Gamma_{V_{1,p}}, E_{V_{1,p}}\,; s_{1,p}, \psi_{1,p})
    \sim (V_{2,p},\Gamma_{V_{2,p}}, E_{V_{2,p}}\,; s_{2,p}, \psi_{2,p})$,
 if
  \begin{itemize}
   \item[(1)]
    $\dimm V_{1,p} - \rank E_{V_{1,p}}
      = \dimm V_{2,p} - \rank E_{V_{2,p}} =: d\,$;

   \item[(2)]
     their exists another Kuranishi neighborhood-in-${\cal C}$
      $(V_p,\Gamma_{V_p}, E_{V_p}\,; s_p, \psi_p)$ of $p$
     such that
     \begin{itemize}
      \item[$\cdot$]
       $\dimm V_p - \rank E_{V_p}=d$,

      \item[$\cdot$]
      there exists
       a group homomorphism
        $h_i:\Gamma_{V_{i,p}}\rightarrow \Gamma_{V_p}$ and
       an $h_i$-equivariant vector-bundle embedding
        $\hat{\phi}_i/\phi_i:
          (E_{V_{i,p}}|_{V_{i,p}^{\flat}})/V_{i.p}^{\flat}
                                    \rightarrow E_{V_p}/V_p$
        of the restriction of $E_{V_{i,p}}$ to a neighborhood
        $V_{i,p}^{\flat}$ of $\psi_{i,p}^{-1}(p)$ in $V_{i,p}$
       so that
        $\hat{\phi}_i\circ s_{i,p} = s_p\circ \phi_i$
         on $V_{i,p}^{\flat}$  and
        $\psi_{i,p}=\psi_p\circ \phi_i$
         on $s_{i,p}^{-1}(0)\cap V_{i,p}^{\flat}\,$;
       $i=1$, $2$.
     \end{itemize}

    \item[]
    (When this happens, we say that
     $(V_p,\Gamma_{V_p}, E_{V_p}\,; s_p, \psi_p)$ {\it dominates}
      $(V_{i,p},\Gamma_{V_{i,p}}, E_{V_{i,p}}\,; s_{i,p}, \psi_{i,p})$
        or
     $(V_{i,p},\Gamma_{V_{i,p}}, E_{V_{i,p}}\,; s_{i,p}, \psi_{i,p})$
      is {\it subordinate to}
      $(V_p,\Gamma_{V_p}, E_{V_p}\,; s_p, \psi_p)$
     via $(h_i, \phi_i, \hat{\phi}_i)$, $i=1,2$.)
  \end{itemize}
} % end-definition

\bigskip

The following definition of Kuranishi structure is modified from
 [F-O-O-O: A.2.1.5 - A.2.1.11],
 [Liu(C): Definition 6.3], and [Th: Sec.~13.2].
It is based on the original definition of orbifolds ([Sat] and [Th])
 and the notion of a ``good coordinate system"
 ([F-O: Definition 6.1]) extracted from a Kuranishi structure
 that is originally defined in a functorially-more-natural and
 closer-to-stack way in [F-O: Definition 5.3]
  (if [F-O: Definition 5.2] is replaced by the equivalence relation
   $\sim$ in Definition 5.1.1 above).
    % Definition [Kuranishi neighborhood-in-${\cal C}$]

\bigskip

\noindent
{\bf Definition 5.1.2 [Kuranishi structure-in-${\cal C}$].} {\rm
 Let $M$ be a Hausdorff topological space and
  ${\cal C}$ be a category of topology/geometry.
 A {\it Kuranishi structure-in-${\cal C}$} $\,{\cal K}$ on $M$
  consists of the following data/assignment:
  \begin{itemize}
   \item[(1)]
    a system
     $$
      {\frak N}^{(0)}\;
      :=\; \{(V_p, \Gamma_{V_p}, E_{V_p}\,; s_p, \psi_p)\}_{p\in M}
     $$
    of Kuranishi neighborhoods-in-${\cal C}$, one for each $p\in M\,$;

   \item[(2)]
    a system
     $$
      {\frak N}^{(1)}\;
      :=\; \{(V_{qp}; h_{qp}, \phi_{qp}, \hat{\phi}_{qp})\}_{p,q}
     $$
    of $4$-tuple {\it transition data}
    $(V_{qp}, h_{qp}, \phi_{qp}, \hat{\phi}_{qp})$,
     one each pair $(p,q)$ with $p\in M$ and $q\in \psi_p(s_p^{-1}(0))$,
    such that
    \begin{itemize}
     \item[$\cdot$] ({\it transition function}):
      $V_{qp}$ is an open neighborhood of $\psi_q^{-1}(q)$ in $V_q$,
      $h_{qp}:\Gamma_{V_q}\rightarrow \Gamma_{V_p}$
       is an injective group homomorphism,
      $\hat{\phi}_{qp}/\phi_{qp}:
        (E_{V_q}|_{V_{qp}})/V_{qp} \rightarrow E_{V_p}/V_p$
       is an $h_{qp}$-equivariant vector-bundle embedding
       such that
        $(V_p, \Gamma_{V_p}, E_{V_p}\,; s_p, \psi_p)$
        dominates the restriction of
        $(V_q, \Gamma_{V_q}, E_{V_q}\,; s_q, \psi_q)$ to $V_{pq}$
        via $(h_{qp}, \phi_{qp}, \hat{\phi}_{pq})$;

     \medskip
     \item[$\cdot$] ({\it orbifold cocycle condition}):
      if $r\in \psi_q(s_q^{-1}(0)\cap V_{qp})$,
       then there exists a $\gamma\in\Gamma_{V_p}$
       such that
         $\phi_{qp}\circ \phi_{rq} = \gamma\,\phi_{rp}$
          on a neighborhood $V_{rqp}$ of $\psi_r^{-1}(r)$ in $V_r$,
         $\hat{\phi}_{qp}\circ \hat{\phi}_{rq}
            = \gamma\,\hat{\phi}_{rp}$
          over $V_{rqp}$,  and
         $h_{qp}\circ h_{rq}(g) = \gamma\cdot h_{rp}(g)\cdot\gamma^{-1}$
          for each $g\in \Gamma_{V_r}$.
    \end{itemize}
  \end{itemize}
 If furthermore ${\cal C}$ allows a well-defined notion of
  dimensions to its objects and we require that
   $\dimm V_p - \rank E_{V_p}$ be a constant $d$ independent of $p$
   in the above data,
 then we say that the Kuranishi structure-in-${\cal C}$ $\,{\cal K}$
  on $M$ has {\it virtual dimension} $d$.

 Two Kuranishi structures-in-${\cal C}$
  $$
   {\cal K}_1\;
    =\; \left(
          {\frak N}^{(0)}_1
           = \{ (V_{1,p}, \Gamma_{V_{1,p}}, E_{V_{1,p}}\,;
                 s_{1,p}, \psi_{1,p}) \}_{p\in M}\;,\;
          {\frak N}^{(1)}_1
           = \{( V_{1,qp};
                 h_{1,qp}, \phi_{1,qp}, \hat{\phi}_{1,qp}) \}_{p,q}
        \right)\,,
  $$
  $$
   {\cal K}_2\;
    =\; \left(
          {\frak N}^{(0)}_2
           = \{ (V_{2,p}, \Gamma_{V_{2,p}}, E_{V_{2,p}}\,;
                 s_{2,p}, \psi_{2,p}) \}_{p\in M}\;,\;
          {\frak N}^{(1)}_2
           = \{( V_{2,qp};
                 h_{2,qp}, \phi_{2,qp}, \hat{\phi}_{2,qp}) \}_{p,q}
        \right)
  $$
  on $M$ are said to be {\it equivalent},
  in notation ${\cal K}_1 \sim {\cal K}_2$,
 if there exist another Kuranishi structure-in-${\cal C}$ on $M$
  $$
   {\cal K}\;
    =\; \left(
          {\frak N}^{(0)}
           = \{ (V_p, \Gamma_{V_p}, E_{V_p}\,;
                                    s_p, \psi_p) \}_{p\in M}\;,\;
          {\frak N}^{(1)}
           = \{( V_{qp};h_{qp}, \phi_{qp}, \hat{\phi}_{qp}) \}_{p,q}
        \right)
  $$
  and a system of triples of
   $(\mbox{group}\,,\,\mbox{space}\,,\,\mbox{bundle})$-embedding
  $$
   \left\{ \left(
            h_{i,p}:\Gamma_{V_{i,p}}\rightarrow \Gamma_{V_p}\;,\;
            \phi_{i,p}:V_{i,p}^{\flat}\rightarrow V_p\;,\;
            \hat{\phi}_{i,p}: E_{V_{i,p}}|_{V_{i,p}^{\flat}}
               \rightarrow E_{V_p}
           \right)
    \right\}_{p\in M}\,,
  $$
   where
    $V_{i,p}^{\flat}$ is a neighborhood of
     $\psi_{i,p}^{-1}(p)$ in $V_{i,p}$  and
    $\hat{\phi}_{i,p}$ covers $\phi_{i,p}$,
  such that
  \begin{itemize}
   \item[$\cdot$] ({\it morphism between Kuranishi neighborhoods})\\
    $(V_p, \Gamma_{V_p}, E_{V_p}\,; s_p, \psi_p)$  dominates
    $(V_{i,p}, \Gamma_{V_{i,p}}, E_{V_{i,p}}\,; s_{i,p}, \psi_{i,p})$
    via $(h_{i,p},\phi_{i,p},\hat{\phi}_{i,p})\,$;

   \item[$\cdot$] ({\it compatibility with gluing})\hspace{1em}
     \parbox[t]{24em}{$h_{i,p}\circ h_{i,qp} = h_{qp}\circ h_{i,q}\,$,
      $\; \phi_{i,p}\circ \phi_{i,qp} = \phi_{qp} \circ \phi_{i,q}$
       on $V_{i,qp}$,\\
      $\hat{\phi}_{i,p}\circ \hat{\phi}_{i,qp}
        = \hat{\phi}_{qp} \circ \hat{\phi}_{i,q}$
       over $V_{i,qp}\,$;}
  \end{itemize}
  $i=1$, $2$.
} % end-definition

\bigskip

\noindent
{\it Remark 5.1.3 $[$orbifold cocycle condition$]$.}
 Though we are not generally looking at a space locally modelled
  on some ${\Bbb R}^n$ modulo faithful finite group actions
  as in the definition of an orbifold,
 the fact that
  all the maps $h_{qp}$ on Kuranishi neighborhoods are
   regarded as being defined up to composition with elements
   in the structure finite group $\Gamma_{V_p}$ and
  the morphism $h_{qp}$ of the structure groups are defined
   up to a conjugation in $\Gamma_{V_p}$
  remain to hold in the definition of Kuranishi
  structure-in-${\cal C}$.
 The expression of the compatibility of gluings via
  the transitions functions $\{(\phi_{qp},\hat{\phi}_{qp})\}_{p,q}$
  in terms of the orbifold cocycle condition,
  rather than the ordinary cocycle condition,
  reflects particularly this fact.
 It is in such form that the setting re-phrases the gluing
  in a Deligne-Mumford stack.

\bigskip

We should remark that
a Hausdorff topological space with a Kuranishi structure is
 a topological analogue to a Deligne-Mumford moduli stack with
 a perfect tangent-obstruction complex ([B-F] and [L-T1]) and
 a coarse moduli space.

\bigskip

\noindent
{\bf Example/Definition 5.1.4 [Kuranishi structure with corners].}
 Let ${\cal C}$ be the category of smooth manifolds with corners,
  locally modelled on open sets in some
  ${\Bbb R}^{n_1}\times ({\Bbb R}_{\ge 0})^{n_2}$,
  or more generally,
   ${\Bbb R}^{n_1}\times(\,\mbox{cone in ${\Bbb R}^{n_2}$)}$,
  ($n_1$ and $n_2$ are allowed to vary).
 This gives the notion of {\it Kuranishi structures with corners}
  in [F-O-O-O: Sec.~A.2] and [Liu(C): Sec.~6.1].

\bigskip

\noindent
{\bf Example/Definition 5.1.5 [family Kuranishi structure].}
 Let
  $M$ be a Hausdorff topological space fibered over
   a base Hausdorff topological space $B$,
   in notation $\pi:M\rightarrow B$ or $M/B$, and
  ${\cal C}$ be a category of topological spaces all of
   whose objects and morphisms are over $B$
   (as in the category of schemes over a base scheme
    in algebraic geometry).
 A ({\it family}) {\it Kuranishi structure-in-${\cal C}$ on $M/B$}
  is a Kuranishi structure-in-${\cal C}$ ${\cal K}$ on $M$,
  for which all the data in Definition~5.1.1 and Definition 5.1.2
  are over $B$.
 By construction, there is a natural morphism
  $\tilde{\pi}: {\cal K}\rightarrow B$,
  which restricts to the defining map $V_p\rightarrow B$
   on each Kuranishi neighborhood $V_p$.
 The {\it fiber} ${\cal K}_b:=\tilde{\pi}^{-1}(b)$ of ${\cal K}$
  over $b\in B$ gives a Kuranishi structure-in-${\cal C}_b$ on
  the fiber $M_b:=\pi^{-1}(b)$,
  where ${\cal C}_b$ is the category whose objects and morphisms
   are from taking the restriction of objects and morphisms
   in ${\cal C}$ to over $b$.
 We will denote such ${\cal K}$ on $M$ by ${\cal K}/B$ on $M/B$
  when the family notion is emphasized.
 If, furthermore, there exists an open dense subset $B_0$ of $B$
  and a category $C^{\prime}$ of topology/geometry such that
  each ${\cal K}_b$ is a Kuranishi structure-in-${\cal C}^{\prime}$
  on $M_b$ for $b\in B_0$,
 then we say that
  ${\cal K}/B$ has {\it general fibers in ${\cal C}^{\prime}$} and
  $M/B$ has general fibers
  with a Kuranishi structure-in-${\cal C}^{\prime}$.

\bigskip

\begin{flushleft}
{\bf Morphisms and fibered product.}
\end{flushleft}
In any category of geometry, once the geometric objects are
 defined, the notion of {\it morphisms} and {\it fibered products}
 between them have to be defined accordingly/compatibly as well
 since these two are the foundation of many other notions and
 constructions.
We will postpone them until Sec.~7.1, where we will define these
 two notions in a way that works for
 the specific type of topological spaces-with-a-Kuranishi-structure
 from the current moduli problem.
We won't need them until then.

\bigskip

\begin{flushleft}
{\bf The category ${\cal C}_{\spsccw}$.}
\end{flushleft}
We now describe the category ${\cal C}_{\spsccw}$ over the complex
 line ${\Bbb C}$ in which our Kuranishi structure will model.
An object in ${\cal C}_{\spsccw}$ is a specific kind of
 {\it s}tratified {\it p}iecewise-{\it s}mooth-with-{\it c}orners
 topological space with complex {\it CW}-complex singularities
 and is fibered over ${\Bbb C}$ with smooth-with-corner fibers
 except at $0$, constructed as follows.

First, we introduce
 a complex stratified space $\Xi_{(\vec{s}_0,\,\ldots,\,\vec{s}_k)}$
 over the complex line ${\Bbb C}$.
Let
 $\vec{s}_i=(s_{i1},\,\ldots,\, s_{i, I_i})\in {\Bbb N}^{I_i}$,
 $\vec{\mu}_i =(\mu_{i1},\,\ldots,\,\mu_{i, I_i})\in {\Bbb C}^{I_i}$,
  for $i=0,\,\ldots,\,k$, and
 $\vec{\lambda}=(\lambda_0,\,\ldots,\,\lambda_k)\in B[k]:={\Bbb C}^{k+1}$.
As an affine variety, $\Xi_{(\vec{s}_0,\,\ldots,\,\vec{s}_k)}$
 is defined as the subvariety in
 ${\Bbb C}^{\,(I_0\,+\,\cdots\,+I_k) + (k+1)}\,$:
 $$
  \Xi_{(\vec{s}_0,\,\ldots,\,\vec{s}_k)}\;
   =\; \{(\vec{\mu}_0,\,\ldots,\,\vec{\mu}_k; \vec{\lambda})\,:\;
           {\mu_{ij}}^{s_{ij}}\,=\,\lambda_i,\;
           i=0,\,\ldots,\, k, \; j=1,\,\ldots,\, I_i\, \}\,.
 $$
It has complex dimension
 $\dimm_{\scriptsizeBbb C}
   (\Xi_{(\vec{s}_0,\,\ldots,\,\vec{s}_k)}) = k+1$.
The projection
 $$
  {\Bbb C}^{\,(I_0\,+\,\cdots\,+I_k) + (k+1)}\;
    \longrightarrow\; {\Bbb C}^{k+1}\,, \hspace{1em}
  (\vec{\mu}_0,\,\ldots,\,\vec{\mu}_k; \vec{\lambda})\;
    \longmapsto\; \vec{\lambda}\,,
 $$
 induces a finite flat morphism from
 $\Xi_{(\vec{s}_0,\,\ldots,\,\vec{s}_k)}$ onto ${\Bbb C}^{k+1}$
 of degree $\prod_{i=0}^k\prod_{j=1}^{I_i}\,s_{ij}$,
 and is \'{e}tale over the complement
 $\{\vec{\lambda}\,:\;\lambda_i\ne 0,\; i=0,\,\ldots,\,k\}$
 of coordinate subspaces in $B[k]$.
After the post-composition with the flat morphism
 ${\mathbf p}[k]:B[k]\rightarrow {\Bbb C}$,
 $\vec{\lambda}\mapsto \lambda_0\,\cdots\,\lambda_k$ (cf.\ Sec.~1.1.1),
one has a flat morphism
 $p:\Xi_{(\vec{s}_0,\,\ldots,\,\vec{s}_k)}\rightarrow {\Bbb C}$
 that is smooth over ${\Bbb C}-\{0\}$.

From the system of defining equations,
 $\Xi_{(\vec{s}_0,\,\ldots,\,\vec{s}_k)}$ is the product
 $\prod_{i=0}^k\,\Xi_{\vec{s}_i}$,
 where
 $$
  \Xi_{\vec{s}_i}\; =\;
   \{ (\vec{\mu}_i;\lambda_i)\,:\;
       {\mu_{ij}}^{s_{ij}}\,=\,\lambda_i,\; j=1,\,\ldots,\, I_i\, \}\,.
 $$
$\Xi_{\vec{s}_i}$ is the fibered product of the morphisms
 $f_j: {\Bbb C}\rightarrow {\Bbb C}$, $z\rightarrow z^{s_{ij}}$,
 $j=1,\,\ldots,\,I_i$.
Its ($C^0$-)topology is thus the gluing $\vee_{n_i}{\Bbb C}$
 of $n_i$ copies of ${\Bbb C}$'s at the origin,
 where $n_i$ is the number of orbits in the group
 $({\Bbb Z}/s_{i1}{\Bbb Z})
    \oplus\,\cdots\,\oplus({\Bbb Z}/s_{i,I_i}{\Bbb Z})$
 under the action generated by the translation
 $(e_1,\,\ldots,\,e_{I_i})\mapsto (e_1+1,\,\ldots,\,e_{I_i}+1)$.
It follows that
\begin{itemize}
 \item[$\cdot$]
 the ($C^0$-)topology of $\Xi_{(\vec{s}_0,\,\ldots,\,\vec{s}_k)}$
  is the product $\prod_{i=0}^k\,(\vee_{n_i}{\Bbb C})$,
  which is a gluing of $n_0\,\cdots\,n_k$ copies of
  ${\Bbb C}^{k+1}$ along proper coordinate subspaces, and

 \item[$\cdot$]
  $p^{-1}(t)$, $t\ne 0$, is a disjoint union of
   $n_0\,\cdots\,n_k$ copies of $({\Bbb C}^{\times})^k$
   with each $({\Bbb C}^{\times})^k$ \'{e}tale over
   ${\mathbf p}[k]^{-1}(t)\subset B[k]$.
\end{itemize}

Denote by $H_I$ the coordinate subspace of $B[k]$
 whose points have coordinates $\lambda_i=0$ for $i\in I$.
It follows from the topology of $\Xi_{(\vec{s}_0,\,\ldots,\,\vec{s}_k)}$
 that
 $p^{-1}(0)\,
   =\, {\mathbf p}[k]^{-1}(H_{\{0\}})\,\cup\,\cdots\,
        \cup\, {\mathbf p}[k]^{-1}(H_{\{k\}})$
 has $n_0\,\cdots\,n_k\,(\frac{1}{n_0}\,+\,\cdots\,+\,\frac{1}{n_k})$
 irreducible components, with
 $n_0\,\cdots\,n_{i-1}\,n_{i+1}\,\cdots\,n_k$ of them
 contained in ${\mathbf p}[k]^{-1}(H_{\{i\}})$.
Each of these irreducible components has ($C^0$-)topology
 isomorphic to ${\Bbb C}^k$.
Let $[{\mathbf p}[k]^{-1}(H_{\{i\}})]_0$ be the formal sum of
 the subvarieties of $\Xi_{(\vec{s}_0,\,\ldots,\,\vec{s}_k)}$
 that appear as irreducible components of $p^{-1}(0)$.
It follows from the defining equation, $\lambda_i=0$, of
 ${\mathbf p}[k]^{-1}(H_{\{i\}})$ in
 $\Xi_{(\vec{s}_0,\,\ldots,\,\vec{s}_k)}$ that
 $$
  [p^{-1}(t)]\,,\; t\ne 0\,,\;
  =\; [p^{-1}(0)]\;
  =\; \sum_{i=0}^k\,
      (s_{i1}\,\cdots\,s_{i,I_i})\,[{\mathbf p}[k]^{-1}(H_{\{i\}})]_0
 $$
in the Chow group $A_k(\Xi_{(\vec{s}_0,\,\ldots,\,\vec{s}_k)})$.

The composition of the projection map with $p$
 $$
  \Xi_{(\vec{s}_0,\,\ldots,\,\vec{s}_k)}
   \times {\Bbb R}^{n_1}\times ({\Bbb R}_{\ge 0})^{n_2}\;
   \longrightarrow\; \Xi_{(\vec{s}_0,\,\ldots,\,\vec{s}_k)}\;
   \stackrel{p}{\longrightarrow}\; {\Bbb C}
 $$
 gives a flat\footnote{For non-algebraic-geometers: the
                                   flatness of the fibration
                                   $p:\Xi_{(\vec{s}_0,\,\ldots,\,\vec{s}_k)}
                                     \rightarrow {\footnotesizeBbb C}$
                                    is in the sense of morphisms of
                                    schemes over the ground field
                                    ${\footnotesizeBbb C}$.
                                   The fiber subscheme $p^{-1}(0)$ over
                                    $0\in{\footnotesizeBbb C}$
                                    is non-reduced;
                                   each of the irreducible components
                                    of $p^{-1}(0)$ carries a multiplicity
                                    in the sense of [Fu].}
 fibration of
 $\Xi_{(\vec{s}_0,\,\ldots,\,\vec{s}_k)}
     \times {\Bbb R}^{n_1}\times ({\Bbb R}_{\ge 0})^{n_2}$
 over ${\Bbb C}$.
Let ${\cal C}_{\spsccw}$ be the category of
 Hausdorff topological spaces fibered over ${\Bbb C}$
  that are locally modelled on an open set in
   $\Xi_{(\vec{s}_0,\,\ldots,\,\vec{s}_k)}
     \times {\Bbb R}^{n_1}\times ({\Bbb R}_{\ge 0})^{n_2}$
    as stratified piecewise-smooth-with-corner spaces,
   with the gluing maps isomorphisms over ${\Bbb C}$.
Here, $k$, $(\vec{s}_0,\,\ldots,\,\vec{s}_k)$, $n_1$, $n_2$
 are all allowed to vary.

\bigskip

We can now state the main theorem of the current work,
 which gives the foundation of the degeneration axiom and
 the gluing axiom of open Gromov-Witten invariants.
Its proof takes Sec.~5.3 - Sec.~5.4.

\bigskip

\noindent
{\bf Theorem 5.1.6 [family Kuranishi structure on
        $\overline{\cal M}_{(g,h),(n,\vec{m})}
                  (W/B,L\,|\,[\beta],\vec{\gamma},\mu)$].}
{\it
 There is a family Kuranishi structure ${\cal K}$ on
  $\overline{\cal M}_{(g,h),(n,\vec{m})}
    (W/B,L\,|\,[\beta],\vec{\gamma},\mu)$ over $B$
  that is modelled in ${\cal C}_{\spsccw}/{\Bbb C}$,
  $($recall that $B\subset {\Bbb C}$$)$.
 ${\cal K}/B$
  is fiberwise of the same virtual dimension
   $$
    \vdim^{\fiber}
     \overline{\cal M}_{\tinybullet}(W/B,L\,|\,\tinybullet\;)/B\;
    :=\; \mu + (N-3)(2-2g-h)
         + 2n + (m_1 +\,\cdots\, + m_h)\,,
   $$
   where $2N$ is the dimension of $X$ $($as a fiber of $W/B$$)$.
 The family Kuranishi neighborhood-in-${\cal C}_{\spsccw}$
  $(V_{\rho},\Gamma_{V_{\rho}},E_{V_{\rho}};s_{\rho},\psi_{\rho})$
  at $\rho=[f:(\Sigma,\partial\Sigma)\rightarrow (Y_{[k]},L_{[k]})]$
  has $V_{\rho}/B$ \footnote{Here we use $V_{\rho}/B$ to indicate that
                                     there is built-in map
                                      $V_{\rho}\rightarrow B$.
                                     The map is not necessarily
                                      surjective.}
  isomorphic to a neighborhood of the origin in the total space
  of the flat fibration
  $$
   \left.
    \left(\,
      \Xi_{(\vec{s}_0,\,\ldots,\,\vec{s}_k)}
      \times {\Bbb R}^{n_1} \times ({\Bbb R}_{\ge 0})^{n_2}\,
     \right)
   \right/{\Bbb C}\,,
  $$
  where
   \begin{itemize}
    \item[$\cdot$]
     $\vec{s}_i$ is the contact order of $f$ along $D_i$
     at the ordered set of distinguished nodes in $f^{-1}(D_i)$,
     $i=0,\,\ldots\,,\,k\,$,
    $\;($and recall that
         $\dimm \Xi_{(\vec{s}_0,\,\ldots,\,\vec{s}_k)}=2k+2$$)$;

    \item[$\cdot$]
     $n_1\;=\;$
      $\vdim^{\fiber}
        \overline{\cal M}_{\tinybullet}(W/B,L\,|\,\tinybullet\;)/B\,
       + \dimm E_{\rho} - (2k+n_2)\,$; and

    \item[$\cdot$]
     $n_2\;=\;$
      the total number of boundary nodes and free marked points
      that land on $\partial\Sigma$.
   \end{itemize}
 The homeomorphism-type $\{Y_{[k^{\prime}]}\}_{0\le k^{\prime}\le k}$
  of the targets of maps gives a $\Gamma_{V_{\rho}}$-invariant
  stratification
  $\{S_{k^{\prime}}\}_{0\le k^{\prime}\le k}$
  on the fiber $V_{\rho;0}$ of $V_{\rho}/B$ over $0\in B$;
 each connected component of $S_{k^{\prime}}$
  is a manifold of codimension $2k^{\prime}$ in $V_{\rho;0}$.
 This stratification coincides with the induced stratification
  on $V_{\rho;0}$ from the stratification\footnote{Recall the map
                             $\Xi_{(\vec{s}_0,\,\ldots,\,\vec{s}_k)}
                               \rightarrow {\footnotesizeBbb C}^{k+1}$.
                             The coordinate-subspace stratification
                              of ${\footnotesizeBbb C}^{k+1}$ induces
                              a stratification on
                              $\Xi_{(\vec{s}_0,\,\ldots,\,\vec{s}_k)}$.}
  of $\Xi_{(\vec{s}_0,\,\ldots,\,\vec{s}_k)}$.
 %%%%%%%%%%%%%%%%%%%%%%%%%%%%%
 %
 % ${\cal K}$ is fiberwise orientable
 %   if $L$ is spin or if $h=1$ and $L$ is relative spin
 %   {\rm (}i.e.\ $L$ is orientable and
 %            $\,w_2(T_{\ast}L)=c|_{L}$
 %             for some $c\in H^2(X;{\Bbb Z}/2{\Bbb Z})\,${\rm )}.
 %
 %%%%%%%%%%%%%%%%%%%%%%%%%%%%%
} % end-theorem

% \bigskip
%
% The construction follows mainly the construction in [MD-S1] and [F-O]
%  for the moduli space of stable maps from closed Riemann surfaces
%  $\Sigma$ to $(X,\omega,J)$.
% The modification of [F-O] to take care of the boundary nodes in
%  the case of open stable maps with $\partial\Sigma$ mapped to
%  a Lagrangian submanifold $L\subset X$ is given in [Liu(C)].
% Another modification to [F-O] is required in our case to take care
%  of the nodes that are mapped to the singular locus of $Y[k]$
%  in Sec.~???. (Such nodes are all interior nodes of $\Sigma$.)
% The related discussion in [L-R] (see also [I-P2]),
%  in the symplectic infinite-end form and in the language of
%  virtual neighborhoods of [Ru],
%  can be converted to the Fukaya-Ono format to take care of this.
% Insights on how pieces from stable maps with different $Y[k]$ targets
%  should glue together are hidden in [I-P2], in yet another symplectic
%  setting.
% They can be incorporated into the Fukaya-Ono format as well.
% The algebro-geometric construction parallel to those of [L-R]
%  and [I-P1, I-P2] is given in [Li1, Li2].
% There is also the work [E-G-H] from the symplectic-field-theory aspect.
% Both provide further insights to the construction
%  in the symplectic/almost-complex category.

\bigskip

\subsection{Local transversality and
            locally regular almost-complex structures.}

There are three types of local transversality issues in
 our moduli problem that have to be understood
 before one can choose a good obstruction space to work on:
 (for a {\it fixed} $J$)
 \begin{quote}
  \parbox[t]{2.7em}{(T1)}
   local surjectivity of $D_f\bar{\partial}_J$,

  \parbox[t]{2.7em}{(T2)}
   local transversality of evaluation maps, and

  \parbox[t]{2.7em}{(T3)}
   \parbox[t]{30em}{local transversality
    of the contact order condition along $D$ and
    local transversality of the pre-deformability conditions
     at a distinguished node with a specified contact order.}
 \end{quote}
Global such issues have been discussed in related symplectic
 Gromov-Witten theories,
 e.g.\ [MD-S1; MD-S3], [R-T1; R-T2], and
       [I-P1; I-P2] (particularly for Item (3)).
In dealing with transversality issues, it is a standard procedure by now
 that one first show the sought-for transversality properties on
 the related universality moduli space ${\cal U}$$\overline{\cal M}$
 of extended tuples $(J, f:\Sigma\rightarrow X)$
 (or $(J,\nu,f:\Sigma\rightarrow X)$ where $\nu$ is an additional
     perturbation in [I-P1], [R-T1], [R-T2])
 that contains a choice of an almost-complex structure $J$ and
 a $J$-holomorphic map $f$ of a fixed class.
One shows that ${\cal U}$$\overline{\cal M}$ is a smooth Banach manifold
 and then apply the Sard-Smale Theorem to the fibration of
 ${\cal U}$$\overline{\cal M}$ over the Banach manifold ${\cal J}$
 of allowed almost-complex structures to obtain the sought-for
 transversality property for the fiber moduli space $\overline{\cal M}^J$
 over a regular value $J\in {\cal J}$ of the fibration
 ${\cal U}$$\overline{\cal M}\rightarrow {\cal J}$.
A good feature for such a setting is that the moduli space
 $\overline{\cal M}^J$ for $J$ regular is a smooth orbifold of
 the correct dimension as expected from deformation theory.
However, it can happen that no regular $J$'s are integrable.

In our current moduli problem, approached along [F-O],
the effect of allowing $J$ to vary to obtain
 a sought-for transversality property is absorbed into a choice of
 a large enough subspace $E$ in
 $L^p(\Sigma; \Lambda^{0,1}\Sigma\otimes_J f^{\ast}T_{\ast}X)$
 so that its preimage
  $(D_f\bar{\partial}_J)^{-1}(E)$ in
  $W^{1,p}(\Sigma,\partial\Sigma;
        f^{\ast}T_{\ast}X, (f|_{\partial\Sigma})^{\ast}T_{\ast}L)$
 can fit into the related transversality statement.
This is because the infinitesimal deformations of $J$ give rise
 to elements in
 $L^p(\Sigma; \Lambda^{0,1}\Sigma\otimes_J f^{\ast}T_{\ast}X)$
 as well, after the pre-composition with $df\circ j$;
 and, similarly, for the additional $\nu$ in [I-P1], [R-T1], [R-T2].
The larger-than-expected dimension and the possibly-worse singularities
 of $\overline{\cal M}^J$ for a fixed $J$ that is not regular now have
 to be compensated in the construction of Kuranishi structure.
However, in doing so, we may retain a good $J$ to work on,
The latter can be important for other parts in the theory,
 cf.\ Example 5.2.3.
      % Example [complex structure]

With these highlights in mind, we now give the precise respective
 statement of Conditions (T1), (T2), and (T3) in the setting of
 Kuranishi structures.
The domain unit disc or half unit-disc in the following discussion
 is considered {\it fixed}.

\bigskip

\noindent
(1) {\it Local surjectivity of $D_f\bar{\partial}_J$.}
This condition says that:
 \begin{quote}
  (T1) \hspace{1em}
   The map
    $$
     D_f\bar{\partial}_J\; :\;
      W^{1,p}(\Sigma,\partial\Sigma;
        f^{\ast}T_{\ast}X, (f|_{\partial\Sigma})^{\ast}T_{\ast}L)\;
      \longrightarrow\;
      L^p(\Sigma; \Lambda^{0,1}\Sigma\otimes_J f^{\ast}T_{\ast}X)
    $$
    is {\it surjective} for any non-constant $J$-holomorphic maps on
    the unit disc $f:D^2:= \{z\in{\Bbb C}:|z|\le 1\}\rightarrow X$ or
    on a half unit disc
    $f: (D^2_+\,,\, \partial_0D^2_+)
         := (\{z\in{\Bbb C}:|z|\le 1,\im(z)\ge 0\},[-1,1])
        \rightarrow (X,L)$.
 \end{quote}

Any almost-complex structure $J$ that is $C^1$ close to a complex
 structure has this property, cf.~Example 5.2.3.
  % Example [complex structure]

\bigskip

\noindent
(2) {\it Local transversality of evaluation maps.}
This condition says that:
 \begin{quote}
  (T2) \hspace{1em}
   Given a $J$-holomorphic map on the marked disc
    $f:(D^2; 0)\rightarrow X$ (resp.\ on the marked half unit disc
    $f:(D^2_+\,,\,\partial_0D^2_+;0)\rightarrow (X,L)$),
   there exists a (finite dimensional) subspace
    $E\subset
     L^p(\Sigma; \Lambda^{0,1}\Sigma\otimes_J f^{\ast}T_{\ast}X)$
   such that the differential of the evaluation map $\ev$ associated
    to the marked point
    $$
     D_f\ev\; :\;  (D_f\bar{\partial}_J)^{-1}(E)\;
        \longrightarrow\; T_{f(0)}X
    $$
   (resp.
    $D_f\ev: (D_f\bar{\partial}_J)^{-1}(E)\rightarrow T_{f(0)}L$)
    is surjective.
 \end{quote}

This is the local Kuranishi statement for [MD-S1: Lemma 6.1.2].
Note that, in the above expression, $D_f\ev$ is defined on
 the whole
 $W^{1,p}(\Sigma,\partial\Sigma;
    f^{\ast}T_{\ast}X, (f|_{\partial\Sigma})^{\ast}T_{\ast}L)$
 by $(D_f\ev)(\xi)=\xi(0)$, where $\Sigma=(D^2,0)$ or $(D^2_+,0)$.

\bigskip

\noindent
(3) {\it Local transversality of the contact order and
         the pre-deformability condition.}
To describe these conditions in the Kuranishi setting, we have to
 introduce the objects from [I-P1] (with a notation change:
 $V$ there = $D$ here):
 %%
 % \marginpar{\raggedright\tiny $\bullet$
 % The role of symplectic structure here is to be clarified.}
 %
 \begin{itemize}
  \item[$\cdot$]
   ${\cal J}^D\,:\;$
   \parbox[t]{33.4em}{the space of pairs $(J^{\prime},\nu^{\prime})$
    where
     $J^{\prime}$ is an admissible almost complex structure on
      the relative pair $(X;D)$ and
     $\nu^{\prime}$ is an element in
       $\Hom(\pi_2^{\ast}T_{\ast}D^2,\pi_1^{\ast}X)$
       (of the lifted bundles on $X\times D^2$)
     that is anti-$J^{\prime}$-linear: $\nu\circ j=-J\circ \nu$,
   (the set of all such $\nu$ will be denoted by
     $\Hom^J(\pi_2^{\ast} T_{\ast}D^2,\pi_1^{\ast}T_{\ast}X)$);}

  \item[$\cdot$]
   $\cal U$${\cal M}\,:\;$
   \parbox[t]{33.4em}{the universal moduli space of
    $(J^{\prime},\nu^{\prime})$-holomorphic maps
    (i.e.\ $(f^{\prime},\phi^{\prime}):D^2\rightarrow X\times D^2$
    such that $\bar{\partial}_{J^{\prime}}f^{\prime}=\nu^{\prime}$)
    for some $(J^{\prime},\nu^{\prime})$.}
 \end{itemize}
Let $(J,0)\in {\cal J}^D$ and $f:D^2\rightarrow (X,D)$ be
 a $J$-holomorphic disc in $X$ with $f^{-1}(D)=s\cdot(0)$.
(we set $\phi=\Id_{D^2}$ for such $f$ by convention.)
Then, [I-P1: Lemma 3.4] implies that
 there is a {\it divisor map} $\divisor$ from a neighborhood of
  $[f]\in {\cal U}$$\cal M$ to the space $\Div^s(D^2)\subset {\Bbb C}^s$
  of degree $s$ divisors on the unit disc $D^2$,
 defined by $f^{\prime}\mapsto {f^{\prime}}^{-1}(D)$.
Let $\End^J(T_{\ast}X)$ be the space of anti-$J$-linear endomorphisms
 of $T_{\ast}X$.
Then, there is a map
 $$
  \begin{array}{cccl}
   T_{(J,0)}{\cal J}^D\;=\;
    \End^J(T_{\ast}X)
    \oplus \Hom^J(\pi_2^{\ast} T_{\ast}D^2, \pi_1^{\ast}T_{\ast}X)
    & \longrightarrow
    & L^p(D^2; \Lambda^{0,1}D^2\otimes_J f^{\ast}T_{\ast}X) \\[.6ex]
   (\delta J\,,\,\delta\nu)
    & \longmapsto
    & \mbox{\large $\frac{1}{2}$}
      (\delta J)\circ df\circ j\,
      -\, \delta\nu    &,
  \end{array}
 $$
 where
  $(\delta J, \delta\nu)$ denotes an infinitesimal deformation
  of $(J,0)$.
Denote the image of the above map by ${\cal H}$.
Then $D_f\divisor$ is defined on the subspace
 $(D_f\bar{\partial}_J)^{-1}({\cal H})$  of
 $W^{1,p}(D^2,f^{\ast}T_{\ast}X)$.
Recall the holomorphic coordinate $z$ on $D^2$ and fix a complex
 normal coordinate to $D$ around $f(0)$ in $D\subset X$
 that is compatible with $J|_{f(0)}$.
For $\xi$ in the subspace $V_0 := (D_f\divisor)^{-1}(0)$ of
 $(D_f\bar{\partial}_J)^{-1}({\cal H})$,
 let $\xi^n$ be its normal component with respect to the normal
 coordinate to $D$.
Then there is a linear $s(0)$-jet-at-$0$ map
$$
 \begin{array}{ccccc}
  \jet^{s(0)}_0   & :
   & V_0
   & \longrightarrow     & {\Bbb C}        \\[.6ex]
  && \xi      & \longmapsto
   &  d^{\,s(0)}\xi^n(0)/{dz}^{s(0)}\,.
 \end{array}
$$

With these preparations, the local transversality of contact order
 condition says that:
 \begin{quote}
  (T3.1) \hspace{1em}
   Given a $J$-holomorphic map $f:D^2\rightarrow (X;D)$
    such that $f^{-1}(D)$ is a divisor $s\cdot(0)$ on $D^2$,
   there exists a (finite dimensional) subspace
     $E\subset {\cal H}\subset
      L^p(D^2; \Lambda^{0,1}D^2\otimes_J f^{\ast}T_{\ast}X)$
    such that
    \begin{itemize}
     \item[$\cdot$]
      $D_f\divisor :
        (D_f\bar{\partial}_J)^{-1}(E) \rightarrow
          T_{s\cdot (0)}\Div^s(D^2) \simeq T_{{\mathbf 0}}{\Bbb C}^s$
      is surjective;

     \item[$\cdot$]
      $\jet^{s(0)}_0$ on the subspace
      $(D_f\divisor)^{-1}(0)$ of $(D_f\bar{\partial}_J)^{-1}(E)$
      is also surjective.
    \end{itemize}
 \end{quote}

This is the local Kuranishi statement for the combination of
 the related part of [I-P1: proof of Lemma 4.2] and [I-P1: Lemma 3.4].

\medskip

For the local transversality of the pre-deformability condition,
 consider first the {\it fixed} unit disc $D^2$ with the marked
 point $0$ and
restrict the above discussion to maps with $0$ sent to $D$ in $X$.
Denote the related universal moduli space by
 ${\cal U}$${\cal M}^0(X;D)$ and let
 $f:(D^2,0)\rightarrow (X,D)$ be a $J$-holomorphic map with
 $f^{-1}(D)=s\cdot(0)$.
Then there are
 the evaluation map $\ev_0:{\cal U}$${\cal M}^0(X;D)\rightarrow D$
  associated to the marked point $0$ and
 the divisor map $\divisor_0$ from a neighborhood of
  $[f]\in {\cal U}$${\cal M}^0(X;D)$ to the space
  $\Div^{s-1}(D^2)\subset {\Bbb C}^{s-1}$,
  defined by $f^{\prime}\mapsto {f^{\prime}}^{-1}(D)-(0)$.
Their differential, $D_f\ev_0$ and $D_f\divisor_0$, are both defined
 on the subspace $(D_f\bar{\partial}_J)^{-1}({\cal H})$  of
 $W^{1,p}(D^2,f^{\ast}T_{\ast}X)$,
 where ${\cal H}$ is from the previous discussion.
Again, for the complex coordinate $z$ on $D^2$ and
 a fixed normal coordinate to $D$ in $X$, has the
 $s(0)$-jet-at-$0$ map $\jet^{s(0)}_0$ from the subspace
 $(D_f\divisor_0)^{-1}(0)$ of $(D_f\bar{\partial}_J)^{-1}({\cal H})$
 to ${\Bbb C}$.

Next consider a pre-deformable $J$-holomorphic map
 $$
  f\,=\, f_1\cup f_2\; :\;  \Sigma\, :=\, D^2_1\cup_0 D^2_2\;
   \longrightarrow\;
    Y\,=\,Y_1\cup_D Y_2
 $$
 of contact order $s$ along $D$ at the distinguished node $0$.
Define
 \begin{eqnarray*}
  \lefteqn{
   W^{1,p}(\Sigma; f^{\ast}T_{\ast}Y)  } \\[.6ex]
   && :=\; \{\, (\xi_1,\xi_2) \in
              W^{1,p}(D^2_1; f_1^{\ast}T_{\ast}Y_1)
              \oplus W^{1,p}(D^2_2; f_2^{\ast}T_{\ast}Y_2)\,:\,
             \xi_1(0)=\xi_2(0)\,\}
 \end{eqnarray*}
 and
 $$
  L^p(\Sigma; \Lambda^{0,1}\Sigma\otimes_J f^{\ast}T_{\ast}Y)\;
  :=\; L^p(D^2_1; \Lambda^{0,1}D^2_1\otimes_J f_1^{\ast}T_{\ast}Y_1)
       \oplus
       L^p(D^2_2; \Lambda^{0,1}D^2_2\otimes_J f_2^{\ast}T_{\ast}Y_2)\,.
 $$
Gluing of evaluation maps and their differential defines
 $$
  D_f\ev_0\; :\;   W^{1,p}(\Sigma; f^{\ast}T_{\ast}Y)\;
    \longrightarrow\; T_{f(0)}D\,.
 $$
Let
 ${\cal H}_i \subset
     L^p(D^2_i; \Lambda^{0,1}D^2_i\otimes_J f_i^{\ast}T_{\ast}Y_i)$
 be the subspace that encodes the infinitesimal deformation of
 $(J,\nu)$ as in the previous discussion and set
${\cal H}
 := {\cal H}_1\oplus {\cal H}_2 \subset
  L^p(\Sigma; \Lambda^{0,1}\Sigma\otimes_J f^{\ast}T_{\ast}Y)$.
Then, gluing of the divisor map $\divisor_{0,i} $, $i=1,\,2$,
       and their differential gives
 $D_f\divisor_0:
   (D_f\bar{\partial}_J)^{-1}({\cal H}) \rightarrow
    T_{(s-1)\cdot(0)}\Div^{s-1}(D^2_1)
                   \oplus T_{(s-1)\cdot(0)}\Div^{s-1}(D^2_2)$.
Again, recall the complex coordinates $z_1$ and $z_2$ on $D^2_1$
 and $D^2_2$ and fixed normal coordinates on $(Y_1,D)$ and $(Y_2,D)$.
 that is compatible with $J|_{f(0)}$.
Then, for $\xi=(\xi_1,\xi_2)$ in the subspace
 $V^{\pd} := (D_f\divisor)^{-1}(0)$ of
 $(D_f\bar{\partial}_J)^{-1}({\cal H})$,
 let $\xi^n=(\xi_1^n,\xi_2^n)$ be its normal component
  with respect to the normal coordinate to $D$.
Then there is a linear $s(0)$-jet-at-$0$ map
 $$
  \begin{array}{ccccc}
   \jet^{s(0)}_0   & :
    & V^{\pd}
    & \longrightarrow     & {\Bbb C}^{\,2}        \\[1ex]
   && \xi      & \longmapsto
    & \left(  d^{\,s(0)}\xi_1^n(0)/dz_1^{s(0)}\,,\,
              d^{\,s(0)}\xi_2^n(0)/dz_2^{s(q)}
     \right)\,.
  \end{array}
 $$

In terms of these, the local transversality of the
 pre-deformability condition says that:
 \begin{quote}
  (T3.2) \hspace{1em}
   Given a pre-deformable $J$-holomorphic map
    $f=f_1\cup f_2: \Sigma:=D^2_1\cup_0 D^2_2\rightarrow
     Y=Y_1\cup_DY_2$ of contact order $s$ along $D$
     at the distinguished node $0$,
   let ${\cal H}:={\cal H}_1\oplus {\cal H}_2$ be the subspace of
    $L^p(\Sigma; \Lambda^{0,1}\Sigma\otimes_J f^{\ast}T_{\ast}Y)
     := L^p(D^2_1; \Lambda^{0,1}D^2_1\otimes_J f_1^{\ast}T_{\ast}Y_1)
        \oplus
        L^p(D^2_2; \Lambda^{0,1}D^2_2\otimes_J f_2^{\ast}T_{\ast}Y_2)$
    from the previous discussion on which
    $D\divisor_0:=D\divisor_{0,1}\oplus D\divisor_{0,2}$ is defined.
   Then, there exists a (finite dimensional) subspace $E\subset {\cal H}$
    such that
    \begin{itemize}
     \item[$\cdot$]
      $D_f\divisor_0:
       (D_f\bar{\partial}_J)^{-1}(E) \rightarrow
       T_{(s-1)\cdot(0)}\Div^{s-1}(D^2)
                    \oplus T_{(s-1)\cdot(0)}\Div^{s-1}(D^2)$
     is surjective;

    \item[$\cdot$]
     denote the subspace $(D_f\divisor_0)^{-1}(0)$ of
      $(D_f\bar{\partial}_J)^{-1}(E)$ by
      $(D_f\bar{\partial}_J)^{-1}(E)^{\pd}$, then
     $\;D_f\ev_0\oplus \jet^{s(0)}_0:
       (D_f\bar{\partial}_J)^{-1}(E)^{\pd}
        \rightarrow T_{f(0)}D\oplus {\Bbb C}^2\;$
      is surjective.
  \end{itemize}
 \end{quote}

This is the local Kuranishi statement for the combination of
 [I-P2: Lemma 3.5] and [I-P1: Lemma 4.2].

Note that, in both (T3.1) and (T3.2),
 though the map $\jet_0^{s(0)}$ depends on the choice
  of a local coordinate around $0$ in the domain and
  a local normal coordinate to $D$ around $f(0)$ in the target,
 the surjectivity condition stated is independent of the choices.
  of such coordinates.

\bigskip

\noindent
{\bf Definition 5.2.1
     [(strongly) locally regular almost-complex structure].}
{\rm
 An almost-complex structure $J$ on $(X,L)$ with $L$
  a maximal totally real submanifold
  (resp.\ on $(X,L;D)$ with $D$ a $\codimm_{\scriptsizeBbb R}$-$2$
   almost-complex submanifold)
 is called {\it locally regular}
 if the transversality conditions $(T1)$, $(T2)$
  (resp.\ in addition, $(T3)$) hold
 for sufficiently small holomorphic discs and half-discs in $X$.
 Such $J$ is called {\it strongly locally regular}
  if, in addition, $E$ in Condition (T2) and Condition $(T3)$,
   can be chosen to be supported in a compact set away from
   the marked point and the distinguished node respectively.
} % end-definition

\bigskip

\noindent
{\it Remark 5.2.2.}
  Condition (T1) is said to be true for all smooth $J$ in
   [F-O: (12.7.3)] and [Liu(C): proof of Lemma 6.18].
  The proof of [MD-S1: Lemma 6.1.2] can be adapted to show that
   Condition (T2) always holds and $E$ can be chosen to be supported
   in a compact set away from the marked point.
  The proof of [I-P1: Lemma 4.2] can be adapted to show that
   Condition (T3.1) also always holds.
  Since the domain $D^2$ is unstable, the perturbation $\nu$ in
   [I-P1: proof of Lemma 4.2] can be set to be $0$.
  The argument in the proof of [I-P1: Lemma 4.2] implies then
   that the $E$ in Condition (T3.1) can be chosen to be $0$.
  Similarly for the case of $Y=Y_1\cup_D Y_2$ and Condition (T3.2).

\bigskip

\noindent
{\bf Example 5.2.3 [complex structure].}
 Let $(X, L; D)$ be a complex manifold $(X,J)$ with
  a maximal totally real submanifold $L$ and
  a smooth divisor $D$.
 Then the local study of [Sie1], [Sik], [Ve]
  implies that Condition (T1) is satisfied and
  a right inverse $Q$ of $D_f\bar{\partial}_J$ is given
   as an singular integral operator.
 Conditions (T2) and (T3) can be directly checked by
  constructing a family of local holomorphic discs or half-discs
  whose associated deformation vectors map surjectively
  to $T_{f(0)}X$, $T_{f(0)}L$, and $T_{s(0)}\Div^s(D^2)$
  respectively, e.g.\ using the local pseudo-automorphism group
  action on $X$ around $f(0)$.
 One can also choose $E$ in Conditions (T2) and (T3) to be $0$
  as long as the holomorphic disc or half-disc is small enough.
 This shows directly that the complex $X$ is strongly locally regular.
 Similarly for a complex manifold-divisor relative pair $(Y;D)$
  and the singular complex space $Y=Y_1\cup_D Y_2$.

\bigskip

\noindent
{\bf Assumption.} {\it
 From now on, we assume that the fixed smooth ($C^{\infty}$)
  almost-complex structure on targets of types $X$, $W[k]/B[k]$,
  $(Y[k];D[k])$, $k\in{\Bbb Z}_{\ge 0}$,
  are all strongly locally regular.
} % end-assumption

\bigskip

\subsection{Construction of family Kuranishi neighborhoods.}

The foundation of the construction is the following two facts,
 applied in a continuous way to Banach-space fibers of a family
 over a finite dimensional base.

\bigskip

\noindent
{\bf Proposition 5.3.0.1 [Newton-Picard iteration].}
([MD-S3: Proposition A.3.4].)
{\it
 Let
  $X$ and $Y$ be Banach spaces,
  $U\subset X$ be an open set, and
  $f:U\rightarrow Y$ be a continuous differentiable map.
 Let $x_0\in U$ be such that $D := df(x_0):X\rightarrow Y$
  is surjective and has a bounded linear right inverse
  $Q:Y\rightarrow X$.
 Choose positive constants $\delta$ and $c$ such that
  $\|Q\|\le c$, $B_{\delta}(x_0;X)\subset U$, and
  $$
   \|x-x_0\|\;<\; \delta\hspace{1em}\Longrightarrow\hspace{1em}
   \|df(x)-D\|\le \frac{1}{2c}\,.
  $$
 Suppose that $x_1\in X$  satisfies
  $$
   \|f(x_1)\|\;<\; \frac{\delta}{4c}\,, \hspace{1em}
   \|x_1-x_0\|\; \le \; \frac{\delta}{8}\,.
  $$
 Then there exists a unique $x\in X$ such that
  $$
   f(x)\;=\;0\,,\;\hspace{1em}
   x-x_1\;\in\; \Image Q\,, \hspace{1em}
   \|x-x_0\|\; \le \; \delta\,.
  $$
 Moreover, $\|x-x_1\|\;\le \; 2c\,\|f(x_1)\|$.
} % end-proposition

\bigskip

\noindent
{\bf Theorem 5.3.0.2 [implicit function theorem].}
([MD-S3: Theorem A.3.3].)
{\it
 Let $X$ and $Y$ be Banach spaces, $U\subset X$ be an open set,
  and $l$ be a positive integer.
 If $f:U\rightarrow Y$ is of class $C^l$ and
    $y$ is a regular value of $f$ $($i.e.\ $df(x)$ surjective
     with a right inverse for every $x\in f^{-1}(y)$$)$.
 then
  ${\cal M}:= f^{-1}(y)\subset X$ is a $C^l$ Banach manifold and
  $T_x{\cal M}=\Ker\, df(x)$ for every $x\in {\cal M}$.
} % end-theorem

\bigskip

With notations therein, Proposition 5.3.0.1 and Theorem 5.3.0.2
 together imply that, for $x_0\in{\cal M}$, there is a homeomorphism
 from a neighborhood of $0\in T_{x_0}{\cal M}$ to a neighborhood of
 $x_0\in{\cal M}$.
We now resume our study and notations.

\bigskip

The construction of a Kuranishi neighborhood involves
 the construction of a (continuous) family of objects and maps that
 fit into Proposition 5.3.0.1  % Proposition [Newton-Picard iteration]
 and Theorem 5.3.0.2.  % Theorem [implicit function theorem]
Relevant techniques and results in
 (closed)
   [MD-S1: Sec.~3.3, Appendix A], [MD-S3: Sec.~3.5, Chapter 10],
   [F-O: Sec.~12 - Sec.~14];
 (closed relative and closed degeneration)
  [I-P1: Sec.3, Sec.~4, Sec.~6, Sec.~7], [I-P2: Sec.~5 - Sec.~9],
  [L-R: Sec.~4]; and
 (open) [Liu(C): Sec.~6.4]
 for various related symplectic Gromov-Witten theories
 will be adapted and used to construct a family Kuranishi
 neighborhood-in-${\cal C}_{\spsccw}$ $V_{\rho}$ at each
 $\rho = [f: (\Sigma, \partial\Sigma)\rightarrow (Y[k],L)]
    \in \overline{\cal M}_{(g,h),(n,\vec{m})}
                 (W/B,L\,|\,[\beta],\vec{\gamma},\mu)$
 for an open Gromov-Witten theory of the degeneration family $W/B$.
The following diagram/flow-chart outlines the construction:

\bigskip

\noindent
$\circ$ {\it step} (1)
$$
\begin{array}{c}
 \fbox{\parbox{43ex}{\small
   choice of a {\it saturated obstruction space} $E_\rho$
   at $\rho$ }}
\end{array}
$$

\medskip

\noindent
$\circ$ {\it step} (2)
$$  % step 2.1
\begin{array}{c}
 \hspace{-2em}\Downarrow\hspace{1ex}
  \begin{array}{l}
   \mbox{\footnotesize {\tiny $\bullet$}
    linearized $\;(J,E_{\rho})$-stability condition}
  \end{array}                                  \\[1.6ex]
 \fbox{\parbox{55ex}{\small
   $\Ker(\pi_{E_{\rho}}\circ D_{\!f}\bar{\partial}_J)^{\pd}\;$
   in $\;W^{1,p}( \Sigma, \partial\Sigma; f^{\ast}T_{\ast}Y_{[k]},
                   (f|_{\partial\Sigma})^{\ast}T_{\ast}L_{[k]} )$ }}
\end{array}
$$
$$  % step 2.2
\begin{array}{c}
 \hspace{2.8em}\Downarrow\hspace{1ex}
  \begin{array}{l}
   \mbox{\footnotesize {\tiny $\bullet$}
    upper semi-continuity of $\;${\it index}$\,(D\bar{\partial}_J)\;$
    w.r.t.\ $\;B[k]$}
  \end{array}                                  \\[1.6ex]
 \fbox{\parbox{66ex}{\small
   the product space
    $\;\Def(\Sigma)\times B[k]
        \times \Ker(\pi_{E_{\rho}}
                     \circ D_{\!f}\bar{\partial}_J)^{\pd}\;$
    is large enough       }}
\end{array}
$$
$$  % step 2.3
\begin{array}{c}
 \hspace{2.2ex}\Downarrow\hspace{1ex}
  \begin{array}{l}
   \mbox{\footnotesize {\tiny $\bullet$}
     system of algebraic equations for         }\\[-.2ex]
   \hspace{1.2ex}\mbox{\footnotesize
     {\it target-deformation-driven deformations} of $\Sigma$ }
  \end{array}                                               \\[3ex]
 \fbox{\parbox{59ex}{\small
   algebraic subset $\;\widetilde{V}_{\rho}\;$ of
    $\;\Def(\Sigma)\times B[k]
         \times \Ker(\pi_{E_{\rho}}
                     \circ D_{\!f}\bar{\partial}_J)^{\pd}\,$,\\
    which projects to a constructible subset
     $\pi_{\mbox{\tiny $\bullet$}}(\widetilde{V}_{\rho})$
     in $\Def(\Sigma)\times B[k]$}}
\end{array}
$$

\medskip

\noindent
$\circ$ {\it step} (3)
\nopagebreak
\vspace{-12pt}
\nopagebreak
$$
\begin{array}{c}
 \hspace{8.6em}\Downarrow\hspace{1ex}
  \begin{array}{l}
   \mbox{\footnotesize {\tiny $\bullet$}
    piecewise-continuous section
    $\pi_{\mbox{\tiny $\bullet$}}(\widetilde{V}_{\rho})
      \rightarrow \widetilde{V}_{\rho}$
    with image closure $\Theta_{\rho}$  }      \\[.6ex]
   \mbox{\footnotesize {\tiny $\bullet$}
    gluing construction around three types -
     ordinary interior, }\\[-.2ex]
    \hspace{1.6ex}\mbox{\footnotesize
     boundary, and distinguished interior -
      of nodes on $\Sigma$ }                   \\[.6ex]
   \mbox{\footnotesize {\tiny $\bullet$}
    exponential-map construction}
  \end{array}                                               \\[7ex]
 \fbox{\parbox{82ex}{\small
   piecewise-continuous-$\pi_{\mbox{\tiny $\bullet$}}
                                 (\widetilde{V}_{\rho})$-family,
    which extends to a {\it continuous-$\widetilde{V}_{\rho}$-family},
    of pre-deformable approximate-$J$-stable $C^{\infty}$ maps
    $\;h_{\raisebox{.6ex}{\rm\tiny $\,$approx$\,,$}
                            \,\mbox{\LARGE $\cdot$}}\;$
    from $\;\Sigma_{\mbox{\LARGE $\cdot$}}\;$
    to fibers of $\;W[k]/B[k]$                       }}
\end{array}
$$

\medskip

\noindent
$\circ$ {\it step} (4)
\nopagebreak
\vspace{-12pt}
\nopagebreak
$$
\begin{array}{c}
 \hspace{6.7em}\Downarrow\hspace{1ex}
  \begin{array}{l}
   \mbox{\footnotesize {\tiny $\bullet$}
    $E_{\rho}\,$ induces a trivialized obstruction bundle
     $E_{\widetilde{V}_{\rho}}$ over $\,\widetilde{V}_{\rho}$                          } \\[-.2ex]
    \hspace{1.6ex}\mbox{\footnotesize
     with fiber
      $\;E_{\mbox{\LARGE $\cdot$}}\,\subset\,
        L^p(\Sigma_{\mbox{\LARGE $\cdot$}};
           \Lambda^{0,1}\Sigma_{\mbox{\LARGE $\cdot$}}
            \otimes_J h_{\mbox{\LARGE $\cdot$}}^{\ast}
              T_{\ast}(W[k]_{\mbox{\LARGE $\cdot$}}))$ } \\[.6ex]
   \mbox{\footnotesize {\tiny $\bullet$}
     construction of
      a $\pi_{\mbox{\tiny $\bullet$}}(\widetilde{V}_{\rho})$-family of
      uniformly bounded }  \\[-.2ex]
   \hspace{1.6ex}\mbox{\footnotesize
      right inverse $\;Q_{\mbox{\Large $\cdot$}}$  to
      $\;\pi_{E_{\mbox{\LARGE .}}}\!\circ
           D_{h_{\mbox{\LARGE .}}}\!\bar{\partial}_J$ } \\[.6ex]
   \mbox{\footnotesize {\tiny $\bullet$}
     Proposition 5.3.0.1 + Theorem 5.3.0.2$\,$: } \\[-.2ex]
   \hspace{1.6ex}\mbox{\footnotesize
     Newton's iteration method to deform approximate solutions}\\[-.2ex]
   \hspace{1.6ex}\mbox{\footnotesize
    to exact solutions to the $(J,E)$-holomorphy equation}
  \end{array}  \\[11ex]
 \fbox{\parbox{73ex}{\small
   $\,\widetilde{V}_{\rho}$-family of (exact)
    {\it $\,(J,E)$-stable maps}$\,$ $\;f_{\mbox{\LARGE $\cdot$}}\;$
    from $\;\Sigma_{\mbox{\LARGE $\cdot$}}\;$
    to fibers of $\;W[k]/B[k]$ }}
\end{array}
$$

\medskip

\noindent
$\circ$ {\it step} (5) [{\sl rigidification}]
$$ % step 5.1
\begin{array}{c}
 \hspace{2.7em}\Downarrow\hspace{1ex}
  \begin{array}{l}
   \mbox{\footnotesize {\tiny $\bullet$}
   the $J$-holomorphy of the
    $\Aut(\Sigma)\times {\footnotesizeBbb G}_m[k]$-action }
  \end{array}                                             \\[1.6ex]
 \fbox{\parbox{79ex}{\small
  a maximal subset $\;V_{\rho}\;$ in $\;\widetilde{V}_{\rho}\;$
   through $\;\rho\,$, transverse to
   the $\Aut(\Sigma)\times{\smallBbb G}_m[k]$-orbit of $\;\rho$ }}\\[1ex]
 \mbox{\footnotesize
  (This converts `$\,$maps to fibers of $\;W[k]/B[k]\,$'
   to `$\,$maps to fibers of $\,\widehat{W}/\widehat{B}\,$'.)}
  \hspace{14ex}
\end{array}
$$
$$ % step 5.2
\begin{array}{c}
 \hspace{8.6em}\Downarrow\hspace{1ex}
  \begin{array}{l}
   \mbox{\footnotesize {\tiny $\bullet$}
    Kursnishi map $\;s_{\rho}:V_{\rho}\rightarrow E_{V_{\rho}}\;$
    from the $\;\bar{\partial}_J$-operator              }\\[.6ex]
   \mbox{\footnotesize  {\tiny $\bullet$}
    stability of $\rho\,$, $\;\Gamma_{V_{\rho}}=\Aut(\rho)$ }\\[.6ex]
   \mbox{\footnotesize {\tiny $\bullet$}
    $\psi:s_{\rho}^{-1}(0)\rightarrow U_{\rho}\;$:
     orbifold quotient map to a neighborhood of $\;\rho$  }
  \end{array}                                                 \\[6ex]
 \fbox{\parbox{83ex}{\small
  $\:V_{\rho}/B\,$:
    a Kuranishi neighborhood-in-${\cal C}_{\mbox{\tiny\rm spsccw}}\;$
    of $\;\rho\;$ on
    $\;\overline{\cal M}_{(g,h),(n,\vec{m})}
      (W/B,L\,|\,[\beta],\vec{\gamma},\mu)\,/B\;$  }}
 \end{array}
$$

\bigskip

\noindent
Step (4) is the analytical core in the construction.
The algebraic system in Step (2),
 the distinguished nodes in Step (3), and
 the conversion in Step (5) from `maps to fibers of $W[k]/B[k]$'
  back to `maps to fibers of $W^+/B^+$'
 are the main substeps for which the singularity of
  the degenerate fiber $W_0=Y_1\cup_DY_2$ plays a role.

\bigskip

Throughout this subsection, we let
 $\rho\; =\;
   (\,\Sigma,\, \dot{\partial}\Sigma;\,
       \vec{p},\,\vec{p}_1,\,\ldots,\,\vec{p}_h;\, f)$
  be a stable map to the central fiber $(Y_{[k]},L_{[k]})$
   of $W[k]/B[k]$,  and
 $\rho_{(i)}:= (\,\Sigma_{(i)},\, (\dot{\partial}\Sigma)_{(i)};\,
    \vec{p}_{(i)},\,
    \vec{p}_{1,(i)},\,\ldots,\,\vec{p}_{h, (i)};\, f_{(i)})$
  be the associated submap to the irreducible component $\Delta_i$
  of $Y_{[k]}$, for $i=0,\,\ldots,\,k+1$.
 (By construction, $(\dot{\partial}\Sigma)_{(i)}$, $\vec{p}_{j,(i)}$
  can be non-empty only for $i=0$ and $k+1$.)
We denote the labelled-bordered Riemann surface with marked points
 $(\Sigma,\, \dot{\partial}\Sigma;\,
             \vec{p},\,\vec{p}_1,\,\ldots,\,\vec{p}_h)$
 also simply by $\Sigma$.
The corresponding point of $\rho$ in
 $\overline{\cal M}_{(g,h),(n,\vec{m})}
                  (W/B,L\,|\,[\beta],\vec{\gamma},\mu)$
 will also be denoted by $\rho$.
Let $\Lambda_i=f^{-1}(D_i)$ and
 $\Lambda=\disjointunion_{i=0}^k\,\Lambda_i$ be the set of
 distinguished nodes on $\Sigma$ under $f$.
Let ${\mathbf s}=(\vec{s}_0,\,\cdots,\,\vec{s}_k)$
  be the tuple of contact orders of $f$ at $\Lambda$.
Both $\Aut(\rho)$ and $\Aut(f)$ mean the same.
Denote by $\Aut(\rho)^{\domain}$ (resp.\ $\Aut(\rho)^{\target}$)
 the subgroup of $\Aut(\Sigma)$ (resp.\ ${\Bbb G}_m[k]$)
 that consists of $\alpha$ (resp.\ $\beta$) such that
  there is an $(\alpha,\beta)$ in $\Aut(f)$.
These groups are all finite.
With a re-adjustment, we assume that
 the auxiliary K\"{a}hler metric on ${\cal C}/\Def(\Sigma)$
  is $\Aut(\rho)^{\domain}$-invariant  and
 the symplectic and, hence, the metric structure on $W[k]$
  are $\Aut(\rho)^{\target}$-invariant.

\bigskip

\subsubsection{Choice of obstruction space $E_{\rho}$ of
    $\,\overline{\cal M}_{(g,h),(n,\vec{m})}
         (W/B,L\,|\,[\beta],\vec{\gamma},\mu)$ at $\rho$.}

\begin{flushleft}
{\bf The index of the linearized operator $D_{\!f}\bar{\partial}_J$
     of $\bar{\partial}_J$ at $f$.}
\end{flushleft}
The fiber of the $\check{W}^{1,p}$-tangent-obstruction fibration complex
 $$
  T^1_{\check{\cal W}^{1,p}_{\tinybullet}
          ((\widehat{W},\widehat{L})/\widehat{B}\,|\,\tinybullet)
            /\widetilde{\cal M}_{\tinybullet}}
    \left.\rule{0ex}{2ex}\right|\,
     _{\overline{\cal M}_{\tinybullet}(W/B,L\,|\,\tinybullet)}\;
  \stackrel{D\bar{\partial}_J}{\longrightarrow}\;
   T^2_{\check{\cal W}^{1,p}_{\tinybullet}
          ((\widehat{W},\widehat{L})/\widehat{B}\,|\,\tinybullet)
            /\widetilde{\cal M}_{\tinybullet}}
    \left.\rule{0ex}{2ex}\right|\,
     _{\overline{\cal M}_{\tinybullet}(W/B,L\,|\,\tinybullet)}
 $$
 at $\rho$ has a $C^l$-, $C^{\infty}$-, and $W^{1,p}$-parallel
 as follows:
  (by convention, $\partial\Sigma_{(i)}=\emptyset ={L_{[k]}}_{(i)}$
   for $i=1,\,\ldots,\,k$)
 \begin{eqnarray*}
  \lefteqn{
   C^l(\Sigma, \partial \Sigma;
     f^{\ast}T_{\ast}Y_{[k]},
     (f|_{\partial\Sigma})^{\ast}T_{\ast}L_{[k]}) } \\[.6ex]
   && :=\; \left\{\, (\xi_{(i)})_{i=0}^{k+1}
           \in \oplus_{i=0}^{k+1}\,
           C^l(\,\Sigma_{(i)}, \partial \Sigma_{(i)};
                f_{(i)}^{\ast}T_{\ast}\Delta_{(i)},
                (f_{(i)}|_{\partial\Sigma_{(i)}})^{\ast}
                             T_{\ast}{L_{[k]}}_{(i)}\,) \right.  \\
   && \left.\rule{0em}{1em} \hspace{12em}
        :\; \xi_{(j)}|_{\Lambda_j}\,=\, \xi_{(j+1)}|_{\Lambda_j}\,
              \in (f|_{\Lambda_j})^{\ast}T_{\ast}D_j\,,
                                                 j=0,\,\ldots,\,k\,
       \right\}\,,
 \end{eqnarray*}
 $$
  C^l(\Sigma; \Lambda^{0,1}\Sigma\otimes_J f^{\ast}T_{\ast}Y_{[k]})\;
  :=\; \oplus_{i=0}^{k+1}\,C^l(\Sigma_{(i)}\,,\,
        \Lambda^{0,1}\Sigma_{(i)}
                 \otimes_J f_{(i)}^{\ast}T_{\ast}\Delta_{(i)})\,,
 $$
 and
 $$
  D_{\!f}\bar{\partial}_J\;:\;
   C^{\infty}(\Sigma, \partial \Sigma;
    f^{\ast}T_{\ast}Y_{[k]},
    (f|_{\partial\Sigma})^{\ast}T_{\ast}L_{[k]})\;
   \longrightarrow\;
   C^{\infty}
    (\Sigma; \Lambda^{0,1}\Sigma\otimes_J f^{\ast}T_{\ast}Y_{[k]})\,,
 $$
 $$
  D_{\!f}\bar{\partial}_J\;:\;
   W^{1,p}(\Sigma, \partial \Sigma;
    f^{\ast}T_{\ast}Y_{[k]},
    (f|_{\partial\Sigma})^{\ast}T_{\ast}L_{[k]})\;
   \longrightarrow\;
   L^p(\Sigma; \Lambda^{0,1}\Sigma\otimes_J f^{\ast}T_{\ast}Y_{[k]})\,.
 $$
 For $\nabla$ the Levi-Civita connection of the metric
  on $Y_{[k]}$ induced by $(\omega,J)$,
 the linearization $D\bar{\partial}_J$ of $\bar{\partial}_J$ is given
  by
  $$
   (D_{\!f}\bar{\partial}_J)(\xi)\;
    =\; \frac{1}{2}\left(\, \rule{0ex}{1.2em}
        \nabla\xi\,\circ\, df\:+\; J\circ \nabla\xi\circ df\circ j\:
        +\: \nabla_{\xi}J\circ df\circ j
       \,\right)\,,
  $$
  on the irreducible components of $\Sigma$ for which $f$ is not constant,
  cf.\ [Liu(C): Proposition 6.12];
  see also [MD-S1: Eq.~(3.2) and Remark 3.3.1].
For an irreducible component of $\Sigma$ on which $f$ is a constant map.
the related bundles,
 $(f^{\ast}T_{\ast}Y_{[k]},
    (f|_{\partial\Sigma})^{\ast}T_{\ast}L_{[k]})$ and
 $\Lambda^{0,1}\Sigma\otimes_J f^{\ast}T_{\ast}Y_{[k]}$,
 on that component are of the respective forms,
  ${\cal O}_{\Sigma}\otimes_{\scriptsizeBbb C} {\Bbb C}^m$ and
  $\Lambda^{0,1}\Sigma\otimes_{\scriptsizeBbb C}{\Bbb C}^m$,
 and have the canonical holomorphic structure from
  the complex structure on $\Sigma$.
$D_f\bar{\partial}_J$ for such component is the restriction to
 that component of the operator
  $\bar{\partial}:
   C^{\infty}(\Sigma,
    {\cal O}_{\Sigma}\otimes_{\scriptsizeBbb C}{\Bbb C}^m)
   \rightarrow
   C^{\infty}(\Sigma,
    \Lambda^{0,1}\Sigma\otimes_{\scriptsizeBbb C}{\Bbb C}^m)$
 associated to the canonical holomorphic structure.
The following lemma should be compared to [I-P2: Lemma 7.2]
 and [L-R: Theorem 5.1].

\bigskip

\noindent
{\bf Lemma 5.3.1.1
     [index of $D_{\!f}\bar{\partial}_J$ for rigid target].}
{\it
 Let
  $f: (\Sigma, \partial\Sigma) \rightarrow (Y_{[k]},L_{[k]})$
  be a stable map to the specified expanded target space as above.
 Then the restriction
  $$
   D_{\!f}\bar{\partial}_J\,:\,
    C^{\infty}(\Sigma, \partial \Sigma;
     f^{\ast}T_{\ast}Y_{[k]},
     (f|_{\partial\Sigma})^{\ast}T_{\ast}L_{[k]})\;
    \longrightarrow\;
    C^{\infty}(\Sigma;
    \Lambda^{0,1}\Sigma \otimes_J f^{\ast}T_{\ast}Y_{[k]})
  $$
  is a Fredholm operator of index
  $$
   \ind(D_{\!f}\bar{\partial}_J)\; =\;
    \mu(f)\, +\, \dimm Y \cdot (1-\tilde{g})\,
     -\, 2\,\sum_{i=0}^k l(\vec{s}_i)\,
     +\, 4\,\sum_{i=0}^k \degree\vec{s}_i\,,
  $$
  where $\tilde{g}$ is the arithmetic genus of
  $\Sigma_{\scriptsizeBbb C}$.
} % end-lemma

\bigskip

\noindent
{\it Proof.}
 Let
  $f=\cup_{i=0}^{k+1}\,f_{(i)}\,:\,
   \Sigma=\cup_{i=0}^{k+1}\,\Sigma_{(i)}\,
    \rightarrow\, Y_{[k]}=\cup_{i=0}^{k+1}\Delta_i$
  be the decomposition of $f$ into submaps.
 Then, it follows from the Riemann-Roch Theorem
  (e.g.\ [F-O: Lemma 12.2], [Liu(C): Lemma 6.13],
         and [MD-S3: Appendix C])
 that
  each of
   $$
    D_{f_{(i)}}\bar{\partial}_J\,:\,
     C^{\infty}(\Sigma_{(i)}, \partial \Sigma_{(i)};
      f_{(i)}^{\ast}T{\Delta_i},
      ({f_{(i)}}|_{\partial\Sigma})^{\ast}TL)\;
     \longrightarrow\;
     C^{\infty}(\Sigma_{(i)};
      \Lambda^{0,1}\Sigma_{(i)} \otimes_J f_{(i)}^{\ast}T\Delta_i)\,,
   $$
    for $i=0,\, k+1$, and
   $$
    D_{f_{(i)}}\bar{\partial}_J\,:\,
     C^{\infty}(\Sigma_{(i)}; f_{(i)}^{\ast}T\Delta_i)\;
     \longrightarrow\;
     C^{\infty}(\Sigma_{(i)};
      \Lambda^{0,1}\Sigma_{(i)} \otimes_J f_{(i)}^{\ast}T\Delta_i)\,,
   $$
    for $i=1,\,\ldots,\,k$,
   is a Fredholm operator of index
   $$
    \ind(D_{f_{(i)}}\bar{\partial}_J)\;
     =\; \mu(f_{(i)})\, +\, \dimm Y\,\tilde{\chi_i}/2\,,
    \hspace{1em}\mbox{for $i=0,\,k+1$}\,,
   $$
  and
   $$
    \ind(D_{f_{(i)}}\bar{\partial}_J)\;
     =\; -2 K_{\Delta_i}\cdot\beta_i + \dimm Y \chi_i/2\,,
    \hspace{1em}\mbox{for $i=1,\,\ldots,\,k$}\,.
   $$
 This implies, in particular, that $D_{\!f}\bar{\partial}_J$
  is Fredholm.

 The matching condition along $TD_i$ at each distinguished node
  of $\Sigma$ imply that
  \begin{eqnarray*}
   \lefteqn{
    C^{\infty}(\Sigma, \partial \Sigma;
     f^{\ast}T{Y_{[k]}}, (f|_{\partial\Sigma})^{\ast}TL_{[k]})  }\\[.6ex]
   && \hookrightarrow \;
      \oplus_{i=0, k+1}\,
       C^{\infty}(\Sigma_{(i)}, \partial \Sigma_{(i)};
           f_{(i)}^{\ast}T{\Delta_i},
           ({f_{(i)}}|_{\partial\Sigma})^{\ast}TL)\;
      \bigoplus
      \oplus_{i=1}^k\,
      C^{\infty}(\Sigma_{(i)}; f_{(i)}^{\ast}T\Delta_i)
  \end{eqnarray*}
  has codimension
  $$
   \sum_{i=0}^{k}\,l(\vec{s}_i)\,(\dimm Y + 2)\,.
  $$
 Denote the quotient vector space of this inclusion by $V$, then
  one has the following short exact sequence of $2$-term complexes:
  $$
   \begin{array}{cccccccccl}
    0  & \longrightarrow  & C^0
       & \longrightarrow  & \oplus_{i=0}^{k+1}\,C^0_{(i)}
       & \longrightarrow  & V
       & \longrightarrow  & 0\\
    &  & \mbox{\scriptsize $D_{\!f}\bar{\partial}_J$}
         \downarrow \hspace{2em}
      && \mbox{\scriptsize
         $\oplus_i D_{f_{(i)}}\bar{\partial}_J$}\downarrow\hspace{3.4em}
      && \downarrow \\[1ex]
    0  & \longrightarrow  & C^1
       & \longrightarrow  & \oplus_{i=0}^{k+1}\,C^1_{(i)}
       & \longrightarrow  & 0
       & \longrightarrow  & 0 &,
   \end{array}
  $$
  where
  $$
   \begin{array}{lcll}
    C^0         & =
     & C^{\infty}(\Sigma, \partial \Sigma; f^{\ast}T{Y_{[k]}},
         (f|_{\partial\Sigma})^{\ast}TL_{[k]})\,,               \\[.6ex]
    C^1         & =
     & C^{\infty}(\Sigma;
          \Lambda^{0,1}\Sigma \otimes_J f^{\ast}TY_{[k]})\,,    \\[1ex]
    C^0_{(i)}   & =
     & \left\{
         \begin{array}{ll}
          C^{\infty}(\Sigma_{(i)}, \partial \Sigma_{(i)};
             f_{(i)}^{\ast}T{\Delta_i},
                  ({f_{(i)}}|_{\partial\Sigma})^{\ast}TL)\,,
            & \mbox{for $i=0, k+1$}\,,                  \\[.6ex]
          C^{\infty}(\Sigma_{(i)}; f_{(i)}^{\ast}T\Delta_i)\,,
            & \mbox{for $i=1,\,\ldots,\,k$}\,,
        \end{array}
       \right.                                                  \\[3ex]
    C^1_{(i)}   & =
     & C^{\infty}(\Sigma_{(i)};
        \Lambda^{0,1}\Sigma_{(i)} \otimes_J f_{(i)}^{\ast}T\Delta_i)\,.
   \end{array}
  $$
 The Snake Lemma, together with
  the additivity property
   of (relative) Maslov index under joining of submaps
    (Definition 3.1.2) and
   of the Euler characteristic of Riemann surfaces under gluing along
    boundaries from removing small discs around distinguished nodes,
  implies then
  $$
    \ind(D_{\!f}\bar{\partial}_J)\;
     =\; \mu(f)\, +\, \dimm Y \cdot (1-\tilde{g})\,
                  -\, 2\,\sum_{i=0}^k l(\vec{s_i})\,
                  +\, 4\,\sum_{i=0}^k \degree\vec{s}_i\,.
  $$

\noindent\hspace{15cm}$\Box$

\bigskip

\noindent
{\it Remark 5.3.1.2 $[$class independence$]$.} ([MD-S1: Remark.~3.2.3].)
 Lemma 5.3.1.1 holds also for
  $$
   D_{\!f}\bar{\partial}_J\,:\,
    W^{l,p}( \Sigma, \partial\Sigma; f^{\ast}T_{\ast}Y_{[k]},
                    (f|_{\partial\Sigma})^{\ast}T_{\ast}L_{[k]} )\;
    \longrightarrow\;
     W^{l-1,p}(\Sigma,\Lambda^{0,1}\Sigma
                      \otimes_J f^{\ast}T_{\ast}Y_{[k]})\,,
  $$
  $$
   D_{\!f}\bar{\partial}_J\,:\,
    \check{W}^{l,p}( \Sigma, \partial\Sigma; f^{\ast}T_{\ast}Y_{[k]},
                    (f|_{\partial\Sigma})^{\ast}T_{\ast}L_{[k]} )\;
    \longrightarrow\;
     \check{W}^{l-1,p}(\Sigma,\Lambda^{0,1}\Sigma
                      \otimes_J f^{\ast}T_{\ast}Y_{[k]})\,,
  $$
 and
  $$
   D_{\!f}\bar{\partial}_J\,:\,
    C^l( \Sigma, \partial\Sigma; f^{\ast}T_{\ast}Y_{[k]},
                   (f|_{\partial\Sigma})^{\ast}T_{\ast}L_{[k]} )\;
    \longrightarrow\;
     C^{l-1}(\Sigma,\Lambda^{0,1}\Sigma
                    \otimes_J f^{\ast}T_{\ast}Y_{[k]})\,.
  $$
 We have taken $J$ to be of class $C^{\infty}$ on each irreducible
  component of $Y_{[k]}$. Thus, elliptic regularity implies that
 $\Ker(D_{\!f}\bar{\partial}_J)$ always lies in
   $C^{\infty}( \Sigma, \partial\Sigma; f^{\ast}T_{\ast}Y_{[k]},
                     (f|_{\partial\Sigma})^{\ast}T_{\ast}L_{[k]} )$,
 independent of the choice of the space on which
  $D_{\!f}\bar{\partial}_J$ is defined.

\bigskip

\begin{flushleft}
{\bf Existence of a saturated obstruction space $E_{\rho}$ at $\rho$.}
\end{flushleft}
For a small enough neighborhood
 $U_{\Lambda^+}=
   \left( \amalg_{q\in\Lambda} U_q \right)
   \amalg \left( \amalg_{p_i}U_{p_i} \right)
   \amalg \left( \amalg_{q_{ij}}U_{q_{ij}} \right)$
 of the set
 $\Lambda^+:=\Lambda\cup{\mathbf p}\cup\cup_{j=1}^{\,h}{\mathbf q}_j$
 of the distinguished nodes and the marked points
 on $\Sigma$,
recall from Sec.~5.2 (with $0$ there replaced by $q$ here)
 the associated subspace
 ${\cal H}_q$ in
 $L^p(\Sigma; \Lambda^{0,1}\Sigma\otimes_J f^{\ast}T_{\ast}Y_{[k]})
     |_{U_q}$, $q\in \Lambda$,
 such that $D_f\divisor_q$ is defined on
  $(D_{\!f}\bar{\partial}_J)^{-1}({\cal H}_q)$ with values in
  $T_{(s(q)-1)\cdot(q)}\Div^{s(q)-1}(U_{q,1})
   \oplus T_{(s(q)-1)\cdot(q)}\Div^{s(q)-1}(U_{q,2})$,
 where $s(q)$ is the contact order of $f$ along
  the singular locus of $Y_{[k]}$ at $q$.

\bigskip

\noindent
{\bf Definition 5.3.1.3 [admissible subspace].} {\rm
 A subspace $V$ in
  $W^{1,p}(\Sigma, \partial\Sigma; f^{\ast}T_{\ast}Y_{[k]},
               (f|_{\partial\Sigma})^{\ast}T_{\ast}L_{[k]})$
  is called {\it admissible}
 if there exists such an $U_{\Lambda^+}$ so that
  $V|_{U_q} \subset (D_{\!f}\bar{\partial}_J)^{-1}({\cal H}_q)$
  for all $q\in \Lambda$.
} % end-definition

\bigskip

As all the maps $D_f\divisor_q$, $\jet_q^{s(q)}$,
 $D_f\ev_q$, $D_f\ev_{p_i}$, and $D_f\ev_{q_{ij}}$ depend only on
 a jet at the specified point in $\Lambda^+$,
they extend canonically to maps on an admissible subspace of
 $W^{1,p}(\Sigma, \partial\Sigma; f^{\ast}T_{\ast}Y_{[k]},
                   (f|_{\partial\Sigma})^{\ast}T_{\ast}L_{[k]})$
 by pre-composition with the restriction-to-$U_{\Lambda^+}$ map.

\bigskip

\noindent
{\bf Definition 5.3.1.4 [saturated/pre-deformable subspace].} {\rm
 A subspace $V$ in
  $W^{1,p}(\Sigma, \partial\Sigma;\\ f^{\ast}T_{\ast}Y_{[k]},
                   (f|_{\partial\Sigma})^{\ast}T_{\ast}L_{[k]})$
  is said to be {\it saturated}$\,$
 if
  \begin{itemize}
   \item[$(1)$]
    $V$ is admissible;

   \item[$(2)$]
   the map
   \begin{eqnarray*}
    \lefteqn{
      \left( \oplus_{q\in\Lambda} D_f\divisor_q \right)
       \mbox{$\bigoplus$} \left( \oplus_{p_i} D_f\ev_{p_i} \right)
      \mbox{$\bigoplus$} \left( \oplus_{q_{ij}} D_f\ev_{q_{ij}} \right)\;
     :\;  V\;\longrightarrow } \\[.6ex]
    &&
     \left( \mbox{$\bigoplus$}_{q\in\Lambda}
        \left( T_{(s(q)-1)\cdot(q)}\Div^{s(q)-1}(U_{q,1})
               \oplus T_{(s(q)-1)\cdot(q)}\Div^{s(q)-1}(U_{q,2})
        \right)
      \right)  \\[.6ex]
    && \hspace{18em}
    \mbox{$\bigoplus$} \left( \oplus_{p_i} T_{f(p_i)}Y_{[k]} \right)
    \mbox{$\bigoplus$} \left( \oplus_{q_{ij}} T_{f(q_{ij})}L \right)\;
   \end{eqnarray*}
   is surjective;

   \item[$(3)$]
   let $V^{\pd}$ be the subspace
    $(\oplus_{q\in\Lambda} D_f\divisor_q)^{-1}({\mathbf 0})$ in $V$,
   then the linear map
    $$
     \oplus_{q\in\Lambda}\,( D_f\ev_q \oplus \jet_q^{s(q)})\; :\;
      V^{\pd}\; \longrightarrow\;
      \oplus_{q\in\Lambda}\,(T_{f(q)}D\oplus {\Bbb C}^2)
    $$
    is surjective, where we have identified $D_i$, $i=0,\,\ldots,\,k$,
    canonically with $D$.
  \end{itemize}
 In the above statement, $V^{\pd}$ is called the
  {\it pre-deformable subspace} of $V$.

 A subspace $E$ of
   $L^p(\Sigma;
      \Lambda^{0,1}\Sigma\otimes_J f^{\ast}T_{\ast}Y_{[k]})$
   is said to be {\it saturated}
 if $(D_{\!f}\bar{\partial}_J)^{-1}(E)\subset
        W^{1,p}(\Sigma, \partial \Sigma;\\   f^{\ast}T_{\ast}Y_{[k]},
                      (f|_{\partial\Sigma})^{\ast}T_{\ast}L_{[k]})$
  is saturated.
} % end-definition

\bigskip

\noindent
{\bf Definition/Lemma 5.3.1.5 [saturated obstruction space].} {\rm
 Denote by $\Image(D_{\!f}\bar{\partial}_J)$ the image of
 $D_{\!f}\bar{\partial}_J$,
  $(D_{\!f}\bar{\partial}_J)
   ( W^{1,\,p}( \Sigma, \partial\Sigma; f^{\ast}T_{\ast}Y_{[k]},
                 (f|_{\partial\Sigma})^{\ast}T_{\ast}L_{[k]}) )$, in
  $L^p(\Sigma;
    \Lambda^{0,1}\Sigma\otimes_J f^{\ast}T_{\ast}Y_{[k]})$.
 Then {\it
 there exists a subspace $E_\rho$ of
  $L^p(\Sigma; \Lambda^{0,1}\Sigma\otimes_J f^{\ast}T_{\ast}Y_{[k]})$
 such that
 \begin{itemize}
  \item[{\rm (1)}]
   $\Image(D_{\!f}\bar{\partial}_J) + E_{\rho}
    = L^p(\Sigma;
          \Lambda^{0,1}\Sigma\otimes_J f^{\ast}T_{\ast}Y_{[k]})$,

  \item[{\rm (2)}]
   $E_{\rho}$ is finite-dimensional, complex linear,
    and $\Aut(\rho)$-invariant,

  \item[{\rm (3)}]
   $E_{\rho}$ consists of smooth sections supported in
    a compact subset of $\Sigma$ disjoint from the set of
    all $\,($three types of$\,)$ nodes on $\Sigma$,

  \item[{\rm (4)}]
   $(D_{\!f}\bar{\partial}_J)^{-1}(E_{\rho})$ is a saturated
    subspace of
    $W^{1,p}( \Sigma, \partial\Sigma; f^{\ast}T_{\ast}Y_{[k]},
                    (f|_{\partial\Sigma})^{\ast}T_{\ast}L_{[k]}) )$.
 \end{itemize}} % end-italic
 $E_{\rho}$ is called a {\it saturated obstruction space} of
  $\,\overline{\cal M}_{(g,h),(n,\vec{m})}
            (W/B,L\,|\,[\beta],\vec{\gamma},\mu)$ at $\rho\,$.
} % end-lemma

\bigskip

\noindent
{\it Proof.}
 Since $J$ is strongly locally regular and
  the vector spaces above is constructed from a gluing of
  the ordinary case for smooth target spaces, the existence of
  $E_{\rho}^{\prime}$ with Properties (1), (2), and (3) follows
  the same argument as in [F-O: 12.7] and [Liu(C): Lemma 6.18].
 It remains now to enlarge $E_{\rho}^{\prime}$ to incorporate
  Property (4).

 As $J$ is locally strongly regular, there exist
  finite-dimensional subspaces
  $V_{U_q}$ (resp.\ $V_{U_{p_i}}$, $V_{U_{q_{ij}}}$)
  in $(D_{f|_{U_q}}\!\circ\bar{\partial}_J)^{-1}({\cal H}_q)$
  (resp.\
    $C^{\infty}(U_p,(f|_{U_p})^{\ast}T_{\ast}Y_{[k]})$,
    $C^{\infty}(U_{q_{ij}} \partial_0U_{ij};
      (f|_{U_{q_{ij}}})^{\ast}T_{\ast}Y_{[k]}, \\
        (f|_{\partial_0U_{q_{ij}}})^{\ast}T_{\ast}L )$)
  such that, for all $q,\, p_i,\, q_{ij}\in\Lambda^+$,
   (a) the restriction of $D_f\divisor_q$
       (resp.\ $D_f\ev_{p_i}$, $D_f\ev_{q_{ij}}$)
       thereon is surjective;
   (b) the restriction of $D_f\ev_q\oplus\jet_q^{s(q)}$
       on the local pre-deformable subspace ${V_{U_q}}^{\pd}$
       is surjective, and
   (c) $D_{f|_{U_q}}(V_{U_q})$
        (resp.\ $D_{f|_{U_{p_i}}}(V_{U_{p_i}})$,
             $D_{f|_{U_{q_{ij}}}}(V_{U_{q_{ij}}})$)
       is supported in the complement of a small neighborhood
       of $q$ (resp.\ $p_i$, $q_{ij}$).
 One can extend $V_{U_q}$, $V_{U_{p_i}}$, $V_{U_{q_{ij}}}$
  to subspaces $V_q$, $V_{p_i}$, $V_{q_{ij}}$ in
  $C^{\infty}(\Sigma, \partial\Sigma; f^{\ast}T_{\ast}Y_{[k]},\\
                   (f|_{\partial\Sigma})^{\ast}T_{\ast}L_{[k]})$
  so that the summation
   $V := (\sum_{q\in\Lambda}V_q)
           +(\sum_{p_i}V_{p_i})+(\sum_{q_{ij}}V_{q_{ij}})$
   in $C^{\infty}(\Sigma, \partial\Sigma;\\ f^{\ast}T_{\ast}Y_{[k]},
                        (f|_{\partial\Sigma})^{\ast}T_{\ast}L_{[k]})$
  is a direct sum.
 Then the image $(D_{\!f}\bar{\partial}_J)(V)$
  is a finite-dimensional saturated subspace of
   $L^p(\Sigma;
        \Lambda^{0,1}\Sigma\otimes_J f^{\ast}T_{\ast}Y_{[k]})$
  that satisfies Condition (3).
 Let $E_{\rho}$ be the span of
  $E_{\rho}^{\prime}+D_{\!f}\bar{\partial}_J(V)$ and its image
   under the complex rotation and the $\Aut(\rho)$-action.
 Then $E_{\rho}$ satisfies Properties (1), (2), (3), (4).

\noindent\hspace{15cm}$\Box$

\bigskip

Let $E_{\rho}$ be a such obstruction space at $\rho$.
Property (4) of $E_{\rho}$ implies that
 $(D_{\!f}\bar{\partial}_J)^{-1}(E_{\rho})^{\pd}$
 has (real) codimension
 $\,4\sum_{i=0}^k \left(\degree\vec{s}_i\,-\,l(\vec{s}_i)\right)$
 in $(D_{\!f}\bar{\partial}_J)^{-1}(E_{\rho})$.
It follows thus from Lemma 5.3.1.1
   % Lemma [index of $D_{\!f}\bar{\partial}_J$ for rigid target]
 that:

\bigskip

\noindent
{\bf Corollary 5.3.1.6 [pre-deformable subspace of
        $(D_{\!f}\bar{\partial}_J)^{-1}(E_{\rho})$].}
{\it
 $$
  \begin{array}{ccll}
   \dimm (D_{\!f}\bar{\partial}_J)^{-1}(E_{\rho})^{\pd}
    & = & \mu(f)\, +\, \dimm Y \cdot (1-\tilde{g})\,
          +\, 2\,\sum_{i=0}^k l(\vec{s}_i)\,
          +\, \dimm E_{\rho}                         \\[1ex]
   & =  & \mu(f)\, +\, \dimm Y \cdot (1-\tilde{g})\,
          +\, 2\,|\Lambda|\,  +\, \dimm E_{\rho}     & .
  \end{array}
 $$
} % end-corollary

\bigskip

\noindent
{\bf Definition 5.3.1.7 [pre-deformable index].} {\it
 We define the {\it pre-deformable index} of $D_{\!f}\bar{\partial}_J$
  to be
  $$
    \ind^{\pd}(D_{\!f}\bar{\partial}_J)\;
     :=\; \mu(f)\, + \, \dimm Y \cdot (1-\tilde{g}) \,
           +\,  2\, |\Lambda|\,.
  $$
} % end-definition

\bigskip

\noindent
{\it Remark 5.3.1.8 $[$ fixed vs.\ non-fixed $($domain, target$)$$\,]$.}
 While there is no local obstruction to extending stable map $f$
  from a fixed nodal curve $\Sigma$ to a fixed transverse nodal
  target $Y_{[k]}$,
 there remain obstructions when extending such maps
  to a partial smoothing of $Y_{[k]}$, enforcing a deformation of
  the domain as well.
 In algebro-geometric/holomorphic setting, such obstructions are
  encoded in the cohomology
  $H^0( \Sigma,
        f^{\ast}\Extsheaf^1(\Omega_{Y_{[k]}},{\cal O}_{Y_{[k]}}) )$.
  % ?????? \simeq {\Bbb C}^{|\Lambda|}$
  % (???, or $2\,\sum_{i=0}^k\degree\vec{s}_i\,$?).
  % %
  % \marginpar{\raggedright\tiny\vspace{-2em} $\bullet$ To be completed.}
  %
 The existence of such obstructions is reflected in the dropping of
  $\ind^{\pd} (D_{\mbox{\LARGE $\cdot$}}\bar{\partial}_J)$
  %%%%%%%%%
  % by, e.g.,\ $2|\Lambda_i|$ (???, or $2\,\degree\vec{s}_i$ ?)
  %%%%%%%%%
  when $f$ is deformed to a nearby stable map to $Y_{[k-1]}$
   that smoothes $D_i$,
 cf.\ Definition 5.3.1.7, % Definition [pre-deformable index]

\bigskip

\noindent
{\bf Definition 5.3.1.9 [$(J,E_{\rho})$-stable map].} {\rm
 Given $E_{\rho}$ in Definition/Lemma 5.3.1.5,
   % Definition/Lemma [saturated obstruction space]
 a map $h:(\Sigma, \partial\Sigma)\rightarrow (Y_{[k]}, L)$
  is called {\it $(J,E_{\rho})$-stable}
 if it satisfies the {\it perturbed $J$-holomorphy equations}
  $\bar{\partial}_J\,h \in E_{\rho}$,
  is {\it pre-deformable} at the distinguished nodes, and
  has a finite $\Aut(h)$.
} % end-definition

\bigskip

For later use, we introduce the quotient map
 $$
  \pi_{E_{\rho}}\; :\;
   L^p(\Sigma;
     \Lambda^{0,1}\Sigma\otimes_J f^{\ast}T_{\ast}Y_{[k]})\;
    \longrightarrow\;
    L^p(\Sigma;\\
    \Lambda^{0,1}\Sigma\otimes_J f^{\ast}T_{\ast}Y_{[k]})/E_{\rho}
 $$
  and
 denote $(D_{\!f}\bar{\partial}_J)^{-1}(E_{\rho})^{\pd}$ also as
 $\Ker(\pi_{E_{\rho}}\!\circ D_{\!f}\bar{\partial}_J)^{\pd}$.
With respect to the holomorphic coordinates around $\Lambda$ in
 $\Sigma$ and normal coordinates to $\cup_{i=0}^k D_i$ around
 $f(\Lambda)$ in $Y_{[k]}$ that defines $\jet_q^{s(q)}$, $q\in\Lambda$,
one thus has the linear map
 $$
  \begin{array}{cccccl}
   \jet_{\Lambda}^{\mathbf s} & :
    & \Ker(\pi_{E_{\rho}}\!\circ D_{\!f}\bar{\partial}_J)^{\pd}
    & \longrightarrow      & {\Bbb C}^{\,2\,|\Lambda|}      \\[1ex]
   && \xi & \longmapsto
    & \left(\,\jet_q^{s(q)}\left((\xi|_{U_q})^n\right)\,
      \right)_{q\in\Lambda}\,.
 \end{array}
$$
For $q\in\Lambda$, suppose that with respect to
 the fixed local coordinates $f|_{U_q}$ is given by
 $$
  f(z_{q,i})
    =\; \left(\, f(q)+O\left(|z_{q,i}| \right)\,,\,
           a_{q,i} z_{q,i}^{s(q)} + O\left(|z_{q,i}|^{s(q)+1} \right)
         \,\right)\,,\hspace{1em} i=1,\,2\,.
 $$
Define the {\it shift-product map}
 $\shiftproduct_q:{\Bbb C}^2 \rightarrow {\Bbb C}$,
 $(\,\mbox{\LARGE $\cdot$}_1\,,\,\mbox{\LARGE $\cdot$}_2) \mapsto
     (\,a_{q,1}+\,\mbox{\LARGE $\cdot$}_1)
     (\,a_{q,2}+\,\mbox{\LARGE $\cdot$}_2)$.
Then, the image of a small enough neighborhood of $0$ in
 $\Ker(\pi_{E_{\rho}}\!\circ D_{\!f}\bar{\partial}_J)^{\pd}$
 under the composition $\shiftproduct_q\circ\jet_q^s$ lies in a
 simply-connected neighborhood of $a_{q,1}a_{q,2}$ in ${\Bbb C}-\{0\}$.
For $f(q)\in D_i$, $\shiftproduct_q\circ\jet_q^{s(q)}$ is
 a (nonlinear) map from
 $\Ker(\pi_{E_{\rho}}\!\circ D_{\!f}\bar{\partial}_J)^{\pd}$ to
 $\left(
   (T_q^{\ast}\Sigma_{(i)})^{\otimes s(q)}\otimes T_{f(q)}\Delta_i
  \right)
  \bigotimes
  \left(
   (T_q^{\ast}\Sigma_{(i+1)})^{\otimes s(q)}\otimes T_{f(q)}\Delta_{i+1}
  \right)$.
Define the nonlinear map
 $$
  \begin{array}{cccccl}
   \shiftproduct_{\Lambda}\circ\jet_{\Lambda}^{\mathbf s} & :
    & \Ker(\pi_{E_{\rho}}\!\circ D_{\!f}\bar{\partial}_J)^{\pd}
    & \longrightarrow   & {\Bbb C}^{\,|\Lambda|} \\[1ex]
   && \xi  & \longmapsto
    & \left(\,
       \shiftproduct_q\circ\jet_q^{s(q)}\left( (\xi|_{U_q})^n\right)
       \,\right)_{q\in\Lambda}\,.
 \end{array}
$$
Property (4) of $E_{\rho}$ in Definition/Lemma 5.3.1.5
      % Definition/Lemma [saturated obstruction space]
 implies that the map
 $\shiftproduct_{\Lambda}\circ\jet_{\Lambda}^{\mathbf s}$
 is a bundle map over a small enough neighborhood of
 $\shiftproduct_{\Lambda}\circ\jet_{\Lambda}^{\mathbf s}(0)$
 in $({\Bbb C}-\{0\})^{|\Lambda|}$
 with fiber of (real) dimension
  $\mu(f)\, + \, \dimm Y \cdot (1-\tilde{g})\,+\,\dimm E_{\rho}$.

\bigskip

\subsubsection{The algebraic subset $\widetilde{V}_{\rho}$ in
    $\Def(\Sigma)\times B[k]  \times
       \Ker(\pi_{E_{\rho}}\circ D_{\!f}\bar{\partial}_J)^{\pd}\,$.}

A family Kuranishi neighborhood $V_{\rho}$ of
 $\rho=[f]\in \overline{\cal M}_{(g,h),(n,\vec{m})}
                  (W/B,L\,|\,[\beta],\vec{\gamma},\mu)$ over $B$
 is to be obtained from an enlarged deformation theory of
 the underlying moduli problem.
Corollary 5.3.1.6
   % Corollary [pre-deformable subspace of
   %            $\Ker(\pi_{E_{\rho}}\!\circ D_{\!f}\bar{\partial}_J)$]
 implies that
  $\ind^{\pd} D_h\bar{\partial}_J$
  is piecewise-constant and upper semi-continuous with respect to
  the stratification of $B[k]$ when $h$ runs over $J$-stable maps
  of the given combinatorial type from deformed $\Sigma$ to fibers
  of $W[k]/B[k]$.
This hints that the product space
 $\Def(\Sigma) \times B[k]  \times
    \Ker(\pi_{E_{\rho}}\circ D_{\!f}\bar{\partial}_J)^{\pd}$
 is large enough to accommodate all the new maps to appear
 in a candidate Kuranishi neighborhood $V_{\rho}$ of $\rho$.
We describe in this subsubsection an algebraic subset
 $\widetilde{V}_{\rho}$ in
 $\Def(\Sigma) \times B[k]  \times
   \Ker(\pi_{E_{\rho}}\circ D_{\!f}\bar{\partial}_J)^{\pd}$,
  characterized by the deformation theory of maps at
  the distinguished nodes,
 that will finally give $V_{\rho}$.

As observed in [I-P2]
 (see also [Li1] and [Gr-V] in algebro-geometric category),
when $Y_{[k]}$ is partially smoothed to some
 $Y_{[k^{\prime}]}$ with $D_i$ being smoothed, all the nodes in
 $\Sigma_i$ in $\Sigma$ have to be simultaneously smoothed in order
 for there to exist a $J$-holomorphic map $f^{\prime}$ from the new
 $\Sigma^{\prime}$ to $Y_{[k^{\prime}]}$ that is close to $f$.
Thus, $V_{\rho}$ should only come from a subset of a locus
 $\widetilde{V}_{\rho}$ in
 $\Def(\Sigma)\times B[k]  \times
       \Ker(\pi_{E_{\rho}}\circ D_{\!f}\bar{\partial}_J)^{\pd}$
 that is characterized by such target-space-driven deformations
 of the domain.
With the higher-order terms omitted, the germ of a target-space-driven
 deformation of $\Sigma$ at a distinguished node is modelled on
 the family of maps from ${\Bbb C}^2/\,{\Bbb C}$ to
 ${\Bbb C}^2/\,{\Bbb C}$ given by
 $$
  \hspace{4em}
  \begin{array}{cccccl}
   {\Bbb C}^2  &  & \longrightarrow  & {\Bbb C}^2 \\[-1ex]
  & \mbox{\small $(z_1,z_2)$} &
    & & \mbox{\small $(w_1,w_2)=(a_1\,{z_1}^s\,,\,a_2\,{z_2}^s)$}  \\
  \downarrow &&& \downarrow \\[.6ex]
   {\Bbb C} &  & \longrightarrow  & {\Bbb C} \\[-1ex]
  & \mbox{\small $\mu=z_1z_2$}  &
    & & \mbox{\small $\lambda\, =\, w_1w_2 \,=\, a_1a_2\,\mu^s$}  &,
  \end{array} \hspace{2em}
   \mbox{\small $a_1$, $a_2\;\in {\smallBbb C}-\{0\}$}\,,
 $$
 where the local target (resp.\ domain) deformations
  are parameterized by $\lambda$ (resp.\ $\mu\,$).
Such constraints from deformation theory at distinguished nodes
 select a subset $\widetilde{V}_{\rho}$ in
 $\Def(\Sigma)\times B[k]  \times
       \Ker(\pi_{E_{\rho}}\circ D_{\!f}\bar{\partial}_J)^{\pd}$
 described as follows.

Fix a factorization
 $$
   \Def(\Sigma)\;
    =\; \Def(\Sigma;\Lambda)\,
            \times\, H_{\rho,\domain}^{(\smooth,\Lambda)}\,,
 $$
 where
  $H_{\rho,\domain}^{(\smooth,\Lambda)}$ is the space of
   local smoothing of distinguished nodes in $\Lambda$ and
  $\Def(\Sigma;\Lambda)$ consists of deformations of $\Sigma$
   that keep $\Lambda$ as nodes.
$H_{\rho,\domain}^{(\smooth,\Lambda)}$ is a neighborhood of
 $0\in {\Bbb C}^{\,|\Lambda|}$, with coordinates
 $(\vec{\lambda}_0,\,\cdots,\,\vec{\lambda}_k)$ with $0$
 corresponding to no smoothing of nodes in $\Lambda$.
Let $H_{\rho,\map}^{(\loc,\Lambda)}\subset {\Bbb C}^{\,|\Lambda|}$
 be a neighborhood of
 $\shiftproduct_{\Lambda}\circ\jet_{\Lambda}^{\mathbf s}(0)$ in
 $({\Bbb C}-\{0\})^{|\Lambda|}$,
 with coordinates $\vec{a}=(\vec{a}_0,\,\cdots,\,\vec{a}_k)$,
 over which the map
  $\shiftproduct_{\Lambda}\circ\jet_{\Lambda}^{\mathbf s}
   :\Ker(\pi_{E_{\rho}}\circ D_{\!f}\bar{\partial}_J)^{\pd}
     \rightarrow {\Bbb C}^{\,|\Lambda|}$ is a bundle map
 (of fiber dimension
  $\mu(f)\, + \, \dimm Y \cdot (1-\tilde{g})\,+\,\dimm E_{\rho}$).
\begin{quote}
 {\bf [Choice].}
 From now on in the construction, we will assume that the local chart
  around $q\in\Lambda$ is chosen so that $a_{q,1}=a_{q.2}$
  in the normal form expression of $f_{\rho}$ around $q$,
  cf.\ Sec.~5.2.
 Fix a section from $H_{\rho,\map}^{(\loc,\Lambda)}$ to
  $\Ker(\pi_{E_{\rho}}\circ D_{\!f}\bar{\partial}_J)^{\pd}$
  so that the condition $a_{q,1}=a_{q.2}$ is preserved for
  all $q\in\Lambda$.
 The value in
  $\Ker(\pi_{E_{\rho}}\circ D_{\!f}\bar{\partial}_J)^{\pd}$
  of this section for $\vec{a}\in H_{\rho,\map}^{(\loc,\Lambda)}$
  will be denoted $\xi_{\vec{a}}$.
\end{quote}
This gives a trivialization
 $$
  \Ker(\pi_{E_{\rho}}\circ D_{\!f}\bar{\partial}_J)^{\pd}_{\Lambda}\;
   :=\; (\shiftproduct_{\Lambda}\circ\jet_{\Lambda}^{\mathbf s})^{-1}
         (H_{\rho,\map}^{(\loc,\Lambda)})\;
   \simeq\; H_{\rho,\map}^{(\loc,\Lambda)}
             \times H_{\rho,\map}^{\,(0,\Lambda)}\,.
 $$
(By convention we fix coordinates on $H_{\rho,\map}^{\,(0,\Lambda)}$
  so that the afore-mentioned section has image
  $H_{\rho,\map}^{(\loc,\Lambda)}\times \{0\}$  in
  $H_{\rho,\map}^{(\loc,\Lambda)}\times H_{\rho,\map}^{\,(0,\Lambda)}$.)
Combining the two, one has a decomposition of the relevant open
 neighborhood of the origin of
 $\Def(\Sigma)\times B[k]  \times
        \Ker(\pi_{E_{\rho}}\circ D_{\!f}\bar{\partial}_J)^{\pd}\,$:
 \begin{eqnarray*}
  \lefteqn{
   \Def(\Sigma)\times B[k] \times
   \Ker(\pi_{E_{\rho}}
          \circ D_{\!f}\bar{\partial}_J)^{\pd}_{\Lambda} } \\[.6ex]
   && \simeq\;
           \Def(\Sigma;\Lambda)\,
           \times\, \left(\,\rule{0ex}{2.2ex}
            H_{\rho,\domain}^{(\smooth,\Lambda)}\,
             \times\, B[k] \,
             \times\, H_{\rho,\map}^{(\loc,\Lambda)}\, \right)\,
           \times\, H_{\rho,\map}^{\,(0,\Lambda)} \\[.6ex]
   && \subset\;
         \Def(\Sigma;\Lambda)\,
          \times\, \left(\,\rule{0ex}{2.2ex}
           {\Bbb C}^{\,|\Lambda|}\times{\Bbb C}^{\,k+1}
                  \times{\Bbb C}^{\,|\Lambda|}\,  \right)\,
          \times\, H_{\rho,\map}^{\,(0,\Lambda)}\,.
 \end{eqnarray*}
 The product
  ${\Bbb C}^{\,|\Lambda|}\times{\Bbb C}^{\,k+1}
            \times{\Bbb C}^{\,|\Lambda|}$
  has coordinates
   $(\,\vec{\mu}_0,\,\cdots,\,\vec{\mu}_k\,;\,
            \vec{\lambda}\,;\,\vec{a}_0,\,\cdots,\,\vec{a}_k\,)\,$
   with
    $$
     \vec{\mu}_i=(\mu_{i1},\,\cdots,\,\mu_{i,|\Lambda_i|})\,,
       \hspace{1em}
     \vec{\lambda}=(\lambda_0,\,\cdots,\,\lambda_k)\,,
       \hspace{1em}\mbox{and}\hspace{1em}
     \vec{a}_i=(a_{i1},\,\cdots,\,a_{i,|\Lambda_i|})
    $$
 that correspond to the deformations of domain, target,
  and maps respectively around $\Lambda$.

Compare this with the basic deformation model above,
one concludes that in terms of these coordinates,
the subset $\widetilde{V}_{\rho}$ of
 $\Def(\Sigma)\times B[k] \times
   \Ker(\pi_{E_{\rho}}\circ D_{\!f}\bar{\partial}_J)^{\pd}_{\Lambda}\,$
 is described by a system of algebraic equations on the
 $\left(\rule{0ex}{2ex}\right.\!
    H_{\rho,\domain}^{(\smooth,\Lambda)}\, \times\, B[k]\,
                \times\,H_{\rho,\map}^{(\loc,\Lambda)}
   \!\left.\rule{0ex}{2ex}\right)$-factor$\,$:
 $$
  \begin{array}{ccl}
   \widetilde{V}_{\rho}
    & = & \left\{\,
       (\,\cdots\,;\,
        \vec{\mu}_0,\,\cdots,\,\vec{\mu}_k\,;\,
        \vec{\lambda}\,;\,\vec{a}_0,\,\cdots,\,\vec{a}_k\,;\,\cdots\,)\,
       \left|\,
        \begin{array}{l}
         {\mu_{ij}}^{s_{ij}}\;=\; \lambda_{\,i}/a_{ij}\,, \\[.6ex]
          i\,=\, 0,\,\ldots,\,k\,;\; j\,=\, 1,\,\ldots,\,|\Lambda_i|\,
        \end{array}
       \right.
      \right\}  \\[3ex]
   & =: & \Def(\Sigma;\Lambda)\,
          \times\,\overline{V}_{\rho}\,
          \times\, H_{\rho,\map}^{(0,\Lambda)}\,.
  \end{array}
 $$
As each $a_{ij}$ takes values in a simply-connected domain in
 ${\Bbb C}-\{0\}$,
 $$
  \widetilde{V}_{\rho}\;
   \simeq\;  \Def(\Sigma;\Lambda)
              \times \Xi_{\mathbf s}
             \times H_{\rho,\map}^{(\loc,\Lambda)}
             \times H_{\rho,\map}^{\,(0,\Lambda)}\;
   =\; \Def(\Sigma;\Lambda)
        \times \Xi_{\mathbf s}
        \times \Ker(\pi_{E_{\rho}}
                    \circ D_{\!f}\bar{\partial}_J)^{\pd}_{\Lambda}
 $$
 in the category of piecewise-smooth stratified spaces,
 where $\Xi_{\mathbf s}$ is defined in Sec.~5.1.

The projection map from
 $\Def(\Sigma)\times B[k]  \times
    \Ker(\pi_{E_{\rho}}\circ D_{\!f}\bar{\partial}_J)^{\pd}\,$
 to $B[k]$ restricts to a morphism
 $\pi_{B[k]}:\widetilde{V}_{\rho}\rightarrow B[k]$
 of constant fiber dimension
  $\mu(f)+\dimm Y\cdot(1-\tilde{g})+\dimm\Def(\Sigma)+\dimm E_{\rho}$.
The restriction of $\widetilde{V}_{\rho}$ over each stratum
 of $B[k]$ can be made a trivial bundle under $\pi_{B[k]}$.
On the other hand, the projection map from
 $\Def(\Sigma)\times B[k]  \times
    \Ker(\pi_{E_{\rho}}\circ D_{\!f}\bar{\partial}_J)^{\pd}\,$
 to $\Def(\Sigma)\times B[k]$ restricts to a morphism
 $\pi_{\scriptsizeDef(\Sigma)\times B[k]}:
  \widetilde{V}_{\rho}\rightarrow \Def(\Sigma)\times B[k]$
  whose image is only a constructible subset in a neighborhood of
   $0\in Def(\Sigma)\times B[k]$ and
  whose fiber dimensions is given by the upper semi-continuous
   function $\ind^{\pd}(D_{\bullet}\bar{\partial}_J) + \dimm E_{\rho}$.

\bigskip

\noindent
{\bf Definition/Convention 5.3.2.1 [linear/nonlinear coordinates on
   $\Ker(\pi_{E_{\rho}}\!\circ D_{\!f}\bar{\partial}_J)^{\pd}$].} {\rm
  Coordinates of
   $\Ker(\pi_{E_{\rho}}\!\circ D_{\!f}\bar{\partial}_J)^{\pd}$
   as a subset of a vector space will be called
   {\it linear coordinates} on
   $\Ker(\pi_{E_{\rho}}\!\circ D_{\!f}\bar{\partial}_J)^{\pd}$.
  Those from the isomorphism with
   $H_{\rho,\map}^{(\loc,\Lambda)}
     \times H_{\rho,\map}^{\,(0,\Lambda)}$
   will be called {\it nonlinear coordinates}.
  Unless otherwise mentioned, we adopt by convention
   the nonlinear coordinates for
    $\Ker(\pi_{E_{\rho}}\!\circ D_{\!f}\bar{\partial}_J)^{\pd}$
   (particularly when written as coordinates from the factorization)
  except the {\it origin}
   $0\in \Ker(\pi_{E_{\rho}}\!\circ D_{\!f}\bar{\partial}_J)^{\pd}$.
} % end-definition

\bigskip

Finally, $\Aut(\rho)$ acts on
 $\Def(\Sigma)\times B[k]\times
  \Ker(\pi_{E_{\rho}}\circ D_{\!f}\bar{\partial}_J)^{\pd}$.
Shrinking if necessary, we take $\widetilde{V}_{\rho}$ to be
 $\Aut(\rho)$-invariant in the above construction.

\bigskip

\subsubsection{A $\widetilde{V}_{\rho}\,$-family of
 approximate-$J$-stable $C^{\infty}$ maps to fibers of $W[k]/B[k]$.}

We construct in this subsubsection an
 $\Aut(\rho)$-invariant $\widetilde{V}_{\rho}$-family
 of approximate-$J$-holomorphic $C^{\infty}$ maps
 $h_{\approxi,\,\raisebox{-.6ex}{\LARGE $\cdot$}}$
 from deformed $\Sigma$ to fibers of $W[k]/B[k]$
 by gluing maps around nodes of $\Sigma$.
 %%%%%%%%%%
 % \marginpar{\raggedright\tiny\vspace{-3em} $\bullet$
 %    Describe as tersely as possible. \newline\newline
 %    {\bf CAUTION.} {\it [F-O] constructs
 %     a full family of approximate $J$-stable $C^{\infty}$ maps
 %     while [Liu(C)] seems to construct a family only
 %      over the space of smoothings of all the nodes.}  \newline
 %    Think more carefully. }
 %%%%%%%%%%
Such construction is given in [MD-S1: Appendix~A] and in
 [F-O], [Liu(C)], [Liu(G)], [R-T1], [R-T2], [Sal], and [I-P2], [L-R]
 for various extensions.

To separate the effect from various types of deformations involved,
 the factorization
 $\Def(\Sigma) = \Def(\Sigma;\Lambda)
                  \times H_{\rho,\domain}^{(\smooth,\Lambda)}$
 is refined to
 $$
  \begin{array}{l}
   \Def(\Sigma)\quad (\,=:\; H_{\rho,\domain}\,) \\[.6ex]
    \hspace{2em}
    =\; \left(
         H_{\rho,\domain}^{(\deform,\Sigma)}
         \times H_{\rho,\domain}^{(\smooth,\oin)}
         \times H_{\rho,\domain}^{(\smooth,\bn)}\right)
         \times H_{\rho,\domain}^{(\smooth,\Lambda)}\,,
  \end{array}
 $$
 where
  $H_{\rho,\domain}^{(\deform,\Sigma)}$
   consists of deformations of the {\it complex structure}
   on $\Sigma$ (as a bordered Riemann surface with marked points)
    without changing the topology of $\Sigma$,
  $H_{\rho,\domain}^{(\smooth,\oin)}$
   consists of local deformations of $\Sigma$ that
   smooth some {\it ordinary interior nodes} of $\Sigma$, and
  $H_{\rho,\domain}^{(\smooth,\bn)}$
   consists of local deformations of $\Sigma$ that smooth
   some {\it boundary nodes} of $\Sigma$.
For $\Sigma$ of genus $g$, $h$ holes, $n_{oin}$ ordinary interior nodes,
 $|\Lambda|$ distinguished interior nodes, $n_{bn}$ boundary nodes,
 $n$ ordinary marked points, and $|\vec{m}|$ boundary marked points,
$H_{\rho,\domain}$ is parameterized by a neighborhood of
 ${\mathbf 0}$ in the $4$-factor product space
 (with coordinates $(\zeta, \vec{t}, \vec{t^{\prime}},\vec{\mu})$)
 $$
  \begin{array}{c}
  (\,{\Bbb C}^{3g-3+h-n_{oin}-|\Lambda|+n^{\prime}+d_c}  \times
     \overline{\Bbb H}^{n^{\prime\prime}} \times
     {\Bbb R}^{h-n_{bn}+|\vec{m}|+d_b}\,)      \times
     {\Bbb C}^{n_{oin}} \times {{\Bbb R}_{\ge 0}}^{n_{bn}}  \times
     {\Bbb C}^{|\Lambda|}\,,
      \hspace{1em}n\doteq n^{\prime}+n^{\prime\prime}\,,
  \end{array}
 $$
 with respect to the above decomposition.
Let ${\cal C}/\Def(\Sigma)$ be the universal curve over $\Def(\Sigma)$,
 with the fiber labelled-bordered Riemann surface-with-marked-points
  over $(\zeta,\vec{t},\vec{t}^{\prime},\vec{\mu})\in\Def(\Sigma)$
  denoted by $\Sigma_{(\zeta,\vec{t},\vec{t}^{\prime},\vec{\mu})}$.
With a fixed local model chart at each node of $\Sigma$,
   cf.\ Definiton 2.1,
        % Definition [prestable labelled-bordered Riemann surface]
 a fixed $\varepsilon>0$ small, and
 the assumption that
  $\|(\vec{t},\vec{t}^{\prime},\vec{\mu})\|\ll \varepsilon$,
following the same construction as in the case of $W[k]/B[k]$,
there is a {\it $\varepsilon$-neck-trunk decomposition}\footnote{Cf.\
                             the thick-thin decompsoition
                             in terms of hyperbolic geometry.}
 of ${\cal C}/\Def(\Sigma)$ and gluing maps\footnote{Cf.\ the maps
                             $I_{\vec{\lambda}}:
                              Y_{[k]} -
                               \cup_{i=0}^k\,N_{\sqrt{|\lambda_i|}}(D_i)
                               \rightarrow W[k]_{\vec{\lambda}}$
                              and
                             $I_{\vec{\lambda},\varepsilon}:
                              Y_{[k]} -
                               \cup_{i=0}^k\, N_{|\lambda_i|/\varepsilon}
                                              (D_i)
                              \rightarrow W[k]_{\vec{\lambda}}$
                              defined in Sec.~1.1.1 by cut-and-glue.}:
 $$
  \begin{array}{lllll}
   I_{(0,\vec{t},\vec{t}^{\prime},\vec{\mu})}  & :
    & \Sigma- \cup_{q:\node}\, N_{\sqrt{|t_q|}}(q)
    & \longrightarrow
    & \Sigma_{(0,\vec{t},\vec{t}^{\prime},\vec{\mu})}\,, \\[.6ex]
   I_{(0,\vec{t},\vec{t}^{\prime},\vec{\mu}),\varepsilon}  & :
    & \Sigma- \cup_{q:\node}\, N_{|t_q|/\varepsilon}(q)
    & \longrightarrow
    & \Sigma_{(0,\vec{t},\vec{t}^{\prime},\vec{\mu})}\,,
  \end{array}
 $$
 where
  $t_q$ is the entry of $(\vec{t},\vec{t}^{\prime},\vec{\mu})$
   associated to the node $q$,  and
  $N_{(\cdots)}(q)$ is the $(\cdots)$-neighborhood of $q$
   in the local model of node $q$.
We also have a fixed family of diffeomorphisms
 $\Sigma_{(\zeta,\vec{t},\vec{t}^{\prime},\vec{\mu})}
  \simeq \Sigma_{(0,\vec{t},\vec{t}^{\prime},\vec{\mu})}$.
The combination of the two defines the gluing maps
 $$
  \begin{array}{lllll}
   I_{(\zeta,\vec{t},\vec{t}^{\prime},\vec{\mu})}  & :
    & \Sigma- \cup_{q:\node}\, N_{\sqrt{|t_q|}}(q)
    & \longrightarrow
    & \Sigma_{(\zeta,\vec{t},\vec{t}^{\prime},\vec{\mu})}\,, \\[.6ex]
   I_{(\zeta,\vec{t},\vec{t}^{\prime},\vec{\mu}),\varepsilon}  & :
    & \Sigma- \cup_{q:\node}\, N_{|t_q|/\varepsilon}(q)
    & \longrightarrow
    & \Sigma_{(\zeta,\vec{t},\vec{t}^{\prime},\vec{\mu})}\,.
  \end{array}
 $$
These maps satisfy the $\Aut(\rho)$-conjugation property that
 $$
  \begin{array}{ccl}
   \alpha  \circ
    I_{(\zeta,\vec{t},\vec{t}^{\prime},\vec{\mu})} \circ \alpha^{-1}
     & = & I_{\alpha\cdot (\zeta,\vec{t}, \vec{t}^{\prime},\vec{\mu})}\,,
                                                             \\[.6ex]
   \alpha  \circ
    I_{(\zeta,\vec{t},\vec{t}^{\prime},\vec{\mu}),\varepsilon}  \circ
     \alpha^{-1}
     & =
     & I_{\alpha\cdot
          (\zeta,\vec{t}, \vec{t}^{\prime},\vec{\mu}),\varepsilon}
  \end{array}
 $$
 for $\alpha\in\Aut(\rho)^{\domain}$ acting on ${\cal C}/\Def(\Sigma)$.
The $\varepsilon$-neck region of the fiber
 $\Sigma_{(\zeta,\vec{t},\vec{t}^{\prime},\vec{\mu})}$
 of ${\cal C}/\Def(\Sigma)$ will be denoted by
 $\Neck_{\varepsilon\,,\,(\zeta,\vec{t},\vec{t}^{\prime},\vec{\mu})}$.
It is a disjoint union of annuli/strips,
 of the form
  $$
   \{(z_1,z_2)\in{\Bbb C}^2\;:\;
      z_1z_2=t_q\,,\; |z_1|<\varepsilon\,,\; |z_2|<\varepsilon\}
  $$
 (resp.\
  $$
   \begin{array}{c}
    \{(z_1,z_2)\; :\; z_1z_2=t_q\,,\;  |z_1|<\varepsilon\,,\;
                      |z_2|<\varepsilon\}/
         (z_1,z_2)\sim(\overline{z_2},\overline{z_1})\,, % type E
                                                         \\[.6ex]
    \{(z_1,z_2)\; :\; z_1z_2=t_q\,,\;  |z_1|<\varepsilon\,,\;
                      |z_2|<\varepsilon\}/
         (z_1,z_2)\sim(\overline{z_1},\overline{z_2})\;) % type H
   \end{array}
  $$
 in $\Sigma_{(\zeta,\vec{t},\vec{t}^{\prime},\vec{\mu})}$,
 associated to smoothed interior
  (resp.\ type-E boundary, type-H boundary) nodes $q$ of $\Sigma$.

To homogenize the notation, we write interchangeably
 $H_{\rho,\target}:= B[k]$
  for the deformations of the target $Y_{[k]}$, and
 $H_{\rho,\map} :=
  \Ker(\pi_{E_{\rho}}\circ D_{\!f}\bar{\partial}_J)^{\pd}_{\Lambda}$
  as the deformation space of $f$ with
   the fixed domain $\Sigma$ and rigid target $Y_{[k]}$.
The coordinates for $H_{\rho,\target}\times H_{\rho,\map}$ will be
 denoted by $(\vec{\lambda},\vec{a},\xi)$ with respect to its
 decomposition as
 $B[k] \times
       H_{\rho,\map}^{(\loc,\Lambda)}\times H_{\rho,\map}^{(0,\Lambda)}$,
 (cf.\ Definition/Convention 5.3.2.1).
   % Definition/Convention [linear/nonlinear coordinates on
   %    $\Ker(\pi_{E_{\rho}}\!\circ D_{\!f}\bar{\partial}_J)^{\pd}$]
Recall then
  $$
   \begin{array}{cccl}
    \widetilde{V}_{\rho}
     & \subset
     & H_{\rho}\; :=\;
       H_{\rho,\domain}\times H_{\rho,\target}\times H_{\rho,\map} \\[.6ex]
    \mbox{\scriptsize $\pi_{\mbox{\tiny {\it Def}$(\Sigma)\times B[k]$}}$}
      \downarrow \hspace{5.2em} &
     & \hspace{4em} \downarrow
       \mbox{\scriptsize
             $\pi_{\mbox{\tiny {\it Def}$(\Sigma)\times B[k]$}}$}  \\[.6ex]
    \pi_{\scriptsizeDef(\Sigma)\times B[k]}(\widetilde{V}_{\rho})
     & \subset & \Def(\Sigma)\times B[k] \hspace{2em}   & .
   \end{array}
  $$
We will use the product coordinates
 $(\zeta, \vec{t}, \vec{t^{\prime}},\vec{\mu},\vec{\lambda},\vec{a},\xi)$
 of $H_{\rho}$ for the algebraic subset $\widetilde{V}_{\rho}$,
 with $\vec{\lambda}$ being the redundant coordinates expressible in terms
 of $(\vec{\mu},\vec{a})$, (cf.\ Sec.~5.3.2).

The intersections
  $$
   \begin{array}{lcl}
   \Theta_{\rho,0}
    & :=
    & \left( H_{\rho,\domain}\times H_{\rho,\target}
        \times
        \{ \shiftproduct_{\Lambda}\circ\jet_{\Lambda}^{\mathbf s}(0) \}
        \times \{0\} \right)
      \cap \widetilde{V}_{\rho}\,, \\[1.6ex]
   \Theta_{\rho}
    & :=
    & \left(\rule{0ex}{2ex}\right.
        H_{\rho,\domain}\times H_{\rho,\target}
        \times H_{\rho,\map}^{(\loc,\Lambda)} \times \{0\}
       \left.\rule{0em}{2ex}\right)
      \cap \widetilde{V}_{\rho} \\[1ex]
    & \simeq
    & \Def(\Sigma;\Lambda)\times \overline{V}_{\rho}
   \end{array}
  $$
 in $H_{\rho}$ are both connected constructible subsets of $H_{\rho}$.
$\Theta_{\rho,0}$ is a deformation retract of
 $\Theta_{\rho}$ and, hence,
$\pi_{\scriptsizeDef(\Sigma)\times B[k]}(\Theta_{\rho,0})$
 is a deformation retract of
 $\pi_{\scriptsizeDef(\Sigma)\times B[k]}(\widetilde{V}_{\rho})$.
The restriction of $\pi_{\scriptsizeDef(\Sigma)\times B[k]}$ to
 $\Theta_{\rho,0}$ is one-to-one.
Its inverse defines a (continuous) section
 $$
  S_0\; :\;
   \pi_{\scriptsizeDef(\Sigma)\times B[k]}(\Theta_{\rho,0})\;
   \longrightarrow\;
   \widetilde{V}_{\rho}|
       _{\pi_{\tinyDef(\Sigma)\times B[k]}(\Theta_{\rho,0})}
 $$
 with image $\Theta_{\rho,0}$.
The restriction of $\pi_{\scriptsizeDef(\Sigma)\times B[k]}$ on
 $\Theta_{\rho}$ is one-to-one only on an open dense subset
 (i.e.\ the subset described by $\lambda_i\ne 0$, $i=0,\,\ldots,\,k$).
It follows that $S_0$ extends uniquely to
 a {\it piecewise-continuous} section
 $$
  S\; :\;
   \pi_{\scriptsizeDef(\Sigma)\times B[k]}(\widetilde{V}_{\rho})\;
      \longrightarrow\;  \widetilde{V}_{\rho}\,,
 $$
 whose image has closure $\Theta_{\rho}$ in $\widetilde{V}_{\rho}$.
Both $S_0$ and $S$ are $\Aut(\rho)$-equivariant.
As $\widetilde{V}_{\rho}$ is a bundle over $\Theta_{\rho}$
 with fiber $H_{\rho,\map}^{(0,\Lambda)}$, this says in particular that,
while it is not possible to make all the ingredients in
 the relative construction (of $\widetilde{V}_{\rho}$-family of maps)
 continuous with respect to
 $\pi_{\scriptsizeDef(\Sigma)\times B[k]}(\widetilde{V}_{\rho})$
 in $\Def(\Sigma)\times B[k]$,
we have to ensure their extendibility and continuity over
 $\Theta_{\rho}$.
We now proceed to construct a
 piecewise-continuous-$\pi_{\scriptsizeDef(\Sigma)
                     \times B[k]}(\widetilde{V}_{\rho})$-family
 of approximate-$J$-holomorphic maps that extends to a
 continuous-$\Theta_{\rho}$-family of approximate-$J$-holomorphic maps.

Fix a rotation-invariant smooth cutoff function
 $\beta_1:{\Bbb C}\rightarrow [0,1]$ such that
 $$
  \beta_1(z)\;
   =\; \left\{ \begin{array}{cl}
                1  & \mbox{if $|z|\ge 2$}\,,\\[.6ex]
                0  & \mbox{if $|z|\le 1$}\,,
               \end{array} \right.
   \hspace{1em}\mbox{and}\hspace{1em} |\nabla\beta_1\,|\le 2\,,
 $$
([MD-S1:~Lemma A.1.1]).
Then the local model of our approximate-$J$-stable maps
 around a smoothed node is given as follows
 for a fixed $\varepsilon>0$ small and
       $|t|,|t^{\prime}|, |\mu|, |\lambda| \ll \varepsilon$.
(Cf.\ [MD-S: Sec.~A.2], [F-O: (12.13)], [Liu(C): Sec.~6.4.1], and
      [L-R:  Sec.~4.1].)

\bigskip

\noindent $(a)$
$H_{\rho,\domain}^{(\smooth,\oin)}\,:\,$
The local model of the deformation/smoothing of an ordinary interior
 node $q$ of $\Sigma$ is given by
 $$
  \begin{array}{cccl}
   B_{\varepsilon}
    := \{(z_1,z_2)\in {\Bbb C}^2\,|\, |z_1|, |z_2|\le \varepsilon \}
    & \longrightarrow           & {\Bbb C}        \\[.6ex]
   (z_1,z_2)   & \longmapsto    & z_1z_2    & .
  \end{array}
 $$
Let $A_t$ be the fiber over $t\in {\Bbb C}$ and
 $f_{\rho}|_{A_0}=f_1\cup f_2$;
then, for $t\ne 0$, define $h_t:A_t\rightarrow Y_{[k]}$ by
$$
 h_t(z,\frac{t}{z})\;
  =\; \exp_{f(q)} \left(
        \beta_1\left(\frac{z}{|t|^{1/4}}\right) \exp_{f(q)}^{-1}(f_1(z))\,
     +\,\beta_1\left(\frac{|t|^{3/4}}{z}\right)
           \exp_{f(q)}^{-1}\left(f_2\left(\frac{t}{z}\right)\right)
             \right)\,.
$$

\bigskip

\noindent $(b)$
$H_{\rho,\domain}^{(\smooth,\bn)}\,:\,$
The local model of the deformation/smoothing of the two types of
 boundary node $q$ of $\Sigma$ is given respectively by
 $$
  \parbox{5em}{(type E)}
  \begin{array}{cccl}
   B_{\varepsilon}/\!\sim_{\rm E}\;
    :=\; \{(z_1,z_2)\in {\Bbb C}^2\,|\, |z_1|, |z_2|\le \varepsilon \}
          /(z_1, z_2)\sim (\overline{z_2},\overline{z_1})
    & \longrightarrow           & {\Bbb R}_{\ge 0}        \\[.6ex]
   (z_1,z_2)   & \longmapsto    & z_1z_2    & ,
  \end{array}
 $$
 $$
  \parbox{5em}{(type H)}
  \begin{array}{cccl}
   B_{\varepsilon}/\!\sim_{\rm H}\;
    :=\; \{(z_1,z_2)\in {\Bbb C}^2\,|\, |z_1|, |z_2|\le \varepsilon \}
          /(z_1, z_2)\sim (\overline{z_1},\overline{z_2})
    & \longrightarrow           & {\Bbb R}_{\ge 0}        \\[.6ex]
   (z_1,z_2)   & \longmapsto    & z_1z_2    & .
  \end{array}
 $$

For $q$ of type E, let $A^{\prime}_{t^{\prime}}$ be the fiber
 over $t^{\prime}\in {\Bbb R}_{\ge 0}$ and $f_{\rho}|_{A^{\prime}_0}=f$;
then, for $t^{\prime}>0$,
 define $h_{t^{\prime}}:A^{\prime}_{t^{\prime}}\rightarrow Y_{[k]}$ by
 $$
  h_{t^{\prime}}(z,\frac{t^{\prime}}{z})\;
   =\; \exp_{f(q)}  \left(
         \beta_1\left(\frac{z}{|t^{\prime}|^{1/4}}\right)
           \exp_{f(q)}^{-1}(f(z))\,
                     \right)\,.
 $$

For $q$ of type H, let $A^{\prime}_{t^{\prime}}$ be the fiber over
 $t^{\prime}\in {\Bbb R}_{\ge 0}$ and
 $f_{\rho}|_{A^{\prime}_0}=f_1\cup f_2$;
then, for $t^{\prime}>0$,
 define $h_{t^{\prime}}:A^{\prime}_{t^{\prime}}\rightarrow Y_{[k]}$ by
 $$
  h_{t^{\prime}}(z,\frac{t^{\prime}}{z})\;
   =\; \exp_{f(q)}  \left(
         \beta_1\left(\frac{z}{|t^{\prime}|^{1/4}}\right)
           \exp_{f(q)}^{-1}(f_1(z))\,
      +\,\beta_1\left(\frac{|t^{\prime}|^{3/4}}{z}\right)
           \exp_{f(q)}^{-1}
             \left(f_2\left(\frac{t^{\prime}}{z}\right)\right)
                    \right)\,.
 $$

\bigskip

\noindent $(c)$
$\overline{V}_{\rho} \subset
 H_{\rho,\domain}^{(\smooth,\Lambda)}\times H_{\rho,\target}
  \times H_{\rho,\map}^{(\loc,\Lambda)}\,:\,$
Let $q\in \Lambda$ be a distinguished node of contact order $s$.
Recall the fixed local coordinates around $f(q)$.
Denote by $(f_1^D, f_1^N)\cup (f_2^D, f_2^N)$
 the restriction of $f$ around $q$ with the expression
 in terms of the coordinates on $D$ and the normal coordinate
  to $D$ around $f(q)$.
Recall also the local model in Sec.~5.3.2 (cf.\ [I-P2])
 $$
  \hspace{4em}
  \begin{array}{clccll}
   B_{\varepsilon}
   &  & \longrightarrow  & {\Bbb C}^2 \\[-1ex]
  & \mbox{\small $(z_1,z_2)$} &
    & & \mbox{\small $(w_1,w_2)=(a_1\,{z_1}^s\,,\,a_2\,{z_2}^s)$}  \\
  \downarrow &&& \downarrow \\[.6ex]
   {\Bbb C} &  & \longrightarrow  & {\Bbb C} \\[-1ex]
  & \mbox{\small $\mu=z_1z_2$}  &
    & & \mbox{\small $\lambda\, =\,w_1w_2 \,=\,a_1a_2\,\mu^s\,
                                =\, a\,\mu^s$}  &,
  \end{array} \hspace{2em}
   \mbox{\small $a_1$, $a_2\;\in {\smallBbb C}-\{0\}$}\,,
 $$
that links
 the deformation/smoothing (here parameterized by $\mu$) of the node $q$,
 the deformation (here parameterized by $\lambda$) of $Y_{[k]}$
  along the $D_i$ that contains $f(q)$, and
 the product (here parameterized by $a$) of the lowest-order
  pre-deformable deformations of the normal-to-$D$ component of
  the germ of $f$ on the two branches of $\Sigma$ at $q$.

Let $(\mu,\lambda,a)$ be the relevant coordinates
 in the coordinates of
  $H_{\rho,\domain}^{(\smooth,\Lambda)} \times H_{\rho,\target}
    \times H_{\rho,\map}^{(\loc,\Lambda)}$  with $\lambda= a\mu^s$.
Define
 $h_{(\mu,\lambda,a)}=(h_{(\mu,\lambda,a)}^D,h_{(\mu,\lambda,a)}^N)$
 as follows:

 \bigskip
 \noindent
 $\cdot$ For $\lambda=0\,$:
 recall the $\xi_{a}$ in
  $\Ker(\pi_{E_{\rho}}\circ D_{\!f}\bar{\partial}_J)^{\pd}$
  associated to $a$ and define
 $$
  h_{(0,0,a)}(\,\cdot\,)\; =\; \exp_{f(\,\cdot\,)}(\xi_a(\,\cdot\,))\,.
 $$

 \bigskip
 \noindent
 $\cdot$ For $\lambda\ne 0\,$:
 express $h_a:=h_{(0,0,a)}$ above as
  $h_{a,1}\cup h_{a,2}
      =(h_{a,1}^D,h_{a,1}^N)\cup (h_{a,2}^D,h_{a,2}^N)$ and define
 $$
 \begin{array}{lll}
  h_{(\mu,\lambda,a)}^D(z,\frac{\mu}{z})
   & = &  \exp_{h_a(q)} \left(
            \beta_1\left(\frac{z}{|\mu|^{1/4}}\right)
            \exp_{h_a(q)}^{-1}(h_{a,1}^D(z))\,
         +\, \beta_1\left(\frac{|\mu|^{3/4}}{z}\right)
             \exp_{h_a(q)}^{-1}\left(h_{a,2}^D
               \left(\frac{\mu}{z}\right)\right)
                  \right)\,,                        \\[3ex]
  h_{(\mu,\lambda,a)}^N(z,\frac{\mu}{z})
   & = & \left\{
   \begin{array}{ll}
    \beta_1\left(\frac{z}{|\mu|^{1/4}}\right)\, h_{a,1}^N(z)\,
     +\,\beta_1\left(\frac{|\mu|^{3/4}}{z}\right)\, \sqrt{a}\,z^s
     & \mbox{for $|\mu|^{1/2}\le |z|\le \varepsilon$}\,,\\[1.6ex]
    \beta_1\left(\frac{z}{|\mu|^{1/4}}\right)\,
       \sqrt{a}\,\left(\frac{\mu}{z}\right)^s\,
     +\,\beta_1\left(\frac{|\mu|^{3/4}}{z}\right)\,h_{a,2}^N(\frac{\mu}{z})
     & \mbox{for $|\mu|^{1/2}\le |\mu/z|\le \varepsilon$}\,,
   \end{array} \right.
 \end{array}
 $$

 where $\sqrt{a}$ is chosen so that $\sqrt{a_f}$ fits
  the normal-form expression of $f$ at $q$.

\bigskip

\noindent
This describes what happens on a smoothed neighborhood of $q$
 with all irrelevant indices of the coordinates of
 $H_{\rho,\domain}^{(\smooth,\Lambda)}\times H_{\rho,\target}
   \times H_{\rho,\map}^{(\loc,\Lambda)}$
 suppressed.
The substitutions
 $\mu\rightarrow \mu_{ij}$, $s\rightarrow s_{ij}$,
  $\lambda\rightarrow\lambda_i$, $a\rightarrow a_{ij}$,
  for $i=0,\,\ldots,\,k$, $j=1,\,\ldots,\,|\Lambda_i|$
 to the above expression recover the complete
 $\overline{V}_{\rho}$-family
 of maps from $\Lambda_{0,\vec{t},\vec{t}^{\prime}}$
  to the fiber $W[k]_{\vec{\lambda}}$ of $W[k]/B[k]\,$.

\bigskip

By construction, these maps are defined on disjoint subsets
 of $\Sigma_{(0,\vec{t},\vec{t}^{\prime},\vec{\mu})}$ and
 coincide with $f$ on their intersection with a compact
 subset $K_{\varepsilon_{\varepsilon_-}}$ of $\Sigma$
 by removing a small $\varepsilon_-$-neighborhood of all the nodes,
  with $\varepsilon_-$ slightly less than $\varepsilon$.
As $W[k]_{\vec{\lambda}}$ are obtained from gluing truncated
 $Y_{[k]}$ around $\vec{\lambda}$-specified $D_i$'s
 (cf.\ Sec.~1.1.1),
they can be combined with and extended by $f|_{K_{\varepsilon_-}}$
 to a map from $\Sigma_{(0,\vec{t},\vec{t}^{\prime},\vec{\mu})}$
 to $W[k]_{\vec{\lambda}}$.
In this way, one obtains a
 (continuous-)$H_{\rho,\domain}^{(\smooth,\oin)} \times
  H_{\rho,\domain}^{(\smooth,\bn)}\times \overline{V}_{\rho}$-family
 of maps
 $$
  h_{\approxi,(0,\vec{t},\vec{t}^{\prime},
                     \vec{\mu},\vec{\lambda},\vec{a},\vec{0})}\; :\;
   \Sigma_{(0,\vec{t},\vec{t}^{\prime},\vec{\mu})}\;
                        \longrightarrow\; W[k]_{\vec{\lambda}}\,.
 $$
For $\zeta\in H_{\rho,\domain}^{(\deform,\Sigma)}$,
 one defines
 $h_{\approxi,(\zeta,\vec{t},\vec{t}^{\prime},
                 \vec{\mu},\vec{\lambda},\vec{a},\vec{0})}
  : \Sigma_{(\zeta,\vec{t},\vec{t}^{\prime},\vec{\mu})}
                           \rightarrow W[k]_{\vec{\lambda}}$
 by setting
 $$
   h_{\approxi,(\zeta,\vec{t},\vec{t}^{\prime},
                            \vec{\mu},\vec{\lambda},\vec{a},\vec{0})}\;
   =\; h_{\approxi,(0,\vec{t},\vec{t}^{\prime},
                            \vec{\mu},\vec{\lambda},\vec{a},\vec{0})}.
 $$
(In all the discussion, though the $\vec{\lambda}$-label is
 determined uniquely by $(\vec{\mu},\vec{a})$, we keep it
 in the notation to remind us of the change of the target.)
To summarize:

\bigskip

\noindent
{\bf Lemma 5.3.3.1 [pre-deformable $\Theta_{\rho}$-family].} {\it
 $h_{\approxi,(\zeta,\vec{t},\vec{t}^{\prime},
                   \vec{\mu},\vec{\lambda},\vec{a},\vec{0})}$,
  $(\zeta,\vec{t},\vec{t}^{\prime},
       \vec{\mu},\vec{\lambda},\vec{a})\in\Theta_{\rho}$,
 defines a {\rm(}continuous-{\rm )}$\Theta_{\rho}$-family of
  $C^{\infty}$ maps of the same contact order and pre-deformability
  behavior as $f$ at un-smoothed distinguished nodes of $\Sigma$.
} % end-lemma

\bigskip

To keep the relative-to-$(\mbox{domain},\,\mbox{target})$-construction
 picture manifest,
one should think of this $\Theta_{\rho}$-family of maps as
 an extension/completion-at-$[f]$ of the corresponding
 (piecewise-continuous-)$\pi_{\scriptsizeDef(\Sigma)\times B[k]}
                                       (\widetilde{V}_{\rho})$-family
 of maps via the open-dense embedding
 $S:\pi_{\scriptsizeDef(\Sigma)\times B[k]}(\widetilde{V}_{\rho})
    \hookrightarrow \Theta_{\rho}$.
The $\Theta_{\rho}$-family of maps can be extended further to
 a (continuous-)$\widetilde{V}_{\rho}$-family of maps by defining
 first
 $$
  h_{\approxi,(0,\vec{0},\vec{0}^{\prime},
              \vec{0},\vec{0},\vec{a},\vec{b}) }(\,\cdot\,)\;
  =\; \exp_{f(\,\cdot\,)}\xi_{(\vec{a},\vec{b})}(\,\cdot\,)\,.
 $$
This is a $H_{\rho,\map}$-family of $C^{\infty}$ maps from $\Sigma$
 to $Y_{[k]}$ for which pre-deformability at each distinguished node
 remains hold with the same order.
Repeating then the above construction that deforms the map
  at the three types of nodes with $f$ replaced by
  $h_{\approxi,(0,\vec{0},\vec{0}^{\prime},
              \vec{0},\vec{0},\vec{a},\vec{b})}$.
For $(\zeta,\vec{t},\vec{t}^{\prime},
       \vec{\mu},\vec{\lambda},\vec{a},\vec{0})\in \Theta_{\rho}$
 this gives $h_{\approxi, (\zeta,\vec{t},\vec{t}^{\prime},
                          \vec{\mu},\vec{\lambda},\vec{a},\vec{0}))}$
 as constructed above.

\bigskip

\noindent
{\bf Lemma 5.3.3.2
     [$\widetilde{V}_{\rho}$-family of pre-deformable
                approximate-$J$-holomorphic maps].}
{\it
 Assume that $\|\zeta\|$, $\|\vec{t}\|$, $\|\vec{t}^{\prime}\|$,
   $\|\vec{\mu}\|$, $\|\vec{a}-\vec{a}_f\|$, $\|\vec{b}\|$
  are all sufficiently small
  {\rm (}say, bounded above uniformly by an $\varepsilon\ll 1${\rm )},
 then
 $$
  \left\|\,
   \bar{\partial}_J h_{\approxi,(\zeta,\vec{t},\vec{t}^{\prime},
                          \vec{\mu},\vec{\lambda},\vec{a},\vec{b})}\,
   \right\|
    _{L^{\raisebox{.3ex}{\tiny $p$}}(\Sigma
             _{(\zeta,\vec{t},\vec{t}^{\prime},\vec{\mu})})}\;
  \le\;
   C\, \left(\, \rule{0ex}{2.4ex}
           \|\zeta\| + \|\vec{t}\| + \|\vec{t}^{\prime}\|
            + \|\vec{\mu}\| + \|\vec{a}-\vec{a}_f\| + \|\vec{b}\|\,
        \right)^{\frac{1}{2p}}\,,
 $$
 where $C$ is a constant that depends only on
  $\varepsilon$, $f$, $\nabla f$, $J$, $\nabla J$,
  the norm of the differential of
   $\shiftproduct_{\Lambda}\circ\jet_{\Lambda}^{\mathbf s}$,
  and the norm of the differential of the exponential map and
   its inverse along $f$.
 Thus, $h_{\approxi,(\,\cdot\,)}$ gives
  a {\rm (}continuous{\rm )} $\Aut(\rho)$-invariant
  $\widetilde{V}_{\rho}$-family of pre-deformable
  approximate-$J$-holomorphic maps.
} % end-lemma

\bigskip

\noindent
{\it Proof.}
 The approximate $J$-holomorphy property follows from
  [MD-S1: Lemma A.4.3], [F-O: Lemma 12.14, Lemma 12.15],
  [Liu(C): Lemma 6.22], and [L-R: Lemma 4.6].
 Here we have assumed that
  $\|\zeta\|$, $\|\vec{t}\|$, $\|\vec{t}^{\prime}\|$, $\|\vec{\mu}\|$,
  $\|\vec{a}-\vec{a}_f\|$, $\|\vec{b}\|$ are all sufficiently small
 so that the combination of all the estimates in ibidem is bounded
  above by the right-hand side of the inequality above.
 The $\Aut(\rho)^{\domain}$-invariance of
  the domain decomposition involved,
 the $\Aut(\rho)^{\target}$-invariance of the metric on $W[k]$, and
  the cutoff function chosen imply that
  the gluing construction is $\Aut(\rho)$-invariant.
 This implies that the $\widetilde{V}_{\rho}$-family of maps
  as constructed is $\Aut(\rho)$-invariant.

\noindent\hspace{15cm}$\Box$

\bigskip

\noindent
{\bf Notation 5.3.3.3.} {\rm
 We will assume that $\|\zeta\|$, $\|\vec{t}\|$, $\|\vec{t}^{\prime}\|$,
  $\|\vec{\mu}\|$, $\|\vec{a}-\vec{a}_f\|$, $\|\vec{b}\|$
  are all sufficiently small so that
 $h_{\approxi,(\zeta,\vec{t},\vec{t}^{\prime},
       \vec{\mu},\vec{\lambda},\vec{a},\vec{b})}(\cdot)
  = \exp_{h_{\tinyapproxi,S(\zeta,\vec{t},\vec{t}^{\prime},
                            \vec{\mu},\vec{\lambda}) }(\,\cdot\,) }
      \xi_{(\zeta,\vec{t},\vec{t}^{\prime},
             \vec{\mu},\vec{\lambda},\vec{a},\vec{b})}(\,\cdot\,)$
  for a unique

  \vspace{-1ex}
  {\small
  $$
   \begin{array}{l}
    \xi_{(\zeta,\vec{t},\vec{t}^{\prime},
        \vec{\mu},\vec{\lambda},\vec{a},\vec{b})}
      \\ [1ex]
    \hspace{1ex}\in\;
    W^{1,p}
    \left(
         \Sigma_{(\zeta,\vec{t},\vec{t}^{\prime},\vec{\mu})},
         \partial\Sigma_{(\zeta,\vec{t},\vec{t}^{\prime},\vec{\mu})};
         h_{\approxi, S(\zeta,\vec{t},\vec{t}^{\prime},
                                            \vec{\mu},\vec{\lambda})}
           ^{\ast}T_{\ast}(W[k]_{\vec{\lambda}}),
        (h_{\approxi,S(\zeta,\vec{t},\vec{t}^{\prime},
                                       \vec{\mu},\vec{\lambda})}
         |_{\partial\Sigma_{(\zeta,\vec{t},\vec{t}^{\prime},\vec{\mu})}})
        ^{\ast} T_{\ast}L
                \right)\,.
   \end{array}
  $$
  } % end-small
} % end-notation

This expression renders the $\widetilde{V}_{\rho}$-family of maps
 a continuous extension of the $\Theta_{\rho}$-family of maps
 by the exponential-map construction
 along the $H_{\rho,\map}^{(0,\Lambda)}$-factor directions;
this helps making the later relative construction over
 $\pi_{\scriptsizeDef(\Sigma)\times B[k]}(\widetilde{V}_{\rho})$
 manifest.

\bigskip

\subsubsection{The $\widetilde{V}_{\rho}\,$-family of
    (exact) $(J,E_{\mbox{\LARGE $\cdot$}})$-stable maps
    $f_{\mbox{\LARGE $\cdot$}}$ to fibers of $W[k]/B[k]$.}

In this subsubsection, we extend
 the $\Aut(\rho)$-invariant obstruction space $E_{\rho}$ at $\rho$
 step by step to trivialized $\Aut(\rho)$-equivariant auxiliary
 obstruction bundles
  $E_{S(\pi_{\tinyDef(\Sigma)\times B[k]}(\widetilde{V}_{\rho}))}^{\aux}$
   over\\
   $S(\pi_{\scriptsizeDef(\Sigma)\times B[k]}(\widetilde{V}_{\rho}))$,
  $E_{\Theta_{\rho}}$ over $\Theta_{\rho}$, and
  $E_{\widetilde{V}_{\rho}}$ over $\widetilde{V}_{\rho}$.
We then deform the $\Aut(\rho)$-invariant $\widetilde{V}_{\rho}$-family
 of approximate-$J$-stable $C^{\infty}$ maps in Sec.~5.3.3
 to a (continuous) $\Aut(\rho)$-invariant $\widetilde{V}_{\rho}$-family
 of $(J,E_{\mbox{\tiny $\bullet$}})$-stable maps.
The major step is a construction of a
 $\pi_{\scriptsizeDef(\Sigma)\times B[k]}(\widetilde{V}_{\rho})$-family
 of right inverses
 (of $\pi_{E_{\,\cdot}^{\tinyaux}}\!\circ
        D_{h_{\tinyapproxi,\,\cdot\,}}\bar{\partial}_J$)

 {\small
 $$
  \begin{array}{c}
   Q_{(\zeta,\vec{t},\vec{t}^{\prime},\vec{\mu},\vec{\lambda})}\;
    :\; \left.
       L^p \left(
       \Sigma_{(\zeta,\vec{t},\vec{t}^{\prime},\vec{\mu})};
      \Lambda^{0,1}\Sigma_{(\zeta,\vec{t},\vec{t}^{\prime},\vec{\mu})}
        \otimes_J
          h_{\approxi,S(\zeta,\vec{t},\vec{t}^{\prime},
                                         \vec{\mu},\vec{\lambda})}
         ^{\ast}T_{\ast}W[k]_{\vec{\lambda}}  \right)
       \right/\!
        \raisebox{-.6ex}{$E_{S(\zeta,\vec{t},\vec{t}^{\prime},
                               \vec{\mu},\vec{\lambda})}^{\aux}$}\;
      \longrightarrow   \hspace{3em}\\[1ex]
   W^{1,p}  \left(
     \Sigma_{(\zeta,\vec{t},\vec{t}^{\prime},\vec{\mu})},
     \partial\Sigma_{(\zeta,\vec{t},\vec{t}^{\prime},\vec{\mu})};
     h_{\approxi, S(\zeta,\vec{t},\vec{t}^{\prime},
                                        \vec{\mu},\vec{\lambda})}
       ^{\ast}T_{\ast}W[k]_{\vec{\lambda}},
    (h_{\approxi,S(\zeta,\vec{t},\vec{t}^{\prime},
                                   \vec{\mu},\vec{\lambda})}
       |_{\partial\Sigma_{(\zeta,\vec{t},\vec{t}^{\prime},\vec{\mu})}})
      ^{\ast} T_{\ast}L
            \right)
  \end{array}
 $$
 } %end-small

 \noindent
 to deform the $\widetilde{V}_{\rho}$-family of
  approximate-$J$-stable $C^{\infty}$ maps
  $$
    h_{\approxi,(\zeta,\vec{t},\vec{t}^{\prime},
                    \vec{\mu},\vec{\lambda},\vec{a},\vec{b})}\; :\;
    ( \Sigma_{(\zeta,\vec{t},\vec{t}^{\prime},\vec{\mu})},
      \partial \Sigma_{(\zeta,\vec{t},\vec{t}^{\prime},\vec{\mu})} )\;
     \longrightarrow\;  ( W[k]_{\vec{\lambda}}, L )
  $$
 recursively to a $\widetilde{V}_{\rho}$-family of
  $(J,E_{\mbox{\tiny $\bullet$}})$-stable maps
 $$
  f_{(\zeta,\vec{t},\vec{t}^{\prime},
                   \vec{\mu},\vec{\lambda},\vec{a},\vec{b})}\; :\;
    ( \Sigma_{(\zeta,\vec{t},\vec{t}^{\prime},\vec{\mu})},
      \partial \Sigma_{(\zeta,\vec{t},\vec{t}^{\prime},\vec{\mu})} )\,
     \rightarrow\,  ( W[k]_{\vec{\lambda}}, L )\,.
 $$
Such construction is provided by [MD-S1: Sec.~3.3 and A.4] and
 its extensions to various situations in [F-O], [Liu(C)], and [L-R].
The discussion below follows these four works with mild necessary
 modifications to fit our overall presentation and notations.

Throughout the discussion, we assume that
 $\|\zeta\|$, $\|\vec{t}\|$, $\|\vec{t}^{\prime}\|$,
    $\|\vec{\mu}\|$, $\|\vec{a}-\vec{a}_f\|$, $\|\vec{b}\|$,
 and, hence, $\widetilde{V}_{\rho}$ are all sufficiently small
 so that statements in the construction hold.

\bigskip

\begin{flushleft}
{\bf The $\Aut(\rho)$-equivariant auxiliary bundle
     $E_{S(\pi_{\tinyDef(\Sigma)\times B[k]}(\widetilde{V}_{\rho}))}
      ^{\aux}$
     over
     $S(\pi_{\scriptsizeDef(\Sigma)\times B[k]}(\widetilde{V}_{\rho}))$.}
\end{flushleft}
Introduce first the following operations.
Let $P_{\bullet,\bullet^{\prime}}$ be the parallel transport
 from point $\bullet$ to point $\bullet^{\prime}$ along the minimal
 geodesic on a $W[k]_{\vec{\lambda}}$ for
 $\bullet,\bullet^{\prime}\in W[k]_{\vec{\lambda}}$
 of distance $<$ the injective radius of $W[k]_{\vec{\lambda}}$ and
$P_{\bullet,\bullet^{\prime}}^{\prime}$ be its $J$-linear part.
For
 $\eta \in
  L^p ( \Sigma_{(\zeta,\vec{t},\vec{t}^{\prime},\vec{\mu})};
        \Lambda^{0,1}\Sigma_{(\zeta,\vec{t},\vec{t}^{\prime},\vec{\mu})}
        \otimes_J
          h_{\approxi,S(\zeta,\vec{t},\vec{t}^{\prime},
                                      \vec{\mu},\vec{\lambda})}
         ^{\ast}T_{\ast}(W[k]_{\vec{\lambda}}) )$,
 define
  $$
   P_{(\zeta,\vec{t},\vec{t}^{\prime},\vec{\mu},
                     \vec{\lambda},\vec{a},\vec{b})}^{\prime} \eta\;
   \in\;
   L^p ( \Sigma_{(\zeta,\vec{t},\vec{t}^{\prime},\vec{\mu})};
         \Lambda^{0,1}\Sigma_{(0,\vec{t},\vec{t}^{\prime},\vec{\mu})}
         \otimes_J
          h_{\approxi,(\zeta,\vec{t},\vec{t}^{\prime},
                          \vec{\mu},\vec{\lambda},\vec{a},\vec{b})}
         ^{\ast}T_{\ast}(W[k]_{\vec{\lambda}}) )
  $$
  by
  $$
   (P_{(\zeta,\vec{t},\vec{t}^{\prime},\vec{\mu},
                  \vec{\lambda},\vec{a},\vec{b})}^{\prime} \eta) (x)\;
   =\;
      P_{h_{\tinyapproxi,S(\zeta,\vec{t},\vec{t}^{\prime},\vec{\mu},
                                         \vec{\lambda})}(x)\,,\,
         h_{\tinyapproxi,(\zeta,\vec{t},\vec{t}^{\prime},\vec{\mu},
            \vec{\lambda},\vec{a},\vec{b})} (x)
            }^{\prime}\, \eta(x)\,,
    \hspace{1em}
    x\in \Sigma_{(\zeta,\vec{t},\vec{t}^{\prime},\vec{\mu})}\,.
  $$
This is the $J$-linear parallel transport along the geodesic
 determined by
 $\xi_{(\zeta,\vec{t},\vec{t}^{\prime},
       \vec{\mu},\vec{\lambda},\vec{a},\vec{b})}(x)$
 in Notation 5.3.3.3.
Recall also the gluing maps for domain curves and targets spaces:
 (with $\varepsilon>0$ small and fixed, and
       $\|(\vec{t},\vec{t}^{\prime},\vec{\lambda})\|\ll\varepsilon$)
  $$
   \begin{array}{lllll}
    I_{(\zeta,\vec{t},\vec{t}^{\prime},\vec{\mu})}
     & : & \Sigma - \cup_{q:\node}\, N_{\sqrt{|t_q|}}(q)
     & \longrightarrow
     & \Sigma_{(\zeta,\vec{t},\vec{t}^{\prime},\vec{\mu})}\,,  \\[.6ex]
   I_{(\zeta,\vec{t},\vec{t}^{\prime},\vec{\mu}),\varepsilon}
    & :  & \Sigma- \cup_{q:\node}\, N_{|t_q|/\varepsilon}(q)
    & \longrightarrow
    & \Sigma_{(\zeta,\vec{t},\vec{t}^{\prime},\vec{\mu})}\,,   \\[.6ex]
   I_{\vec{\lambda}}
    & :  & Y_{[k]} - \cup_{\,i=0}^{\,k}\,N_{\sqrt{|\lambda_i|}}(D_i)
    & \longrightarrow   & W[k]_{\vec{\lambda}}\,,              \\[.6ex]
   I_{\vec{\lambda},\varepsilon}
    & :  & Y_{[k]} - \cup_{i=0}^k\,N_{|\lambda_i|/\varepsilon}(D_i)
    & \longrightarrow   & W[k]_{\vec{\lambda}}
   \end{array}
  $$
 and conjugation properties:
  $$
   \begin{array}{ccl}
    \alpha  \circ
     I_{(\zeta,\vec{t},\vec{t}^{\prime},\vec{\mu})} \circ \alpha^{-1}
      & = & I_{\alpha\cdot (\zeta,\vec{t}, \vec{t}^{\prime},\vec{\mu})}\,,
                                                             \\[.6ex]
    \alpha  \circ
     I_{(\zeta,\vec{t},\vec{t}^{\prime},\vec{\mu}),\varepsilon}  \circ
      \alpha^{-1}
      & =
      & I_{\alpha\cdot
           (\zeta,\vec{t}, \vec{t}^{\prime},\vec{\mu}),\varepsilon}\,,
                                                             \\[.6ex]
    \beta\circ I_{\vec{\lambda}}\circ\beta^{-1}
      & = & I_{\beta\cdot\vec{\lambda}}\,,                   \\[.6ex]
    \beta\circ I_{\vec{\lambda},\varepsilon}\circ\beta^{-1}
      & = & I_{\beta\cdot\vec{\lambda},\varepsilon}
  \end{array}
 $$
 for
  $\alpha\in\Aut(\rho)^{\domain}$ acting on ${\cal C}/\Def(\Sigma)$ and
  $\beta\in\Aut(\rho)^{\target}$ action on $W[k]/B[k]$.

Since $E_{\rho}$ is supported in a compact subset in the complement
 of all three types of nodes of $\Sigma$,
it can be canonically realized as a subspace in
 $L^p
  ( \Sigma_{(0,\vec{t},\vec{t}^{\prime},\vec{\mu})};
   \Lambda^{0,1}\Sigma_{(0,\vec{t},\vec{t}^{\prime},\vec{\mu})}
    \otimes_J
     h_{\approxi,S(0,\vec{t},\vec{t}^{\prime},\vec{\mu},\vec{\lambda})}
                            ^{\ast}T_{\ast}(\\  W[k]_{\vec{\lambda}})  )$
 via the composition
 $I_{\vec{\lambda}\,\ast}  \circ
  I_{(0,\vec{t},\vec{t}^{\prime},\vec{\mu})}^{\;-1\;\ast}  \circ
  P_{(0,\vec{t},\vec{t}^{\prime},\vec{\mu},
                \vec{\lambda},\vec{a},\vec{b})}^{\prime}$
 on $E_{\rho}$.
Define
 $E_{S(0,\vec{t},\vec{t}^{\prime},\vec{\mu},\vec{\lambda})}^{\aux}$
 to be this subspace in
 $L^p
  ( \Sigma_{(0,\vec{t},\vec{t}^{\prime},\vec{\mu})};
     \Lambda^{0,1}\Sigma_{(0,\vec{t},\vec{t}^{\prime},\vec{\mu})}
      \otimes_J
      h_{\approxi,S(0,\vec{t},\vec{t}^{\prime},\vec{\mu},\vec{\lambda})}
                      ^{\ast}  T_{\ast}(W[k]_{\vec{\lambda}})  )$.
To extend the above along
 the $H_{\rho,\domain}^{(\deform,\Sigma)}$-factor,
note that
 $E_{S(0,\vec{t},\vec{t}^{\prime},\vec{\mu},\vec{\lambda})}^{\aux}$
 is canonically a subspace of
 $L^p
  ( \Sigma_{(\zeta,\vec{t},\vec{t}^{\prime},\vec{\mu})};
     \Omega^1_{\scriptsizeBbb C}
      \Sigma_{(\zeta,\vec{t},\vec{t}^{\prime},\vec{\mu})}
      \otimes_J
      h_{\approxi,S(\zeta,\vec{t},\vec{t}^{\prime},
                          \vec{\mu},\vec{\lambda})}^{\ast}
                              T_{\ast}(W[k]_{\vec{\lambda}})  )$
 via the composition
  $\Lambda^{0,1}\Sigma_{(0,\vec{t},\vec{t}^{\prime},\vec{\mu})}
    \hookrightarrow
    \Omega^1_{\scriptsizeBbb C}
     \Sigma_{(0,\vec{t},\vec{t}^{\prime},\vec{\mu})}
    \stackrel{\sim}{\rightarrow}
     \Omega^1_{\scriptsizeBbb C}
      \Sigma_{(\zeta,\vec{t},\vec{t}^{\prime},\vec{\mu})}$
 of the canonical inclusion and the fixed isomorphism
  from the fixed
   $\Sigma_{(0,\vec{t},\vec{t}^{\prime},\vec{\mu})}
    \simeq \Sigma_{(\zeta,\vec{t},\vec{t}^{\prime},\vec{\mu})}$.
The restriction to
 $E_{S(0,\vec{t},\vec{t}^{\prime},\vec{\mu},\vec{\lambda})}^{\aux}$
of the following projection map
 $$
  \begin{array}{l}
   P^{0,1}_{\zeta}\;:\;
   L^p
   ( \Sigma_{(\zeta,\vec{t},\vec{t}^{\prime},\vec{\mu})};
     \Omega^1_{\scriptsizeBbb C}
      \Sigma_{(\zeta,\vec{t},\vec{t}^{\prime},\vec{\mu})}
      \otimes_J
      h_{\approxi,S(\zeta,\vec{t},\vec{t}^{\prime},
                       \vec{\mu},\vec{\lambda})}^{\ast}
                          T_{\ast}(W[k]_{\vec{\lambda}}) )\\[2ex]
  \hspace{5em}\longrightarrow\;
  L^p
   ( \Sigma_{(\zeta,\vec{t},\vec{t}^{\prime},\vec{\mu})};
     \Lambda^{0,1}\Sigma_{(\zeta,\vec{t},\vec{t}^{\prime},\vec{\mu})}
      \otimes_J
      h_{\approxi,S(\zeta,\vec{t},\vec{t}^{\prime},
                   \vec{\mu},\vec{\lambda})}^{\ast}
                                 T_{\ast}(W[k]_{\vec{\lambda}}) )
  \end{array}
 $$
 induced by the projection map
 $\Omega^1_{\scriptsizeBbb C}
   \Sigma_{(\zeta,\vec{t},\vec{t}^{\prime},\vec{\mu})}
  \rightarrow
  \Lambda^{0,1}
   \Sigma_{(\zeta,\vec{t},\vec{t}^{\prime},\vec{\mu})}$
 is injective for $\|\zeta\|$ sufficiently small.
Define
 $E_{S(\zeta,\vec{t},\vec{t}^{\prime},\vec{\mu},\vec{\lambda})}^{\aux}$
 to be the image of
 $E_{S(0,\vec{t},\vec{t}^{\prime},\vec{\mu},\vec{\lambda})}^{\aux}$
 in
 $L^p (
     \Sigma_{(\zeta,\vec{t},\vec{t}^{\prime},\vec{\mu})};
      \Lambda^{0,1}\Sigma_{(\zeta,\vec{t},\vec{t}^{\prime},\vec{\mu})}
        \otimes_J
          h_{\approxi,S(\zeta,\vec{t},\vec{t}^{\prime},
                            \vec{\mu},\vec{\lambda})}^{\ast}\\
                                       T_{\ast}(W[k]_{\vec{\lambda}}) )$
 under this projection.
This gives a trivialized vector bundle
 $E_{S(\pi_{\tinyDef(\Sigma)\times B[k]}(\widetilde{V}_{\rho}))}^{\aux}$
 over $S(\pi_{\scriptsizeDef(\Sigma)\times B[k]}(\widetilde{V}_{\rho}))$.
One can further define
 $$
  E_{(\zeta,\vec{t},\vec{t}^{\prime},\vec{\mu},
                    \vec{\lambda},\vec{a},\vec{b})}^{\aux}\;
  :=\; \left\{\, P_{(\zeta,\vec{t},\vec{t}^{\prime},\vec{\mu},
                 \vec{\lambda},\vec{a},\vec{b})}^{\prime}\, \eta\;
  :\; \eta \in
      E_{S(\zeta,\vec{t},\vec{t}^{\prime},\vec{\mu},
                         \vec{\lambda})}^{\aux}\, \right\}
 $$
to extend
 $E_{S(\pi_{\tinyDef(\Sigma)\times B[k]}(\widetilde{V}_{\rho}))}^{\aux}$
 to a trivialized vector bundle $E_{\widetilde{V}_{\rho}}^{\aux}$
 over $\widetilde{V}_{\rho}$,
 with specified isomorphisms of fibers to $E_{\rho}$.
In particular,
 $E_{S(\pi_{\tinyDef(\Sigma)\times B[k]}(\widetilde{V}_{\rho}))}^{\aux}$
 extends to
 $E_{\Theta_{\rho}}
   := E_{\widetilde{V}_{\rho}}^{\aux}|_{\Theta_{\rho}}$
  over $\Theta_{\rho}$.
The various group-invariance and conjugation properties of
 the objects and maps used in the construction implies that
 these trivialized bundles are $\Aut(\rho)$-equivariant.

\bigskip

\noindent
{\bf Definition 5.3.4.1 [auxiliary obstruction bundle].} {\rm
 We will call the $\Aut(\rho)$-equivariant trivialized bundle
  $E_{S(\pi_{\tinyDef(\Sigma)\times B[k]}(\widetilde{V}_{\rho}))}^{\aux}$
  (resp.\ $E_{\Theta_{\rho}}^{\aux}$, $E_{\widetilde{V}_{\rho}}^{\aux}$)
  as constructed above
 the {\it auxiliary obstruction bundle} over
  $S(\pi_{\scriptsizeDef(\Sigma)\times B[k]}(\widetilde{V}_{\rho}))$
  (resp.\ $\Theta_{\rho}$, $\widetilde{V}_{\rho}$)
  induced by $E_{\rho}$ at $\rho$.
} % end-definition

\bigskip

\begin{flushleft}
{\bf
 $\pi_{\scriptsizeDef(\Sigma)\times B[k]}(\widetilde{V}_{\rho})$-family
 of right inverse $Q_{\mbox{\tiny $\bullet$}}$  of
 $\pi_{E_{\cdot}^{\tinyaux}}
   \circ D_{\mbox{\tiny $\bullet$}}\bar{\partial}_J$
 from approximate one.}
\end{flushleft}
Let $\Ker(\pi_{E_{\rho}}\!\circ D_{\!f}\bar{\partial}_J)^{\perp}$
 be the $L^2$-orthogonal complement of
 $\Ker(\pi_{E_{\rho}}\!\circ D_{\!f}\bar{\partial}_J)$ in
 $W^{1,p}( \Sigma, \partial\Sigma;\\
           f^{\ast}T_{\ast}Y_{[k]},
                   (f|_{\partial\Sigma})^{\ast}T_{\ast}L )$.
This space is $\Aut(\rho)$-invariant, as the metric on $\Sigma$
 and $W[k]$ are respectively $\Aut(\rho)^{\domain}$- and
 $\Aut(\rho)^{\target}$-invariant.
Then
 $$
  \pi_{E_{\rho}}\circ D_{\!f}\bar{\partial}_J\; :\;
   \Ker(\pi_{E_{\rho}}\!\circ D_{\!f}\bar{\partial}_J)^{\perp}\;
    \longrightarrow\;
    \left.
      L^p(\Sigma; \Lambda^{0,1}\Sigma\otimes_J f^{\ast}T_{\ast}Y_{[k]})
    \right/\!\raisebox{-.6ex}{$E_{\rho}$}
 $$
 is an isomorphism and its inverse
 $$
  Q_{\rho}\;:\;
  \left.
   L^p(\Sigma; \Lambda^{0,1}\Sigma\otimes_J f^{\ast}T_{\ast}Y_{[k]})
  \right/\!\raisebox{-.6ex}{$E_{\rho}$}\;
   \longrightarrow\;
   \Ker(\pi_{E_{\rho}}\!\circ D_{\!f}\bar{\partial}_J)^{\perp}\;
 $$
 is a bounded operator.
This defines $Q_{\rho}$ as a right inverse of
 $$
  \pi_{E_{\rho}}\circ D_{\!f}\bar{\partial}_J\; :\;
  W^{1,p}( \Sigma, \partial\Sigma; f^{\ast}T_{\ast}Y_{[k]},
                   (f|_{\partial\Sigma})^{\ast}T_{\ast}L_{[k]} )\;
    \longrightarrow\;
    \left.
      L^p(\Sigma; \Lambda^{0,1}\Sigma\otimes_J f^{\ast}T_{\ast}Y_{[k]})
    \right/\!\raisebox{-.6ex}{$E_{\rho}$}\,.
 $$
We now proceed to construct first a suitable
 $\pi_{\scriptsizeDef(\Sigma)\times B[k]}(\widetilde{V}_{\rho})$-family
 of approximate right inverse
 $Q^{\prime}_{S(\zeta,\vec{t},\vec{t}^{\prime},
                              \vec{\mu},\vec{\lambda})}$
 of
 $$
  \begin{array}{l}
   \pi_{E_{S(\zeta,\vec{t},\vec{t}^{\prime},
                \vec{\mu},\vec{\lambda})}^{\tinyaux}}
                \circ D_{h_{\tinyapproxi,
                         S(\zeta,\vec{t},\vec{t}^{\prime},
                           \vec{\mu},\vec{\lambda})}}
                                 \bar{\partial}_J\;  :  \\[1.6ex]
   \hspace{1em}
   W^{1,p}(
      \Sigma_{(\zeta,\vec{t},\vec{t}^{\prime},\vec{\mu})},
      \partial \Sigma_{(\zeta,\vec{t},\vec{t}^{\prime},\vec{\mu})};
      h_{\approxi,S(\zeta,\vec{t},\vec{t}^{\prime},
                          \vec{\mu},\vec{\lambda})}^{\ast}
      T_{\ast}(W[k]_{\vec{\lambda}}),
       ( h_{\approxi,S(\zeta,\vec{t},\vec{t}^{\prime},
                             \vec{\mu},\vec{\lambda})}
        |_{\partial\Sigma_{(\zeta,\vec{t},\vec{t}^{\prime},
                                          \vec{\mu})}} )^{\ast}
      T_{\ast}L )                                       \\[1.6ex]
   \hspace{3em}
   \longrightarrow\;
   \left.
    L^p(
     \Sigma_{(\zeta,\vec{t},\vec{t}^{\prime},\vec{\mu})};
     \Lambda^{0,1} \Sigma_{(\zeta,\vec{t},\vec{t}^{\prime},\vec{\mu})}
      \otimes_J  h_{\approxi,S(\zeta,\vec{t},\vec{t}^{\prime},
                                     \vec{\mu},\vec{\lambda})}^{\ast}
                  T_{\ast}(W[k]_{\vec{\lambda}}) ) \right/
       \raisebox{-.6ex}{$E_{S(\zeta,\vec{t},\vec{t}^{\prime},
                                    \vec{\mu},\vec{\lambda})}^{\aux}$}
  \end{array}
 $$
 by passing to $Q_{\rho}$ at $\rho$.

The combination of $I_{(0,\vec{t},\vec{t}^{\prime},\vec{\mu})}$ and
 $I_{\vec{\lambda}}$ on domains and targets induces a map
 $$
  \begin{array}{c}
   I_{(0,\vec{t},\vec{t}^{\prime},\vec{\mu},\vec{\lambda})}^{\ast}\; :
    L^p(
      \Sigma_{(0,\vec{t},\vec{t}^{\prime},\vec{\mu})};
      \Lambda^{0,1} \Sigma_{(0,\vec{t},\vec{t}^{\prime},\vec{\mu})}
       \otimes_J  h_{\approxi,S(0,\vec{t},\vec{t}^{\prime},
                                  \vec{\mu},\vec{\lambda})}^{\ast}
                          T_{\ast}(W[k]_{\vec{\lambda}}) )   \\[1.6ex]
   \longrightarrow\;
   L^p(\Sigma; \Lambda^{0,1}\Sigma\otimes_J f^{\ast}T_{\ast}Y_{[k]})
  \end{array}
 $$
 by first using $I_{\vec{\lambda}}$ to turn an
  $$
   \eta\; \in\;
     L^p(
       \Sigma_{(0,\vec{t},\vec{t}^{\prime},\vec{\mu})};
       \Lambda^{0,1} \Sigma_{(0,\vec{t},\vec{t}^{\prime},\vec{\mu})}
        \otimes_J  h_{\approxi,S(0,\vec{t},\vec{t}^{\prime},
                                   \vec{\mu},\vec{\lambda})}^{\ast}
                           T_{\ast}(W[k]_{\vec{\lambda}}) )
  $$
  to an element
  $$
   \eta^{\prime}\; =\; I_{\vec{\lambda}}^{\ast}\,\eta\; \in\;
      L^p(
        \Sigma_{(0,\vec{t},\vec{t}^{\prime},\vec{\mu})};
        \Lambda^{0,1} \Sigma_{(0,\vec{t},\vec{t}^{\prime},\vec{\mu})}
         \otimes_J  (I_{\vec{\lambda}}^{-1} \circ
                     h_{\approxi,S(0,\vec{t},\vec{t}^{\prime},
                                   \vec{\mu},\vec{\lambda})} )^{\ast}
                                             T_{\ast}Y_{[k]} )
  $$
 and then using parallel transport
  $P_{( I_{\vec{\lambda}}^{-1} \circ
        h_{\tinyapproxi,
         S(0,\vec{t},\vec{t}^{\prime},\vec{\mu},\vec{\lambda})} )
          (I_{(0,\vec{t},\vec{t}^{\prime},\vec{\mu})}(x))\,,\,
        f(x)}$ on $Y_{[k]}$
  for\\  $x \in I_{(0,\vec{t},\vec{t}^{\prime},\vec{\mu})}^{-1}
              (\Sigma_{(0,\vec{t},\vec{t}^{\prime},\vec{\mu})})
             \subset \Sigma$
  to move $\eta^{\prime}$ to an element
  $$
   \eta^{\prime\prime}\;
   =\; P_{\mbox{\tiny $\bullet$},\mbox{\tiny $\bullet$}}(\eta^{\prime})\;
   =:\;
     I_{(0,\vec{t},\vec{t}^{\prime},\vec{\mu},\vec{\lambda})}^{\ast}
      (\eta)\; \in\;
     L^p(\Sigma; \Lambda^{0,1}\Sigma\otimes_J f^{\ast}T_{\ast}Y_{[k]})\,.
  $$
The composition of
$$
 \begin{array}{ll}
   & L^p(
     \Sigma_{(\zeta,\vec{t},\vec{t}^{\prime},\vec{\mu})};
     \Lambda^{0,1} \Sigma_{(\zeta,\vec{t},\vec{t}^{\prime},\vec{\mu})}
      \otimes_J  h_{\approxi,S(\zeta,\vec{t},\vec{t}^{\prime},
                                     \vec{\mu},\vec{\lambda})}^{\ast}
                  T_{\ast}(W[k]_{\vec{\lambda}}) )             \\
  \stackrel{P^{0,1}_0\hspace{2ex}}{\longrightarrow}
   & L^p(
     \Sigma_{(0,\vec{t},\vec{t}^{\prime},\vec{\mu})};
     \Lambda^{0,1} \Sigma_{(0,\vec{t},\vec{t}^{\prime},\vec{\mu})}
      \otimes_J  h_{\approxi,S(0,\vec{t},\vec{t}^{\prime},
                                 \vec{\mu},\vec{\lambda})}^{\ast}
                  T_{\ast}(W[k]_{\vec{\lambda}}) )             \\
  \stackrel{I_{(0, \vec{t},\vec{t}^{\prime},\vec{\mu},
                   \vec{\lambda})}^{\ast}}{\longrightarrow\hspace{1.2em}}
   & L^p(\Sigma;\Lambda^{0,1}\Sigma\otimes_J f^{\ast}T_{\ast}Y_{[k]})
 \end{array}
$$
gives the map
$$
 \begin{array}{c}
  \IP_{S(\zeta,\vec{t},\vec{t}^{\prime},
                       \vec{\mu},\vec{\lambda})}\; :\;
  L^p(
    \Sigma_{(\zeta,\vec{t},\vec{t}^{\prime},\vec{\mu})};
    \Lambda^{0,1} \Sigma_{(\zeta,\vec{t},\vec{t}^{\prime},\vec{\mu})}
     \otimes_J  h_{\approxi,S(\zeta,\vec{t},\vec{t}^{\prime},
                                    \vec{\mu},\vec{\lambda})}^{\ast}
                 T_{\ast}(W[k]_{\vec{\lambda}}) )           \\[2ex]
  \longrightarrow\;
    L^p(\Sigma;\Lambda^{0,1}\Sigma\otimes_J f^{\ast}T_{\ast}Y_{[k]})\,.
 \end{array}
$$

Fix a $C^{\infty}$ rotation-invariant cutoff function
 $\beta_{\delta}:{\Bbb C}\rightarrow {\Bbb R}$
 as in [MD-S1: Lemma A.1.1] with the following properties:
 (mixed with the presentation of [F-O])
 $$
  \beta_{\delta}(z)\;=\; \left\{
   \begin{array}{ll}
    1   & \mbox{if $\,|z|\; \le\; \delta\; (<1)$} \\[.6ex]
    0   & \mbox{if $\,|z| \ge 1-o$} \\
        & \hspace{1em}\mbox{for some $0< o \ll 1-\delta$}
   \end{array}           \right.
  \hspace{1em}\mbox{and}\hspace{1em}
  \int_{|z|\le 1} |\nabla\beta_{\delta}(z)|^2\;
     \le\; \frac{4\,\pi}{|\log\delta|}\,;
 $$
and recall that $\varepsilon>0$ is small and fixed, and
 $\|(\vec{t},\vec{t}^{\prime},\vec{\mu},\vec{\lambda})\|\ll\varepsilon$.
Define the map.

{\small
 $$
  \begin{array}{l}
  \Glue_{S(\zeta,\vec{t},\vec{t}^{\prime},
                         \vec{\mu},\vec{\lambda})}\,:
    W^{1,p}( \Sigma, \partial\Sigma; f^{\ast}T_{\ast}Y_{[k]},
                       (f|_{\partial\Sigma})^{\ast}T_{\ast}L )\;
     \longrightarrow \\[2ex]
    \hspace{1em}
    W^{1,p}(
      \Sigma_{(\zeta,\vec{t},\vec{t}^{\prime},\vec{\mu})},
      \partial \Sigma_{(\zeta,\vec{t},\vec{t}^{\prime},\vec{\mu})};
      h_{\approxi,S(\zeta,\vec{t},\vec{t}^{\prime},
                                  \vec{\mu},\vec{\lambda})}^{\ast}
      T_{\ast}(W[k]_{\vec{\lambda}}),
       ( h_{\approxi,S(\zeta,\vec{t},\vec{t}^{\prime},
                    \vec{\mu},\vec{\lambda})}
         |_{\partial\Sigma_{(\zeta,\vec{t},\vec{t}^{\prime},\vec{\mu})}} )
       ^{\ast} T_{\ast}L )
  \end{array}
 $$
} % end-small

\noindent
by gluing locally defined bundle-valued fields to a continuous field
 as follows.
Let\\
 $\xi \in
  W^{1,p}( \Sigma, \partial\Sigma; f^{\ast}T_{\ast}Y_{[k]},
                       (f|_{\partial\Sigma})^{\ast}T_{\ast}L )$  and
recall the gluing construction of the maps
 $h_{\approxi,\,\mbox{\tiny $\bullet$}\,}$ in Sec.~5.3.3.

\begin{itemize}
 \item[(o)]
  For $x$ in
   $\Sigma_{(\zeta,\vec{t},\vec{t}^{\prime},\vec{\mu})}
    -\Neck
      _{\varepsilon\,,\,(\zeta,\vec{t},\vec{t}^{\prime},\vec{\mu})}$,
  define

 \vspace{-1ex}
 \item[]
 {\small
  $$
   (\Glue_{S(\zeta,\vec{t},\vec{t}^{\prime},
                   \vec{\mu},\vec{\lambda})}(\xi))(x)\;
   =\; \left(\rule{0ex}{2ex}\right.
         {I_{\vec{\lambda}}}_{\ast} \circ
         P_{ f(I_{(\zeta,\vec{t},\vec{t}^{\prime},\vec{\mu}) }^{-1}(x))
             \,,\,
             I_{\vec{\lambda}}^{-1}(h_{\tinyapproxi,
                S(\zeta,\vec{t},\vec{t}^{\prime},
                        \vec{\mu},\vec{\lambda})}(x)) }
        \left.\rule{0ex}{2ex}\right)\,
       \xi( {I_{(\zeta,\vec{t},\vec{t}^{\prime},\vec{\mu})}}^{-1} (x) )\,.
  $$ } % end-small

 \item[$(a)$]
  For $x$ in the
   {\it annulus $A_t$ from a smoothed ordinary interior node $q$},
  let $x=(z,\frac{t}{z})\in{\Bbb C}^2$ in the local model
  in Sec.~5.3.3,
  where
   $t$ is an entry of $\vec{t}$ involved and
   $\xi=\xi_1\cup \xi_2$ on the two irreducible components
    the neighborhood of $q=(0,0)$ in
    $\{(z_1,z_2):z_1z_2=0, |z_1|<\varepsilon, |z_2|<\varepsilon\}
     \subset \Sigma$,
    with $\xi_1=\xi_1(z_1)$ and $\xi_2=\xi_2(z_2)$.
  Define

  \vspace{-1em}
  \item[]
  {\small
  \begin{eqnarray*}
   \lefteqn{
   (\Glue_{S(\zeta,\vec{t},\vec{t}^{\prime},
                   \vec{\mu},\vec{\lambda})} (\xi))(x) } \\[1ex]
   && =\; \left\{
   \begin{array}{l}
    \left(\rule{0ex}{2ex}\right.
      {I_{\vec{\lambda}}}_{\ast} \circ
      P_{f_1(I_{(\zeta,\vec{t},\vec{t}^{\prime},\vec{\mu})}^{-1}(x))\,,\,
         I_{\vec{\lambda}}^{-1}(
            h_{\tinyapproxi,
            S(\zeta,\vec{t},\vec{t}^{\prime},
                    \vec{\mu},\vec{\lambda})}(x) )}
      \left.\rule{0ex}{2ex}\right)\,
      \xi_1(z)                        \\[1.6ex]
    \hspace{1em}
    +\, \left( 1-\beta_{\delta}(\frac{|t|^{1/2}}{z}) \right)
        \left(
         \left(\rule{0ex}{2ex}\right.
           {I_{\vec{\lambda}}}_{\ast} \circ
           P_{f_2(
               I_{(\zeta,\vec{t},\vec{t}^{\prime},\vec{\mu}) }^{-1}(x))
              \,,\,
              I_{\vec{\lambda}}^{-1}(
                h_{\tinyapproxi,
                  S(\zeta,\vec{t},\vec{t}^{\prime},
                          \vec{\mu},\vec{\lambda})}(x))}
          \left.\rule{0ex}{2ex}\right)\,
             \xi_2(\frac{t}{z})\,    \right.    \\[1ex]
        \hspace{12em}\left.\rule{0ex}{3ex}
         -\, \left(\rule{0ex}{2ex}\right.
               {I_{\vec{\lambda}}}_{\ast} \circ
               P_{f(q)\,,\,
                  I_{\vec{\lambda}}^{-1}(
                   h_{\tinyapproxi,
                   S(\zeta,\vec{t},\vec{t}^{\prime},
                           \vec{\mu},\vec{\lambda})}(x) )}
                \left.\rule{0ex}{2ex}\right)\, \xi(q)
        \right)                                        \\[2.6ex]
    \hspace{25em}
     \mbox{for $\;|t|^{1/2}\le |z|\le \varepsilon$}\,,      \\[2ex]
    \left(\rule{0ex}{2ex}\right.
      {I_{\vec{\lambda}}}_{\ast} \circ
      P_{f_2(
         I_{(\zeta,\vec{t},\vec{t}^{\prime},\vec{\mu}) }^{-1}(x) )
         \,,\,
         I_{\vec{\lambda}}^{-1}(
          h_{\tinyapproxi,
             S(\zeta,\vec{t},\vec{t}^{\prime},
                     \vec{\mu},\vec{\lambda})}(x) )}
       \left.\rule{0ex}{2ex}\right)\,
      \xi_2(\frac{t}{z})                            \\[1.6ex]
     \hspace{1em}
    +\, \left( 1-\beta_{\delta}(\frac{z}{|t|^{1/2}}) \right)
        \left(
          \left.\rule{0ex}{2ex}\right(
           {I_{\vec{\lambda}}}_{\ast} \circ
           P_{f_1(
                I_{(\zeta,\vec{t},\vec{t}^{\prime},\vec{\mu}) }^{-1}(x))
              \,,\,
             I_{\vec{\lambda}}^{-1}(
              h_{\tinyapproxi,
                S(\zeta,\vec{t},\vec{t}^{\prime},
                        \vec{\mu},\vec{\lambda})}(x) )}
            \left.\rule{0ex}{2ex}\right)\,
          \xi_1(z)\,                       \right.  \\[1ex]
         \hspace{12em}\left.\rule{0ex}{3ex}
         -\, \left.\rule{0ex}{2ex}\right(
               {I_{\vec{\lambda}}}_{\ast} \circ
               P_{f(q)\,,\,
                  I_{\vec{\lambda}}^{-1}(
                   h_{\tinyapproxi,
                    S(\zeta,\vec{t},\vec{t}^{\prime},
                            \vec{\mu},\vec{\lambda})}(x) )}
              \left.\rule{0ex}{2ex}\right)\,
             \xi(q)
        \right)                                       \\[2.6ex]
    \hspace{25em}
     \mbox{for $\;|t|^{1/2}\le |t/z|\le \varepsilon$}\,.
   \end{array}
          \right.
  \end{eqnarray*}
  } % end-small
\end{itemize}

\begin{itemize}
 \item[$(b)$]
  For $x$ in the {\it annulus $A^{\prime}_{t^{\prime}}$, $t^{\prime}>0$,
   from smoothing a type E boundary node $q$},
  let $x=(z,\frac{t^{\prime}}{z})$ in the local model in Sec.~5.3.3,
   where $t^{\prime}$ is an entry of $\vec{t}^{\prime}$ involved, and
  define

  \vspace{-1em}
  \item[]
  {\small
  \begin{eqnarray*}
   \lefteqn{
   (\Glue_{S(\zeta,\vec{t},\vec{t}^{\prime},
                   \vec{\mu},\vec{\lambda})} (\xi))(x) }    \\[1ex]
   && =\;
    \left.\rule{0ex}{2ex}\right(
      {I_{\vec{\lambda}}}_{\ast} \circ
      P_{f(I_{(\zeta,\vec{t},\vec{t}^{\prime},\vec{\mu}) }^{-1}(x))
         \,,\,
         I_{\vec{\lambda}}^{-1}(
          h_{\tinyapproxi,
            S(\zeta,\vec{t},\vec{t}^{\prime},
                    \vec{\mu},\vec{\lambda})}(x) )}
           \left.\rule{0ex}{2ex}\right)\,
    \xi(z)                                       \\[1.6ex]
   && \hspace{2em}
    -\, \left( 1-\beta_{\delta}(\frac{|t^{\prime}|^{1/2}}{z}) \right)\,
         \left.\rule{0ex}{2ex}\right(
           {I_{\vec{\lambda}}}_{\ast} \circ
           P_{f(q)\,,\,
               I_{\vec{\lambda}}^{-1}(
                h_{\tinyapproxi,
                   S(\zeta,\vec{t},\vec{t}^{\prime},
                           \vec{\mu},\vec{\lambda})}(x) )}
               \left.\rule{0ex}{2ex}\right)\,
         \xi(q)                                  \\[1ex]
   && \hspace{26em}
     \mbox{for $\;|t^{\prime}|^{1/2}\le |z|\le \varepsilon$}\,.
  \end{eqnarray*}
  } % end-small
\end{itemize}

\begin{itemize}
 \item[]
  For $x$ in the {\it band $A^{\prime}_{t^{\prime}}$, $t^{\prime}>0$,
   from smoothing a type H boundary node $q$},
  let
   $x=(z,\frac{t^{\prime}}{z})$ in the local model in Sec.~5.3.3,
    where $t^{\prime}$ is an entry of $\vec{t}^{\prime}$ involved,
   $\xi=\xi_1\cup \xi_2$ on the two irreducible components
    the neighborhood of $q=(0,0)$ in
    $\{(z_1,z_2):z_1z_2=0, |z_1|<\varepsilon, |z_2|<\varepsilon\}/
      (z_1,z_2)\sim(\overline{z_1},\overline{z_2})
     \subset \Sigma$,
    with $\xi_1=\xi_1(z_1)$ and $\xi_2=\xi_2(z_2)$, and
  define

  \vspace{-1em}
  \item[]
  {\small
  \begin{eqnarray*}
   \lefteqn{
   (\Glue_{S(\zeta,\vec{t},\vec{t}^{\prime},
                   \vec{\mu},\vec{\lambda})} (\xi))(x) }   \\[1ex]
   && =\; \left\{
   \begin{array}{l}
    \left(\rule{0ex}{2ex}\right.
      {I_{\vec{\lambda}}}_{\ast} \circ
      P_{f_1(I_{(\zeta,\vec{t},\vec{t}^{\prime},\vec{\mu}) }^{-1}(x))
          \,,\,
         I_{\vec{\lambda}}^{-1}(
          h_{\tinyapproxi,
            S(\zeta,\vec{t},\vec{t}^{\prime},
                            \vec{\mu},\vec{\lambda})}(x))}
           \left.\rule{0ex}{2ex}\right)\,
     \xi_1(z)                                     \\[1.6ex]
    \hspace{1em}
    +\, \left( 1-\beta_{\delta}(\frac{|t^{\prime}|^{1/2}}{z}) \right)
        \left(
         \left(\rule{0ex}{2ex}\right.
           {I_{\vec{\lambda}}}_{\ast} \circ
           P_{f_2(I_{(\zeta,\vec{t},\vec{t}^{\prime},\vec{\mu}) }^{-1}(x))
                \,,\,
              I_{\vec{\lambda}}^{-1}(
               h_{\tinyapproxi,
                S(\zeta,\vec{t},\vec{t}^{\prime},
                                \vec{\mu},\vec{\lambda})}(x) )}
               \left.\rule{0ex}{2ex}\right)\,
             \xi_2(\frac{t^{\prime}}{z})\,    \right.    \\[1ex]
        \hspace{12em}\left.\rule{0ex}{3ex}
         -\, \left(\rule{0ex}{2ex}\right.
               {I_{\vec{\lambda}}}_{\ast} \circ
               P_{f(q)\,,\,
                  I_{\vec{\lambda}}^{-1}(
                    h_{\tinyapproxi,
                     S(\zeta,\vec{t},\vec{t}^{\prime},
                             \vec{\mu},\vec{\lambda})}(x) )}
                 \left.\rule{0ex}{2ex}\right)\,
              \xi(q)
        \right)                                         \\[2.6ex]
    \hspace{25em}
     \mbox{for $\;|t^{\prime}|^{1/2}\le |z|\le \varepsilon$}\,,  \\[2ex]
    \left(\rule{0ex}{2ex}\right.
      {I_{\vec{\lambda}}}_{\ast} \circ
      P_{f_2(I_{(\zeta,\vec{t},\vec{t}^{\prime},\vec{\mu}) }^{-1}(x))
           \,,\,
         I_{\vec{\lambda}}^{-1}(
           h_{\tinyapproxi,
            S(\zeta,\vec{t},\vec{t}^{\prime},
                            \vec{\mu},\vec{\lambda})}(x) )}
        \left.\rule{0ex}{2ex}\right)\,
      \xi_2(\frac{t^{\prime}}{z})                  \\[1.6ex]
     \hspace{1em}
    +\, \left( 1-\beta_{\delta}(\frac{z}{|t^{\prime}|^{1/2}}) \right)
        \left(
         \left(\rule{0ex}{2ex}\right.
           {I_{\vec{\lambda}}}_{\ast} \circ
           P_{f_1(
               I_{(\zeta,\vec{t},\vec{t}^{\prime},\vec{\mu}) }^{-1}(z))
               \,,\,
               I_{\vec{\lambda}}^{-1}(
                h_{\tinyapproxi,
                 S(\zeta,\vec{t},\vec{t}^{\prime},
                         \vec{\mu},\vec{\lambda})}(x) )}
             \left.\rule{0ex}{2ex}\right)\,
          \xi_1(z)\, \right.                 \\[1ex]
         \hspace{12em}\left.\rule{0ex}{3ex}
         -\, \left(\rule{0ex}{2ex}\right.
               {I_{\vec{\lambda}}}_{\ast} \circ
               P_{f(q)\,,\,
                  I_{\vec{\lambda}}^{-1}(
                   h_{\tinyapproxi,
                    S(\zeta,\vec{t},\vec{t}^{\prime},
                            \vec{\mu},\vec{\lambda})}(x) )}
                \left.\rule{0ex}{2ex}\right)\,
               \xi(q)
        \right)                                      \\[2.6ex]
    \hspace{25em}
     \mbox{for $\;|t^{\prime}|^{1/2}\le |t/z|\le \varepsilon$}\,.
   \end{array}
          \right.
  \end{eqnarray*}
  } % end-small
\end{itemize}

\begin{itemize}
 \item[$(c)$]
  For $x$ in the {\it annulus $A_{\mu}$ from a smoothed distinguished
   interior node $q$},
  let
   $x=(z,\frac{\mu}{z})$ in the local model in Sec.~5.3.3,
    where $\mu$ is an entry of $\vec{\mu}$ involved.
  Suppose that $f(q)\in D_i\subset Y_{[k],\sing}$;
  then denote the restriction of $I_{\vec{\lambda},\varepsilon}$
   to $\Delta_i$ (resp.\ $\Delta_{i+1}$)
   by $I_{\vec{\lambda},\varepsilon}^{f(q),1}$
    (resp.\ $I_{\vec{\lambda},\varepsilon}^{f(q),2}$).
 Define
\end{itemize}

\vspace{-1ex}
{\small
\begin{eqnarray*}
 \lefteqn{
 (\Glue_{S(\zeta,\vec{t},\vec{t}^{\prime},
                 \vec{\mu},\vec{\lambda})}(\xi))(x) } \\[1.6ex]
 && =\; \left\{
 \begin{array}{l}
  \left(\rule{0ex}{2ex}\right.
    {I_{\vec{\lambda},\varepsilon}^{f(q),1}}_{\ast} \circ
      P_{f_1(I_{(\vec{\zeta,\varepsilon},
                      \vec{t},\vec{t}^{\prime},\vec{\mu})}^{-1}(x))
         \,,\,
         (I_{\vec{\lambda},\varepsilon}^{f(q),1})^{-1}(
            h_{\tinyapproxi,
               S(\zeta,\vec{t},\vec{t}^{\prime},
                       \vec{\mu},\vec{\lambda})}(x) )}
    \left.\rule{0ex}{2ex}\right)\,
  \xi_1(z)          \\[1ex]
  \hspace{1em}
  +\, \left( 1-\beta_{\delta}(\frac{|\mu|^{1/2}}{z}) \right)
      \left(\rule{0ex}{3ex}\right.
       -\,\frac{1}{2}\,
           \left(\rule{0ex}{2ex}\right.
            {I_{\vec{\lambda},\varepsilon}^{f(q),1}}_{\ast} \circ
            P_{f(q)\,,\,
               (I_{\vec{\lambda,\varepsilon}}^{f(q),1})^{-1}
                 (h_{\tinyapproxi,
                     S(\zeta,\vec{t},\vec{t}^{\prime},
                             \vec{\mu},\vec{\lambda})}(x)) }
            \left.\rule{0ex}{2ex}\right)\,
          \xi(q)                            \\[1ex]
      \hspace{6em}
       +\,
          \left(\rule{0ex}{2ex}\right.
           {I_{\vec{\lambda},\varepsilon}^{f(q),2}}_{\ast}\circ
           P_{f_2(I_{(\vec{\zeta},\vec{t},\vec{t}^{\prime},
                                          \vec{\mu})}^{-1}(x))
                \,,\,
              (I_{\vec{\lambda},\varepsilon}^{f(q),2})^{-1}(
                h_{\tinyapproxi,
                   S(\zeta,\vec{t},\vec{t}^{\prime},
                           \vec{\mu},\vec{\lambda})}(x) )}
              \left.\rule{0ex}{2ex}\right)\,
          \xi_2(\frac{\mu}{z})     \\[1ex]
      \hspace{10em}
       -\,\frac{1}{2}\,
           \left(\rule{0ex}{2ex}\right.
            {I_{\vec{\lambda},\varepsilon}^{f(q),2}}_{\ast} \circ
            P^2_{f(q)\,,\,
                 (I^2_{\vec{\lambda},\varepsilon})
                                  ^{\raisebox{.2ex}{\tiny $-1$}}
            (h_{\tinyapproxi,
                S(\zeta,\vec{t},\vec{t}^{\prime},
                        \vec{\mu},\vec{\lambda})}(x)) }\,
             \left.\rule{0ex}{2ex}\right)\,
            \xi(q)
      \left.\rule{0ex}{3ex}\right)                     \\[3ex]
  \hspace{25em}
   \mbox{for $\;|\mu|^{1/2}\le |z|\le \varepsilon_1$}\,,     \\[2ex]
  \left(\rule{0ex}{2ex}\right.
    {I_{\vec{\lambda},\varepsilon}^{f(q),2}}_{\ast} \circ
      P_{f_2(I_{(\vec{\zeta},\vec{t},\vec{t}^{\prime},\vec{\mu})}^{-1}(x))
         \,,\,
         (I_{\vec{\lambda},\varepsilon}^{f(q),2})^{-1}(
            h_{\tinyapproxi,
               S(\zeta,\vec{t},\vec{t}^{\prime},
                       \vec{\mu},\vec{\lambda})}(x) )}
    \left.\rule{0ex}{2ex}\right)\,
   \xi_2(\frac{\mu}{z})                \\[1ex]
  \hspace{1em}
  +\, \left( 1-\beta_{\delta}(\frac{z}{|\mu|^{1/2}}) \right)
      \left(\rule{0ex}{3ex}\right.
       -\,\frac{1}{2}\,
           \left(\rule{0ex}{2ex}\right.
            {I_{\vec{\lambda},\varepsilon}^{f(q),2}}_{\ast} \circ
            P_{f(q)\,,\,
               (I_{\vec{\lambda},\varepsilon}^{f(q),2})^{-1}
                 (h_{\tinyapproxi,
                     S(\zeta,\vec{t},\vec{t}^{\prime},
                             \vec{\mu},\vec{\lambda})}(x)) }
              \left.\rule{0ex}{2ex}\right)\,
          \xi(q)                            \\[1ex]
      \hspace{6em}
       +\,
          \left(\rule{0ex}{2ex}\right.
          {I_{\vec{\lambda},\varepsilon}^{f(q),1}}_{\ast} \circ
          P_{f_1(I_{(\vec{\zeta},\vec{t},\vec{t}^{\prime},
                                         \vec{\mu})}^{-1}(x))
               \,,\,
             (I_{\vec{\lambda},\varepsilon}^{f(q),1})^{-1}(
                h_{\tinyapproxi,
                   S(\zeta,\vec{t},\vec{t}^{\prime},
                           \vec{\mu},\vec{\lambda})}(x) )}
         \left.\rule{0ex}{2ex}\right)\,
        \xi_1(z)   \\[1ex]
     \hspace{10em}
       -\,\frac{1}{2}\,
           \left(\rule{0ex}{2ex}\right.
            {I_{\vec{\lambda},\varepsilon}^{f(q),1}}_{\ast} \circ
            P_{f(q)\,,\,
                 (I_{\vec{\lambda},\varepsilon}^{f(q),1})^{-1}
                   (h_{\tinyapproxi,
                       S(\zeta,\vec{t},\vec{t}^{\prime},
                               \vec{\mu},\vec{\lambda})}(x)) }
               \left.\rule{0ex}{2ex}\right)
          \xi(q)
      \left.\rule{0ex}{3ex}\right)                                              \\[3ex]
  \hspace{25em}
   \mbox{for $\;|\mu|^{1/2}\le |\mu/z|\le \varepsilon_1$}\,.
 \end{array}
        \right.
\end{eqnarray*}
} % end-small

Then the composition

{\small
$$
 \begin{array}{l}
  Q^{\prime}_{S(\zeta,\vec{t},\vec{t}^{\prime},\vec{\mu},\vec{\lambda})}\;
   :=\; \Glue_{S(\zeta,\vec{t},\vec{t}^{\prime},\vec{\mu},\vec{\lambda})}\,
        \circ\, Q_{\rho}\, \circ\,
        \IP_{S(\zeta,\vec{t},\vec{t}^{\prime},\vec{\mu},\vec{\lambda})}
                                               \\[2ex]
  \hspace{1em}:\;
   \left.
    L^p(
     \Sigma_{(\zeta,\vec{t},\vec{t}^{\prime},\vec{\mu})};
     \Lambda^{0,1} \Sigma_{(\zeta,\vec{t},\vec{t}^{\prime},\vec{\mu})}
      \otimes_J  h_{\approxi,S(\zeta,\vec{t},\vec{t}^{\prime},
                                      \vec{\mu},\vec{\lambda})}^{\ast}
                  T_{\ast}(W[k]_{\vec{\lambda}}) ) \right/
       \raisebox{-.6ex}{$E_{S(\zeta,\vec{t},\vec{t}^{\prime},
                                      \vec{\mu},\vec{\lambda})}^{\aux}$}
       \hspace{1ex} \longrightarrow                        \\[2ex]
  \hspace{2em}W^{1,p}(
      \Sigma_{(\zeta,\vec{t},\vec{t}^{\prime},\vec{\mu})},
      \partial \Sigma_{(\zeta,\vec{t},\vec{t}^{\prime},\vec{\mu})};
      h_{\approxi,S(\zeta,\vec{t},\vec{t}^{\prime},
                          \vec{\mu},\vec{\lambda})}^{\ast}
       (T_{\ast}W[k]_{\vec{\lambda}}),
       ( h_{\approxi,S(\zeta,\vec{t},\vec{t}^{\prime},
                                     \vec{\mu},\vec{\lambda})}
         |_{\partial\Sigma_{(\zeta,\vec{t},\vec{t}^{\prime},
                                           \vec{\mu}}} )^{\ast}
      T_{\ast}L )\,,
 \end{array}
$$
} % end-small

\noindent
where we regard $Q_{\rho}$ as a linear map on
 $L^p(\Sigma;\Lambda^{0,1}\Sigma\otimes_J f^{\ast}T_{\ast}Y_{[k]})$
 that is $0$ on $E_{\rho}$,
has the following property:

\bigskip

\noindent
{\bf Lemma 5.3.4.2 [approximate right inverse].} {\it
 $Q^{\prime}_{S(\zeta,\vec{t},\vec{t}^{\prime},
                              \vec{\mu},\vec{\lambda})}$
  is an approximate right inverse of
  $ \pi_{E_{S(\zeta,\vec{t},\vec{t}^{\prime},
                           \vec{\mu},\vec{\lambda})}^{\tinyaux}}
    \circ D_{h_{\tinyapproxi,
             S(\zeta,\vec{t},\vec{t}^{\prime},\vec{\mu},\vec{\lambda})}}
           \bar{\partial}_J$
  in the sense that
 $$
  \left\|\, \left(\pi_{E_{S(\zeta,\vec{t},\vec{t}^{\prime},
                                  \vec{\mu},\vec{\lambda})}^{\tinyaux}}
               \circ D_{h_{\tinyapproxi,
                           S(\zeta,\vec{t},\vec{t}^{\prime},
                                           \vec{\mu},\vec{\lambda})}}
                       \bar{\partial}_J\right)
       \circ\, Q^{\prime}_{S(\zeta,\vec{t},\vec{t}^{\prime},
                                   \vec{\mu},\vec{\lambda})}(\eta)\,
       -\,\eta\,\right\|_{L^{\raisebox{.3ex}{\tiny $p$}}}\;
  \le\;\frac{1}{2}\,\|\eta\|_{L^{\raisebox{.3ex}{\tiny $p$}}}
 $$
 for $\|\zeta\|$, $\|\vec{t}\|$, $\|\vec{t}^{\prime}\|$, $\|\vec{\mu}\|$
     small enough.
} % end-lemma

\bigskip

\noindent
{\it Proof.}
See [MD-S1: Lemma A.4.2], [F-O: Lemma 13.11],
    [Liu(C): Proposition 6.30], [L-R: proof of Lemma 4.8].

\noindent\hspace{15cm}$\Box$

\bigskip

Recall the universal approximate-$J$-holomorphic map
 $h_{\approxi}:{\cal C}/\Theta_{\rho} \rightarrow W[k]/B[k]$
 associated to the family
 $h_{\approxi,(\zeta,\vec{t},\vec{t}^{\prime},\vec{\mu},
               \vec{\lambda},\vec{a},\vec{0})}$,
 $(\zeta,\vec{t},\vec{t}^{\prime},\vec{\mu},
               \vec{\lambda},\vec{a},\vec{0}) \in \Theta_{\rho}$.
The following definition is inspired by the built-in family-treatment
 in the study of moduli problems in algebraic geometry and
 the fact that a $W^{k,p}$ Sobolev space is the completion of
     the related $C^{\infty}$ space with the $W^{k,p}$ norm:

\bigskip

\noindent
{\bf Definition 5.3.4.3
      [continuous-$\pi_{\scriptsizeDef(\Sigma)\times B[k]}
                    (\widetilde{V}_{\rho})$-family  of operators].}
{\rm
 A collection of linear operators

 \vspace{-.6ex}
 {\small
 $$
 \begin{array}{l}
  O_{S(\zeta,\vec{t},\vec{t}^{\prime},\vec{\mu},\vec{\lambda})}\; :\;
   \left.
    L^p(
     \Sigma_{(\zeta,\vec{t},\vec{t}^{\prime},\vec{\mu})};
     \Lambda^{0,1} \Sigma_{(\zeta,\vec{t},\vec{t}^{\prime},\vec{\mu})}
      \otimes_J  h_{\approxi,S(\zeta,\vec{t},\vec{t}^{\prime},
                                      \vec{\mu},\vec{\lambda})}^{\ast}
                  T_{\ast}(W[k]_{\vec{\lambda}}) ) \right/
       \raisebox{-.6ex}{$E_{S(\zeta,\vec{t},\vec{t}^{\prime},
                                      \vec{\mu},\vec{\lambda})}^{\aux}$}
       \hspace{1ex} \longrightarrow                        \\[2ex]
  \hspace{2em}W^{1,p}(
      \Sigma_{(\zeta,\vec{t},\vec{t}^{\prime},\vec{\mu})},
      \partial \Sigma_{(\zeta,\vec{t},\vec{t}^{\prime},\vec{\mu})};
      h_{\approxi,S(\zeta,\vec{t},\vec{t}^{\prime},
                          \vec{\mu},\vec{\lambda})}^{\ast}
       (T_{\ast}W[k]_{\vec{\lambda}}),
       ( h_{\approxi,S(\zeta,\vec{t},\vec{t}^{\prime},
                                     \vec{\mu},\vec{\lambda})}
         |_{\partial\Sigma_{(\zeta,\vec{t},\vec{t}^{\prime},
                                           \vec{\mu}}} )^{\ast}
      T_{\ast}L )\,,
 \end{array}
 $$
 } % end-small

 \noindent
 over $\pi_{\scriptsizeDef(\Sigma)\times B[k]}(\widetilde{V}_{\rho})$
 are said to form
 a {\it continuous-$\pi_{\scriptsizeDef(\Sigma)\times B[k]}
                   (\widetilde{V}_{\rho})$-family} of operators
 if the collection can be enlarged to a collection of linear operators

 \vspace{-.6ex}
 {\small
 $$
  \begin{array}{l}
   O_{(\zeta,\vec{t},\vec{t}^{\prime},\vec{\mu},\vec{\lambda},
                                                \vec{a},\vec{0})}\;
   :\;
    \left.
     L^p(
      \Sigma_{(\zeta,\vec{t},\vec{t}^{\prime},\vec{\mu})};
      \Lambda^{0,1} \Sigma_{(\zeta,\vec{t},\vec{t}^{\prime},\vec{\mu})}
       \otimes_J  h_{\approxi,(\zeta,\vec{t},\vec{t}^{\prime},
                       \vec{\mu},\vec{\lambda},\vec{a},\vec{0})}^{\ast}
                  T_{\ast}(W[k]_{\vec{\lambda}}) ) \right/
        \raisebox{-.6ex}{$E_{(\zeta,\vec{t},\vec{t}^{\prime},\vec{\mu},
                              \vec{\lambda},\vec{a},\vec{0})}^{\aux}$}
        \hspace{1ex} \longrightarrow                        \\[2ex]
   \hspace{2em}W^{1,p}(
       \Sigma_{(\zeta,\vec{t},\vec{t}^{\prime},\vec{\mu})},
       \partial \Sigma_{(\zeta,\vec{t},\vec{t}^{\prime},\vec{\mu})};
       h_{\approxi,(\zeta,\vec{t},\vec{t}^{\prime},\vec{\mu},
                    \vec{\lambda},\vec{a},\vec{0})}^{\ast}
        (T_{\ast}W[k]_{\vec{\lambda}}),
       ( h_{\approxi,(\zeta,\vec{t},\vec{t}^{\prime},
                      \vec{\mu},\vec{\lambda},\vec{a},\vec{0})}
          |_{\partial\Sigma_{(\zeta,\vec{t},\vec{t}^{\prime},
                                            \vec{\mu}}} )^{\ast}
       T_{\ast}L )\,,
  \end{array}
 $$
 } % end-small

 \noindent
 over $\Theta_{\rho}$ such that, for all
  $\eta\in
   C^{\infty}({\cal C};
       \Lambda^{0,1}_{{\cal C}/\Theta_{\rho}} \otimes
       h_{\approxi}^{\ast} T_{W[k]/B[k]}/
         C^{\infty}(E_{\Theta_{\rho}}^{\aux})$
  with
  $\eta|_{C_{(\zeta,\vec{t},\vec{t}^{\prime},\vec{\mu},
          \vec{\lambda},\vec{a},\vec{0})}}\\  \in
    \left.
     C^{\infty}(
      \Sigma_{(\zeta,\vec{t},\vec{t}^{\prime},\vec{\mu})};
      \Lambda^{0,1} \Sigma_{(\zeta,\vec{t},\vec{t}^{\prime},\vec{\mu})}
       \otimes_J  h_{\approxi,(\zeta,\vec{t},\vec{t}^{\prime},
                       \vec{\mu},\vec{\lambda},\vec{a},\vec{0})}^{\ast}
                  T_{\ast}(W[k]_{\vec{\lambda}}) ) \right/
        \raisebox{-.6ex}{$E_{(\zeta,\vec{t},\vec{t}^{\prime},\vec{\mu},
                              \vec{\lambda},\vec{a},\vec{0})}^{\aux}$}$,
 there exists a $\xi\in C^0({\cal C}; h_{\approxi}^{\ast}T_{W[k]/B[k]})$
  such that
  $\xi|_{C_{(\zeta,\vec{t},\vec{t}^{\prime},\vec{\mu},
          \vec{\lambda},\vec{a},\vec{0})}}
   = O_{(\zeta,\vec{t},\vec{t}^{\prime},\vec{\mu},
             \vec{\lambda},\vec{a},\vec{0})}
      (\eta|_{C_{(\zeta,\vec{t},\vec{t}^{\prime},\vec{\mu},
         \vec{\lambda},\vec{a},\vec{0})}})$.
} % end-definition

\bigskip

\noindent
{\bf Proposition 5.3.4.4
     [continuous-$\pi_{\scriptsizeDef(\Sigma)\times B[k]}
                     (\widetilde{V}_{\rho})$-family
                      of right inverse].} {\it
 For
  $\|\zeta\|$, $\|\vec{t}\|$, $\|\vec{t}^{\prime}\|$, $\|\vec{\mu}\|$
  small enough,
 there exist a constant $c$ and right inverses
  $Q_{S(\zeta,\vec{t},\vec{t}^{\prime},\vec{\mu},\vec{\lambda})}$  of
  $\pi_{E_{S(\zeta,\vec{t},\vec{t}^{\prime},\vec{\mu},\vec{\lambda})}
         ^{\tinyaux}}
           \circ D_{h_{\tinyapproxi,
                       S(\zeta,\vec{t},\vec{t}^{\prime},
                                       \vec{\mu},\vec{\lambda})}}
                       \bar{\partial}_J$
 such that
  their operator norm is uniformly bounded by $c$ and that
  they form a continuous-$\pi_{\scriptsizeDef(\Sigma)\times B[k]}
                                 (\widetilde{V}_{\rho})$-family
  of linear operators in the sense of Definition 5.3.4.3.
     % Definition [continuous-$$\pi_{\scriptsizeDef(\Sigma)\times B[k]}
     %              (\widetilde{V}_{\rho})$-family of operators]
} % end-corollary

\bigskip

\noindent{\it Proof.}
 Lemma 5.3.4.2
  % Lemma [approximate right inverse]
  implies that
  $\left(\rule{0ex}{2ex}\right.
     \pi_{E_{S(\zeta,\vec{t},\vec{t}^{\prime},
                   \vec{\mu},\vec{\lambda})}^{\tinyaux}}
               \circ D_{h_{\tinyapproxi,
                           S(\zeta,\vec{t},\vec{t}^{\prime},
                             \vec{\mu},\vec{\lambda})}}
                      \bar{\partial}_J
     \left.\rule{0ex}{2ex}\right)
       \circ\, Q^{\prime}_{S(\zeta,\vec{t},\vec{t}^{\prime},
                                   \vec{\mu},\vec{\lambda})}$
  is invertible.
 A right inverse of
  $\pi_{E_{S(\zeta,\vec{t},\vec{t}^{\prime},
                   \vec{\mu},\vec{\lambda})}^{\tinyaux}}
           \circ D_{h_{\tinyapproxi,
                       S(\zeta,\vec{t},\vec{t}^{\prime},
                               \vec{\mu},\vec{\lambda})}}
                  \bar{\partial}_J$
  is thus given by
  $$
   Q_{S(\zeta,\vec{t},\vec{t}^{\prime},\vec{\mu},\vec{\lambda})}\;
   =\; Q^{\prime}_{S(\zeta,\vec{t},\vec{t}^{\prime},
                                   \vec{\mu},\vec{\lambda})}\,
       \circ\,
       \left(
         \left(\rule{0ex}{2ex}\right.
          \pi_{E_{S(\zeta,\vec{t},\vec{t}^{\prime},
                         \vec{\mu},\vec{\lambda})}^{\tinyaux}}
                  \circ D_{h_{\tinyapproxi,
                              S(\zeta,\vec{t},\vec{t}^{\prime},
                                      \vec{\mu},\vec{\lambda})}}
                         \bar{\partial}_J
          \left.\rule{0ex}{2ex}\right)
          \circ\, Q^{\prime}_{S(\zeta,\vec{t},\vec{t}^{\prime},
                                      \vec{\mu},\vec{\lambda})}
       \right)^{-1}\,,
  $$

 \noindent
 (cf.\ [MD-S1: Sec.~3.3], [F-O: (13.2)], [Liu(C): Corollary 6.31],
       and [L-R: (4.31)]).

 To see that they form
  a continuous-$\pi_{\scriptsizeDef(\Sigma)\times B[k]}
                                 (\widetilde{V}_{\rho})$-family,
  recall from Notation 5.3.3.3 that
   $h_{\approxi,(\zeta,\vec{t},\vec{t}^{\prime},
       \vec{\mu},\vec{\lambda},\vec{a},\vec{b})}(\cdot)
  = \exp_{h_{\tinyapproxi,S(\zeta,\vec{t},\vec{t}^{\prime},
                            \vec{\mu},\vec{\lambda}) }(\,\cdot\,) }
      \xi_{(\zeta,\vec{t},\vec{t}^{\prime},
             \vec{\mu},\vec{\lambda},\vec{a},\vec{b})}(\,\cdot\,)$
  for a unique

  \vspace{-1ex}
  {\small
  $$
   \begin{array}{l}
    \xi_{(\zeta,\vec{t},\vec{t}^{\prime},
        \vec{\mu},\vec{\lambda},\vec{a},\vec{b})}
      \\ [1ex]
    \hspace{1ex}\in\;
    W^{1,p}
    \left(
         \Sigma_{(\zeta,\vec{t},\vec{t}^{\prime},\vec{\mu})},
         \partial\Sigma_{(\zeta,\vec{t},\vec{t}^{\prime},\vec{\mu})};
         h_{\approxi, S(\zeta,\vec{t},\vec{t}^{\prime},
                                            \vec{\mu},\vec{\lambda})}
           ^{\ast}T_{\ast}(W[k]_{\vec{\lambda}}),
        (h_{\approxi,S(\zeta,\vec{t},\vec{t}^{\prime},
                                       \vec{\mu},\vec{\lambda})}
         |_{\partial\Sigma_{(\zeta,\vec{t},\vec{t}^{\prime},\vec{\mu})}})
        ^{\ast} T_{\ast}L
                \right)\,.
   \end{array}
  $$
  } % end-small

 \noindent
 Using the parallel transport
  along the geodesic determined
  by $\xi_{(\zeta,\vec{t},\vec{t}^{\prime},
       \vec{\mu},\vec{\lambda},\vec{a},\vec{b})}(x)$,
 one can extend the collections of operators
   $\IP_{S(\zeta,\vec{t},\vec{t}^{\prime},\vec{\mu},\vec{\lambda})}$,
   $\Glue_{S(\zeta,\vec{t},\vec{t}^{\prime},\vec{\mu},\vec{\lambda})}$
  to the collections of operators
   $\IP_{(\zeta,\vec{t},\vec{t}^{\prime},\vec{\mu},\vec{\lambda},
          \vec{a},\vec{b})}$,
   $\Glue_{(\zeta,\vec{t},\vec{t}^{\prime},\vec{\mu},\vec{\lambda},
           \vec{a},\vec{b})}$
   over $\widetilde{V}_{\rho}$, and, hence,
  in particular over $\Theta_{\rho}$.
 The collections of operators

  {\small
  $$
   \begin{array}{l}
    Q^{\prime}_{(\zeta,\vec{t},\vec{t}^{\prime},\vec{\mu},\vec{\lambda},
                 \vec{a},\vec{0})}\;
    :=\;
     \Glue_{(\zeta,\vec{t},\vec{t}^{\prime},\vec{\mu},\vec{\lambda},
               \vec{a},\vec{0})}
       \circ Q_{\rho} \circ
       \IP_{(\zeta,\vec{t},\vec{t}^{\prime},\vec{\mu},\vec{\lambda},
             \vec{a},\vec{0})}\,, \\[.6ex]
    Q_{(\zeta,\vec{t},\vec{t}^{\prime},\vec{\mu},\vec{\lambda},
        \vec{a},\vec{0})}\;
     :=\;
     Q^{\prime}_{(\zeta,\vec{t},\vec{t}^{\prime},
                   \vec{\mu},\vec{\lambda},\vec{a},\vec{0})}\,
       \circ\,
       \left(
         \left(\rule{0ex}{2ex}\right.
          \pi_{E_{(\zeta,\vec{t},\vec{t}^{\prime},
                   \vec{\mu},\vec{\lambda},\vec{a},\vec{0})}^{\tinyaux}}
                  \circ D_{h_{\tinyapproxi,
                              (\zeta,\vec{t},\vec{t}^{\prime},
                               \vec{\mu},\vec{\lambda},\vec{a},\vec{0})}}
                         \bar{\partial}_J
          \left.\rule{0ex}{2ex}\right)
          \circ\, Q^{\prime}_{(\zeta,\vec{t},\vec{t}^{\prime},
                               \vec{\mu},\vec{\lambda},\vec{a},\vec{0})}
       \right)^{-1}
   \end{array}
  $$
  } % end-small

  \noindent
  extend the collections
  $Q^{\prime}_{S(\zeta,\vec{t},\vec{t}^{\prime},\vec{\mu},\vec{\lambda})}$,
  $Q_{S(\zeta,\vec{t},\vec{t}^{\prime},\vec{\mu},\vec{\lambda})}$.
  %%%%%%
  % \marginpar{\raggedright\tiny $\bullet$ Think more carefully.}
  % Note that
  %  $\IP_{(\zeta,\vec{t},\vec{t}^{\prime},\vec{\mu},\vec{\lambda},
  %               \vec{a},\vec{0})}$
  %  takes a $C^{\infty}$ section $\eta^{\prime}$ to
  %  a piecewise-$C^{\infty}$ section $\eta^{\prime\prime}$
  %  with the discontinuity occurs in a small neighborhood of nodes,
  %  elliptic regularity implies that
  %  $Q_{\rho}$ sends $\eta^{\prime\prime}$ to
  %  a piecewise-$C^{\infty}$ section $\xi^{\prime\prime}$
  %  with the singularity locus supported in the above discontinuity locus,
  %  $\Glue_{(\zeta,\vec{t},\vec{t}^{\prime},\vec{\mu},\vec{\lambda},
  %                 \vec{a},\vec{0})}$
  %  in the end smoothes out the singularity locus and sends
  %  $\xi^{\prime\prime}$ to a $C^{\infty}$ $\xi^{\prime}$.
  % It follows that
  %  $Q^{\prime}_{(\zeta,\vec{t},\vec{t}^{\prime},\vec{\mu},\vec{\lambda},
  %                \vec{a},\vec{0})}$ and , hence,
  %  $Q_{(\zeta,\vec{t},\vec{t}^{\prime},\vec{\mu},\vec{\lambda},
  %             \vec{a},\vec{0})}$
  %  send $C^{\infty}$ sections to $C^{\infty}$ sections.
  %%%%%%
 The explicit expressions of
  $\IP_{\mbox{\tiny $\bullet$}}$ and $\Glue_{\mbox{\tiny $\bullet$}}$
  imply that
 these operators over $\Theta_{\rho}$ together satisfy the
  continuous-family behavior required in Definition 5.3.4.3.
    %  Definition [continuous-$\pi_{\scriptsizeDef(\Sigma)\times B[k]}
    %              (\widetilde{V}_{\rho})$-family  of operators].

\noindent\hspace{15cm}$\Box$

\bigskip

\begin{flushleft}
{\bf Newton-Picard iteration:
     deforming $h_{\approxi,\mbox{\tiny $\bullet$}}$ to
      a $(J,E_{\mbox{\tiny $\bullet$}})$-holomorphic map
      $f_{\mbox{\tiny $\bullet$}}$.}
\end{flushleft}
Once one realized that the continuity of the relative construction
 has to be made over $\Theta_{\rho}$,
 not directly over
  $\pi_{\scriptsizeDef(\Sigma)\times B[k]}(\widetilde{V}_{\rho})
   \subset \Def(\Sigma)\times B[k]$, and
 has constructed the ingredients accordingly,
the rest of the discussion is similar to those in
 [MD-S1: proof of Theorem 3.3.4],
 [F-O: pp.~987-988] (directly on the maps), and
 [Liu(C): proof of Proposition 6.32];
see also [I-P2: Sec.~9] and [L-R: proof of Proposition 4.10].
We give a sketch below to conclude the discussion.

Beginning with the $\widetilde{V}_{\rho}$-family of maps
 $h_{\approxi,(\zeta,\vec{t},\vec{t}^{\prime},
               \vec{\mu},\vec{\lambda},\vec{a},\vec{b})}$,
define a sequence of $\widetilde{V}_{\rho}$-family of maps
 as follows:
 \begin{itemize}
  \item[$\cdot$]
   Set
    $$
     \begin{array}{lcl}
      h_{1,(\zeta,\vec{t},\vec{t}^{\prime},\vec{\mu},
                          \vec{\lambda},\vec{a},\vec{b})}
       & = & h_{\approxi,(\zeta,\vec{t},\vec{t}^{\prime},
                          \vec{\mu},\vec{\lambda},\vec{a},\vec{b})}\;
             =\; \exp_{h_{\tinyapproxi,S(\zeta,\vec{t},\vec{t}^{\prime},
                                         \vec{\mu},\vec{\lambda}) } }
                  \xi_{(\zeta,\vec{t},\vec{t}^{\prime},
                    \vec{\mu},\vec{\lambda},\vec{a},\vec{b})}\,,\\[.6ex]
      \xi_{1, (\zeta,\vec{t},\vec{t}^{\prime},
               \vec{\mu},\vec{\lambda},\vec{a},\vec{b})}
       & = & \xi_{(\zeta,\vec{t},\vec{t}^{\prime},
                   \vec{\mu},\vec{\lambda},\vec{a},\vec{b})}\,.
     \end{array}
    $$

  \item[$\cdot$]
   Suppose that
    $h_{n,(\zeta,\vec{t},\vec{t}^{\prime},\vec{\mu},
                       \vec{\lambda},\vec{a},\vec{b})}
    = \exp_{\,h_{\tinyapproxi,S(\zeta,\vec{t},\vec{t}^{\prime},
                            \vec{\mu},\vec{\lambda}) } }
      \xi_{n, (\zeta,\vec{t},\vec{t}^{\prime},
                  \vec{\mu},\vec{\lambda},\vec{a},\vec{b})}$
    is defined,
   let
    $$
     \xi_{n+1, (\zeta,\vec{t},\vec{t}^{\prime},
                    \vec{\mu},\vec{\lambda},\vec{a},\vec{b})}\;
     =\; \xi_{n, (\zeta,\vec{t},\vec{t}^{\prime},\vec{\mu}
                              \vec{\lambda},\vec{a},\vec{b})}\,
         -\,
          Q_{S(\zeta,\vec{t},\vec{t}^{\prime},\vec{\mu},\vec{\lambda})}
          \circ \pi_{E_{S(\zeta,\vec{t},\vec{t}^{\prime},
                           \vec{\mu},\vec{\lambda})}^{\tinyaux}}
          \circ P_n \circ
          (\bar{\partial}_J
             h_{n,(\zeta,\vec{t},\vec{t}^{\prime},\vec{\mu},
                                 \vec{\lambda},\vec{a},\vec{b})})\,,
    $$
     where
      $$
       \begin{array}{l}
        P_n\; :\;
         L^p(
          \Sigma_{(\zeta,\vec{t},\vec{t}^{\prime},\vec{\mu})};
           \Lambda^{0,1}
            \Sigma_{(\zeta,\vec{t},\vec{t}^{\prime},\vec{\mu})}
           \otimes_J
            h_{n, (\zeta,\vec{t},\vec{t}^{\prime},
                   \vec{\mu},\vec{\lambda},\vec{a},\vec{b})}^{\ast}
            T_{\ast}(W[k]_{\vec{\lambda}}) )                    \\[1ex]
        \hspace{4em}\longrightarrow\;
         L^p(
         \Sigma_{(\zeta,\vec{t},\vec{t}^{\prime},\vec{\mu})};
          \Lambda^{0,1}
           \Sigma_{(\zeta,\vec{t},\vec{t}^{\prime},\vec{\mu})}
          \otimes_J  h_{\approxi,S(\zeta,\vec{t},\vec{t}^{\prime},
                                   \vec{\mu},\vec{\lambda})}^{\ast}
                      T_{\ast}(W[k]_{\vec{\lambda}}) )
       \end{array}
      $$
      is the map induced by the parallel transport along the geodesics
      determined by
      $\xi_{n,(\zeta,\vec{t},\vec{t}^{\prime},\vec{\mu}
                     \vec{\lambda},\vec{a},\vec{b})}$,  and
    define
    $$
     h_{n+1,(\zeta,\vec{t},\vec{t}^{\prime},\vec{\mu},
                       \vec{\lambda},\vec{a},\vec{b})}\;
     =\; \exp_{\,h_{\tinyapproxi,S(\zeta,\vec{t},\vec{t}^{\prime},
                               \vec{\mu},\vec{\lambda}) } }
          \xi_{n+1, (\zeta,\vec{t},\vec{t}^{\prime},
                     \vec{\mu},\vec{\lambda},\vec{a},\vec{b})}\,.
    $$
 \end{itemize}

The series
 $$
  -\,\sum_{n=1}^{\infty}
      \pi_{E_{S(\zeta,\vec{t},\vec{t}^{\prime},\vec{\mu},\vec{\lambda})}
         ^{\tinyaux}}
       \circ P_n \circ
       (\bar{\partial}_J
          h_{n,(\zeta,\vec{t},\vec{t}^{\prime},\vec{\mu},
                              \vec{\lambda},\vec{a},\vec{b})})
 $$
 converges to an
 $$
  \eta_{(\zeta,\vec{t},\vec{t}^{\prime},\vec{\mu},
                       \vec{\lambda},\vec{a},\vec{b})}\;
  \in\;
    \left.
     L^p(
      \Sigma_{(\zeta,\vec{t},\vec{t}^{\prime},\vec{\mu})};
      \Lambda^{0,1} \Sigma_{(\zeta,\vec{t},\vec{t}^{\prime},\vec{\mu})}
       \otimes_J  h_{\approxi,S(\zeta,\vec{t},\vec{t}^{\prime},
                                       \vec{\mu},\vec{\lambda})}^{\ast}
                   T_{\ast}(W[k]_{\vec{\lambda}}) ) \right/
        \raisebox{-.6ex}{$E_{S(\zeta,\vec{t},\vec{t}^{\prime},
                               \vec{\mu},\vec{\lambda})}^{\aux}$}
 $$
and the sequence of maps
 $h_{n,(\zeta,\vec{t},\vec{t}^{\prime},\vec{\mu},
                      \vec{\lambda},\vec{a},\vec{b})}$,
  $n=1,\,\ldots,\,\infty$,
 converge both uniformly and with respect to the $W^{1,p}$-topology
 (as the $W^{1,p}$-norm dominates the $C^0$-norm for $p\gg 0$) to
 $$
  f_{(\zeta,\vec{t},\vec{t}^{\prime},\vec{\mu},
                    \vec{\lambda},\vec{a},\vec{b})}\;
  =\; \exp_{\,h_{\tinyapproxi,S(\zeta,\vec{t},\vec{t}^{\prime},
                                    \vec{\mu},\vec{\lambda}) } }
          \left(
            \xi_{(\zeta,\vec{t},\vec{t}^{\prime},
              \vec{\mu},\vec{\lambda},\vec{a},\vec{b})}\,
           +\, Q_{S(\zeta,\vec{t},\vec{t}^{\prime},\vec{\mu},\vec{\lambda})}
                \eta_{(\zeta,\vec{t},\vec{t}^{\prime},\vec{\mu},
                       \vec{\lambda},\vec{a},\vec{b})}
          \right)\,.
 $$
This gives rise to a continuous-$\widetilde{V}_{\rho}$-family of maps.

Define the trivialized
 {\it obstruction bundle $E_{\widetilde{V}_{\rho}}$
      over $\widetilde{V}_{\rho}$}
 by setting its fiber
 $$
  E_{(\zeta,\vec{t},\vec{t}^{\prime},\vec{\mu},
                                     \vec{\lambda},\vec{a},\vec{b})}\;
  \subset\;
  L^p(
   \Sigma_{(\zeta,\vec{t},\vec{t}^{\prime},\vec{\mu})};
   \Lambda^{0,1} \Sigma_{(\zeta,\vec{t},\vec{t}^{\prime},\vec{\mu})}
    \otimes_J  f_{(\zeta,\vec{t},\vec{t}^{\prime},
                                   \vec{\mu},\vec{\lambda})}^{\ast}
                T_{\ast}(W[k]_{\vec{\lambda}}) )
  $$
 at  $(\zeta,\vec{t},\vec{t}^{\prime},\vec{\mu},
                                 \vec{\lambda},\vec{a},\vec{b})$
 to be the parallel transport of
 $E_{S(\zeta,\vec{t},\vec{t}^{\prime},\vec{\mu},\vec{\lambda})}^{\aux}$
 along the geodesics determined by
 $\xi_{(\zeta,\vec{t},\vec{t}^{\prime},
        \vec{\mu},\vec{\lambda},\vec{a},\vec{b})}\,
   +\, Q_{S(\zeta,\vec{t},\vec{t}^{\prime},\vec{\mu},\vec{\lambda})}
        \eta_{(\zeta,\vec{t},\vec{t}^{\prime},\vec{\mu},
               \vec{\lambda},\vec{a},\vec{b})}$.
Let
 $$
  \begin{array}{l}
   \pi_{E_{(\zeta,\vec{t},\vec{t}^{\prime},\vec{\mu},
                          \vec{\lambda},\vec{a},\vec{b})}}\; :\;
   L^p (
    \Sigma_{(\zeta,\vec{t},\vec{t}^{\prime},\vec{\mu})};
    \Lambda^{0,1}\Sigma_{(\zeta,\vec{t},\vec{t}^{\prime},\vec{\mu})}
      \otimes_J
        f_{(\zeta,\vec{t},\vec{t}^{\prime},\vec{\mu},
                          \vec{\lambda},\vec{a},\vec{b})}^{\ast}
        T_{\ast}W[k]_{\vec{\lambda}}  )   \\[1ex]
   \hspace{2em}\longrightarrow\;
     \left.
       L^p (
       \Sigma_{(\zeta,\vec{t},\vec{t}^{\prime},\vec{\mu})};
      \Lambda^{0,1}\Sigma_{(\zeta,\vec{t},\vec{t}^{\prime},\vec{\mu})}
        \otimes_J
          f_{(\zeta,\vec{t},\vec{t}^{\prime},\vec{\mu},
                            \vec{\lambda},\vec{a},\vec{b})}^{\ast}
          T_{\ast}W[k]_{\vec{\lambda}}  )
       \right/\!
        \raisebox{-.6ex}{$E_{(\zeta,\vec{t},\vec{t}^{\prime},\vec{\mu},
                              \vec{\lambda},\vec{a},\vec{b})}$}\;
  \end{array}
 $$
 be the quotient map; then, by construction,
 $$
  \pi_{E_{(\zeta,\vec{t},\vec{t}^{\prime},\vec{\mu},
                         \vec{\lambda},\vec{a},\vec{b})}}
   \circ \bar{\partial}_J\,
   f_{\zeta,\vec{t},\vec{t}^{\prime},\vec{\mu},
                    \vec{\lambda},\vec{a},\vec{b}}\;
   =\;0\,.
 $$
In other words, the collection of maps
 $f_{(\zeta,\vec{t},\vec{t}^{\prime},\vec{\mu},
                    \vec{\lambda},\vec{a},\vec{b})}$
 form a continuous-$\widetilde{V}_{\rho}$-family of
 $(J,E_{\mbox{\tiny $\bullet$}})$-holomorphic maps.

\bigskip

\noindent
{\bf Proposition 5.3.4.5
     [$\Aut(\rho)$-equivariant pre-deformable family].}
{\it
 The bundle $\widetilde{E}_{\widetilde{V}_{\rho}}$ is
  $\Aut(\rho)$-equivariant over $\widetilde{V}_{\rho}$ and
 the collection of maps
  $f_{(\zeta,\vec{t},\vec{t}^{\prime},
                   \vec{\mu},\vec{\lambda},\vec{a},\vec{b})}$
 form a $\Aut(\rho)$-equivariant
  continuous-$\widetilde{V}_{\rho}$-family of pre-deformable
  $(J,E_{\mbox{\tiny $\bullet$}})$-holomorphic maps.
} % end-proposition

\bigskip

\noindent{\it Proof.}
 The $\Aut(\rho)$-equivariance of the family
  $f_{(\zeta,\vec{t},\vec{t}^{\prime},
             \vec{\mu},\vec{\lambda},\vec{a},\vec{b})}$
  follows from the $\Aut(\rho)$-invariance of
  the family of maps
   $h_{\approxi,(\zeta,\vec{t},\vec{t}^{\prime},
                 \vec{\mu},\vec{\lambda},\vec{a},\vec{0})}$,
  the bundle $E_{S(\pi_{\tinyDef(\Sigma)\times B[k]}
                           (\widetilde{V}_{\rho}))}^{\aux}$, and
  the family of operators
   $Q_{\approxi,(\zeta,\vec{t},\vec{t}^{\prime},
                 \vec{\mu},\vec{\lambda},\vec{a},\vec{0})}$.
 The $\Aut(\rho)$-equivariance of $E_{\widetilde{V}_{\rho}}$
  follows then from the $\Aut(\rho)$-equivariant of the family of maps
  $f_{(\zeta,\vec{t},\vec{t}^{\prime},
       \vec{\mu},\vec{\lambda},\vec{a},\vec{b})}$.
 It remains to prove the pre-deformability of
  $f_{(\zeta,\vec{t},\vec{t}^{\prime},
             \vec{\mu},\vec{\lambda},\vec{a},\vec{b})}$.

 Note first that, by construction, elements of
  $E_{(\zeta,\vec{t},\vec{t}^{\prime},
              \vec{\mu},\vec{\lambda},\vec{a},\vec{b})}$
  are supported in a compact subset in the complement of
  the union of $\varepsilon$-neighborhood of nodes and
  the annuli or bands on
  $\Sigma_{(\zeta,\vec{t},\vec{t}^{\prime},\vec{\mu})}$
  from smoothing related nodes of $\Sigma$.
 This implies in particular that
  $f_{(\zeta,\vec{t},\vec{t}^{\prime},
              \vec{\mu},\vec{\lambda},\vec{a},\vec{b})}$
  is $J$-holomorphic in the $\varepsilon$-neighborhood of nodes
  and hence it makes sense to talk about pre-deformability of
  $f_{(\zeta,\vec{t},\vec{t}^{\prime},
             \vec{\mu},\vec{\lambda},\vec{a},\vec{b})}$
  at distinguished nodes.
 Furthermore, as the universal map
  $F:{\cal C}/\widetilde{V}_{\rho}\rightarrow W[k]/B[k]$
  associated to the family of maps
  $f_{(\zeta,\vec{t},\vec{t}^{\prime},
              \vec{\mu},\vec{\lambda},\vec{a},\vec{b})}$
  is continuous with the central $f$ pre-deformable,
  there can be no mass falling into the locus
   $(W[k]/B[k])_{\sing}$ of singularities of the fibers of $W[k]/B[k]$.
 In other words, $F$ is a family of flat maps in the sense of
  [I-P2: Definition 3.1]
  as long as
   $(\zeta,\vec{t},\vec{t}^{\prime},
           \vec{\mu},\vec{\lambda},\vec{a}-\vec{a}_f,\vec{b})$
   is sufficiently small, a condition that is already incorporated
   implicitly into the definition of $\widetilde{V}_{\rho}$.
 As the fibers of $F$ over an open-dense subset of
  $\widetilde{V}_{\rho}$ are maps from smooth domains(-with-boundary)
  to smooth fibers of $W[k]/B[k]$,
 it follows from [I-P2: Lemma 3.3] that
 this above flatness property implies that the fibers
  $f_{(\zeta,\vec{t},\vec{t}^{\prime},
              \vec{\mu},\vec{\lambda},\vec{a},\vec{b})}$
  of $F$ over $\widetilde{V}_{\rho}$ must be all pre-deformable
  for $(\zeta,\vec{t},\vec{t}^{\prime},
        \vec{\mu},\vec{\lambda},\vec{a}-\vec{a}_f,\vec{b})\in
        \widetilde{V}_{\rho}$).
 This concludes the proof.

\noindent\hspace{15cm}$\Box$

%%%%%%%%%%%%%%%%%%%%%%%%%%%%%%
% \bigskip
%
% \noindent
% {\it Remark ??? $[$slight variation of the construction$]$.}
% %
% \marginpar{\raggedright\tiny $\bullet$ To be completed.}
%
% \bigskip
%
% \noindent
% {\it Remark ??? $[$slightly deformed $\widetilde{V}_{\rho}$-family$]$.}
% %
% \marginpar{\raggedright\tiny $\bullet$ To be completed.}
%
% \bigskip
%
% \bigskip
% \noindent $\circ$
% Another $\widetilde{V}_{\rho}$-family of maps can be constructed
%  by applying Proposition ???
%   % Proposition
%   % [deforming approximate- to exact $(J,E)$-holomorphic map]
%  to the subfamily $\widetilde{V}_{\rho}^0$ of
%   approximate-$(J,E)$-holomorphic maps
%   $h_{\zeta,\vec{t},\vec{t}^{\prime},
%                     \vec{\mu},\vec{\lambda},\vec{a},0}$.
%
% \bigskip
% \noindent $\circ$
% The new $\widetilde{V}_{\rho}$-family of maps
%  $f^{\prime}_{\zeta,\vec{t},\vec{t}^{\prime},
%                         \vec{\mu},\vec{\lambda},\vec{a},\xi}$
%  are $(J,E)$-holomorphic respect to the restriction of $E$ to
%  $\widetilde{V}_{\rho}^0$.
% This gives a different map from $\widetilde{V}_{\rho}$ to
%  $W^{1,p}_{(g,h),(n,\vec{m})}
%    (W[k]_{\vec{\lambda}},L\,|\,[\beta],\vec{\gamma},\mu)$
%  and is closer to the construction in [Liu(C)].
%
% \bigskip
%
% \noindent $\circ$
% Both contains a common locus in a neighborhood of $\rho$
%  on the subspace $???$ of\\
%  $W^{1,p}_{(g,h),(n,\vec{m})}
%     (W[k]_{\vec{\lambda}},L\,|\,[\beta],\vec{\gamma},\mu)$
%  that consists of $J$-holomorphic maps.
%%%%%%%%%%%%%%%%%%%%%%%%%%%%%%%%

\bigskip

\subsubsection{Rigidification:
    a Kuranishi neighborhood-in-${\cal C}_{\mbox{\tiny\rm spsccw}}$
     $V_{\rho}/B$ of $\rho$ on\\
      $\;\overline{\cal M}_{(g,h),(n,\vec{m})}
       (W/B,L\,|\,[\beta],\vec{\gamma},\mu)\,/B$. }

How the $\widetilde{V}_{\rho}$-family of
 $(J,E_{\tinybullet})$-holomorphic maps
 $f_{(\zeta,\vec{t},\vec{t}^{\prime},
           \vec{\mu},\vec{\lambda},\vec{a},\vec{b})}$
 gives rise to a Kuranishi neighborhood
 $V_{\rho}/B$ of $\rho$ on
  $\;\overline{\cal M}_{(g,h),(n,\vec{m})}
          (W/B,L\,|\,[\beta],\vec{\gamma},\mu)\,/B$
 is explained in this subsubsection.

\bigskip

\begin{flushleft}
{\bf A stratified subset $V_{\rho}/B$ of $\widetilde{V}_{\rho}/B$
     from the rigidification of $\Aut(\Sigma)\times {\Bbb G}_m[k]$.}
\end{flushleft}
$(\widetilde{V}_{\rho},E_{\widetilde{V}_{\rho}})$ and
 the associated $\widetilde{V}_{\rho}$-family of
 $(J,E_{\tinybullet})$-holomorphic maps from deformed $\Sigma$
 to fibers of $W[k]/B[k]$ are only $\Aut(\rho)$-equivariant.
However,
 the equivariant approximate product pseudo-action of
  $\Aut(\Sigma)\times{\Bbb G}_m[k]$ on
  $({\cal C}/\Def(\Sigma))\times W[k]/B[k]$
  remains to induce an equivalence relation on $\widetilde{V}_{\rho}$,
  defined by setting
  $(\zeta_1,\vec{t}_1,\vec{t}^{\prime}_1,
            \vec{\mu}_1,\vec{\lambda}_1,\vec{a}_1,\vec{b}_1)
    \sim (\zeta_2,\vec{t}_2,\vec{t}^{\prime}_2,
                  \vec{\mu}_2,\vec{\lambda}_2,\vec{a}_2,\vec{b}_2)$
  if there exists a pair
     $(\alpha,\beta)\in\Aut(\Sigma)\times{\Bbb G}_m[k]$
     such that
      $\beta
         \circ f_{(\zeta_1,\vec{t}_1,\vec{t}^{\prime}_1,
                   \vec{\mu}_1,\vec{\lambda}_1,\vec{a}_1,\vec{b}_1)}
         \circ \alpha^{-1}
        = f_{(\zeta_2,\vec{t}_2,\vec{t}^{\prime}_2,
            \vec{\mu}_2,\vec{\lambda}_2,\vec{a}_2,\vec{b}_2)}$.
Denote the $\sim$-equivalence class of $\rho$ by $O_{\rho}$;
then one has:
%%%%%%%%%%%%%%%%%%%%%%%%%%%
% Let $\tilde{s}_{\rho}$ be the section of $E_{\widetilde{V}_{\rho}}$
%  induced from the operator $s_{\bar{\partial}_J}$, (cf. Sec.~???);
% then $\tilde{s}_{\rho}^{-1}(0)$ parameterizes those
%  $(J,E_{\tinybullet})$-holomorphic maps in
%  the $\widetilde{V}_{\rho}$-family that are $J$-holomorphic.
% The equivariant approximate pseudo-$\Aut(\Sigma)$-action on
%  ${\cal C}/\Def(\Sigma)$ induces an equivalence $\sim$
%  on $\widetilde{V}_{\rho}$
%  ?????????????
%   approximate pseudo-$\Aut(\Sigma)$-action
%   on $\widetilde{s}_{\rho}^{-1}(0)/B[k]$ via
%  pre-compositions with elements in $\Aut(\Sigma)$ whenever defined,
% while the equivariant pseudo-${\Bbb G}_m[k]$-action on $W[k]/B[k]$ induces
%  a pseudo-${\Bbb G}_m[k]$-action on $\widetilde{s}^{-1}(0)/B[k]$ via
%  post-compositions with elements in ${\Bbb G}_m[k]$ whenever defined.
%  %%%%%%%%%%%%%
%  % (Here a {\it pseudo-group$\,(G)$-action} on a space $M$ is
%  %  a morphism from an {\it open subset} $A\subset G\times M$ to $M$
%  %  that satisfies the usual group action law whenever the map
%  %  and the compositions are defined.)
%  %%%%%%%%%%%%%
%
% \bigskip
%
% Since these approximate pseudo-group-actions preserve $J$-holomorphy,
%   the $\Aut(\Sigma)\times{\Bbb G}_m[k]$-orbit $O_{\rho}$ of $\rho$
%   is contained in $\widetilde{s}_{\rho}^{-1}(0)$.
%%%%%%%%%%%%%%%%%%%%%%%%%%%%%%%%%

\bigskip

\noindent
{\bf Lemma 5.3.5.1 [$O_{\rho}$ maximal].} {\it
 $O_{\rho}$ is a maximal equivalence class at $\rho$  in the sense that
 a small enough neighborhood of $\rho$ in $O_{\rho}$ is homeomorphic to
 a neighborhood of the identity element in
 $\Aut(\Sigma)\times{\Bbb G}_m[k]$.
} % end-lemma

\bigskip

\noindent
{\it Proof.}
 This is a consequence of transversality at $\rho$.
 The fiber $\widetilde{V}_0$ of
 $\widetilde{V}_{\rho}/(\Def(\Sigma)\times B[k])$
  over $(\vec{0},\vec{0})$ is embedded in the Banach manifold
  $W^{1,p}(\Sigma, Y_{[k]})$ of $W^{1,p}$-maps
  from $\Sigma$ to (the rigid) $Y_{[k]}$.
 The latter is (approximate-pseudo-)acted upon by
  $\Aut(\Sigma)\times{\Bbb G}_m[k]$.
 Under this embedding,
  $E_{\widetilde{V}_0}
   := \widetilde{E}_{\widetilde{V}_{\rho}}|_{\widetilde{V}_0}$
  is embedded in the $L^p$ -obstruction bundle
  $T^2_{W^{1,p}(\Sigma,Y_{[k]})}$ of $W^{1,p}(\Sigma,Y_{[k]})$,
  whose fiber at $\rho$ is precisely
  $L^p(\Sigma; \Lambda^{0,1}\Sigma\otimes_J f^{\ast}T_{\ast}Y_{[k]})$.
 The operator $\bar{\partial}_J$ defines a section
  $s_{\bar{\partial}_J}$ of $T^2_{W^{1,p}(\Sigma,Y_{[k]})}$,
  whose linearization at $\rho$ gives precisely the map
  $$
   D_{\!f}\bar{\partial}_J\,:\,
    W^{l,p}( \Sigma, \partial\Sigma; f^{\ast}T_{\ast}Y_{[k]},
                    (f|_{\partial\Sigma})^{\ast}T_{\ast}L_{[k]} )\;
    \longrightarrow\;
     W^{l-1,p}(\Sigma,\Lambda^{0,1}\Sigma
                      \otimes_J f^{\ast}T_{\ast}Y_{[k]})\,,
  $$
  where
   $W^{l,p}( \Sigma, \partial\Sigma; f^{\ast}T_{\ast}Y_{[k]},
                   (f|_{\partial\Sigma})^{\ast}T_{\ast}L_{[k]})$
   is now regarded as the fiber of the tangent bundle
   $T^1_{W^{1,p}(\Sigma,Y_{[k]})}$ of $W^{1,p}(\Sigma,Y_{[k]})$
   at $\rho$.
 Extend $E_{\widetilde{V}_0}$ to a subbundle $E_U$ of
  $T^2_{W^{1,p}(\Sigma,Y_{[k]})}$ over a neighborhood $U$ of
   $\rho$ in $W^{1,p}(\Sigma,Y_{[k]})$  and
 let $\pi_{E_U}$ be the quotient map
  $\pi_{E_U}:T^2_U\rightarrow T^2_U/E_U$ over $U$.
 Then
  the saturatedness of $E_{\rho}$,
  Lemma 5.3.1.1,
     % Lemma [index of $D_f\bar{\partial}_J$ for rigid target],
  Corollary 5.3.1.6,
     % Corollary [pre-deformability subspace of
     %            $(D_f\bar{\partial}_J)^{-1}(E_{\rho})$]
   and the Implicit Function Theorem (Theorem 5.3.0.2)
              % Theorem [implicit function theorem]
  together imply that
   the pre-deformability condition on
    $(J,E_U)$-holomorphic $W^{1,p}$-maps is a transverse condition on
    $(\pi_{E_{\tinybullet}}\circ s_{\bar{\partial}_J})^{-1}(0)$
    at $\rho$  and that
   the space of pre-deformable $(J,E_U)$-holomorphic $W^{1,p}$-maps near
    $\rho$ coincides with a neighborhood of $\rho$ in $\widetilde{V}_0$.
 The equivalence relation $\sim$ on $\widetilde{V}_{\rho}$ is
  the restriction of the equivalence relation on
  $W^{1,p}(\Sigma,Y_{[k]})$ defined by the
  $(\Aut(\Sigma)\times{\Bbb G}_m[k])$-orbits on $W^{1,p}(\Sigma,Y_{[k]})$.
 All these together imply that the intersection
  $((\Aut(\Sigma)\times{\Bbb G}_m[k])\cdot\rho)\cap \widetilde{V}_0$
  in $W^{1,p}(\Sigma,Y_{[k]})$ coincides with
  $O_{\rho}\subset \widetilde{V}_0$ around $\rho$.
 This concludes the lemma.

\noindent\hspace{15cm}$\Box$

\bigskip

With respect to the embedding
 $\widetilde{V}_{\rho}
  \subset \Def(\Sigma)\times B[k] \times
           \Ker(\pi_{E_{\rho}}\circ D_{\!f}\bar{\partial}_J)^{\pd}$
 in Sec.~5.3.2,
 $T_fO_{\rho}$ lies in the subspace
  $\{0\}\times\{0\}\times
   \left( \Ker(D_{\!f}\bar{\partial}_J) \cap
     \Ker(\pi_{E_{\rho}}\circ D_{\!f}\bar{\partial}_J)^{\pd} \right)$.
We will denote the quotient space
 $\Ker(\pi_{E_{\rho}}\circ D_{\!f}\bar{\partial}_J)^{\pd}
   /T_fO_{\rho}$ by $H_{\rho,\map}^{\rigidified}$.

By a combination of
 the same center-of-mass construction in [Sie1: Sec.~5.3]
  that rigidifies the approximate pseudo-$\Aut(\Sigma)$-action  and
 the same construction in Sec.~4.2
  that rigidifies the ${\Bbb G}_m[k]$-action,
there exists a $\Aut(\rho)$-equivariant rigidifying map
 $$
  R_{\rho}\; :\;
   \widetilde{V}_{\rho}\; \longrightarrow\;
   {\Bbb R}^{a+a^{\prime}}\times {\Bbb C}^{b+b^{\prime}+k}\,,
 $$
 where
  $a$ (resp.\ $b$) is the total number of unstable
   disc-components (resp.\ sphere-components) of $\Sigma$ and
  $a^{\prime}$ (resp.\ $b^{\prime}$) is the total number of
   unstable disc-components (resp.\ sphere-components) of $\Sigma$
   that has only one special point.
Let
 $$
   V_{\rho}\;=\;
    \mbox{a small enough $\Aut(\rho)$-invariant open neighborhood
     of $\rho$ in $R_{\rho}^{-1}(R_{\rho}(0))$}\,,
 $$
then $\Aut(\rho)$ acts on $V_{\rho}$ effectively.
The composition of the standard fibrations $\widetilde{V}_{\rho}/B[k]$
 and $B[k]/B$ induces a standard fibration $V_{\rho}/B$.
The stratified space $\Xi_{\mathbf s}$ induces a stratification
 on $\widetilde{V}_{\rho}$ via the projection map
 $\widetilde{V}_{\rho}\rightarrow \Xi_{\mathbf s}$.
The latter stratification restricts to a stratification on $V_{\rho}$.

\bigskip

\noindent
{\bf Lemma 5.3.5.2 [piecewise-transverse slice at $\rho$].}
{\it
  As a fibered stratified space,
   $V_{\rho}/B$ is isomorphic to
   $(\Def(\Sigma;\Lambda)\times \Xi_{\mathbf s}
       \times H_{\rho,\map;\Lambda}^{\rigidified})/B$.
} % end-lemma

\bigskip

\noindent
{\it Proof.}
 % If the pseudo $\Aut(\Sigma)\times {\Bbb G}_k[m]$-action on
 %  $\widetilde{s}_{\rho}^{-1}(0)$
 %  is also defined in the same way on the whole $\widetilde{V}_{\rho}$,
 % then this lemma is immediate.
 %
 % In general,....
 Embed $\widetilde{V}_{\rho}$ in a singular {\it un-rigidified} chart
  $\widetilde{V}_{\rho}^{\,\prime\,\sharp}$
  in Siebert's construction (cf.\ [Sie1: Sec.~5.2, Sec.~5.3])
  for $\rho$ regarded as a point in
  $\check{\cal W}^{1,p}_{(g,h),(n,\vec{m})}
        (W[k],L[k]\,|\,[\beta],\vec{\gamma},\mu)^{W[k]/B[k]}$;
 then
  $\Aut(\Sigma)\times {\Bbb G}_m[k]$ now does approximate-pseudo-act
   on $\widetilde{V}_{\rho}^{\,\prime\,\sharp}$.
 The rigidifying map
  $R_{\rho}: \widetilde{V}_{\rho}\rightarrow
   {\Bbb R}^{a+a^{\prime}}\times{\Bbb C}^{b+b^{\prime}+k}$
  extends canonically to
  $R_{\rho}^{\prime}:\widetilde{V}_{\rho}^{\prime}
    \rightarrow {\Bbb R}^{a+a^{\prime}}\times{\Bbb C}^{b+b^{\prime}+k}$
  since
   the average-weight functions in [Sie1: Sec.~5.3] and Sec.~4.2 that
   constitutes $R_{\rho}$ are well-defined for $\check{W}^{1,p}$-maps
   from deformed $\Sigma$ to fibers of $W[k]/B[k]$.
 In particular, after shrinking $\widetilde{V}_{\rho}$ if necessary,
  $V_{\rho}=\widetilde{V}_{\rho}
   \cap {R_{\rho}^{\prime}}^{-1}(R_{\rho}(0))$.
 Recall the stratification of $\widetilde{V}_{\rho}/B[k]$
  and $\widetilde{V}_{\rho}^{\,\prime\,\sharp}/B[k]$
  induced from the coordinate-subspace stratification of $B[k]$.
 It follows from the three facts:
  (1) $R_{\rho}$, $R_{\rho}^{\prime}$ are continuous, and
      are continuously differentiable when restricted to each stratum,
  (2) the pseudo-action on $\widetilde{V}_{\rho}^{\prime\sharp}$ of
      a small enough neighborhood of the identity element of
      $\Aut(\Sigma)\times{\Bbb G}_m[k]$ is free, and
  (3) $\widetilde{V}_{\rho}$ contains a whole orbit $O_{\rho}$
      (cf.\ Lemma 5.3.5.1),  %  Lemma [$O_{\rho}$ maximal]
 that, for $\widetilde{V}_{\rho}$ small enough,
  $V_{\rho}$ can be interpreted as a stratified space through $0$
  (i.e.\ $\rho$) in $\widetilde{V}_{\rho}$ that, in each strata, is
  transverse to the span of the $(a+a^{\prime}+2b+2b^{\prime}+2k)$-many
  gradient-flow directions from the component weight functions that
  constitute $R_{\rho}$.
 % Thus, $V_{\rho}$ has the same topology as in the case when
 %  the pseudo $\Aut(\Sigma)\times {\Bbb G}_m[k]$ is defined on
 %  the whole $\widetilde{V}_{\rho}$.
 The lemma then follows.

\noindent\hspace{15cm}$\Box$

\bigskip

\begin{flushleft}
{\bf $V_{\rho}/B$ as a Kuranishi neighborhood of
     $\rho\in \overline{\cal M}_{(g,h),(n,\vec{m})}
                       (W/B,L\,|\,[\beta],\vec{\gamma},\mu)/B$.}
\end{flushleft}
To recapitulate, we have constructed
 \begin{itemize}
  \item[$\cdot$]
   ${\cal C}_{V_{\rho}}/V_{\rho}\,$:\hspace{1em}
    \parbox[t]{30em}{an
     $\Aut(\rho)$-equivariant family ${\cal C}_{V_{\rho}}/V_{\rho}$
     of labelled-bordered Riemann surfaces with marked points over
     $V_{\rho}$;}

  \item[$\cdot$]
   $F_{V_{\rho}}:
     ({\cal C}_{V_{\rho}}, \dot{\partial}{\cal C}_{V_{\rho}})/V_{\rho}
                           \rightarrow (W[k],L[k])/B[k]\,$:\\[.6ex]
    a map over $V_{\rho}\rightarrow B[k]$ that satisfies
    $\beta\circ f_{(\zeta,\vec{t},\vec{t}^{\prime},\vec{\mu},
                    \vec{\lambda},\vec{a},\vec{b})}\circ\alpha^{-1}
     =f_{(\alpha,\beta)\cdot(\zeta,\vec{t},\vec{t}^{\prime},
                 \vec{\mu},\vec{\lambda},\vec{a},\vec{b})}$.
 \end{itemize}
Here
 ${\cal C}_{V_{\rho}}/V_{\rho}$ is the pull-back of the family
  ${\cal C}/\Def(\Sigma)$ to $V_{\rho}$ via
  $V_{\rho}\rightarrow \Def(\Sigma)$ from the construction,
 $\dot{\partial}{\cal C}_{V_{\rho}}$ is the labelled boundary
  of ${\cal C}_{V_{\rho}}$ relative to $V_{\rho}$,  and
 $F_{V_{\rho}}|_{(\zeta,\vec{t},\vec{t}^{\prime},\vec{\mu},
                               \vec{\lambda},\vec{a},\vec{b})}
   =f_{(\zeta,\vec{t},\vec{t}^{\prime},\vec{\mu},
                      \vec{\lambda},\vec{a},\vec{b})}$.
Through the construction, $V_{\rho}$ is equipped with the following data:
 \begin{itemize}
  \item[$\cdot$]
   $\Gamma_{V_{\rho}}=\Aut(\rho)$ that acts on $V_{\rho}$,

  \item[$\cdot$]
   $E_{V_{\rho}}\,$,
    the $\Gamma_{V_{\rho}}$-equivariant bundle on
    $(V_{\rho},\Gamma_{V_{\rho}})$ from the restriction of
    $E_{\widetilde{V}_{\rho}}$ to $V_{\rho}$,

  \item[$\cdot$]
   $s_{\rho}: V_{\rho}\rightarrow E_{\rho}$
    from the operator $\bar{\partial}_J$, and

 \item[$\cdot$]
  $\psi_{\rho}:s_{\rho}^{-1}(0) \rightarrow
    \overline{\cal M}_{(g,h),(n,\vec{m})}
             (W/B,L\,|\,[\beta],\vec{\gamma},\mu)$,
  the map over $B$ that sends each pre-deformable $J$-holomorphic map
   $f_{(\zeta,\vec{t},\vec{t}^{\prime},
             \vec{\mu},\vec{\lambda},\vec{a},\vec{b})}$,
    $(\zeta,\vec{t},\vec{t}^{\prime},
        \vec{\mu},\vec{\lambda},\vec{a},\vec{b})
                        \in s_{\rho}^{-1}(0)\subset V_{\rho}$,
   to its isomorphism class
    $[f_{(\zeta,\vec{t},\vec{t}^{\prime},
                \vec{\mu},\vec{\lambda},\vec{a},\vec{b})}]$  in
   $\overline{\cal M}_{(g,h),(n,\vec{m})}
         (W/B,L\,|\,[\beta],\vec{\gamma},\mu)/B$.
\end{itemize}

\bigskip

\noindent
{\bf Proposition 5.3.5.3 [$V_{\rho}$ Kuranishi neighborhood].} {\it
 The $5$-tuple
  $(V_{\rho}, \Gamma_{V_{\rho}}, E_{V_{\rho}}; s_{\rho},\psi_{\rho})$
  forms a Kuranishi neighborhood-in-${\cal C}_{\spsccw}$ of
  $\rho\in \overline{\cal M}_{(g,h),(n,\vec{m})}
           (W/B,L\,|\,[\beta],\vec{\gamma},\mu)/B$.
} % end-corollary

\bigskip

\noindent{\it Proof.}
 That $V_{\rho}/B$ is an object in the category ${\cal C}_{\spsccw}$
  follows from Lemma 5.3.5.2.
    % Lemma [piecewise-transverse slice at $\rho$]
 Injectivity of the  $\psi_{\rho}$-induced map
  $s_{\rho}^{-1}(0)/\Gamma_{V_{\rho}}\rightarrow
   \overline{\cal M}_{(g,h),(n,\vec{m})}
             (W/B,L\,|\,[\beta],\vec{\gamma},\mu)$
  follows from rigidification.
 To show that the image of $\psi_{\rho}$ contains
  a neighborhood of $\rho$ in
  $\overline{\cal M}_{(g,h),(n,\vec{m})}
                     (W/B,L\,|\,[\beta],\vec{\gamma},\mu)/B$,
 let $\tilde{s}_{\rho}$ be the section of $E_{\widetilde{V}_{\rho}}$
  associated to the operator $\bar{\partial}_J$.
 By construction, $s_{\rho}^{-1}(0)$ is the rigidification of
  $\tilde{s}_{\rho}^{-1}(0)$ by $R_{\rho}$
   with respect to the approximate
   pseudo-$(\Aut(\Sigma)\times{\Bbb G}_m[k])$-action
   on $\widetilde{s}^{-1}(0)$, and
  there is a (continuous) map
  $\tilde{\psi}_{\rho}: \tilde{s}_{\rho}^{-1}(0)
    \rightarrow
      \overline{\cal M}_{(g,h),(n,\vec{m})}
               (W/B,L\,|\,[\beta],\vec{\gamma},\mu)$  over $B$
  that sends each
   $f_{(\zeta,\vec{t},\vec{t}^{\prime},
            \vec{\mu},\vec{\lambda},\vec{a},\vec{b})}$,
   $(\zeta,\vec{t},\vec{t}^{\prime},
       \vec{\mu},\vec{\lambda},\vec{a},\vec{b})
        \in \tilde{s}_{\rho}^{-1}(0)\subset \widetilde{V}_{\rho}$,
   to its isomorphism class
   $[f_{(\zeta,\vec{t},\vec{t}^{\prime},
               \vec{\mu},\vec{\lambda},\vec{a},\vec{b}}]$  in
   $\overline{\cal M}_{(g,h),(n,\vec{m})}
             (W/B,L\,|\,[\beta],\vec{\gamma},\mu)$.
 We will show that the image of $\tilde{\psi}_{\rho}$
  contains a neighborhood of $\rho$ in
  $\overline{\cal M}_{(g,h),(n,\vec{m})}
                     (W/B,L\,|\,[\beta],\vec{\gamma},\mu)/B$.
 This then implies the same for $\psi_{\rho}$.

 Recall the standard piecewise-continuous section
  $S:\pi_{\scriptsizeDef(\Sigma)\times B[k]}(\widetilde{V}_{\rho})
                                    \rightarrow \widetilde{V}_{\rho}$
  of the fibration
  $\pi_{\scriptsizeDef(\Sigma)\times B[k]}:
     \widetilde{V}_{\rho} \rightarrow \Def(\Sigma)\times B[k]$.
 For a fixed
  $(\zeta,\vec{t},\vec{t}^{\prime},\vec{\mu};\vec{\lambda})
   \in \pi_{\scriptsizeDef(\Sigma)\times B[k]}(\widetilde{V}_{\rho})
            \subset \Def(\Sigma)\times B[k]$,
 let
  $$
   W^{1,p}
    \left(
     (\Sigma_{(\zeta, \vec{t},\vec{t}^{\prime},\vec{\mu})},
       \partial\Sigma_{(\zeta,\vec{t},\vec{t}^{\prime},\vec{\mu})})\,,\,
     (W[k]_{\vec{\lambda}},L)
    \right)
  $$
  be the Banach manifold of $W^{1,p}$-maps from
  $(\Sigma_{(\zeta, \vec{t},\vec{t}^{\prime},\vec{\mu})},
     \partial\Sigma_{(\zeta,\vec{t},\vec{t}^{\prime},\vec{\mu})})$
  to (the {\it rigid}) $(W[k]_{\vec{\lambda}},L)$.
 Then the same transversality argument as in the proof of
  Lemma~5.3.5.1  % Lemma [$O_{\rho}$ maximal]
  implies that
   $\widetilde{V}_{\rho}$ contains all pre-deformable
    $(J,E_{\tinybullet})$-holomorphic $W^{1,p}$-maps from
    $(\Sigma_{(\zeta, \vec{t},\vec{t}^{\prime},\vec{\mu})},
       \partial\Sigma_{(\zeta,\vec{t},\vec{t}^{\prime},\vec{\mu})})$
    to $(W[k]_{\vec{\lambda}},L)$
    that are close to
    $f_{S(\zeta,\vec{t},\vec{t}^{\prime},\vec{\mu};\vec{\lambda})}$.

 Let $(\zeta,\vec{t},\vec{t}^{\prime},\vec{\mu};\vec{\lambda})$
  now vary in
  $\pi_{\scriptsizeDef(\Sigma)\times B[k]}(\widetilde{V}_{\rho})
   \subset \Def(\Sigma)\times B[k]$.
 Note that
  the fiber-dimension of $\widetilde{V}_{\rho}$ over
   $\pi_{\scriptsizeDef(\Sigma)\times B[k]}(\widetilde{V}_{\rho})$
   is upper semi-continuous  and that
  there is a well-defined flattening stratification on
   $\pi_{\scriptsizeDef(\Sigma)\times B[k]}(\widetilde{V}_{\rho})$
   so that the restriction of the fibration
    $\widetilde{V}_{\rho}/
      \pi_{\scriptsizeDef(\Sigma)\times B[k]}(\widetilde{V}_{\rho})$
    to each stratum is a bundle whose fibers do not shrink or get pinched
   when moving toward the boundary of the stratum.
  Together with the conclusion of the previous paragraph,
   these imply that $\widetilde{V}_{\rho}$ contains
   all pre-deformable $(J,E_{\tinybullet})$-holomorphic $W^{1,p}$-maps
   that are close to some
   $f_{S(\zeta,\vec{t},\vec{t}^{\prime},\vec{\mu};\vec{\lambda})}$,
   $(\zeta,\vec{t},\vec{t}^{\prime},\vec{\mu};\vec{\lambda})  \in
     (\pi_{\scriptsizeDef(\Sigma)\times B[k]}(\widetilde{V}_{\rho}))$.
 In particular, $\widetilde{V}_{\rho}$ contains all
  ($J$-holomorphic, pre-deformable) stable maps near $f$.
 This concludes the proof.

\noindent\hspace{15cm}$\Box$

\bigskip

\noindent
{\it Remark 5.3.5.4 {\rm [}$E_{\rho}$-dependence of $V_{\rho}${\rm ]}.}
Different choices of $E_{\rho}$ in Definition/Lemma 5.3.1.5
  % Definition/Lemma [saturated obstruction space]
 give rise to different but equivalent family Kuranishi
 neighborhoods of $\rho$ in the sense of Definition 5.1.1.
  % %
  % \marginpar{\raggedright\tiny\vspace{-6em} $\bullet$
  %   Add a definition in Sec.~???. \newline
  %    Cf.\ [Liu(C): Definition 6.2]. \newline\newline
  %   Explain the notion of isomorphic, sub-, ... system of PDE's
  %    to relate Kuranishi neighborhoods $V_{1,\rho}$, $V_{2,\rho}$,
  %     and $V_{\rho}$ from $E_{1,\rho}$, $E_{2,\rho}$, and
  %     $E_{1.,\rho}+E_{2,\rho}$ respectively.
  %   CAUTION that in general $V_{1,\rho}$ and $V_{2,\rho}$ do not
  %    intersect properly in $V_{\rho}$.}
  % %
 E.g.\ taking $E_{1, \rho}+E_{2,\rho}$ creates a third family
  Kuranishi neighborhood of $\rho$ that dominates both
  $V_{1,\rho}$ and $V_{2,\rho}$, as in [Liu(C): Remark 6.34].
 %%%%%%%%%%%%%%%%%%%%%%%
 % Thus, varying $E_{\rho}$ in
 %  $L^p(\Sigma;
 %     \Lambda^{0,1}\Sigma\otimes_J f^{\ast}T_{\ast}Y_{[k]})$
 %  changes $V_{\rho}$ accordingly in the same equivalence
 % class of family Kuranishi neighborhoods of $\rho$.
 % For a technical necessity, we will describe in the next subsection
 %  a particular kind of small variations that move a given $E_{\rho}$
 %  totally into the complement of a specified subset in
 %  $L^p(\Sigma;
 %    \Lambda^{0,1}\Sigma\otimes_J f^{\ast}T_{\ast}Y_{[k]})$
 %  while keeping the corresponding variations of
 %  $V_{\rho}$ and
 %  $U_{\rho}=\psi_{\rho}(s_{\rho}^{-1}(0))$ small.
 %%%%%%%%%%%%%%%
%
% \bigskip
%
% \noindent
% {\it Remark ??? $[$generalization to multi-centered construction$]$.}
%
% \marginpar{\raggedright\tiny $\bullet$
%  To be completed after finishing the next subsection so that
%   we know better what we may still miss.}
%
% \bigskip
% \noindent $\bullet$
% ?????????????????????????.
% \bigskip
%%%%%%%%%%%%%%%%%%%%%%%%%%%%%%%%%%%%

\bigskip

\subsection{Construction of a family Kuranishi structure.}

We now proceed to construct a family Kuranishi structure on
  $\overline{\cal M}_{(g,h),(n,\vec{m})}
      (W/B,L\,|\,[\beta],\vec{\gamma},\mu)/B$
by relating Kuranishi neighborhoods on
 $\overline{\cal M}_{(g,h),(n,\vec{m})}
     (W/B,L\,|\,[\beta],\vec{\gamma},\mu)/B$
 with sub-fibrations of the $\check{L}^p$-obstruction-space fibration
 $T^2_{\check{\cal W}^{1,p}_{\tinybullet}
          ((\widehat{W},\widehat{L})/\widehat{B}\,|\,\tinybullet)
            /\widetilde{\cal M}_{\tinybullet}}$.
The construction connects Fukaya-Ono's construction in [F-O: Sec.~15]
 with Siebert's construction in [Sie1: Sec.~5 - Sec.~6].
%%%%%%%%%%%%%%%%%%%%%%%%%%%%%
 % ?????????????????????????????
% by making a judicious choice of a family Kuranishi neighborhood
% $V_{\rho}$ at each
%  $\rho \in
%    \overline{\cal M}_{(g,h),(n,\vec{m})}
%       (\\  W/B,L\,|\,[\beta],\vec{\gamma},\mu)/B$
% so that the transition data among them follow immediately.
% %
% \marginpar{\raggedright\tiny\vspace{-3em} $\bullet$ {\bf CAUTION.}
%   Synchronizing $E_{\rho}$'s is {\it enough} to make
%    the construction of transition data immediate.
%   Need to deal with choices of slice when $\Sigma_{\rho}$ has
%    unstable components. \newline
%   To be revised. }
% %
% The construction presented here is guided by [F-O: Sec.~15].
% However, the technical discussion ibidem on rigidifying isomorphisms
%  when one compares different representatives of an isomorphism class
%  of maps is replaced by [Sie1: Sec.~5], extended to our case in Sec.~???.
%
% As gluing involves comparing/relating points in different Kuranishi
%  neighborhoods, it is important to keep the following distinction
%  in mind before we start:
%%%%%%%%%%%%%%%%%%%%%%%%%%%%%%%%

The following remark should be kept in mind as it is everywhere
 behind the discussion.

\bigskip

\noindent
{\it Remark 5.4.1 $[$isomorphism class vs.\ representative$]$.}
 A point $\rho$ in the moduli space
  $\overline{\cal M}_{(g,h),(n,\vec{m})}
            (\\  W/B,L\,|\,[\beta],\vec{\gamma},\mu)/B$
  represents an {\it isomorphism class} of maps
   % %
   % \marginpar{\raggedright\tiny $\bullet$
   %    Some definition in early section that defines the moduli space.
   %    \newline To be completed.}
  while a family Kuranishi neighborhood
   $(V_{\rho},\Gamma_{V_{\rho}},E_{V_{\rho}}; s_{\rho},\psi_{\rho})$
   of $\rho$, as constructed in Sec.~5.3, parameterizes a collection
   of {\it maps} that contains a sub-collection,
   namely $s_{\rho}^{-1}(0)$, of
   {\it representatives}
   $f_{(\zeta,\vec{t},\vec{t}^{\prime},\vec{\mu},
                      \vec{\lambda},\vec{a},\vec{b})}$,
   $(\zeta,\vec{t},\vec{t}^{\prime},\vec{\mu},
                   \vec{\lambda},\vec{a},\vec{b})\in s_{\rho}^{-1}(0)$,
   as in Sec.~5.3.5, whose corresponding set of isomorphism classes
   covers a neighborhood,
   namely $U_{\rho}:=\psi_{\rho}(s_{\rho}^{-1}(0))$, of $\rho$ in
   $\overline{\cal M}_{(g,h),(n,\vec{m})}
            (W/B,L\,|\,[\beta],\vec{\gamma},\mu)/B$
   via $\psi_{\rho}$.
 In particular, every $p\in V_{\rho}$ goes with a unique
  representative $h_p:\Sigma_p/\pt\rightarrow W[k_{\rho}]/B[k_{\rho}]$.
 The set of isomorphisms from a representative $f_1$ to another
  representative $f_2$ of $\rho$
  (which may come from two different Kuranishi neighborhoods
   $V_{\rho_1}$ and $V_{\rho_2}$ that cover $\rho$)
  is parameterized by $\Aut(f_1)$ up to a right multiplication
  and by $\Aut(f_2)$ up to a left multiplication.
 By definition,
  $\Aut(f_1)\simeq \Aut(f_2)\simeq \Aut(\rho)=\Gamma_{\rho}$.
 The same distinction holds between
  points on the moduli space
   $\check{W}^{1,p}_{(g,h),(n,\vec{m})}
     ((\widehat{W},\widehat{L})/\widehat{B}\,|\,
                             [\beta],\vec{\gamma},\mu)$  and
  points on its (singular) orbifold local charts.

\bigskip

\begin{flushleft}
{\bf Kuranishi neighborhoods in terms of
     $T^2_{\check{\cal W}^{1,p}_{\tinybullet}
              ((\widehat{W},\widehat{L})/\widehat{B}\,|\,\tinybullet)
                /\widetilde{\cal M}_{\tinybullet}}$.}
\end{flushleft}
Note that the same construction in Sec.~5.3 works also with
 $W^{1,p}$ replaced by $\check{W}^{1,p}$ and
 $L^p$ replaced by $\check{L}^p$.
Let
 $\rho\in \overline{\cal M}_{(g,h),(n,\vec{m})}
                     (W/B,L\,|\,[\beta],\vec{\gamma},\mu)$
 be represented by $f:\Sigma\rightarrow (Y_{[k]},L_{[k]})$,
Then, in terms of the fibration
 $T^2_{\check{\cal W}^{1,p}_{\tinybullet}
          ((\widehat{W},\widehat{L})/\widehat{B}\,|\,\tinybullet)
            /\widetilde{\cal M}_{\tinybullet}}$,
 the construction of a Kuranishi neighborhood
 $(V_{\rho},\Gamma_{V_{\rho}},E_{V_{\rho}};s_{\rho},\psi_{\rho})$
 of $\rho$ in Sec.~5.3 can be deformed and rephrased as follows:
\begin{itemize}
 \item[(1)]
  Choose a saturated obstruction space
   $E_{\rho}
    \subset
    C^{\infty}
     (\Sigma; \Lambda^{0,1}\Sigma\otimes_J f^{\ast}T_{\ast}Y_{[k]})$
   at $\rho$.
  Regard $\rho$ as a point in
   $\check{W}^{1,p}_{(g,h),(n,\vec{m})}
    ((\widehat{W},\widehat{L})/\widehat{B}\,|\,
                                        [\beta],\vec{\gamma},\mu)$
   that is represented also by $f$
   gives an embedding
   $E_{\rho}
    \hookrightarrow
     \check{L}^p
     (\Sigma; \Lambda^{0,1}\Sigma\otimes_J f^{\ast}T_{\ast}Y_{[k]})$.

  \item[$\cdot$]
   Extend $E_{\rho}$ at $\rho$ to a trivialized
    $\Aut(\rho)$-equivariant trivial bundle $E_{\check{V}_{\rho}}$
    over a sufficiently small orbifold local chart
    $\check{V}_{\rho}$ of $\rho$ in
    $\check{W}^{1,p}_{(g,h),(n,\vec{m})}
     ((\widehat{W},\widehat{L})/\widehat{B}\,|\,
                                        [\beta],\vec{\gamma},\mu)$
    such that,
     for all $p\in \check{V}_{\rho}$ with its corresponding
       representative $h_p$ is $J$-holomorphic,
     the fiber $E_{\check{V}_{\rho}}|_p$ is
     a saturated obstruction space
      $\subset  C^{\infty}
                  (\Sigma_p; \Lambda^{0,1}\Sigma_p \otimes_J
                    h^{\ast}T_{\ast}W[k_{\rho}]_{\vec{\lambda}_p})$.
   By construction there is a map
    $$
     E_{\check{V}_{\rho}}\;
      \longrightarrow\;
     T^2_{\check{\cal W}^{1,p}_{\tinybullet}
          ((\widehat{W},\widehat{L})/\widehat{B}\,|\,\tinybullet)
             /\widetilde{\cal M}_{\tinybullet}}
    $$
    as an orbifold sub-fibration;
   we shall think of $(E_{\check{V}_{\rho}},\Aut(\rho))$
    equally as a sub-orbifold of
    $T^2_{\check{\cal W}^{1,p}_{\tinybullet}
         ((\widehat{W},\widehat{L})/\widehat{B}\,|\,\tinybullet)
                               /\widetilde{\cal M}_{\tinybullet}}$.

 \item[(2)]
  Recall the global section
   $$
    s_{\bar{\partial}_J}\,:
    \check{\cal W}^{1,p}_{(g,h),(n,\vec{m})}
     ((\widehat{W},\widehat{L})/\widehat{B}\,
       |\,[\beta],\vec{\gamma},\mu)\;
    \longrightarrow\;
    T^2_{\check{\cal W}^{1,p}_{\tinybullet}
         ((\widehat{W},\widehat{L})/\widehat{B}\,|\,\tinybullet)
           /\widetilde{\cal M}_{\tinybullet}}
   $$
   of  $T^2_{\check{\cal W}^{1,p}_{\tinybullet}
          ((\widehat{W},\widehat{L})/\widehat{B}\,|\,\tinybullet)
            /\widetilde{\cal M}_{\tinybullet}}$
   as a morphism of orbifolds.
  Denote its image sub-orbifold in
   $T^2_{\check{\cal W}^{1,p}_{\tinybullet}
         ((\widehat{W},\widehat{L})/\widehat{B}\,|\,\tinybullet)
           /\widetilde{\cal M}_{\tinybullet}}$
   by $\Image(s_{\bar{\partial}_J})$.
  Let
   $$
    \pi^2\;:\;
       T^2_{\check{\cal W}^{1,p}_{\tinybullet}
       ((\widehat{W},\widehat{L})/\widehat{B}\,|\,\tinybullet)
           /\widetilde{\cal M}_{\tinybullet}}\;
     \longrightarrow\;
     \check{W}^{1,p}_{(g,h),(n,\vec{m})}
       ((\widehat{W},\widehat{L})/\widehat{B}\,|\,
                                    [\beta],\vec{\gamma},\mu)
   $$
   be the fibration orbifold-map.
  Then, on the orbifold local chart $\check{V}_{\rho}$,
   $$
    V_{\rho}\; :=\;
     \pi^2\,\left(\, \Image(s_{\bar{\partial}_J})\,
                         \cap\,E_{\check{V}_{\rho}}  \,\right)
   $$
   is $\Aut(\rho)$-invariant.
  Furthermore, $V_{\rho}$ defines
   a Kuranishi neighborhood-in-${\cal C}_{\spsccw}$
   of $\rho \in \overline{\cal M}_{(g,h),(n,\vec{m})}
                    (W/B,L\,|\,[\beta],\vec{\gamma},\mu)$
   with
    $E_{V_{\rho}}=E_{\check{V}_{\rho}}|_{V_{\rho}}$,
    $\Gamma_{V_{\rho}}=\Aut(\rho)$ now acting on
     $E_{V_{\rho}}/V_{\rho}$ equivariantly,
    $s_{\rho}= s_{\bar{\partial}_J}|_{V_{\rho}}$, and
    $\psi_{\rho}: s_{\rho}^{-1}(0)
      \rightarrow \overline{\cal M}_{(g,h),(n,\vec{m})}
                         (W/B,L\,|\,[\beta],\vec{\gamma},\mu)$
    by sending $p\in s_{\rho}^{-1}(0)$ to $[h_p]$.

  \item[$\cdot$]
   On $\check{V}_{\rho}$ it follows by construction that
    $\bar{\partial}_J h_p\in E_{V_{\rho}}|_p$
     if and only if $p\in V_{\rho}$.
   Thus, $V_{\rho}$ parameterizes all the $(J,E)$-holomorphic
    $\check{W}^{1,p}$-maps near $\rho$.
   Indeed it parameterizes also all the $(J,E)$-holomorphic
    $W^{1,p}$-maps near $\rho$.
\end{itemize}

Deformations of the bundle $E_{\check{V}_{\rho}}$ in
 $T^2_{\check{\cal W}^{1,p}_{\tinybullet}
       ((\widehat{W},\widehat{L})/\widehat{B}\,|\,\tinybullet)
           /\widetilde{\cal M}_{\tinybullet}}$
 as orbifold sub-fibrations
  without violating
   the $C^{\infty}$-class and the saturatedness condition
   on the locus $(s_{\bar{\partial}})^{-1}(0)$
 give rise to Kuranishi neighborhoods-in-${\cal C}_{\spsccw}$
  of $\rho$,
 all of the same actual dimension and
  the same virtual dimension.\footnote{This deformation freedom is
                  crucial in the construction of a Kuranishi structure on
                   $\overline{\cal M}_{(g,h),(n,\vec{m})}
                                  (W/B,L\,|\,[\beta],\vec{\gamma},\mu)/B$.
                  In Sec.~5.3, $E_{\check{V}_{\rho}}$ is constructed
                   via parallel transport from the trivialized trivial
                    bundle on
                    $S_0(\pi_{\scriptsizeDef(\Sigma)\times B[k]}
                                                    (\Theta_{\rho,0}))$.
                  Such parallel transport construction depends on
                   the metric on the fibers of $W[k_{\rho}]/B[k_{\rho}]$.
                  The curvature of the metric makes the bundle
                   $E_{\check{V}_{\rho^{\prime}}}$ constructed from nearby
                   $\rho^{\prime}\in
                       \overline{\cal M}_{(g,h),(n,\vec{m})}
                                   (W/B,L\,|\,[\beta],\vec{\gamma},\mu)$
                   distinct on
                   $\Image(\psi_{\rho})\cap \Image(\psi_{\rho^{\prime}})$.
                  The corresponding $V_{\rho}$ and $V_{\rho}^{\prime}$
                   for such $E_{\check{V}_{\rho}}$ and
                   $E_{\check{V}_{\rho^{\prime}}}$ cannot be glued
                    at the level of the universal map on the universal
                    curve.
                  Furthermore, while the almost-complex structure on
                   fibers of $\widehat{W}/\widehat{B}$ is well-defined,
                   the metric is not.
                  So a deformation to the construction in Sec.~5.3
                   that preserves the $C^{\infty}$-class and
                   the saturatedness condition is indispensable.}

\bigskip

\noindent
{\bf Definition 5.4.2 [saturated obstruction local bundle].} {\rm
 The $\Aut(\rho)$-equivariant bundle $E_{\check{V}_{\rho}}$ on
  $\check{V}_{\rho}$ in the above rephrasing,
  with
   the prescribed properties and
   the orbifold sub-fibration map
    $E_{\check{V}_{\rho}} \rightarrow
      T^2_{\check{\cal W}^{1,p}_{\tinybullet}
          ((\widehat{W},\widehat{L})/\widehat{B}\,|\,\tinybullet)
             /\widetilde{\cal M}_{\tinybullet}}$,
  is called
  a {\it saturated obstruction local bundle}
  on  $\check{W}^{1,p}_{(g,h),(n,\vec{m})}
          ((\widehat{W},\widehat{L})/\widehat{B}\,|\,
                                        [\beta],\vec{\gamma},\mu)$.
 The local orbifold chart $\check{V}_{\rho}$ is called
  the {\it support} of $E_{\check{V}_{\rho}}$.
 The tuple
  $$
   V_{\rho}(E_{\check{V}_{\rho}}) \;  :=\;
    (V_{\rho},\Gamma_{V_{\rho}},E_{V_{\rho}};s_{\rho},\psi_{\rho})\,,
  $$
  (also denoted by $V_{\rho}$ in shorthand),
  in the rephrasing is called {\it the Kuranishi neighborhood of
  $\rho\in \overline{\cal M}_{(g,h),(n,\vec{m})}
                     (W/B,L\,|\,[\beta],\vec{\gamma},\mu)$
  determined by} $E_{\check{V}_{\rho}}$.
} % end-definition

\bigskip

\begin{flushleft}
{\bf Kuranishi structures associated to
                   a fine system of local bundles.}
\end{flushleft}
%
%\noindent
{\bf Definition 5.4.3 [direct-sum/fine system of local bundles].}
{\rm
 A collection $\{E_{\check{V}_{\rho_i}}\}_{i\in I}$,
  $\rho_i \in \overline{\cal M}_{(g,h),(n,\vec{m})}
                     (W/B,L\,|\,[\beta],\vec{\gamma},\mu)$,
  of saturated obstruction local bundles is said to form
  a {\it direct-sum system} for
  $\overline{\cal M}_{(g,h),(n,\vec{m})}
             (W/B,L\,|\,[\beta],\vec{\gamma},\mu)$
 if the following two conditions are satisfied:
  \begin{itemize}
   \item[(1)]
    $\{\Image(\psi_{\rho_i})\}_{i\in I}$
     is a locally finite (open) cover of
     $\overline{\cal M}_{(g,h),(n,\vec{m})}
                (W/B,L\,|\,[\beta],\vec{\gamma},\mu)/B$
     that is finite over a compact subset of $B$,
    (here $\psi_{\rho_i}$ is from the Kuranishi neighborhood data
     associated to $E_{\check{V}_{\rho}}$);

   \item[(2)]
    the span of $\{E_{\check{V}_{\rho_i}}\}_{i\in I}$
     in each vector-space fiber of a fibration local chart of\\
     $T^2_{\check{\cal W}^{1,p}_{\tinybullet}
            ((\widehat{W},\widehat{L})/\widehat{B}\,|\,\tinybullet)
            /\widetilde{\cal M}_{\tinybullet}}$   is a direct sum
     of the related fibers of $E_{\check{V}_{\rho_i}}$'s.
  \end{itemize}
 $\{E_{\check{V}_{\rho_i}}\}_{i\in I}$ is said to be {\it fine}
  if, in addition,
  \begin{itemize}
   \item[(3)]
    there exists an open cover $\{U_{\rho_i}^{\flat}\}_{i\in I}$  of
     $\overline{\cal M}_{(g,h),(n,\vec{m})}
                 (W/B,L\,|\,[\beta],\vec{\gamma},\mu)$
     such that
      $U_{\rho_i}^{\flat}$ is an open neighborhood of $\rho_i$
      with the closure
      $\overline{U_{\rho_i}^{\flat}}$ a compact subset of
      $\Image(\psi_{\rho_i})$.
  \end{itemize}
} % end-definition

\bigskip

\bigskip

\noindent
{\bf Lemma 5.4.4 [existence of fine system].} {\it
 A fine system of saturated obstruction local bundles for
  $\overline{\cal M}_{(g,h),(n,\vec{m})}
               (W/B,L\,|\,[\beta],\vec{\gamma},\mu)$
  exists.
} % end-lemma

\bigskip

\noindent
{\it Proof.}
 Let $E_{\check{\rho}}^{\prime}$ be a saturated obstruction
  local bundle at
  $\rho\in \overline{\cal M}_{(g,h),(n,\vec{m})}
                  (W/B,L\,|\,[\beta],\vec{\gamma},\mu)$
  and $U_{\rho}^{\flat}$ be an open neighborhood of $\rho$
   with the closure
   $\overline{U_{\rho}^{\flat}}$ a compact subset of
   $\Image(\psi_{\rho}^{\prime})$.
 Since $\overline{\cal M}_{(g,h),(n,\vec{m})}
          (W/B,L\,|\,[\beta],\vec{\gamma},\mu)/B$
   is compact over a compact subset of $B$,
  one can choose a subcover $\{U_{\rho_i}\}_{i\in I}$ of
   $\{U_{\rho}\}_{\rho}$ that is locally finite and
   is finite over a compact subset of $B$.
 We may assume that each $E_{\check{V}_{\rho}}^{\prime}$,
  and hence $E_{\check{V}_{\rho_i}}^{\prime}$,
  is constructed as in Sec.~5.3 so that elements in the fiber of
  $E_{\check{V}_{\rho}}^{\prime}$ are sections of sheaves
  supported away from the nodes of bordered Riemann surfaces.
 As this is a locally finite system of trivial bundles,
 the direct-sum condition can be achieved by
  deforming $E_{\check{V}_{\rho_i}}^{\prime}$ inductively
   to another equivariant $E_{\check{V}_{\rho_i}}$ that satisfies
   also the $C^{\infty}$-class and the saturatedness conditions,
   and with the same support,
  as sub-fibrations in
   $T^2_{\check{\cal W}^{1,p}_{\tinybullet}
         ((\widehat{W},\widehat{L})/\widehat{B}\,|\,\tinybullet)
         /\widetilde{\cal M}_{\tinybullet}}$.
 This makes the $\psi_{\rho_i}$ from $E_{\check{V}_{\rho_i}}$
  coincides with the $\psi_{\rho_i}^{\prime}$ from
  $E_{\check{V}_{\rho}}^{\prime}$.
 Thus the cover $\{\Image(\psi_{\rho_i})\}_{i\in I}$  of
   $\overline{\cal M}_{(g,h),(n,\vec{m})}
             (W/B,L\,|\,[\beta],\vec{\gamma},\mu)$
   from the new system $\{E_{\check{V}_{\rho_i}}\}_{i\in I}$
   coincides with $\{\Image(\psi_{\rho_i}^{\prime})\}_{i\in I}$
  and, hence, Condition (2) and Condition (3) in Definition 5.4.3
          % Definition [direct-sum/fine system of local bundles]
   are also satisfied.

\noindent\hspace{15cm}$\Box$

\bigskip

Recall the canonical orbifold-embedding
 $$
  \overline{\cal M}_{(g,h), (n, \vec{m})}
                (W/B,L\,|\,[\beta],\vec{\gamma},\mu)\;
   \hookrightarrow\;
   \check{\cal W}^{1,p}_{(g,h),(n,\vec{m})}
         ((\widehat{W},\widehat{L})/\widehat{B}\,|\,
                                  [\beta],\vec{\gamma},\mu)\,.
 $$
Let ${\cal E}:=\{E_{\check{V}_{\rho_i}}\}_{i\in I}$
 be a fine system of saturated obstruction local bundles for
  $\overline{\cal M}_{(g,h),(n,\vec{m})}
             (W/B,\\  L\,|\,[\beta],\vec{\gamma},\mu)$.
Let ${\mathbf F}({\cal E})$ be the fiberwise linear span
 of the union of the image set of\\
 $E_{\rho_i}\rightarrow
   T^2_{\check{\cal W}^{1,p}_{\tinybullet}
              ((\widehat{W},\widehat{L})/\widehat{B}\,|\,\tinybullet)
              /\widetilde{\cal M}_{\tinybullet}}$
 with the induced subset topology.
The orbifold structure on\\
 $T^2_{\check{\cal W}^{1,p}_{\tinybullet}
             ((\widehat{W},\widehat{L})/\widehat{B}\,|\,\tinybullet)
             /\widetilde{\cal M}_{\tinybullet}}$
 induces an orbifold structure on ${\mathbf F}({\cal E})$
 that fibers over\\
 $\check{\cal W}^{1,p}_{(g,h),(n,\vec{m})}
  ((\widehat{W},\widehat{L})/\widehat{B}
                               |\,[\beta],\vec{\gamma},\mu)$.
This realizes ${\mathbf F}({\cal E})$ as an orbifold sub-fibration
 of \\  $T^2_{\check{\cal W}^{1,p}_{\tinybullet}
              ((\widehat{W},\widehat{L})/\widehat{B}\,|\,\tinybullet)
              /\widetilde{\cal M}_{\tinybullet}}$
 that is mapped to a neighborhood of
 $\overline{\cal M}_{(g,h), (n, \vec{m})}
                (W/B,L\,|\,[\beta],\vec{\gamma},\mu)$
 in  $\check{\cal W}^{1,p}_{(g,h),(n,\vec{m})}
         ((\widehat{W},\widehat{L})/\widehat{B}\,|\,
                                  [\beta],\vec{\gamma},\mu)$
 under
  $$
   \pi^2\; :\;
      T^2_{\check{\cal W}^{1,p}_{\tinybullet}
      ((\widehat{W},\widehat{L})/\widehat{B}\,|\,\tinybullet)
          /\widetilde{\cal M}_{\tinybullet}}\;
    \longrightarrow\;
    \check{W}^{1,p}_{(g,h),(n,\vec{m})}
      ((\widehat{W},\widehat{L})/\widehat{B}\,|\,
                                   [\beta],\vec{\gamma},\mu)\,.
  $$
The map
  $$
    \check{\cal W}^{1,p}_{(g,h),(n,\vec{m})}
          ((\widehat{W},\widehat{L})/\widehat{B}\,|\,
                                   [\beta],\vec{\gamma},\mu)\
    \longrightarrow\; {\Bbb Z}_{\ge 0}\,,
      \hspace{2em}
    p\; \longmapsto\; \dimm ({\mathbf F}({\cal E})|_p)
  $$
 defines the flattening stratification of
 ${\mathbf F}({\cal E})$ on
    $\check{\cal W}^{1,p}_{(g,h),(n,\vec{m})}
          ((\widehat{W},\widehat{L})/\widehat{B}\,|\,
                              [\beta],\vec{\gamma},\mu)$
 by its preimage subsets.
Over each stratum,
 ${\mathbf F}({\cal E})|
   _{\overline{\cal M}_{\tinybullet}(W/B,L\,|\,\tinybullet)}$
 is an orbi-bundle.
For any $I^{\prime}\subset I$, the same construction applied to
 ${\cal E}_{I^{\prime}}
                 =\{E_{\check{V}_{\rho_i}}\}_{i\in I^{\prime}}$
 gives an orbifold sub-fibration
  ${\mathbf F}({\cal E}_{I^{\prime}})$
   in $T^2_{\check{\cal W}^{1,p}_{\tinybullet}
          ((\widehat{W},\widehat{L})/\widehat{B}\,|\,\tinybullet)
         /\widetilde{\cal M}_{\tinybullet}}$.
The flattening stratification of
 ${\mathbf F}({\cal E}_{I^{\prime}})$ is defined similarly.
By construction, ${\mathbf F}({\cal E}_{I^{\prime}})$
 is an orbifold sub-fibration of ${\mathbf F}({\cal E})$.

Recall the locally finite cover $\{ U_{\rho_i}^{\flat}\}_{i\in I}$ of
 $\overline{\cal M}_{(g,h),(n,\vec{m})}
           (W/B,L\,|\,[\beta],\vec{\gamma},\mu)\,/B$.
This induces a stratification
 ${\cal S}:=\{S_{I^{\prime}}\}_{I^{\prime}\subset I}$ of
 $\overline{\cal M}_{(g,h),(n,\vec{m})}
        (W/B,L\,|\,[\beta],\vec{\gamma},\mu)\,/B$
 by setting
 $$
  S_{I^{\prime}}\;  =\;
   (\cap_{i\in I^{\prime}}\,\overline{U_{\rho_i}^{\flat}})\,
    -\,(\cup_{i\in I-I^{\prime}}\,\overline{U_{\rho_i}^{\flat}} )\,.
 $$
Define also the subset
 $S_{I^{\prime}}^{\prime}
  =(\cap_{i\in I^{\prime}}\,\Image(\psi_{\rho_i})\,
     -\,(\cup_{i\in I-I^{\prime}}\,\overline{U_{\rho_i}^{\flat}} )$,
 $I^{\prime}\subset I$.
For $\rho\in S_{I^{\prime}}$, let $\check{V}_{\rho}$
 be an orbifold local chart of $\rho$ in
 $\check{\cal W}^{1,p}_{(g,h),(n,\vec{m})}
           ((\widehat{W},\widehat{L})/\widehat{B}\,|\,
                                    [\beta],\vec{\gamma},\mu)$
 such that
 \begin{itemize}
  \item[$\cdot$]
   the image of $\check{V}_{\rho}$ in
    $\check{\cal W}^{1,p}_{(g,h),(n,\vec{m})}
     ((\widehat{W},\widehat{L})/\widehat{B}
                                  |\,[\beta],\vec{\gamma},\mu)$
    is covered by the union of the image of
    $\check{V}_{\rho_i}$, $i\in I^{\prime}$,

  \item[$\cdot$]
   the $J$-holomorphy locus of $\check{V}_{\rho}$
    is mapped to $S_{I^{\prime}}^{\prime}$.
 \end{itemize}
Denote the image of
 $\check{V}_{\rho}$ in
    $\check{\cal W}^{1,p}_{(g,h),(n,\vec{m})}
          ((\widehat{W},\widehat{L})/\widehat{B}\,|\,
                                   [\beta],\vec{\gamma},\mu)$
 by $\underline{\check{V}_{\rho}}$.
Then, ${\mathbf F}({\cal E}_{I^{\prime}})$ is an orbi-bundle
 when restricted to $\underline{\check{V}_{\rho}}$.
Let $E_{\check{V}_{\rho}}$ be the associated orbi-bundle
 local chart of
 ${\mathbf F}({\cal E}_{I^{\prime}})|_{\underline{\check{V}_{\rho}}}$;
then, by construction,
 $E_{\check{V}_{\rho}}$ is a saturated obstruction local bundle
 on  $\check{\cal W}^{1,p}_{(g,h),(n,\vec{m})}
           ((\widehat{W},\widehat{L})/\widehat{B}\,|\,
                                    [\beta],\vec{\gamma},\mu)$
 in the sense of Definition 5.4.2.
     % Definition [saturated obstruction local bundle]

In this way, one recovers a family
 $\,\{ E_{\check{V}_{\rho}} \}
    _{\rho\in\overline{\cal M}_{\tinybullet}(W/B,L\,|\,\tinybullet)}\,$
 of saturated obstruction local bundles
 from the orbifold sub-fibration ${\mathbf F}({\cal E})$ of
 $T^2_{\check{\cal W}^{1,p}_{\tinybullet}
         ((\widehat{W},\widehat{L})/\widehat{B}\,|\,\tinybullet)
               /\widetilde{\cal M}_{\tinybullet}}$.
Define
 $$
  {\frak N}^{(0)}_{\Kuranishi}({\cal E})\; :=\;
   \left\{\,
    V_{\rho}(E_{\check{V}_{\rho}})\,=\,
      (V_{\rho},\Gamma_{\rho},E_{\rho};s_{\rho},\psi_{\rho})
      \,\right\}
        _{\rho\in \overline{\cal M}_{\bullet}(W/B,L\,|\,\bullet\,)}
 $$
 from Definition 5.4.2.
      % Definition [saturated obstruction local bundle]
This gives the set of family Kuranishi
 neighborhoods-in-${\cal C}_{\spsccw}$ on\\
 $\overline{\cal M}_{(g,h),(n,\vec{m})}
                   (W/B,L\,|\,[\beta],\vec{\gamma},\mu)/B$.
The orbifold fibration transition data of
 $T^2_{\check{\cal W}^{1,p}_{\tinybullet}
         ((\widehat{W},\widehat{L})/\widehat{B}\,|\,\tinybullet)
               /\widetilde{\cal M}_{\tinybullet}}$,
 or of ${\mathbf F}({\cal E})$,
 induces a collection of $4$-tuples
 \begin{eqnarray*}
  \lefteqn{{\frak N}^{(1)}_{\Kuranishi}({\cal E})\; := } \\[.6ex]
   && \left\{ ( V_{\rho}, h_{\rho^{\prime}\rho},
       \phi_{\rho^{\prime}\rho}, \hat{\phi}_{\rho^{\prime}\rho} )\;
      :\; \rho \in \overline{\cal M}_{(g,h),(n,\vec{m})}
               (W/B,L\,|\,[\beta],\vec{\gamma},\mu)\,/B\,,\;
          \rho^{\prime}\in \psi_{\rho}\left(s_{\rho}^{-1}(0)\right)
  \right\}
 \end{eqnarray*}
 that gives the set of transition functions between elements in
 ${\frak N}^{(0)}_{\Kuranishi}({\cal E})$ in the sense of
 Definition 5.1.2.  % Definition [Kuranishi structure-in-${\cal C}$]
We shall call the pair
 $$
   {\cal K}({\cal E})\;  =\;
     \left(\, {\frak N}^{(0)}_{\Kuranishi}({\cal E})\,,\,
              {\frak N}^{(1)}_{\Kuranishi}({\cal E})\, \right)
 $$
 a {\it Kuranishi structure} associated to the fine system ${\cal E}$
 of saturated obstruction local bundles.
We remark that
 the gluing thus constructed is at the level of the universal map
  on the universal curve  and that
 different choices of $\{U_{\rho_i}^{\flat}\}_{i\in I}$
 give equivalent Kuranishi structures.

To summarize:

\bigskip

\noindent
{\bf Proposition 5.4.5 [Kuranishi structure from fine system].} {\it
 A fine system of saturated obstruction local bundles for
  $\overline{\cal M}_{(g,h),(n,\vec{m})}
                    (W/B,L\,|\,[\beta],\vec{\gamma},\mu)/B$
  determine a unique equivalence class of Kuranishi
  structures-in-${\cal C}_{\spsccw}$ on
  $\overline{\cal M}_{(g,h),(n,\vec{m})}
                   (W/B,L\,|\,[\beta],\vec{\gamma},\mu)/B$.
} % end-proposition

%%%%%%%%%%%%%%%%%%%%%%%%%%%%%
%
% \bigskip
%
% \subsection{Family orientations of the family Kuranishi structure.}
%
% \bigskip
%
% \noindent $\bullet$
% ????????????.
%
%%%%%%%%%%%%%%%%%%%%%%%%%%%%%

\bigskip

\section{The moduli space
          $\overline{\cal M}_{(g,h),(n+l(\vec{s}),\vec{m})}
           (Z,L;D\,|\,\beta^{\prime},\vec{\gamma},\mu^{\prime};\vec{s})$
          of relative stable maps  and
         its Kuranishi structure.}

We apply and extend the construction in Sec.~1 - Sec.~3
  to a relative pair $(Z,L;D)$ and its expansions
 to define the moduli space
  $\overline{\cal M}_{(g,h),(n+l(\vec{s}),\vec{m})}
    (Z,L;D\,|\,\beta^{\prime},\vec{\gamma},\mu^{\prime};\vec{s})$
  of relative stable maps of type
  $((g,h),(n+l(\vec{s}),\vec{m})\,|\,
          \beta^{\prime},\vec{\gamma},\mu^{\prime};\vec{s})$,
  from labelled-bordered Riemann surfaces with marked points to
  the fibers of the expanded relative pairs
  $(\widehat{Z},\widehat{L};\widehat{D})/\widehat{A}$
  associated to $(Z,L;D)$; (Sec.~6.1).
The same technique in Sec.~4 and Sec.~5 is used to construct
 a Kuranishi structure thereupon; (Sec.~6.2).
See also
 [I-P1], [L-R] for the symplecto-analytic setting
  in different formats and
 [Li1: Sec.~4], [Li2: Sec.~2], [Gr-V] for the algebro-geometric
  setting.

\bigskip

\subsection{The moduli space
            $\overline{\cal M}_{(g,h),(n+l(\vec{s}),\vec{m})}
            (Z,L;D\,|\,\beta^{\prime},\vec{\gamma},\mu^{\prime};\vec{s})$
            of relative stable maps.}

Let $(Z,L;D)$ be a symplectic pair $(Z;D)$,
  with a compatible almost-complex structure,
 together with a Lagrangian/almost-complex submanifold $L$
  that is disjoint from $D$.
Recall
  the space $(\widehat{Z};\widehat{D})/\widehat{A}$ of expanded
   relative pairs associated to $(Z;D)$ with the quotient topology,
  its standard local charts
   $\varphi[k]:(Z[k];D[k])/A[k]\rightarrow
               (\widehat{Z};\widehat{D})/\widehat{A}$
   with $k\in{\Bbb Z}_{\ge 0}$, and
  the equivariant pseudo-${\Bbb G}_m[k]$-action on $(Z[k];D[k])/A[k]$
 from Sec.~1.2.
Let $L[k]$ be the submanifold $\tilde{\mathbf p}[k]^{-1}(L)$ of
 $Z[k]$ from the map
 $\tilde{\mathbf p}[k]:(Z[k];D[k])/A[k]\rightarrow (Z;D)/\pt$.
Over $A[k]$, $L[k]=A[k]\times L$.
Sec.~1.2 can be made to incorporate $L[k]$.
This gives the space
 $(\widehat{Z},\widehat{L};\widehat{D})/\widehat{A}$.
The central fiber of $(Z[k],L[k];D[k])/A[k]$
 is almost-complex isomorphic to the
 pair-with-a-totally-real-submanifold $(Z_{[k]},L_{[k]};D_{[k]})$.

Recall the definition of the relative Maslov index $\mu^{\rel}(h)$ of
 a smooth map
 $h:(\Sigma,\partial\Sigma)\rightarrow (Z_{[k]},L_{[k]};D_{[k]})$
 from Sec.~3.1.
Note also that
 the monodromies of $(Z[k],L[k];D[k])/A[k]$, $k\in{\Bbb Z}_{\ge 0}$,
  on a smooth fiber,
    which is almost-complex isomorphic to $(Z,L;D)$,
  are relatively isotopic to the identity map with respect to $(L;D)$;
  thus,
 the monodromy
  $(\widehat{Z}[k],\widehat{L}[k];\widehat{D})/\widehat{A}$-action
  on $H_1(L;{\Bbb Z})$, $H_2(Z,L;{\Bbb Z})$, and $H_2(Z,L\cup D;{\Bbb Z})$
  are all trivial.

\bigskip

\begin{flushleft}
{\bf Moduli space of relative stable maps to fibers of
      $(\widehat{Z},\widehat{L};\widehat{D})/\widehat{A}$.}
\end{flushleft}
%
%\noindent
{\bf Definition 6.1.1
     [relative stable map to fibers of $(Z[k],L[k];D[k])/A[k]$].}
{\rm
 Let
  $\beta^{\prime}\in H_2(Z,L;{\Bbb Z})$,
  $\vec{\gamma}
    =(\gamma_1\,,\,\ldots\,,\gamma_h)\in H_1(L;{\Bbb Z})^{\oplus h}$
  such that $\partial\beta=\gamma_1+\,\cdots\,\gamma_h$,
  $\mu^{\prime}\in{\Bbb Z}$,  and
  $\vec{s}
    =(s_1,\,\ldots\,,\,s_l)\in ({\Bbb Z}_{\ge 0})^l$.\footnote{For
                          $\vec{s}$,
                          we define its
                          {\it length}
                           $\;l(\vec{s}):=l\,,\,$
                          {\it degree}
                           $\;\footnotesizedeg(\vec{s})
                               := s_1+\,\cdots\,+s_l\,,\,$ and
                          {\it multiplicity}
                           $\;m(\vec{s}) := s_1\,\cdots\,s_l\,$.}
 A relative map
  $f:(\Sigma,\partial\Sigma)/\pt\rightarrow (Z[k], L[k];D[k])/A[k]$
  from a bordered Riemann surface $\Sigma$ to a fiber
  of $(Z[k],L[k];D[k])/A[k]$ is called {\it prestable} of
   ({\it combinatorial}) {\it type}
   $((g,h),(n+l(\vec{s}),\vec{m})\,
           |\,\beta^{\prime},\vec{\gamma},\mu^{\prime};\vec{s})$
  if
   \begin{itemize}
    \item[$\cdot$]
     $f$ is prestable of type
     $((g,h),(n+l(\vec{s}),\vec{m})\,|\,
        \beta^{\prime},\vec{\gamma},\mu^{\prime}+2\,\degree(\vec{s}))$
     as a map to a fiber of $(Z[k],L[k])/A[k]$,
     cf.\ Definition 3.3.1 (with $[\beta]=\{\beta^{\prime}\}$);
        % Definition [stable map to fibers of $(W[k],L[k])/B[k]$]
     the last $l(\vec{s})$ free marked points on $\Sigma$
      shall be called the {\it distinguished marked points};

    \item[$\cdot$]
     ($f$ is non-degenerate with respect to $D[k]$;)$\;$
     $f^{-1}(D[k]) = s_1\,p_{n+1}+\,\cdots\,+s_l\,p_{n+l(\vec{s})}$,
      where $p_{n+1}\,,\,\ldots\,,\,p_{n+l(\vec{s})}$
      are the distinguished marked points on $\Sigma$;
     (in particular, $\mu^{\rel}(f)=\mu^{\prime}$  and
      all distinguished marked points are smooth interior points
       on $\Sigma$).
   \end{itemize}

 An {\it isomorphism} between two relative prestable maps
    $f_1:\Sigma_1/\pt\rightarrow (Z[k],L[k];D[k])/A[k]$,
    $f_2:\Sigma_2/\pt\rightarrow (Z[k],L[k];D[k])/A[k]$
   of the same type
  is a pair $(\alpha,\beta)$,
  where
   $\alpha:\Sigma_1\rightarrow \Sigma_2$
     is an isomorphism of prestable labelled-bordered Riemann surfaces
      with marked points and
   $\beta\in{\Bbb G}_m[k]$
   such that $f_1\circ \beta = f_2\circ\alpha$.
 The isomorphism class of $f$ is denoted by $[f]$.
 The notion of
   {\it non-degenerate} (resp.\ {\it pre-deformable} )
    relative prestable maps,
   {\it distinguished nodes $q$}, and the {\it contact order} at $q$
   are defined exactly the same as in Definition 3.3.1.
         % Definition [stable map to fibers of $(W[k],L[k])/B[k]$]

 A relative prestable map
   $f:\Sigma/\pt \rightarrow (Z[k].L[k];D[k])/A[k]$
   is called {\it stable}
  if $f$ is pre-deformable and
     its group $\Aut(f)$ of automorphisms is finite.
 The moduli space of isomorphism classes of stable maps to fibers of
  $(Z[k],L[k];D[k])/A[k]$  of type
  $((g,h),(n+l(\vec{s}),\vec{m})\,
       |\,\beta^{\prime},\vec{\gamma},\mu^{\prime};\vec{s})$
  is denoted by
  ${\cal M}_{(g,h),(n+l(\vec{s}),\vec{m})}^{\nonrigid}
   ((Z[k],L[k];D[k])/A[k]\,|\,
               \beta^{\prime},\vec{\gamma},\mu^{\prime};\vec{s})$;
 it is equipped with the {\it $C^{\infty}$-topology},
  defined similarly as in Sec.~3.3.
} % end-definition

\bigskip

The pseudo-embedding
 $\varphi^{\prime}_{k^{\prime},k; I}:
  (Z[k^{\prime}],L[k];D[k])/A[k^{\prime}]
    \hookrightarrow (Z[k],L[k];D[k])/A[k]$,
  $k^{\prime}<k$ and $I\subset\{0\,,\,\ldots\,,\,k-1\}$,
 from Sec.~1.2 induces a {\it pseudo-embedding}
 $$
  \begin{array}{r}
   \varphi^{\prime}_{k^{\prime},k; I}\;  :\;
    {\cal M}_{(g,h), (n+l(\vec{s}), \vec{m})}^{\nonrigid}
             ((Z[k^{\prime}],L[k^{\prime}];D[k^{\prime}])/A[k^{\prime}]\,
                      |\,\beta^{\prime},\vec{\gamma},\mu^{\prime};\vec{s})
                                   \hspace{6em}\\[1ex]
   \hookrightarrow\;
      {\cal M}_{(g,h), (n+l(\vec{s}), \vec{m})}^{\nonrigid}
        ((Z[k],L[k];D[k])/A[k]\,
               |\,\beta^{\prime},\vec{\gamma},\mu^{\prime};\vec{s})\,.
  \end{array}
 $$
Define the set of isomorphism classes of relative stable maps to
 fibers of $(\widehat{Z},\widehat{L};\widehat{D})/\widehat{A}\,$:
 \begin{eqnarray*}
  \lefteqn{\overline{\cal M}_{(g,h),(n+l(\vec{s}),\vec{m})}
   (Z,L;D\,|\,\beta^{\prime},\vec{\gamma},\mu^{\prime}; \vec{s})} \\[.6ex]
  &&
  :=\;
   \left.
     \left( \amalg_{k=0}^{\infty}
      {\cal M}_{(g,h), (n, \vec{m})}^{\nonrigid}
      ((Z[k],L[k];D[k])/A[k]\,|\,
                  \beta^{\prime},\vec{\gamma},\mu^{\prime};\vec{s})
     \right)
   \right/ \raisebox{-.6ex}{$\sim$}\,,
 \end{eqnarray*}
 where the equivalence relation $\sim$ is generated by
  $[f]\sim \varphi^{\prime}_{k^{\prime},k;I}([f^{\prime}])$
  for\\
   $[f]  \in
    {\cal M}_{(g,h),(n+l(\vec{s}),\vec{m})}^{\nonrigid}
     ((Z[k],L[k];D[k])/A[k]\,|\,
            \beta^{\prime},\vec{\gamma},\mu^{\prime};\vec{s})$  and
   $[f^{\prime}]  \in$
    the defining domain of $\varphi_{k^{\prime},k;I}$  on
     ${\cal M}_{(g,h),(n+l(\vec{s}),\vec{m})}^{\nonrigid}
     ((Z[k^{\prime}],L[k^{\prime}];D[k^{\prime}])/A[k^{\prime}]\,
              |\,\beta^{\prime},\vec{\gamma},\mu^{\prime};\vec{s})$.
By construction, there are embeddings of sets
 \begin{eqnarray*}
  \lefteqn{\varphi^{\prime}_{(k)}\,:\,
      {\cal M}_{(g,h),(n+l(\vec{s}),\vec{m})}^{\nonrigid}
         ((Z[k],L[k];D[k])/A[k]\,
          |\,\beta^{\prime},\vec{\gamma},\mu^{\prime};\vec{s})} \\[.6ex]
  && \hookrightarrow\,
      \overline{\cal M}_{(g,h),(n+l(\vec{s}),\vec{m})}
       (Z,L;D\,|\,\beta^{\prime},\vec{\gamma},\mu^{\prime};\vec{s})\,,
      \hspace{1em}
  k\in{\Bbb Z}_{\ge 0}\,.
 \end{eqnarray*}
A subset $U$ of
  $\overline{\cal M}_{(g,h),(n+l(\vec{s}),\vec{m})}
   (Z,L;D\,|\,\beta^{\prime},\vec{\gamma},\mu^{\prime};\vec{s})$
  is said to be {\it open}
 if $U=\cup_{\alpha}U_{\alpha}$ such that
     $U_{\alpha}$ is contained in the image of some
      $\varphi^{\prime}_{(k)}$  and
     ${\varphi^{\prime}_{(k)}}^{-1}(U_{\alpha})$ is open  in\\
     ${\cal M}_{(g,h),(n+l(\vec{s}),\vec{m})}^{\nonrigid}
         ((Z[k],L[k];D[k])/A[k]\,|\,
                \beta^{\prime},\vec{\gamma},\mu^{\prime};\vec{s})$.
This defines the {\it $C^{\infty}$-topology} on the
 moduli space
  $\overline{\cal M}_{(g,h),(n+l(\vec{s}),\vec{m})}
          (Z,L;D\,|\,\beta^{\prime},\vec{\gamma},\mu^{\prime};\vec{s})$
 of relative stable maps to fibers of
 $(\widehat{Z},\widehat{L};\widehat{D})/\widehat{A}$.

\bigskip

\noindent
{\bf Definition 6.1.2 [tautological cover].} {\rm
 By construction,
  $$
   \left\{\, {\cal M}_{(g,h),(n+l(\vec{s}),\vec{m})}^{\nonrigid}
             ((Z[k],L[k];D[k])/A[k]\,|\,
               \beta^{\prime},\vec{\gamma},\mu^{\prime};\vec{s})\,
    \right\}_{k\in{\scriptsizeBbb Z}_{\ge 0}}
  $$
  is an open cover of
  $\overline{\cal M}_{(g,h),(n+l(\vec{s}),\vec{m};\vec{s})}
   (Z,L;D\,|\,\beta^{\prime},\vec{\gamma},\mu^{\prime};\vec{s})$.
 We will call it the {\it tautological cover} of
  $\overline{\cal M}_{(g,h),(n+l(\vec{s}), \vec{m})}
   (Z,L;D\,|\,\beta^{\prime},\vec{\gamma},\mu^{\prime};\vec{s})$.
} % end-definition

\bigskip

Indeed, there exists $k_0$ depending $(Z,L;D)$ and
 $((g,h),(n+l(\vec{s}),\vec{m})|
         \beta^{\prime},\vec{\gamma},\mu^{\prime};\vec{s})$
 such that

 \vspace{-1ex}
 {\scriptsize
 $$
  \begin{array}{l}
  {\cal M}_{(g,h),(n+l(\vec{s}),\vec{m})}^{\tinynonrigid}
     ((Z[k_0],L[k_0];D[k_0])/A[k_0]\,|\,
       \beta^{\prime},\vec{\gamma},\mu^{\prime},\vec{s}) \\[1.6ex]
  \hspace{2em} \supset\;
    {\cal M}_{(g,h),(n+l(\vec{s}),\vec{m})}^{\tinynonrigid}
    ((Z[k_0+1],L[k_0+1];D[k_0+1])/A[k_0+1]\,|\,
     \beta^{\prime},\vec{\gamma},\mu^{\prime},\vec{s})  \\[1.6ex]
  \hspace{2em} \supset\;
    {\cal M}_{(g,h),(n+l(\vec{s}),\vec{m})}^{\tinynonrigid}
    ((Z[k_0+2],L[k_0+2];D[k_0+2])/A[k_0+2]\,|\,
               \beta^{\prime},\vec{\gamma},\mu^{\prime},\vec{s})\;
   \supset\;\cdots\,.
  \end{array}
 $$ } % end-scriptsize

\vspace{-1ex}
\noindent
Thus, the tautological cover of
  $\overline{\cal M}_{(g,h),(n+l(\vec{s}), \vec{m})}
   (Z,L;D\,|\,\beta^{\prime},\vec{\gamma},\mu^{\prime};\vec{s})$
 is finite in effect, cf.~Theorem 6.1.3.
  % Theorem [Hausdorffness and compactness]
The universal maps on the universal curve over each
 ${\cal M}_{(g,h),(n+l(\vec{s}),\vec{m})}^{\nonrigid}
  ((Z[k],L[k];D[k])/A[k]\,|\\[0ex]
         \beta^{\prime},\vec{\gamma},\mu^{\prime};\vec{s})$
 are glued to give the universal map (between spaces with charts)
  $$
   F\; :\; {\cal C}
    /\overline{\cal M}_{(g,h),(n+l(\vec{s}),\vec{m})}
     (Z,L;D\,|\,\beta^{\prime},\vec{\gamma},\mu^{\prime};\vec{s})\;
   \longrightarrow\;
   (\widehat{Z},\widehat{L};\widehat{D})/\widehat{A}\,.
  $$

\bigskip

\begin{flushleft}
{\bf Hausdorffness, finite stratification, and compactness.}
\end{flushleft}
Recall
 the notion of {\it weighted layered $(A_2\rightarrow A_1)$-graphs}
  from Definition 3.3.5,
 the category ${\mathfrak G}(A_2\rightarrow A_1)$ of graphs,  and
       % Definition [weighted layered $(A_2\rightarrow A_1)$-graph]
 how a stable map $f$ to fibers of $(W[k-1],L[k-1])/B[k-1]$
 (now $=(Z[k],L[k])/A[k]$) corresponds a such graph $\tau_{[f]}$.
To encode the contact-order data $\vec{s}$ of relative maps with $D[k]$,
 we add to the objects $\tau$ in ${\mathfrak G}(A_2\rightarrow A_1)$
 the following data:
 \begin{itemize}
  \item[$\cdot$]
   an ordered set $R(\tau)$ of $l$-many {\it roots} $r_i$,
    $i=1,\,\ldots\,,\,l$, that are attached to vertices
   of the largest layer-value;

  \item[$\cdot$]
   an additional {\it weight function}
   $\;\order^{\prime}:R(\tau)\rightarrow {\Bbb Z}_{\ge 0}$,
    $r_i\mapsto s_i$;

  \item[$\cdot$]
   replace $\mu(\tau)$ in Definition 3.3.5
        % Definition [weighted layered $(A_2\rightarrow A_1)$]
    by $\mu^{\prime}(\tau)$, called the
    {\it relative index} of $\tau$.\footnote{Let
                            $\vec{s} =
                             (\footnotesizeord^{\prime}(r_1)\,,\,
                              \cdots\,,\,
                              \footnotesizeord^{\prime}(r_l))$.
                            Then, we will call
                             the quantity
                             $\mu^{\prime}+2\,\footnotesizedeg(\vec{s})$
                             the (absolute) {\it index} of $\tau$ and
                             denote it by $\mu(\tau)$.
                            When $l=0$, $\mu^{\prime}(\tau)=\mu(\tau)$.}
 \end{itemize}
Denote a such graph still by $\tau$ with the same name:
 {\it weighted layered $(A_1\rightarrow A_1)$-graph}.
An {\it isomorphism} $\alpha:\tau_1\rightarrow \tau_2$
 between two such graphs is defined the same as in Definition 3.3.5
  with the index replaced by relative index and
       % Definition [weighted layered $(A_2\rightarrow A_1)$-graph]
  the additional requirement that
   $\alpha$ induces an isomorphism
    $R(\tau_1)\stackrel{\sim}{\rightarrow} R(\tau_2)$
    as ordered weighted sets.
The corresponding new {\it category of graphs} enlarges
  the previous one  and
 will be denoted still by ${\mathfrak G}(A_2\rightarrow A_1)$
 (or simply ${\mathfrak G}$ when $(A_2\rightarrow A_1)$
  is understood).

The notion of {\it genus}, {\it $b$-weight},
 {\it contraction}, and ({\it red-to-blue}) {\it color change} of
 weighted layered $(A_2\rightarrow A_1)$-graphs
 extend to the new ${\mathfrak G}(A_2\rightarrow A_1)$.
The correspondence of a point
 $[f:\Sigma/\pt \rightarrow
     (\widehat{Z},\widehat{L};\widehat{D})/\widehat{A}]
   \in   \overline{\cal M}_{(g,h),(n+l(\vec{s}),\vec{m})}
           (Z,L;D\,|\,\beta^{\prime},\vec{\gamma},\mu;\vec{s})$,
  with target isomorphic to $(Z_{[k]},L_{[k]};D_{[k]})$,
 to an element $\tau_{[f]}\in{\mathfrak G}(A_2\rightarrow A_1)$
 is the same as in Sec.~3.3 with the following addition/modification:

{\footnotesize
\noindent\hspace{1em}
\begin{tabular}{lcl}
  & & \\
 \hspace{2em}
 $f:\Sigma \rightarrow (Z_{[k]}, L;D_{[k]})$
  && $(H_2(Z,L;{\footnotesizeBbb Z})\stackrel{\partial}{\rightarrow}
       H_1(L;{\footnotesizeBbb Z}))$-graph $\tau$ \\[.6ex]
 \hline  \\[-1.6ex]
 \hspace{2em} $\cdots\cdots\cdots\cdots$
    && \hspace{2em}
       $\cdots\cdots\cdots\cdots\cdots\cdots\hspace{15em}$\\[.6ex]
 $i$-th distinguished marked point $p_{n+i}$
   &&  {\it root} $r_i\in R(\tau)$ attached to
       $v\in V(\tau)$ with $\layer(v)=k$  \\
 {\it contact order} $s_i$ of $f$ with $D_{[k]}$ at $p_{n+i}$
   &&  $\order^{\prime}(r_i)$, $r_i\in R(\tau)$  \\
 {\it relative Maslov index} $\mu^{\rm rel}(f)$
   && {\it relative index} $\mu^{\prime}\,$.     \\
 && \\
\end{tabular}
} % end-footnotesize

\noindent
Two relative stable maps
  $f_i:\Sigma_i/\pt \rightarrow (Z[k_i],L[k_i];D[k_i])/A[k_i]$,
    $i=1,\,2$,
  are said to be {\it of the same topological type}
 if $\tau_{[f_1]}$ is isomorphic to $\tau_{[f_2]}$ in
  the category ${\mathfrak G}$.
Degenerations of relative stable maps to fibers of
 $(\widehat{Z},\widehat{L};\widehat{D})/\widehat{A}$
 are reflected contravariantly by compositions of contractions
 and color-changes of their dual graphs.

Same reasons that give
  Proposition 3.3.4,  % Proposition [Hausdorffness]
  Lemma 3.3.7,        % Lemma [finite stratification],
    and
  Theorem 3.3.8       % Theorem [compactness$/B$]
 now imply:

\bigskip

\noindent
{\bf Theorem 6.1.3 [Hausdorffness and compactness].} {\it
 The classification of relative stable maps by their topological types
  gives rise to a finite stratification of
  $\overline{\cal M}_{(g,h),(n+l(\vec{s}),\vec{m})}
            (Z,L;D\,|\,\beta^{\prime},\vec{\gamma},\mu;\vec{s})$,
  with each stratum $S_{\tau}$ labelled by
  a weighted layered $(H_2(Y,L;{\Bbb Z}), H_1(L;{\Bbb Z}))$-graph
  $\tau\in{\mathfrak G}$.
 The moduli space
  $\overline{\cal M}_{(g,h),(n+l(\vec{s}),\vec{m})}
            (Z,L;D\,|\,\beta^{\prime},\vec{\gamma},\mu;\vec{s})$
  of relative stable maps to fibers of
   $(\widehat{Z},\widehat{L};\widehat{D})/\widehat{A}$  of
  combinatorial type
   $((g,h),(n+l(\vec{s}),\vec{m})\,|\,
             \beta^{\prime},\vec{\gamma},\mu^{\prime};\vec{s})$,
  with the $C^{\infty}$-topology,
 is Hausdorff and compact.
} % end-theorem

\bigskip

\noindent
Cf.\ [L-R: Sec.~3.3], [I-P1: Theorem 7.4]; [Li1: Theorem 4.10].

\bigskip

\subsection{A Kuranishi structure for
            $\overline{\cal M}_{(g,h),(n+l(\vec{s}),\vec{m})}
             (Z,L;D\,|\,\beta^{\prime},\vec{\gamma},\mu^{\prime};\vec{s})$.}

Introduce first the following category of topological spaces,
 which is closely related to ${\cal C}_{\spsccw}$:

\bigskip

\noindent
{\bf Definition 6.2.1 [category ${\cal C}_{\spsccw}^{\,\prime}$].}
{\rm
 We define ${\cal C}_{\spsccw}^{\,\prime}$ to be the category
  that has the same objects as ${\cal C}_{\spsccw}$ but with
  the fibrations over the complex line ${\Bbb C}$ removed.
 A morphism in ${\cal C}_{\spsccw}^{\,\prime}$ is a continuous map
  as stratified spaces.
} % end-definition

\bigskip

The same construction in Sec.~4 - Sec.~5 gives a
 Kuranishi structure on
 $\overline{\cal M}_{(g,h),(n+l(\vec{s}),\vec{m})}
   (Z,L;\\  D\,|\,\beta^{\prime},\vec{\gamma},\mu^{\prime};\vec{s})$
 that is modelled in the category ${\cal C}_{\spsccw}^{\,\prime}$.
There are only two major modifications in the discussion:
 \begin{quote}
  \parbox[t]{2.7em}{\hspace{1em}\Large $\cdot$}
   \parbox[t]{32,4em}{(the {\it non-rigidity} of target)$\,$:
    \hspace{1ex}
   while treating $(Z[k],L[k];D[k])/A[k]$ as the $(k-1)$-th expanded
    degeneration of the degeneration $(Z[1],L[1];D[1])/A[1]$,
   it is ${\Bbb G}_m[k]$  -- rather than ${\Bbb G}_m[k-1]$ --
    that acts equivariantly on $(Z[k],L[k];D[k])/A[k]$ and
    that corresponds to choices of the renormalization in
    removing degeneracy/falling-into-$D$;}

  \parbox[t]{2.7em}{(T4)}
   \parbox[t]{32.4em}{(additional {\it transversality})$\,$:
    \hspace{1ex}
   local transversality of the contact-order-$s_i$ condition
    along $D_{[k]}$ at the distinguished marked point
    $p_{n+i}$, for $i=1\,,\,\ldots\,,\,l$;
   cf.\ Conditions (T1) - (T3) in Sec.~5.2.}
 \end{quote}

Let
 $(\,\Sigma,\, \dot{\partial}\Sigma;\,
       \vec{p},\,\vec{p}_1,\,\ldots,\,\vec{p}_h;\, f)$
  be a relative stable map to the central fiber
   $(Z_{[k]},L_{[k]};D_{[k]})$ of $(Z[k],L[k];D[k])/A[k]$
%   and
%  $\rho_{(i)}:= (\,\Sigma_{(i)},\, (\dot{\partial}\Sigma)_{(i)};\,
%     \vec{p}_{(i)},\,
%     \vec{p}_{1,(i)},\,\ldots,\,\vec{p}_{h, (i)};\, f_{(i)})$
%   be the associated submap to the irreducible component $\Delta_i$
%   of $Y_{[k]}$, for $i=0,\,\ldots,\,k+1$.
%  (By construction, $(\dot{\partial}\Sigma)_{(i)}$, $\vec{p}_{j,(i)}$
%  can be non-empty only for $i=0$ and $k+1$.)
% We denote the labelled-bordered Riemann surface with marked points
%  $(\Sigma,\, \dot{\partial}\Sigma;\,
%              \vec{p},\,\vec{p}_1,\,\ldots,\,\vec{p}_h)$
% also simply by $\Sigma$.
 that represents
 $\rho  \in
  \overline{\cal M}_{(g,h),(n+l(\vec{s}),\vec{m})}
           (Z,L;D\,|\,\beta^{\prime},\vec{\gamma},\mu;\vec{s})$.
Recall that
 $Z_{[k]}
  = Z\cup_{D=D_{1,\infty}}\Delta_1
     \cup_{D_{1,0} = D_{2,\infty}}\,\cdots\,
     \cup_{D_{k-1,0}=D_{k,\infty}}\Delta_k$ and
 $D_i := \Delta_i\cap \Delta_{i+1}$ in $Z_{[k]}$
 for $i=1,\,\ldots\,,\,k-1$.
Here we set $\Delta_0=Z$ by convention.
Let $\Lambda_i=f^{-1}(D_i)$ and
 $\Lambda=\disjointunion_{i=0}^{k-1}\,\Lambda_i$ be the set of
 distinguished nodes on $\Sigma$ under $f$.
Let ${\mathbf s}=(\vec{s}_0,\,\cdots\,,\,\vec{s}_{k-1}\,;\,\vec{s}_k)$,
  with $\vec{s}_k =\vec{s}$, be the tuple of contact orders of $f$ at
  $\Lambda\cup \{p_{n+1}\,,\,\cdots\,,\,p_{n+l(\vec{s})}\}$.
Recall the discussion and notations in Sec.~5.2.
The notion of a saturated subspace in
 $W^{1,p}(\Sigma, \partial\Sigma;f^{\ast}T_{\ast}Z_{[k]},
                    (f|_{\partial\Sigma})^{\ast}T_{\ast}L_{[k]})$
 from Definition 5.3.1.4
      % Definition [saturated/pre-deformable subspace]
 now has to be revised to incorporate Condition (T4) as well:

\bigskip

\noindent
{\bf Definition 6.2.2 [saturated/relative pre-deformable subspace].}
{\rm
 A subspace $V$ in
  $W^{1,p}(\Sigma,\\  \partial\Sigma;f^{\ast}T_{\ast}Z_{[k]},
                   (f|_{\partial\Sigma})^{\ast}T_{\ast}L_{[k]})$
  is said to be {\it saturated}$\,$
 if
  \begin{itemize}
   \item[$(1)$]
    $V$ is admissible;

   \item[$(2)$]
   the map
   \begin{eqnarray*}
    \lefteqn{
      \left( \oplus_{q\in\Lambda} D_f\divisor_q \right)
       \mbox{$\bigoplus$}
         \left( \oplus_{i=1}^{n+l(\vec{s})}\, D_f\ev_{p_i} \right)
      \mbox{$\bigoplus$} \left( \oplus_{q_{ij}} D_f\ev_{q_{ij}} \right)
      \mbox{$\bigoplus$}
       \left( \oplus_{i=1}^{l(\vec{s})}\,D_f\divisor_{p_{n+i}}
        \right)\;
     :\;  V\;\longrightarrow } \\[.6ex]
    &&
     \left( \mbox{$\bigoplus$}_{q\in\Lambda}
        \left( T_{(s(q)-1)\cdot(q)}\Div^{s(q)-1}(U_{q,1})
               \oplus T_{(s(q)-1)\cdot(q)}\Div^{s(q)-1}(U_{q,2})
        \right)
      \right)  \\[.6ex]
    && \hspace{6em}
    \mbox{$\bigoplus$} \left( \oplus_{p_i} T_{f(p_i)}Y_{[k]} \right)
    \mbox{$\bigoplus$} \left( \oplus_{q_{ij}} T_{f(q_{ij})}L \right)
   \mbox{$\bigoplus$}_{i=1}^{l(\vec{s})}\,
       T_{(s_i-1)\cdot(p_{n+i})}\Div^{s_i-1}(U_{p_{n+i}})
   \end{eqnarray*}
   is surjective;

   \item[$(3)$]
   let $V^{\relpd}$ be the subspace
    $( (\oplus_{q\in\Lambda} D_f\divisor_q)  \oplus
       (\oplus_{q\in\Lambda} D_f\divisor_q) )^{-1}({\mathbf 0})$
    in $V$,
   then the linear map
     \begin{eqnarray*}
      \lefteqn{\left(
        \oplus_{q\in\Lambda}\,( D_f\ev_q \oplus \jet_q^{s(q)} )\right)
        \mbox{$\bigoplus$}
        \left( \oplus_{i=1}^{l(\vec{s})}\,
                   ( D_f\ev_{p_{n+i}} \oplus \jet_{p_{n+i}}^{s_i} )
         \right)  \; :\hspace{3em}}\\[.6ex]
      && \hspace{3em}
       V^{\relpd}\; \longrightarrow\;
       \left( \oplus_{q\in\Lambda}\,(T_{f(q)}D\oplus {\Bbb C}^2) \right)
        \mbox{$\bigoplus$}
       \left( \oplus_{i=1}^{l(\vec{s})}
             ( T_{f(p_{n+i})}D \oplus {\Bbb C} ) \right)
     \end{eqnarray*}
    is surjective, where we have identified $D_i$, $i=0,\,\ldots,\,k-1$,
    canonically with $D$.
  \end{itemize}
 In the above statement, $V^{\relpd}$ is called the
  {\it relative pre-deformable subspace} of $V$.

 A subspace $E$ of
   $L^p(\Sigma;
      \Lambda^{0,1}\Sigma\otimes_J f^{\ast}T_{\ast}Y_{[k]})$
   is said to be {\it saturated}
 if $(D_{\!f}\bar{\partial}_J)^{-1}(E)\subset
        W^{1,p}(\Sigma,\\  \partial \Sigma;f^{\ast}T_{\ast}Y_{[k]},
                      (f|_{\partial\Sigma})^{\ast}T_{\ast}L_{[k]})$
  is saturated.
} % end-definition

\bigskip

The notion of a saturated obstruction space $E_{\rho}$ in
 $L^p(\Sigma;
     \Lambda^{0,1}\Sigma\otimes_J f^{\ast}T_{\ast}Z_{[k]})$
 for $\rho = [f]  \in
        \overline{\cal M}_{(g,h),(n+l(\vec{s}),\vec{m})}
        (Z,L;D\,|\,\beta^{\prime},\vec{\gamma},\mu^{\prime};\vec{s})$
 from Definition/Lemma 5.3.1.5
    % Definition/Lemma [saturated obstruction space]
 is converted accordingly:

\bigskip

\noindent
{\bf Definition/Lemma 6.2.3 [saturated obstruction space].} {\rm
 Denote by $\Image(D_{\!f}\bar{\partial}_J)$ the image of
 $D_{\!f}\bar{\partial}_J$,
  $(D_{\!f}\bar{\partial}_J)
   ( W^{1,\,p}( \Sigma, \partial\Sigma; f^{\ast}T_{\ast}Z_{[k]},
                 (f|_{\partial\Sigma})^{\ast}T_{\ast}L_{[k]}) )$, in
  $L^p(\Sigma;
    \Lambda^{0,1}\Sigma\otimes_J f^{\ast}T_{\ast}Z_{[k]})$.
 Then {\it
 there exists a subspace $E_\rho$ of
  $L^p(\Sigma; \Lambda^{0,1}\Sigma\otimes_J f^{\ast}T_{\ast}Z_{[k]})$
 such that
 \begin{itemize}
  \item[{\rm (1)}]
   $\Image(D_{\!f}\bar{\partial}_J) + E_{\rho}
    = L^p(\Sigma;
         \Lambda^{0,1}\Sigma\otimes_J f^{\ast}T_{\ast}Z_{[k]})$,

  \item[{\rm (2)}]
   $E_{\rho}$ is finite-dimensional, complex linear,
    and $\Aut(\rho)$-invariant,

  \item[{\rm (3)}]
   $E_{\rho}$ consists of smooth sections supported in
    a compact subset of $\Sigma$ disjoint from the union of
     the set of all $\,($three types of$\,)$ nodes and
     the set $\{p_{n+1}\,,\,\cdots\,,\,p_{n+l(\vec{s})}\}$
             of all distinguished marked points on $\Sigma$,

  \item[{\rm (4)}]
   $(D_{\!f}\bar{\partial}_J)^{-1}(E_{\rho})$ is a saturated
    subspace of
   $W^{1,p}( \Sigma, \partial\Sigma; f^{\ast}T_{\ast}Z_{[k]},
                    (f|_{\partial\Sigma})^{\ast}T_{\ast}L_{[k]}) )$.
 \end{itemize}} % end-italic
 $E_{\rho}$ is called a {\it saturated obstruction space} of
  $\,\overline{\cal M}_{(g,h),(n+l(\vec{s}),\vec{m})}
      (Z,L;D\,|\,\beta^{\prime},\vec{\gamma},\mu;\vec{s})$ at $\rho\,$.
} % end-lemma

\bigskip

The index of
 $$
  D_{\!f}\bar{\partial}_J\,:\,
   W^{1,p}( \Sigma, \partial\Sigma; f^{\ast}T_{\ast}Z_{[k]},
                   (f|_{\partial\Sigma})^{\ast}T_{\ast}L_{[k]} )\;
   \longrightarrow\;
    L^p(\Sigma,\Lambda^{0,1}\Sigma
                     \otimes_J f^{\ast}T_{\ast}Z_{[k]})
 $$
 is given by
 $$
 \begin{array}{ccl}
  \ind(D_{\!f}\bar{\partial}_J)
   & =  & \mu(f)\, +\, \dimm Z \cdot (1-\tilde{g})\,
          -\, 2\,\sum_{i=0}^{k-1} l(\vec{s}_i)\,
          +\, 4\,\sum_{i=0}^{k-1} \degree\vec{s}_i \\[.6ex]
   & =  & \mu^{\rel}(f)\, +\, \dimm Z \cdot (1-\tilde{g})\,
          -\, 2\,\sum_{i=0}^{k-1} l(\vec{s}_i)\,
          +\, 4\,\sum_{i=0}^{k-1} \degree\vec{s}_i\,
          +\, 2\,\degree\vec{s}\,,
 \end{array}
 $$
 where $\tilde{g}$ is the arithmetic genus of
 $\Sigma_{\scriptsizeBbb C}$.

\bigskip

\noindent
{\bf Definition 6.2.4 [relative pre-deformable index].} {\rm
 We define the {\it relative pre-deformable index} of
  $D_{\!f}\bar{\partial}_J$
  % (i.e.\ pre-deformable + correct contact order with $D_{[k]}$)
  to be
  $$
    \ind^{\relpd}(D_{\!f}\bar{\partial}_J)\;
    :=\; \mu^{\rel}(f)\,
          + \, \dimm Z \cdot (1-\tilde{g}) \, +\,  2\, |\Lambda|\,.
  $$
} % end-definition

% \bigskip

\noindent
Note that
 $$
  \dimm (D_{\!f}\bar{\partial}_J)^{-1}(E_{\rho})^{\relpd}\;
   =\; \mu^{\rel}(f)\, +\, \dimm Z \cdot (1-\tilde{g})\,
       +\, 2\,|\Lambda|\,  +\, \dimm E_{\rho}\,.
 $$

The same routine of Sec.~5.3 - Sec.~5.4 now proves that:

\bigskip

\noindent
{\bf Theorem 6.2.5 [Kuranishi structure on
      $\overline{\cal M}_{\tinybullet}(Z,L;D\,|\;\tinybullet\,)$].}
{\it
 The moduli space
  $\overline{\cal M}_{(g,h),(n+l(\vec{s}),\vec{m})}
    (\\  Z,L;D\,|\,\beta^{\prime},\vec{\gamma},\mu^{\prime};\vec{s})$
  of relative stable maps to fibers of
  $(\widehat{Z},\widehat{L};\widehat{D})/\widehat{A}$
  admits a Kuranishi structure ${\cal K}^{\prime}$
  modelled in ${\cal C}_{\spsccw}^{\,\prime}$.
 ${\cal K}^{\prime}$ has the expected dimension
  $$
   \vdim \overline{\cal M}_{\tinybullet}(Z,L;D\,|\;\tinybullet\,)\;
   :=\;  \mu^{\prime} + (N-3)(2-2g-h)
          + 2\,(n+l(\vec{s})) + (m_1 +\,\cdots\, + m_h)\,,
  $$
  where $2N$ is the dimension of $Z$.
 The Kuranishi neighborhood-in-${\cal C}_{\spsccw}^{\,\prime}$
  $(V_{\rho},\Gamma_{V_{\rho}},E_{V_{\rho}};s_{\rho},\psi_{\rho})$ at
  $\rho=[f:(\Sigma,\partial\Sigma)\rightarrow (Z_{[k]},L_{[k]};D_{[k]})]$
  has $V_{\rho}$ isomorphic to a neighborhood of the origin of
  $$
   \Xi_{(\vec{s}_0,\,\ldots,\,\vec{s}_{k-1})}
    \times {\Bbb R}^{n_1} \times ({\Bbb R}_{\ge 0})^{n_2}
  $$
  where
   \begin{itemize}
    \item[$\cdot$]
     $\vec{s}_i$ is the contact order of $f$ along $D_i$
      at the ordered set of distinguished nodes in $f^{-1}(D_i)$,
      $i=0,\,\ldots\,,\,k-1\,$,
     $\;($and recall that
          $\dimm \Xi_{(\vec{s}_0,\,\ldots,\,\vec{s}_{k-1})}=2k$$)$;

    \item[$\cdot$]
     $n_1\;=\;$
      $\vdim \overline{\cal M}_{\tinybullet}(Z,L;D\,|\;\tinybullet\,)
       + \dimm E_{\rho} - (2k+n_2)\,$;  and

    \item[$\cdot$]
     $n_2\;=\;$
      the total number of boundary nodes and free marked points
      that land on $\partial\Sigma$.
   \end{itemize}
 The homeomorphism-type $\{Z_{[k^{\prime}]}\}_{0\le k^{\prime}\le k}$
  of the targets of maps gives a $\Gamma_{V_{\rho}}$-invariant
  stratification
  $\{S_{k^{\prime}}\}_{0\le k^{\prime}\le k}$
  on $V_{\rho}$;
 each connected component of $S_{k^{\prime}}$
  is a manifold of codimension $2k^{\prime}$ in $V_{\rho}$.
 This stratification coincides with the induced stratification
  on $V_{\rho}$ from the stratification\footnote{\rm Which,
                               recall that, is induced by the map
                                $\Xi_{(\vec{s}_0,\,\ldots,\,\vec{s}_{k-1})}
                                 \rightarrow {\footnotesizeBbb C}^k$
                                  and
                                the stratification of
                                 ${\footnotesizeBbb C}^k$
                                 by the coordinate subspaces.}
  of $\Xi_{(\vec{s}_0,\,\ldots,\,\vec{s}_{k-1})}$.
 %%%%%%%%%%%%%%%%%%%%
 %
 % ${\cal K}^{\prime}$ is orientable if $L$ is spin
 % or if $h=1$ and $L$ is relative spin.
 %
 %%%%%%%%%%%%%%%%%%%%
} % end-theorem

\bigskip

\section{Degeneration and gluing of Kuranishi structures and
         axioms of open Gromov-Witten invariants
         under a symplectic cut.}

In this last section of the current work, we derive
 a degeneration-gluing relation of the Kuranishi structure
 of the moduli space of stable maps to $(X,L)$  with
 the Kuranishi structure of the moduli spaces of relative stable maps
  to $(Y_1,L_1;D)$, $(Y_2,L_2;D)$ that occur in a symplectic cut
  $\xi:(X,L)\rightarrow (Y,L)=(Y_1,L_1)\cup_D(Y_2,L_2)$.
This degeneration-gluing relation is insensitive to the real
 codimension-$1$ boundary of the Kuranishi structures involved
 when $L$ is non-empty.
Taking this formula as the foundation, together with
 (a) its reduction to the degeneration/gluing formula of
     virtual fundamental classes and Gromov-Witten invariants
     in closed Gromov-Witten theory when $L$ is empty and
 (b) the deformation-invariance requirement of Gromov-Witten invariants,
we propose a degeneration axiom and a gluing axiom under a symplectic cut
 for open Gromov-Witten invariants of a symplectic/almost-complex manifold
 with a decorated Lagrangian/totally-real submanifold.

\bigskip

\subsection{The degeneration-gluing relations of Kuranishi structures.}

\begin{flushleft}
{\bf Central fiber, layer-structure stratification, and
     descendent Kuranishi structure.}
\end{flushleft}
% \noindent
{\bf Definition 7.1.1
 [category ${\cal C}_{\spsccw,0}$ and
 its descendants ${\cal C}_{\spsccw,0}^{\,(i)}$].}
{\rm
 We define ${\cal C}_{\spsccw,0}$ to be the {\it category of
  weighted stratified spaces} $Q_0$ that occur in the central fiber
  of objects $Q/{\Bbb C}$ in ${\cal C}_{\spsccw}$.
 Here the {\it weight} to an irreducible component of $Q_0$ is given by
  the {\it multiplicity} of that component in terms of the associated
  flat affine fibrations $\Xi_{\mathbf s}/\Spec{\Bbb C}[t]$ of schemes,
  cf.\ footnote 19.   % CAUTION: Double-check the footnote number ?????.
 Define also the {\it depth-$i$ descendant}
  ${\cal C}_{\spsccw,0}^{\,(i)}$ of ${\cal C}_{\spsccw,0}$
  to be the category of stratified spaces locally modelled on
  the central fiber of the fibration
  $(\Xi_{(\vec{s}_0,\,\ldots,\,\vec{s}_i)}  \times
   {\Bbb R}^{n_1}\times ({\Bbb R}_{\ge 0})^{n_2})/{\Bbb C}^{i+1}$
   for some $n_1\,,\, n_2$.
 Note that
  ${\cal C}_{\spsccw,0}^{\,(0)} = {\cal C}_{\spsccw,0}$.
} % end-definition

\bigskip

\noindent
{\bf Definition 7.1.2
 [descendants ${\cal C}_{\spsccw}^{\,\prime,\,(i)}$  of
  ${\cal C}_{\spsccw}^{\,\prime}$].}
{\rm
 We define the
  {\it depth-$i$ descendants ${\cal C}_{\spsccw}^{\,\prime,\,(i)}$  of
  ${\cal C}_{\spsccw}^{\,\prime}$
  to be the category of stratified spaces locally modelled on
  the central fiber of the fibration
  $(\Xi_{(\vec{s}_0,\,\ldots,\,\vec{s}_{i-1})}  \times
   {\Bbb R}^{n_1}\times ({\Bbb R}_{\ge 0})^{n_2})/{\Bbb C}^i$
  for some $n_1\,,\, n_2$.
 Note that
  ${\cal C}_{\spsccw}^{\,\prime,\,(0)}
   = {\cal C}_{\spsccw}^{\,\prime}$.
} % end-definition

\bigskip

\noindent
{\bf Definition 7.1.3 [standard Kuranishi structure].} {\it
 We will call a Kuranishi structure ${\cal K}$ on a moduli space
 ${\cal M}$ of stable maps {\it standard} if ${\cal K}$ is constructed
 via the routine in Sec.~4 - Sec.~5.}

\bigskip

In particular, a standard Kuranishi structure-in-${\cal C}_{\spsccw}$
 ${\cal K}/B$ on
 $\overline{\cal M}_{(g,h),(n,\vec{m})}
                   (W/B,L\,|\,[\beta],\vec{\gamma},\\  \mu)/B$
 is flat over $B$ in the sense that
  each $\lambda\in B$  has a neighborhood $U_{\lambda}$
  over which ${\cal K}$ is equivalent to a standard Kuranishi structure
  $\hat{\cal K}/B$ with the Kuranishi neighborhoods and
  obstruction bundles from $\hat{K}/B$ flat over $U_{\lambda}$.
Indeed, a standard ${\cal K}/B$ constructed through Sec.~4 - Sec.~5
 is already flat over a neighborhood of $0\in B$.
This motivates/implies the following definiton/theorem,
 which is a corollary of
  Proposition 3.3.4,   % Proposition [Hausdorffness]
  Theorem 3.3.8,  and  % Theorem [compactness$/B$]
  Theorem 5.1.6:  % Theorem [family Kuranishi structure on
                  %  $\overline{\cal M}_{(g,h),(n,\vec{m})}
                  %   (W/B,L\,|\,\underline{\beta},\vec{\gamma},\mu)$]

\bigskip

\noindent
{\bf Definition/Theorem 7.1.4 [stable maps to $(Y,L)$].} {\rm
 Recall the symplectic cut $\xi:X\rightarrow Y$ and let
  $\underline{\beta}=\xi_{\ast}([\beta])\in H_2(Y,L;{\Bbb Z})$.
 Define the {\it moduli space}
   $\overline{\cal M}_{(g,h),(n,\vec{m})}
            (Y,L\,|\,\underline{\beta},\vec{\gamma},\mu)$
   of stable maps of type
   $((g,h),(n,\vec{m})\,|\,\underline{\beta},\vec{\gamma},\mu)$
   from labelled-bordered Riemann surfaces to $(Y,L)$
  to be the central fiber of
  $\overline{\cal M}_{(g,h),(n,\vec{m})}
                  (W/B,L\,|\,[\beta],\vec{\gamma},\mu)/B$
  over $0\in B$, with the induced {\it $C^{\infty}$-topology}.
 Then {\it
  $\overline{\cal M}_{(g,h),(n,\vec{m})}
      (Y,L\,|\,\underline{\beta},\vec{\gamma},\mu)$
   is Hausdorff and compact.
  The correspondence
   $$
    \overline{\cal M}_{(g,h),(n,\vec{m})}
         (Y,L\,|\,\underline{\beta},\vec{\gamma},\mu)\;
    \longrightarrow\;
    {\mathfrak G}(H_2(Y,L;{\Bbb Z})\stackrel{\partial}{\rightarrow}
                  H_1(L;{\Bbb Z}))\,, \hspace{2em}
   [f]\;\longmapsto\; \tau_{[f]}
  $$
  gives a finite stratification of
  $\overline{\cal M}_{(g,h),(n,\vec{m})}
        (Y,L\,|\,\underline{\beta},\vec{\gamma},\mu)$ by the
   {\it topological type} of maps.
 The central fiber ${\cal K}_0$ of  a
  standard Kuranishi structure-in-${\cal C}_{\spsccw}$ ${\cal K}/B$
  on  $\overline{\cal M}_{(g,h),(n,\vec{m})}
                  (W/B,L\,|\,[\beta],\vec{\gamma},\mu)\\  /B$
  gives a Kuranishi structure-in-${\cal C}_{\spsccw,0}$  on
  $\overline{\cal M}_{(g,h),(n,\vec{m})}
            (Y,L\,|\,\underline{\beta},\vec{\gamma},\mu)$.}
 We will call a Kuranishi structure on
  $\overline{\cal M}_{(g,h),(n,\vec{m})}
            (Y,L\,|\,\underline{\beta},\vec{\gamma},\mu)$
  thus obtained a {\it standard Kuranishi structure}  on
  $\overline{\cal M}_{(g,h),(n,\vec{m})}
      (Y,L\,|\,\underline{\beta},\vec{\gamma},\mu)$.
 The virtual dimension of ${\cal K}_0$ is the same as
  the virtual dimension of
  $\overline{\cal M}_{(g,h),(n,\vec{m})}
                     (X,L\,|\,[\beta],\vec{\gamma},\mu)$.
} % end-definition/lemma

\bigskip

For a
 $\tau\in{\mathfrak G} :=
   {\mathfrak G}(H_2(Y,L;{\Bbb Z})\stackrel{\partial}{\rightarrow}
                 H_1(L;{\Bbb Z}))$,
 denote the layer map $V(\tau)\rightarrow {\Bbb Z}_{\ge 0}$
 by $\layer_{\tau}$.
Then the composition
 $$
  [f]\;  \longmapsto\;  \tau_{[f]}\;   \longmapsto\;
         \max\{\,0\,,\, |\Image(\layer_{\tau_{[f]}})|-2\,\}
 $$
 gives a correspondence
 $$
  \overline{\cal M}_{(g,h),(n,\vec{m})}
        (Y,L\,|\,\underline{\beta},\vec{\gamma},\mu)\;
   \longrightarrow\; {\Bbb Z}_{\ge 0}\,.
 $$

\bigskip

\noindent
{\bf Definition 7.1.5 [layer-structure stratification].} {\rm
 The (finite) collection of the pre-image of the elements of
  ${\Bbb Z}_{\ge 0}$ under the above correspondence gives,
  by definition, the {\it layer-structure stratification} of
  $\overline{\cal M}_{(g,h),(n,\vec{m})}
            (Y,L\,|\,\underline{\beta},\vec{\gamma},\mu)$.
 The stratum
  ${\cal M}_{(g,h),(n,\vec{m})}^{(i)}
           (Y,L\,|\,\underline{\beta},\vec{\gamma},\mu)$
  associated to $i\in {\Bbb Z}_{\ge 0}$
  is called the stratum of {\it depth $i$}.
 A standard Kuranishi structure ${\cal K}_0$ on
  $\overline{\cal M}_{(g,h),(n,\vec{m})}
            (Y,L\,|\,\underline{\beta},\\  \vec{\gamma},\mu)$
  restricts to a Kuranishi
   structure-in-${\cal C}_{\spsccw,0}^{(i)}$
   ${\cal K}_0^{(i)}$  on
   ${\cal M}_{(g,h),(n,\vec{m})}^{(i)}
          (Y,L\,|\,\underline{\beta},\vec{\gamma},\mu)$
  as follows:
   \begin{itemize}
    \item[$\cdot$]
     Let
      ${\cal K}/B$ be a standard Kuranishi
       structure-in-${\cal C}_{\spsccw}$ on
       $\overline{\cal M}_{(g,h),(n,\vec{m})}
                    (W/B,L\,|\,[\beta],\vec{\gamma},\mu)/B$
       that gives the Kuranishi structure-in-${\cal C}_{\spsccw,0}$
       ${\cal K}_0$,
      $\rho  \in
         {\cal M}_{(g,h),(n,\vec{m})}^{(i)}
           (Y,L\,|\,\underline{\beta},\vec{\gamma},\mu)$,\\
      $(V_{\rho},E_{V_{\rho}},\Gamma_{V_{\rho}};
                             s_{V_{\rho}},\psi_{V_{\rho}})/B$
       be a Kuranishi neighborhood of $\rho$ from ${\cal K}$ with
       $\rho$ treated as a point in
       $\overline{\cal M}_{(g,h),(n,\vec{m})}
                         (W/B,L\,|\,[\beta],\vec{\gamma},\mu)/B$.
     Then, by definition, the Kuranishi
      neighborhood-in-${\cal C}_{\spsccw,0}$ of $\rho$ from
      ${\cal K}_{0}$ is given by
      $$
       (V_{\rho,0},\Gamma_{V_{\rho,0}}=\Gamma_{V_{\rho}},
         E_{V_{\rho,0}}=E_{V_{\rho}}|_{V_{\rho,0}};
         s_{V_{\rho,0}}=s_{V_{\rho}}|_{V_{\rho,0}},
         \psi_{V_{\rho,0}}=\psi_{V_{\rho}}|_{V_{\rho,0}})\,,
      $$
      where
       $V_{\rho,0}$ is the central fiber of $V_{\rho}/B$,
        which is invariant under $\Gamma_{V_{\rho}}$, and
       $\Gamma_{V_{\rho,0}}$ is $\Gamma_{V_{\rho}}$
        that acts on $V_{\rho,0}$.

     \item[$\cdot$]
      By construction $V_{\rho}$ also fibers over $B[i]$.
      Let $V_{\rho,0}^{(i)}$ be the central fiber of $V_{\rho}/B[i]$;
      then
       $V_{\rho,0}^{(i)}$ is $\Gamma_{V_{\rho}}$-invariant and
       the restriction
       $$
        (\, V_{\rho,0}^{(i)}\;,\;
           \Gamma_{V_{\rho,0}^{(i)}}=\Gamma_{V_{\rho}}\;,\;
           E_{V_{\rho,0}^{(i)}}=E_{V_{\rho}}|_{V_{\rho,0}^{(i)}}\;;\;
           s_{V_{\rho,0}^{(i)}}=s_{V_{\rho}}|_{V_{\rho,0}^{(i)}}\;,\;
           \psi_{V_{\rho,0}^{(i)}}=\psi_{V_{\rho}}|_{V_{\rho}^{(i)}}\, )\,.
       $$
       define a Kuranishi neighborhood-in-${\cal C}_{\spsccw,0}^{\,(i)}$
       of  $\rho  \in
                  {\cal M}_{(g,h),(n,\vec{m})}^{(i)}
                     (Y,L\,|\,\underline{\beta},\vec{\gamma},\mu)$.
       The system of transition data in ${\cal K}$ restricts to
        a system of transition data for such system of Kuranishi
        neighborhoods-in-${\cal C}_{\spsccw,0}^{\,(i)}$   for
        ${\cal M}_{(g,h),(n,\vec{m})}^{(i)}
           (Y,L\,|\,\underline{\beta},\vec{\gamma},\mu)$.
        This defines ${\cal K}_0^{(i)}$.

   \end{itemize}
 We shall call such ${\cal K}_0^{(i)}$  a
  {\it standard Kuranishi structure} on
  ${\cal M}_{(g,h),(n,\vec{m})}^{(i)}
         (Y,L\,|\,\underline{\beta},\vec{\gamma},\mu)$.
 Note that in the above description, $V_{\rho}^{(i)}$ has codimension
   $2i$ in $V_{\rho,0}$; thus
  $\vdim {\cal K}_0^{(i)}=\vdim {\cal K}_0-2i$.
 We say that
  the stratum
  ${\cal M}_{(g,h),(n,\vec{m})}^{(i)}
     (Y,L\,|\,\underline{\beta},\vec{\gamma},\mu)$
   has {\it virtual codimension} $2i$ (everywhere)  in
   $\overline{\cal M}_{(g,h),(n,\vec{m})}
               (Y,L\,|\,\underline{\beta},\vec{\gamma},\mu)$
 Note also that ${\cal K}_0^{(0)}={\cal K}_0$.
} % end-definition

\bigskip

Similarly, the composition
 $$
  [f]\;  \longmapsto\;  \tau_{[f]}\;
         \longmapsto\;  |\Image(\layer_{\tau_{[f]}})|
 $$
 gives a correspondence
 $$
  \overline{\cal M}_{(g,h),(n+l(\vec{s}),\vec{m})}
        (Z,L;D\,|\,\beta^{\prime},\vec{\gamma},\mu;\vec{s})\;
   \longrightarrow\; {\Bbb Z}_{\ge 0}\,.
 $$
This defines a {\it layer-structure stratification}
 $$
   \left\{\,
     {\cal M}_{(g,h),(n+l(\vec{s}),\vec{m})}^{(i)}
       (Z,L;D\,|\,\beta^{\prime},\vec{\gamma},\mu^{\prime};\vec{s})\,
   \right\}_{i\in{\scriptsizeBbb Z}_{\ge 0}}
 $$
 of  $\overline{\cal M}_{(g,h),(n+l(\vec{s}),\vec{m})}
       (Z,L;D\,|\,\beta^{\prime},\vec{\gamma},\mu^{\prime};\vec{s})$.
Given a standard Kuranishi structure-in-${\cal C}_{\spsccw}^{\,\prime}$
 ${\cal K}^{\prime}$ for
 $\overline{\cal M}_{(g,h),(n+l(\vec{s}),\vec{m})}
   (Z,L;D\,|\,\beta^{\prime},\vec{\gamma},\mu^{\prime};\vec{s})$,
the same {\it take-central-fiber-then-restrict} construction in
 Definition 7.1.5   % Definition [layer-structure stratification]
 gives a {\it standard
 Kuranishi structure}-in-${\cal C}_{\spsccw}^{\,\prime,\,(i)}$
 ${\cal K}^{\prime,(i)}$ on
 ${\cal M}_{(g,h),(n+l(\vec{s}),\vec{m})}^{(i)}
  (Z,L;D\,|\,\beta^{\prime},\\  \vec{\gamma},\mu^{\prime};\vec{s})$.
The {\it depth-$i$ stratum}
 ${\cal M}_{(g,h),(n+l(\vec{s}),\vec{m})}^{(i)}
  (Z,L;D\,|\,\beta^{\prime},\vec{\gamma},\mu^{\prime};\vec{s})$
 has {\it virtual codimension} $2i$ in
 $\overline{\cal M}_{(g,h),(n+l(\vec{s}),\vec{m})}
  (Z,L;D\,|\,\beta^{\prime},\vec{\gamma},\mu^{\prime};\vec{s})$.

\bigskip

\noindent
{\bf Lemma 7.1.6 [unique equivalence class].} {\it
 Any two standard Kuranishi structures on ${\cal M}$ are equivalent,
  where ${\cal M}$ is any of the following moduli spaces:
  $$
  \begin{array}{lcl}
   \overline{\cal M}_{(g,h),(n,\vec{m})}
         (W/B,L\,|\,[\beta],\vec{\gamma},\mu)/B\,,
    &&  \overline{\cal M}_{(g,h),(n,\vec{m})}
         (W_{\lambda},L\,|\,[\beta],\vec{\gamma},\mu)\,,\,
        \lambda\in B-\{0\}\,,\\[.6ex]
   \overline{\cal M}_{(g,h),(n,\vec{m})}
    (Y,L\,|\,\underline{\beta},\vec{\gamma},\mu)\,,
    &&  {\cal M}_{(g,h),(n,\vec{m})}^{(i)}
          (Y,L\,|\,\underline{\beta},\vec{\gamma},\mu)\,, \\[.6ex]
   \overline{\cal M}_{(g,h),(n+l(\vec{s}),\vec{m})}
    (Z,L;D\,|\,\beta^{\prime},\vec{\gamma},\mu^{\prime};\vec{s})\,,
    && {\cal M}_{(g,h),(n+l(\vec{s}),\vec{m})}^{(i^{\prime})}
        (Z,L;D\,|\,\beta^{\prime},\vec{\gamma},\mu^{\prime};\vec{s})\,.
 \end{array}
 $$
} % end-lemma

%\bigskip

\noindent{\it Proof.}
 For
  ${\cal M} = \overline{\cal M}_{(g,h),(n,\vec{m})}
                       (W/B,L\,|\,[\beta],\vec{\gamma},\mu)/B$,
 let ${\cal K}_i$ be the Kuranishi structure associated to
  a fine system of saturated obstruction local bundles ${\cal E}_i$
  $i=1,\,2$; cf.\ Sec.~5.4.
 Then there exists another fine system ${\cal E}_3$ of saturated
  obstruction local bundles so that both
  ${\mathbf F}({\cal E}_1)$ and ${\mathbf F}({\cal E}_2)$ are
  orbifold sub-fibrations of ${\mathbf F}({\cal E}_3)$.
 The lemma for
  ${\cal M} = \overline{\cal M}_{(g,h),(n,\vec{m})}
                       (W/B,L\,|\,[\beta],\vec{\gamma},\mu)/B$
  follows then from the construction in Sec.~5.4.
 Similarly for all other choices of ${\cal M}$ in the list.

\noindent\hspace{15cm}$\Box$

%%%%%%%%%%%%%%%%%%%%%%%
% \bigskip
%
% \marginpar{\raggedright\tiny $\bullet$
%     {\bf Q.} {\it Is the following true?} Think more carefully. }
% We remark that, for
%  ${\cal M}=\overline{\cal M}_{(g,h),(n,\vec{m})}
%              (Y,L\,|\,\underline{\beta},\vec{\gamma},\mu)$  or
%  $\overline{\cal M}_{(g,h),(n+l(\vec{s}),\vec{m})}
%           (Z,L;D\,|\,\beta^{\prime},\vec{\gamma},\mu;\vec{s})$,
%  denote its depth-$i$ stratum by ${\cal M}^{(i)}$;
% then the closure
%  $\overline{{\cal M}^{(i)}}=\amalg_{j\ge i}{\cal M}^{(j)}$.
%%%%%%%%%%%%%%%%%%%%%%%

\bigskip

\begin{flushleft}
{\bf Topological spaces with a Kuranishi structure:
     morphisms and fibered products.\footnote{Our definitions here are
                                tailored to what we have explicitly,
                                 what we are allowed to do in these cases,
                                 and what we are aiming for.
                                There is still room for further
                                 polishments/generalizations
                                 of these notions/definitions.}}
\end{flushleft}
We digress here to define two fundamental notions that
 we did not truly need until now:
 {\it morphisms} and {\it fibered products} of
  topological spaces with a Kuranishi structure.
These two notions are fundamental in any category of
 spaces/geometries.

\bigskip

\noindent {\bf Definition 7.1.7 [Kuranishi structure: morphism].}
{\rm
 Let $X_i$ be a topological space with a Kuranishi structure
  ${\cal K}_{i,0}$ modelled in a category ${\cal C}$, $i=1,\,2$.
 A {\it morphisms} from $(X_1,{\cal K}_{1,0})$ to $(X_2,{\cal K}_{2,0})$
  is a {\it continuous map}
   $\;\varphi:X_1\rightarrow X_2\;$
   together with a {\it tuple of systems of morphisms}
   $\;\varphi^{\sharp}
     :=(\varphi_{V_{\cdot}}\,,\, \varphi_{\Gamma_{\cdot}}\,,\,
     \varphi_{E_{\cdot}})\,:\, {\cal K}_1\rightarrow{\cal K}_2\,$,
    where ${\cal K}_1\sim{\cal K}_{1,0}$ and
          ${\cal K}_2\sim{\cal K}_{2,0}$,
    consisting of a system of
     continuous maps
      $\varphi_{V_{\cdot}}:V_{x_1}\rightarrow V_{\varphi(x_1)}$,
     group homomorphisms
      $\varphi_{\Gamma_{\cdot}}:\Gamma_{V_{x_1}}\rightarrow
                               \Gamma_{V_{\varphi(x_1)}}$, and
     $\varphi_{\Gamma_{\cdot}}$-equivariant bundle maps
      $\varphi_{E_{\cdot}}: E_{V_{x_1}}\rightarrow E_{V_{\varphi(x_1)}}$
      that covers $\varphi_{V_{\cdot}}$,
    such that\footnote{The necessity of passing to equivalent
                          Kuranishi structures to define morphisms is
                           enforced on us when one considers
                           the simplest case:
                            the notion of {\it embeddings} of a
                            topological space-with-Kuranishi-structure
                            to another.
                          This also makes the definition ring
                           more compatibly with its parallel
                           in algebraic geometry.
                          There one has the notion of two-term
                           locally-free resolutions of
                           a perfect tangent-obstruction complex
                           on the moduli stack in question.
                          Morphisms between such complexes are at
                           the level of derived categories of
                           coherent sheaves on the moduli stacks.
                          In particular, they have to pass to
                           quasi-isomorphisms of chain complexes,
                           rather than directly on the two
                           chain complexes one wants to compare.}
  \begin{itemize}
   \item[$(1)$] [compatibility on each Kuranishi neighborhood]$\,$:
    \begin{itemize}
    \item[]
     $\varphi_{E_{\cdot}}\circ s_p
       = s_{\varphi(p)} \circ\varphi_{V_{\cdot}}\;$  on $V_p\,$,
     \hspace{1em}
     $\varphi \circ \psi_p
       = \psi_{\varphi(p)} \circ\varphi_{V_{\cdot}}\;$
      on $s_p^{-1}(0)\subset V_p$ \hspace{1em}
     for $p\in X_1$;
    \end{itemize}
  \end{itemize}
  \begin{itemize}
   \item[$(2)$] [gluability: compatibility with transition data]$\,$:
    \begin{itemize}
     \item[]
      $\varphi_{V_{\cdot}}(V_{qp})\subset V_{\varphi(q)\varphi(p)}\,$,
       \hspace{3.5em}
      $\varphi_{\Gamma_{\cdot}}\circ h_{qp}
       = h_{\varphi(q)\varphi(p)}\circ \varphi_{\Gamma_{\cdot}}\,$,\\
      $\varphi_{V_{\cdot}}\circ \phi_{qp}
       = \phi_{\varphi(q)\varphi(p)}\circ \varphi_{V_{\cdot}}\,$,
       \hspace{1em}
      $\varphi_{E_{\cdot}}\circ \hat{\phi}_{qp}
       = \hat{\phi}_{\varphi(q)\varphi(p)}\circ \varphi_{E_{\cdot}}\,$.
    \end{itemize}
  \end{itemize}
 For convenience, we will denote a morphism as
  $(\varphi,\varphi^{\sharp}):
   (X_1,{\cal K}_{1,0})\rightarrow (X_2,{\cal K}_{2,0})$
  with it understood that $\varphi^{\sharp}$ is defined
   subject to passing to an equivalent Kuranishi structure.
} % end-definition

\bigskip

\noindent
{\bf Definition/Example 7.1.8 [embedding].}
 A morphism
  $(\varphi,\varphi^{\sharp}):
       (X_1,{\cal K}_1)\rightarrow (X_2,{\cal K}_2)$
  is called an {\it embedding}
  if both $\varphi$ and $\varphi^{\sharp}$ are embeddings.

\bigskip

\noindent
{\bf Definition/Example 7.1.9 [covering map].}
 A morphism
  $(\varphi,\varphi^{\sharp}):
       (X_1,{\cal K}_1)\rightarrow (X_2,{\cal K}_2)$
  is called a {\it covering map}
  if $\varphi$ is a covering map and
     $\varphi^{\sharp}$ is an isomorphism\footnote{By this
                          we mean that all maps in $\varphi^{\sharp}$
                           are isomorphisms.
                          Note that $\varphi^{\sharp}$ alone sees
                           only the local properties of the topology.
                          That maps in $\varphi^{\sharp}$ are all
                           isomorphisms implies only that
                           $\varphi:X_1\rightarrow X_2$
                           is a local isomorphism.}.
 In this case,\\ $\vdim {\cal K}_1=\vdim{\cal K}_2$.

\bigskip

\noindent
{\bf Definition/Example 7.1.10 [virtual bundle map].}
 Given a topological space $S$, we shall regard it
  also as a topological space with the
  {\it trivial Kuranishi structure}
  ${\cal K}^{\trivial}$ that consists of exactly
   one Kuranishi neighborhood
   $(S,\{e\},{\mathbf 0}_S:=S\times\{0\};0,\Id_S)$.
 A morphism
  $(\varphi,\varphi^{\sharp}):
       (X,{\cal K})\rightarrow S=(S,{\cal K}^{\trivial})$,
  is called a {\it virtual bundle map}
 if $\varphi:X\rightarrow S$ is continuous and
   $\varphi^{\sharp}:{\cal K}\rightarrow {\cal K}^{\trivial}$
   is a bundle map\footnote{By this we mean that each
                           $\varphi_{V_{\cdot}}:V_p\rightarrow S$,
                             $p\in X $,
                            in $\varphi^{\sharp}$ is a {\it bundle map}
                            (i.e.\ locally trivial fibration)
                            {\it over a non-empty open subset} of $S$.}.
 Note that, in this case,
  the $\Gamma_{V_p}$-action on $V_p$ leaves each fiber of
   $V_p\rightarrow S$ invariant.

\bigskip

\noindent
{\bf Definition 7.1.11 [Kuranishi structure: fibered product].} {\rm
 Let $S$ be a topological space with
  the trivial Kuranishi structure ${\cal K}^{\trivial}$.
 Given two virtual bundle maps
  $$
   (X_1,{\cal K}_1)\;
    \stackrel{(\varphi_1,\varphi_1^{\sharp})}{\longrightarrow}\;
    S\;\stackrel{(\varphi_2,\varphi_2^{\sharp})}{\longleftarrow}\;
   (X_2,{\cal K}_2)\,,
  $$
 define the {\it fibered product}
   $(X_1\times_S X_2,{\cal K}_1\times_S{\cal K}_2)$
  of $(X_1,{\cal K}_1)$ and $(X_2,{\cal K})$ {\it over $S$}
 to be the topological space
  $$
   X_1\times_S X_2\;
    :=\; (\varphi_1\times \varphi_2)^{-1}(\Delta_S)
   \subset X_1\times X_2\,,
  $$
  where
   $\varphi_1\times\varphi_2:X_1\times X_2\rightarrow S\times S$ and
   $\Delta_S\subset S\times S$ is the diagonal,
 equipped with the following Kuranishi structure:
  \begin{itemize}
   \item[$(1)$]
    [the induced Kuranishi neighborhood
                           at $(p_1,p_2)\in X_1\times_S X_2$]$\,$:
    \begin{itemize}
     \item[$\cdot$]
      define $V_{(p_1,p_2)}:=V_{p_1}\times_S V_{p_2}$ and
      let  $V_{p_1}\stackrel{\pi_1}{\leftarrow}
            V_{p_1}\times_S V_{p_2}
            \stackrel{\pi_2}{\rightarrow} V_{p_2}$
       be the projection maps;

     \item[$\cdot$]
      the diagonal action of
       $\Gamma_{V_{p_1}}\times \Gamma_{V_{p_2}}$ on
       $V_{p_1}\times V_{p_2}$ leaves
       $V_{(p_1,p_2)}=V_{p_1}\times_S V_{p_2}$ invariant,
      define $\Gamma_{V_{(p_1,p_2)}}
                =\Gamma_{V_{p_1}}\times \Gamma_{V_{p_2}}$
        now acting on $V_{(p_1,p_2)}$;

     \item[$\cdot$]
      let  $E_{V_{(p_1,p_2)}}
            := \pi_1^{\ast}E_{V_{p_1}}\oplus \pi_2^{\ast}E_{V_{p_2}}$
        on $V_{(p_1,p_2)}$,
      then the induced $\Gamma_{V_{(p_1,p_2)}}$-action on
        $E_{V_{(p_1,p_2)}}$ is equivariant;

     \item[$\cdot$]
      let  $s_{(p_1,p_2)}
             =(\pi_1^{\ast}s_{p_1}\,,\,\pi_2^{\ast}s_{p_2})$,
      then $s_{(p_1,p_2)}$ is a $\Gamma_{V_{(p_1,p_2)}}$-invariant
       section of $E_{V_{(p_1,p_2)}}$;

     \item[$\cdot$]
      let
       $\psi_{(p_1,p_2)}
         =(\,\psi_{p_1}\times \psi_{p_2}\,:\,
                 s_{p_1}^{-1}(0)\times s_{p_2}^{-1}(0)
              \rightarrow X_1\times X_2\,)|
                                     _{V_{p_1}\times_S V_{p_2}}$,
      then $\psi_{(p_1,p_2)}$ is a map from
       $s_{(p_1,p_2)}^{-1}(0)$ to $X_1\times_S X_2$.
    \end{itemize}
   The $5$-tuple
    $(V_{(p_1,p_2)}\,,\,\Gamma_{V_{(p_1,p_2)}}\,,\,
      E_{V_{(p_1,p_2)}}\,;\, s_{(p_1,p_2)}\,,\, \psi_{(p_1,p_2)})$
    defined above is called the {\it induced Kuranishi neighborhood}
     of $(p_1,p_2)\in X_1\times_S X_2$ from ${\cal K}_1$ and
     ${\cal K}_2$.
   Define ${\mathfrak N}^{(0)}$ to be the system of Kuranishi
    neighborhoods of $X_1\times_S X_2$ thus constructed.
  \end{itemize}
  \begin{itemize}
   \item[$(2)$] [transition data]$\,$:
    \begin{itemize}
     \item[$\cdot$]
      the diagonal product construction defines a canonical
       Kuranishi structure ${\cal K}_1\times{\cal K}_2$ on
       $X_1\times X_2$;
      the Kuranishi neighborhoods for $X_1\times_S X_2$,
        as constructed above,
       are embedded in the Kuranishi neighborhoods in
       ${\cal K}_1\times {\cal K}_2$;
      the canonical transition data in ${\cal K}_1\times {\cal K}_2$
       restricts to a system ${\mathfrak N}^{(1)}$
       of transition data for ${\mathfrak N}^{(0)}$.
    \end{itemize}
  \end{itemize}
 Define ${\cal K}_1\times_S{\cal K}_2
           =({\mathfrak N}^{(0)}\,,\,{\mathfrak N}^{(1)})$.
 When $S=\{\pt\}$, we call $(X_1\times X_2,{\cal K}_1\times {\cal K}_2)$
  the {\it direct product}, or simply the {\it product}, of
  $(X_1,{\cal K}_1)$ and $(X_2,{\cal K}_2)$.
} % end-definition

\bigskip

By construction, there are
 a {\it tautological virtual bundle map}
  $$
   (\varphi_1\times_S\varphi_2\,,\,
    \varphi_1^{\sharp}\times_S\varphi_2^{\sharp})\;:
     (X_1\times_S X_2\,,\, {\cal K}_1\times_S {\cal K}_2)\;
    \longrightarrow\; S\,,
  $$
 an embedding morphism
  $(X_1\times_S X_2,{\cal K}_1\times_S{\cal K}_2)
     \rightarrow (X_1\times X_2, {\cal K}_1\times{\cal K}_2)$,  and
 projection morphisms
  $$
   (X_1,{\cal K}_1)\;
    \stackrel{(\pi_1,\pi_1^{\sharp})}{\longleftarrow}\;
    (X_1\times_S X_2, {\cal K}_1\times_S{\cal K}_2)\;
     \stackrel{(\pi_2,\pi_2^{\sharp})}{\longrightarrow}\;
     (X_2,{\cal K}_2)\,.
  $$
Note that
 $\;\vdim({\cal K}_1\times_S {\cal K}_2)
   = \vdim{\cal K}_1 + \vdim{\cal K}_2 -\dimm S\;$
 when $S$ is a manifold and both ${\cal K}_1$ and ${\cal K}_2$
  are modelled on the category of CW-complexes.

\bigskip

\begin{flushleft}
{\bf The degeneration-gluing relations of Kuranishi structures.}
\end{flushleft}
We are now ready to give the degeneration-gluing relations
 of Kuranishi structures of the several moduli spaces
 that occur in the study.

The following bookkeeping graphs are adapted from [Li1: Sec.~4.2]:

\bigskip

\noindent
{\bf Definition 7.1.12 [admissible weighted graph].} {\rm
 Given a relative pair $(Z,L;D)$ with a symplectic/totally-real
   submanifold,
  an {\it admissible weighted graph} $\Gamma$ for $(Z,L;D)$
  is a graph without edges together with the following data:
 \begin{itemize}
 \item [(1)]
  an ordered collection of {\it hands}, {\it fingers}\footnote{The
                                    order of fingers is lexicographic:
                                    first by the order of the hands
                                    they are attached to and then by
                                    the order within each group that
                                    are attached to the same hand.},
   and {\it legs};
  an ordered collection of {\it weighted roots};
  a {\it relative index function} and two {\it weight functions} on
   the vertex set
   $\mu^{\prime}:V(\Gamma)\rightarrow {\Bbb Z}$,
   $g:V(\Gamma)\rightarrow {\Bbb Z}_{\ge 0}$, and
   $b:V(\Gamma)\rightarrow H_2(Z,L;{\Bbb Z})$;
  a {\it weight function} on the ordered set of hands
   $\gamma: H(\Gamma)\rightarrow H_1(L;{\Bbb Z})$
   such that $\partial b(v)=\sum_{\tinybullet} h_{v,\tinybullet}$,
   where $v\in V(\Gamma)$ and the sum is over the ordered subset
    of hands that are attached to $v$;

  \item[(2)]
   $\Gamma$ is {\it relatively connected} in the sense that either
    $|V(\Gamma)|=1$ or each vertex in $V(\Gamma)$ has at least one root
    attached to it.
 \end{itemize}
} % end-definition

\bigskip

\noindent
{\bf Definition 7.1.13 [admissible quadruple].} {\rm
 Given a gluing $(Y,L)=(Y_1,L_1)\cup_D (Y_2,L_2)$ of relative pairs
  from a symplectic cut,
 let $\Gamma_1$ and $\Gamma_2$ be a pair of admissible weighted graphs
  for $(Y_1,L_1;D)$ and $(Y_2,L_2;D)$ respectively.
 Suppose that $\Gamma_1$ and $\Gamma_2$ have
  identical number $l$ of roots, $h_1$-many and $h_2$-many hands,
  % $(m_{1,1}\,,\,\cdots\,,\,m_{1,h_1})$-many\footnote{Here, $m_{1,j}$
  %                                   is the number of fingers that are
  %                                   attached to the $j$-th hand.}  and
  % $(m_{2,1}\,,\,\cdots\,,\,m_{2,h_2})$-many fingers, and
  $n_1$-many and $n_2$-many legs respectively.
 Let
  $h=h_1+h_2$,
  % $\vec{m}=(m_1\,,\,\cdots\,,\,m_h)$ that is identical to ????,
  $n=n_1+n_2$,
  $I_{\hand}\subset \{1,\,\ldots,\, h\}$ be a set of $h_1$ elements, and
  $I_{\leg}\subset \{1,\,\ldots,\, n\}$ be a set of $n_1$ elements.
 Then $(\Gamma_1,\Gamma_2,I_{\hand}, I_{\leg})$ is called
 an {\it admissible quadruple} if
  the following conditions hold:
  \begin{itemize}
   \item[(1)]
    the weights on the roots of $\Gamma_1$ and $\Gamma_2$ coincide:
     $r_{1,i}=r_{2,i}$, $i=1,\,\ldots,\, l\,$;

   \item[(2)]
    after connecting the $i$-th root of $\Gamma_1$ and the $i$-th root
    of $\Gamma_2$ for all $i$, the resulting new graph with $h$ hands,
    (accompanying fingers), $n$ legs and no roots is connected.
  \end{itemize}
 Re-ordering of roots defines an equivalence relation $\sim$ on the set
  $\Omega$ of admissible quadruples.
 Define $\bar{\Omega}:=\Omega/\sim$.
 Given an admissible quadruple
   $\eta=(\Gamma_1,\Gamma_2,I_{\hand},I_{\leg})$,
  denote by $\Per_r(\eta)$ the set of permutations of the roots
   in $\Gamma_1$ that leaves $\eta$ unchanged.
} % end-definition

\bigskip

Note that $I_{\hand}$ determines the order of the hands
  on the graph from gluing paired roots
 by the unique bijection
  $\{1,\,\cdots\,,\,h_1\}\amalg\{1,\,\cdots\,,\,h_2\}
     \rightarrow \{1,\,\cdots\,,\,h\}$
  such that it preserves the orders of both
   $\{1,\,\cdots\,,\,h_1\}$ and $\{1,\,\cdots\,,\,h_2\}$
  and that the image of $\{1,\,\cdots\,,\,h_1\}$ is $I_{\hand}$.
The order of fingers on the glued graph is then determined
 lexicographically.
Similarly, $I_{\leg}$ determines the order of legs on the glued graph.

Given an admissible quadruple
 $\eta=(\Gamma_1,\Gamma_2,I_{\hand},I_{\leg})$ as above with
 $(Y,L)=(Y_1,L_1)\cup _D (Y_2,L_2)=$ the degenerate fiber $W_0$
 of $W/B$, one has
 \begin{itemize}
  \item[$\cdot$] the {\it genus function}
   $$
    \mbox{$g(\eta)\; :=\;
         l+1-|V(\Gamma_1\,
               \raisebox{.2ex}{\scriptsize $\coprod$}\, \Gamma_2)|
         + \sum_{v\in V(\Gamma_1)\cup V(\Gamma_2)}\,g(v)\;
         \in\; {\Bbb Z}_{\ge 0}\,,$}
   $$

  \item[$\cdot$] the {\it curve-class function}
   $$
    \mbox{$b(\eta)\;  :=\;
     \iota_{1,\ast}\left(\sum_{v\in V(\Gamma_1)}\,b_{\Gamma_1}(v)\right)
      + \iota_{2,\ast}\left(
          \sum_{v\in V(\Gamma_2)}\, b_{\Gamma_2}(v)\right)\;
     \in\; H_2(Y,L;{\Bbb Z})\,,$}
   $$
   where $\iota_i:(Y_i,L_i)\hookrightarrow (Y,L)$, $i=1,\,2\,$,
   are the inclusion maps,

  \item[$\cdot$]
   the total {\it index}
    $\mu(\eta) =     \sum_{v\in V(\Gamma_1)}\mu^{\prime}(v)\,
                 +\, \sum_{v\in V(\Gamma_2)}\mu^{\prime}(v)\,$.
 \end{itemize}
Let
 $\vec{m}(\eta) :=(m_1,\,\cdots,\,m_h)$
  be the tuple of numbers of fingers attached to
   hands $\in H(\Gamma_1)\cup H(\Gamma_2)$ and
 $\vec{\gamma}(\eta)$ be the tuple of values of
  $\gamma_1\cup\gamma_2:
    H(\Gamma_1)\cup H(\Gamma_2)\rightarrow H_1(L;{\Bbb Z})$,
 both with respect to the order on
  $H(\Gamma_1)\cup H(\Gamma_2)$ specified by $I_{\hand}$.
Define the {\it type} of $\eta$ to be
 $$
  |\eta|\; :=\;
   ( (g(\eta),\, h_1+h_2),\, (n_1+n_2,\, \vec{m}(\eta))\,|\,
      b(\eta),\, \vec{\gamma}(\eta),\, \mu(\eta))\,.
 $$

For each $\eta=(\Gamma_1,\Gamma_2,I_{\hand},I_{\leg})$,
  with $l$-many roots,
 such that
 $|\eta| =
  ((g,h),(n,\vec{m})\,|\,\underline{\beta},\vec{\gamma},\mu)$,
 there are five moduli spaces of stable map associated to it:

 \vspace{-1.6ex}
 {\small
 $$
  \begin{array}{l}
   \overline{\cal M}(Y_1,L_1;D\,|\,\Gamma_1)\,,\hspace{2em}
    \overline{\cal M}(Y_2,L_2;D\,|\,\Gamma_2)\,,\hspace{2em}
    \overline{\cal M}(Y_1,L_1;D\,|\,\Gamma_1)   \times_{D^l}
      \overline{\cal M}(Y_2,L_2;D\,|\,\Gamma_2)\,,  \\[1ex]
   \mbox{sub-orbifolds
    $\;\overline{\cal M}((Y_1,L_1;D)\amalg (Y_2,L_2;D)\,|\,\eta)\;$
     and
    $\;\overline{\cal M}(Y,L\,|\,\eta)\;$ of
         $\;\overline{\cal M}_{(g,h),(n,\vec{m})}
           (Y,L\,|\,\underline{\beta},\vec{\gamma},\mu)$}\,.
  \end{array}
 $$
 \normalsize We}} % end-small
 explain each of these spaces and their standard
 Kuranishi structures below.

Let $\Gamma$ be an admissible weighted graph for $(Z,L;D)$.,
The restriction of the data encoded by $\Gamma$ to
 each vertex $v\in |\Gamma|$ specifies a unique type data
 $((g_v,h_v),(n_v,\vec{m}_v)\,|\,
     \underline{\beta}_v,\vec{\gamma}_v,\mu^{\prime}_v;\vec{s}_v)$.
Define the moduli space by the direct product:
 $$
  \overline{\cal M}(Z,L;D\,|\,\Gamma)\;
  :=\; \prod_{v\in |\Gamma|}\,
       \overline{\cal M}_{(g_v,h_v),(n_v,\vec{m}_v)}(Z,L;D\,|\,
       \underline{\beta_v},\vec{\gamma}_v,\mu^{\prime}_v;\vec{s}_v)\,.
 $$
A standard Kuranishi structure
 ${\cal K}^{\prime}_{(Z,L;D\,|\,\Gamma)}$
 on $\overline{\cal M}(Z,L;D\,|\,\Gamma)$ is by definition
 the direct product of
 a standard Kuranishi structure on each moduli-space component
 $\overline{\cal M}_{(g_v,h_v),(n_v,\vec{m}_v)}(Z,L;D\,|\,
   \underline{\beta_v},\\  \vec{\gamma}_v,\mu^{\prime}_v;\vec{s}_v)$.
Let $l$ be the number of roots of $\Gamma$.
Then the saturatedness of the obstruction-space local bundles for each
 $\overline{\cal M}_{(g_v,h_v),(n_v,\vec{m}_v)}(Z,L;D\,|\,
      \underline{\beta_v},\vec{\gamma}_v,\mu^{\prime}_v;\vec{s}_v)$
 implies that there is a virtual bundle map
 $$
  ({\mathbf q}\,,\,{\mathbf q}^{\sharp})\;:\;
   (\overline{\cal M}(Z,L;D\,|\,\Gamma)\,,\,
    {\cal K}^{\prime}_{(Z,L;D\,|\,\Gamma)} )\;
   \longrightarrow\; (D^l,{\cal K}^{\trivial})\,.
 $$

Apply the above to $\Gamma_1$ and $\Gamma_2$
 from the admissible quadruple with $l$-many roots,
 one obtains the fibered-product moduli space
 $\overline{\cal M}(Y_1,L_1;D\,|\,\Gamma_1)
    \times_{D^l}  \overline{\cal M}(Y_2,L_2;D\,|\,\Gamma_2)$
  with a standard Kuranishi structure defined to be
  the fibered product
  ${\cal K}^{\prime}_{(Y_1,L_1;D\,|\,\Gamma_1)}
    \times_{D^l} {\cal K}^{\prime}_{(Y_2,L_2;D_2\,|\,\Gamma_2)}$.

Let
 $$
  \Phi_{\eta}:
   \overline{\cal M}(Y_1,L_1;D\,|\,\Gamma_1)  \times_{D^l}
   \overline{\cal M}(Y_2,L_2;D\,|\,\Gamma_2)\;
    \longrightarrow\;
   \overline{\cal M}_{(g,h),(n,\vec{m})}
      (Y,L\,|\,\underline{\beta},\vec{\gamma},\mu)\,
 $$
  be the gluing orbifold map,
 whose corresponding map at the underlying topological space
  is given by
  $$
   \begin{array}{l}
    (\, f_1:\Sigma_1\rightarrow (Y_{1,[k_1]},L_{1,[k_1]};D_{[k_1]})\,,\,
        f_2:\Sigma_2\rightarrow (Y_{2,[k_2]},L_{2,[k_2]};D_{[k_2]}) \,)
                                                            \\[.6ex]
    \hspace{12em}\longmapsto\;
    f=f_1\cup f_2: \Sigma\rightarrow (Y_{[k]},L_{[k]})\,,
     \hspace{1em}k=k_1+k_2\,,
   \end{array}
  $$
  where
   $\Sigma$ is the gluing $\Sigma_1\cup \Sigma_2$ of $\Sigma_1$ and
     $\Sigma_2$ along their paired distinguished marked points;
   $Y_{[k_1+k_2]}$ is the gluing of
    $(Y_{1,[k_1]},L_{1,[k_1]};D_{[k_1]})$ and
    $(Y_{2,[k_2]},L_{2,[k_2]};D_{[k_2]})$
    by  $D_{[k_1]}\simeq D\simeq D_{[k_2]}$.
Denote the image by
 $$
   \overline{\cal M}((Y_1,L_1;D)\amalg (Y_2,L_2;D)\,|\,\eta)
 $$
 with the induced sub-orbifold structure and
 the $C^{\infty}$-topology from
 $\overline{\cal M}_{(g,h),(n,\vec{m})}
     (Y,L\,|\,\underline{\beta},\vec{\gamma},\mu)$;
then $\Phi_{\eta}$ is an orbifold covering map of pure degree
 $|Per_r(\eta)|$ to
 $\overline{\cal M}((Y_1,L_1;D)\amalg (Y_2,L_2;D)\,|\,\eta)$.
A standard Kuranishi structure on
 $\overline{\cal M}((Y_1,L_1;D)\amalg (Y_2,L_2;D)\,|\,\eta)$
 can be constructed as follows.
Since $\Phi_{\eta}$ is a covering map, a Kuranishi neighborhood
 $(V_{\rho},\Gamma_{V_{\rho}}, E_{V_{\rho}};s_{\rho},\psi_{\rho})$
 of $\rho\in
       \overline{\cal M}((Y_1,L_1;D)\amalg (Y_2,L_2;D)\,|\,\eta)$
 can be taken to be a Kuranishi neighborhood of
 a $\rho^{\prime}\in\Phi_{\eta}^{-1}(\rho)$;
 i.e.\ via the fibered-product construction.
In this way one obtains a system
 ${\mathfrak N}^{(0)}_{((Y_1,L_1;D)\amalg (Y_2,L_2;D)\,|\,\eta)}$
 of Kuranishi neighborhoods on
 $\overline{\cal M}((Y_1,L_1;D)\amalg (Y_2,L_2;D)\,|\,\eta)$.
Assume that all these neighborhoods are small,
then the system $\{\iota_{(k_1,k_2)}\}_{k_1+k_2=k}$
 of almost-complex pseudo-embeddings

 \vspace{-1ex}
 {\small
 $$
  \begin{array}{l}
  \iota_{(k_1,k_2)}\; :\; \\[.6ex]
   \left(\,
    \left( (Y_1[k_1],L_1[k_1];D[k_1])\times A[k_2] \right)
     \cup_{D[k_1]\times A[k_2]\simeq A[k_1]\times D[k_2]}
     \left(A[k_1]\times (Y_2[k_2],L[k_2];D[2])\right)\,\right)
                         /(A[k_1]\times A[k_2])\;  \\[.6ex]
    \longrightarrow\; W[k]/B[k]
  \end{array}
 $$} % end-small

 \noindent
 induces a natural embedding of
 ${\mathfrak N}^{(0)}_{((Y_1,L_1;D)\amalg (Y_2,L_2;D)\,|\,\eta)}$
 into a standard Kuranishi structure ${\cal K}$ on
 $\overline{\cal M}_{(g,h),(n,\vec{m})}
     (Y,L\,|\,\underline{\beta},\vec{\gamma},\mu)$.
Here,
 $D[k_1]\times A[k_2]$ and $A[k_1]\times D[k_2]$ are glued
   via their canonical isomorphisms with
   $D\times A[k_1]\times A[k_2]=D\times A[k]$, and
 the pseudo-embedding
  $A[k]=A[k_1]\times A[k_2]\rightarrow B[k]$ is given by
  $(\vec{\lambda},\vec{\lambda}^{\prime})\mapsto
     (\vec{\lambda}\,,\,0\,,\,\vec{\lambda}^{\prime})$.
% After a re-adjustment of both
%  ${\mathfrak N}^{(0)}_{((Y_1,L_1;D)\amalg (Y_2,L_2;D)\,|\,\eta)}$
%  and ${\cal K}$ by passing to equivalences,  the ......
The transition data from ${\cal K}$ then restricts\footnote{Note
                               that in general one has to pass to an
                                equivalence to make a system of
                                Kuranishi neighborhoods gluable.
                               However, here the system
                                ${\mathfrak N}^{(0)}
                                  _{((Y_1,L_1;D)\amalg (Y_2,L_2;D)\,|\,
                                 \eta)}$
                                is descended from a covering morphism
                                of a Kuranishi structure.
                               We only need to know whether
                                the transition data also descends
                                in our case.
                               The latter is implied by
                                the existence of a natural embedding of
                                ${\mathfrak N}^{(0)}
                                 _{((Y_1,L_1;D)\amalg (Y_2,L_2;D)\,|\,
                                  \eta)}$ into ${\cal K}$, (i.e.\
                                an embedding at the level of
                                universal maps on universal curves).}
 to an transition data on
 ${\mathfrak N}^{(0)}_{((Y_1,L_1;D)\amalg (Y_2,L_2;D)\,|\,\eta)}$.
By construction, one has an embedding morphism
 $$
  \begin{array}{l}
   \left(\,
     \overline{\cal M}((Y_1,L_1;D)\amalg (Y_2,L_2;D)\,|\,\eta)\;,\;
      {\cal K}_{((Y_1,L_1;D)\amalg (Y_2,L_2;D)\,|\,\eta)}\,
    \right) \\[1ex]
   \hspace{3em}\longrightarrow\;
   (\overline{\cal M}_{(g,h),(n,\vec{m})}
       (Y,L\,|\,\underline{\beta},\vec{\gamma},\mu)\,,\,
    {\cal K})\,.
  \end{array}
 $$

With respect to ${\cal K}_{((Y_1,L_1;D)\amalg(Y_2,L_2;D)\,|\,\eta)}$,
 the covering map $\Phi_{\eta}$ lifts to a covering morphism
 $$
  \begin{array}{l}
   (\Phi_{\eta},\Phi_{\eta}^{\sharp})\,:\;
    \left(\,\rule{0ex}{2ex}\right.
     \overline{\cal M}(Y_1,L_1;D\,|\,\Gamma_1)
      \times_{D^l} \overline{\cal M}(Y_2,L_2;D\,|\,\Gamma_2)\;,\;
    {\cal K}^{\prime}_{(Y_1,L_1;D\,|\,\Gamma_1)}
      \times_{D^l} {\cal K}^{\prime}_{(Y_2,L_2;D\,|\,\Gamma_2)}\,
      \left.\rule{0ex}{2ex}\right)                      \\[1ex]
   \hspace{3em}
    \longrightarrow\;
    \left(\,
     \overline{\cal M}((Y_1,L_1;D)\amalg (Y_2,L_2;D)\,|\,\eta)\;,\;
      {\cal K}_{((Y_1,L_1;D)\amalg (Y_2,L_2;D)\,|\,\eta)}\,
    \right)\,.
  \end{array}
 $$
One can check that, with these standard Kuranishi structures,
 $$
  \begin{array}{l}
   \vdim \left(\,\rule{0ex}{2ex}\right.
     \overline{\cal M}(Y_1,L_1;D\,|\,\Gamma_1)
      \times_{D^l} \overline{\cal M}(Y_2,L_2;D\,|\,\Gamma_2)\,
      \left.\rule{0ex}{2ex}\right) \\[1ex]
   \hspace{2em}
   =\; \vdim\overline{\cal M}(Y_1,L_1;D\,|\,\Gamma_1)\,
      +\, \vdim\overline{\cal M}(Y_1,L_1;D\,|\,\Gamma_1)\,
      -\, 2\,l\,(N-1)                                       \\[1ex]
   \hspace{2em}
   =\; \vdim \overline{\cal M}_{(g,h),(n,\vec{m})}
               (Y,L\,|\,\underline{\beta},\vec{\gamma},\mu)\,,
   \end{array}
 $$
 where, recall that, $\dimm Y = 2N$.
This implies that
 $$
   \vdim\,
   \overline{\cal M}((Y_1,L_1;D)\amalg (Y_2,L_2;D)\,|\,\eta)\;
   =\; \vdim \overline{\cal M}_{(g,h),(n,\vec{m})}
               (Y,L\,|\,\underline{\beta},\vec{\gamma},\mu)\,.
 $$

Finally, let
 $$
  \overline{\cal M}(Y,L\,|\,\eta)
 $$
 be the same suborbifold
 $\overline{\cal M}((Y_1,L_1;D)\amalg (Y_2,L_2;D)\,|\,\eta)$
 of  $\overline{\cal M}_{(g,h),(n,\vec{m})}
           (Y,L\,|\,\underline{\beta},\vec{\gamma},\mu)$
 but with a Kuranishi structure constructed as follows.
Consider the defining embedding morphism
 $$
  (\overline{\cal M}_{(g,h),(n,\vec{m})}
      (Y,L\,|\,\underline{\beta},\vec{\gamma},\mu)\,,\,
   {\cal K}_0)\;
  \longrightarrow\;
  (\overline{\cal M}_{(g,h),(n,\vec{m})}
      (W/B,L\,|\,[\beta],\vec{\gamma},\mu)\,,\,
   {\cal K})\,.
 $$
Let ${\cal K} = ({\mathfrak N}^{(0)},{\mathfrak N}^{(1)}$.
For $\rho\in\overline{\cal M}(Y,L\,|\,\eta)$ from gluing
 $f_1:\Sigma_1\rightarrow (Y_{1,[k_1]},L_{1,[k_1]};D_{[k_1]})$
 and $f_2:\Sigma_2\rightarrow (Y_{2,[k_2]},L_{2,[k_2]};D_{[k_2]})$,
one has that
 $V_{\rho}\in{\mathfrak N}^{(0)}$ fibers over $B[k]$,where $k=k_1+k_2$.
Let
 $V_{\rho,\eta}\subset V_{\rho}$ be the preimage of the hyperplane
  $\{\vec{\lambda}=(\lambda_0,\,\cdots\,,\,\lambda_k)\,:\,
      \lambda_{k_1}=0\}\subset B[k]$
 under this fibration with the multiplicity of
 the irreducible components of fibers encoded.
Then $V_{\rho,\eta}$ is $\Gamma_{V_{\rho}}$-invariant.
Define
 $\Gamma_{V_{\rho,\eta}}=\Gamma_{V_{\rho}}$,
  now action on $V_{\rho,\eta}$;
 $E_{V_{\rho,\eta}}=E_{V_{\rho}}|_{V_{\rho,\eta}}$;
 $s_{\rho,\eta}=s_{\rho}|_{V_{\rho,\eta}}$; and
 $\psi_{\rho,\eta}=\psi_{\rho}|_{V_{\rho,\eta}}$.
Then the system ${\mathfrak N}_{\eta}^{(0)}$ of $5$-tuples
 $(V_{\rho,\eta}, \Gamma_{V_{\rho,\eta}}, E_{V_{\rho,\eta}};
   s_{\rho,\eta},\psi_{\rho,\eta})$
 thus constructed defines a system of Kuranishi neighborhood on
 $\overline{\cal M}(Y,L\,|\,\eta)$.
The system ${\mathfrak N}^{(1)}$ of transition data in ${\cal K}$
 restricts to give a system ${\mathfrak N}_{\eta}^{(1)}$ of
 transition data for ${\mathfrak N}_{\eta}^{(0)}$.
The pair
 ${\cal K}_{\eta}
   :=({\mathfrak N}_{\eta}^{(0)},{\mathfrak N}_{\eta}^{(1)})$
 thus defines a Kuranishi structure on
 $\overline{\cal M}(Y,L\,|\,\eta)$.
Kuranishi structures on $\overline{\cal M}(Y,L\,|\,\eta)$
 thus obtained will be called {\it standard} Kuranishi structures
 on $\overline{\cal M}(Y,L\,|\,\eta)$.
By construction, one also has:
 $$
   \vdim\,
   \overline{\cal M}(Y,L\,|\,\eta)\;
   =\; \vdim \overline{\cal M}_{(g,h),(n,\vec{m})}
               (Y,L\,|\,\underline{\beta},\vec{\gamma},\mu)\,.
 $$

The following theorem that relates
  these moduli spaces and their standard Kuranishi structures
 should be compared to [Li2: Corollary 3.13. Lemma 3.14, Theorem 3.15].
It is in effect a re-phrasing of [Li2] in terms of the
 Fukaya-Ono setting and at the level of Kuranishi structures,
 rather than of virtual fundamental classes or chains:

\bigskip

\noindent
{\bf Theorem 7.1.14 [degeneration-gluing$\,$: Kuranishi structure].}
{\it
 Regard $X$ as a fiber $W_{\lambda_0}$ of $W/B$
  over $\lambda_0\in B-\{0\}$.
 Recall the symplectic cut
  $\xi:(X,L)\rightarrow (Y,L)=(Y_1,L_1)\cup_D (Y_2,L_2)$
 Given a type $((g,h),(n,\vec{m})\,|\,[\beta],\vec{\gamma},\mu)$
  of stable maps to $(X,L)$,
 let $\underline{\beta}=\xi_{\ast}([\beta])\in H_2(Y,L;{\Bbb Z})$ and
  $\bar{\Omega}_{((g,h),(n,\vec{m})\,|\,
                        \underline{\beta},\vec{\gamma},\mu)}$
  be the equivalence of admissible quadruples $\eta$ such that
  $|\eta|=((g,h),(n,\vec{m})\,|\,\underline{\beta},\vec{\gamma},\mu)$.
 Then, the following statements hold,
  up to an equivalence of Kuranishi structures:
 \begin{itemize}
  \item[$(1)$]
   A standard Kuranishi structure ${\cal K}_{\lambda_0}$
    on $\overline{\cal M}_{(g,h),(n,\vec{m})}
                  (X,L\,|\,[\beta],\vec{\gamma},\mu)$  and
   a standard Kuranishi structure-in-${\cal C}_{\spsccw,0}$ ${\cal K}_0$
    on  $\overline{\cal M}_{(g,h),(n,\vec{m})}
                  (Y,L\,|\,\underline{\beta},\vec{\gamma},\mu)$
   are related as fibers of a standard Kuranishi
   structure-in-${\cal C}_{\spsccw}$ ${\cal K}/B$, flat over $B$.

  \item[$(2)$]
   There is a decomposition of moduli space
    $$
     \overline{\cal M}_{(g,h),(n,\vec{m})}
             (Y,L\,|\,\underline{\beta},\vec{\gamma},\mu)\;
     =\; \cup_{ \eta \in
                \bar{\Omega}_{( (g,h),(n,\vec{m})\,|\,
                                \underline{\beta},\vec{\gamma},\mu )} }
                \overline{\cal M}(Y,L\,|\,\eta)\,.
    $$
   The two sub-orbifolds
    $\overline{\cal M}((Y_1,L_1;D)\amalg (Y_2,L_2;D)\,|\,\eta)$  and\\
    $\overline{\cal M}(Y,L\,|\,\eta)$ of
         $\overline{\cal M}_{(g,h),(n,\vec{m})}
           (Y,L\,|\,\underline{\beta},\vec{\gamma},\mu)$
    are identical in
    $\overline{\cal M}_{(g,h),(n,\vec{m})}
              (Y,L\,|\,\underline{\beta},\vec{\gamma},\mu)$.

  \item[$(3)$]
   The restriction ${\cal K}_{0,\eta}$ of the Kuranishi structure
    ${\cal K}_0$ on
    $\overline{\cal M}_{(g,h),(n,\vec{m})}
             (Y,L\,|\,\underline{\beta},\vec{\gamma},\mu)$
    to the component $\overline{\cal M}(Y,L\,|\,\eta)$
    is equivalent to the Kuranishi structure ${\cal K}_{\eta}$ on\\
    $\overline{\cal M}((Y_1,L_1;D)\amalg(Y_2,L_2;D)\,|\,\eta)$,
    except that ${\cal K}_{0,\eta}$ carries a multiplicity $m(\eta)$.
   Let $\vec{s}=(s_1,\,\cdots\,,s_l)$ be the weights of
    the ordered roots in $\eta$;
   then $m(\eta)=m(\vec{s}):=s_1\,\cdots\,s_l$.
   In notation ${\cal K}_{0,\eta}=m(\eta)\,{\cal K}_{\eta}$.

  \item[$(4)$]
   Let $\eta=(\Gamma_1,\Gamma_2,I_{\hand},I_{\leg})$
    with $l$-many roots.
   Then, the Kuranishi structure ${\cal K}_{\eta}$  on
    $\overline{\cal M}((Y_1,L_1;D)\amalg(Y_2,L_2;D)\,|\,\eta)$
    is locally equivalent to the Kuranishi structure
    ${\cal K}^{\prime}_1\times_{D^l}{\cal K}^{\prime}_2$
    on  $\overline{\cal M}(Y_1,L_1;D\,|\,\Gamma_1)   \times_{D^l}
                         \overline{\cal M}(Y_2,L_2;D\,|\,\Gamma_2)$
    under the $|\Per_r(\eta)|$-fold covering map $\Phi_{\eta}$.
 \end{itemize}

 We use the following ``formula" to summarize/encapsulate
 $(1)$, $(2)$, $(3)$, and $(4)$:
 $$
  \begin{array}{rcl}
  [{\cal K}_{\lambda}]\;\leftrightarrow\; [{\cal K}_0]
   & = & \cup_{ \eta \in
                \bar{\Omega}_{( (g,h),(n,\vec{m})\,|\,
                        \underline{\beta},\vec{\gamma},\mu )} }\,
           [{\cal K}_{0,\eta}]  \hspace{1em}
         = \hspace{1em}
           \cup_{ \eta \in
               \bar{\Omega}_{( (g,h),(n,\vec{m})\,|\,
                        \underline{\beta},\vec{\gamma},\mu )} }\,
               m(\eta)\,[{\cal K}_{\eta}]          \\[1ex]
   & = &  \cup_{ \eta \in
               \bar{\Omega}_{( (g,h),(n,\vec{m})\,|\,
                        \underline{\beta},\vec{\gamma},\mu )} }\,
          \frac{m(\eta)}{|\smallPer_r(\eta)|}\,
          \Phi_{\eta\,\ast}\,
           [{\cal K}^{\prime}_1\times_{D^l}{\cal K}^{\prime}_2]\,.
  \end{array}
 $$
} % end-theorem

\bigskip

\noindent
{\it Proof.}
 We give only a sketch here and omit the tedious details.
 Statement (1) is by the definition of ${\cal K}_0$.
 Statement (2) follows by considering the topological types of maps.
 The multiplicity $m(\eta)$ in Statement (3) arise from the scheme
  structure of the centra fiber of $\Xi_{\vec{s}}\rightarrow{\Bbb C}$.
 Statement (4) requires a comparison of the fibered product of
  Kuranishi structures and that from a restriction.
 Here, as well as whenever we need to justify the equivalence of
  two standardly constructed Kuranishi structures on a same moduli space
  in question, is where Siebert's work [Sie1] plays roles again and again.
 Associated to a Kuranishi structure ${\cal K}$ is a fine system
  ${\cal E}_{\cal K}$ of saturated obstructed local bundles
  as sub-fibrations in the related $\check{L}^p$-obstruction space
  fibration $T^2_{\check{\cal W}^{1,p}(\,\cdots\,)}$ as in Sec.~5.4;
  and vice versa.
 To construct the equivalence of two given two Kuranishi structures
  ${\cal K}_1$ and ${\cal K}_2$, one constructs
  an appropriate fine system ${\cal E}$ of local bundles
  that contains both ${\cal E}_{{\cal K}_i}$ as sub-fibrations.

\noindent\hspace{15cm}$\Box$

\bigskip

We emphasize that, at the level of Kuranishi structures, the above
 degeneration-gluing relations under a symplectic cut hold for
 {\it both} closed Gromov-Witten theory and open Gromov-Witten theory
 and by {\it the same} reason.

\bigskip

\begin{flushleft}
{\bf Example: Li-Ruan/Li degeneration formula of closed
              Gromov-Witten invariants.}
\end{flushleft}
When $L$ is empty, the domain of maps are closed nodal Riemann surfaces
 and
we resume the moduli space
 $$
  \overline{\cal M}_{g,n}(W_{\lambda},[\beta])\;
   :=\; \overline{\cal M}_{(g,0),(n,0)}(W_{\lambda}\,|\,[\beta])\;
    =\; \amalg_{\beta^{\prime\prime}\in [\beta]}\,
         \overline{\cal M}_{g,n}(W_{\lambda},\beta^{\prime\prime})
 $$
 in closed Gromov-Witten theory.
The notion of admissible quadruples in Definition 7.1.13
        % Definition [admissible quadruple]
 is reduced to the notion of {\it admissible triples}
 $\eta=(\Gamma_1,\Gamma_2,I=I_{\leg})$ (cf.\ [Li1: Definition 4.11])
  and
 its {\it type} is now a triple of the form
  $(\hat{g},\hat{n};\underline{\hat{\beta})}$.
Denote by $\bar{\Omega}_{(g,n,\underline{\beta})}$
 the set of equivalence classes of admissible triples $\eta$
 such that $|\eta|=(g,n;\underline{\beta})$.
Then, Theorem 7.1.14
      % Theorem [degeneration/gluing - Kuranishi structure]
 reduces to
 $$
  \overline{\cal M}_{g,n}(Y,\underline{\beta})\;
  =\; \cup_{ \eta \in \bar{\Omega}_{(g,n;\underline{\beta})} }
             \overline{\cal M}(Y\,|\,\eta)\,.
 $$
 and, in the encapsulated form,
 $$
  \begin{array}{rcl}
  [{\cal K}_{\lambda}]\;\leftrightarrow\; [{\cal K}_0]
   & = & \cup_{ \eta \in\bar{\Omega}_{(g,n;\underline{\beta})} }\,
           [{\cal K}_{0,\eta}]  \hspace{1em}
         = \hspace{1em}
           \cup_{ \eta\in\bar{\Omega}_{(g,n;\underline{\beta})} }\,
               m(\eta)\,[{\cal K}_{\eta}]          \\[1ex]
   & = &  \cup_{ \eta \in\bar{\Omega}_{(g,n;\underline{\beta})} }\,
          \frac{m(\eta)}{|\smallPer_r(\eta)|}\,
          \Phi_{\eta\,\ast}\,
           [{\cal K}^{\prime}_1\times_{D^l}{\cal K}^{\prime}_2]\,.
  \end{array}
 $$

A virtual fundamental class
 $[\overline{\cal M}_{g,n}(W_{\lambda},[\beta])]^{\virt}$
 of the expected dimension and supported in $s_{\rho,\lambda}^{-1}(0)$
  on each Kuranishi neighborhood $V_{\rho;\lambda}$,
  $\rho  \in
   \overline{\cal M}_{(g,0),(n,0)}(W_{\lambda},[\beta])$
 can be constructed\footnote{This step is not trivial.
        It includes a re-doing of [L-T3] and [Sie2] in the Fukaya-Ono
         family setting.
        Readers who are not familiar with ibidem may think of
         a Kuranishi neighborhood $(V,\Gamma,E_V;s,\psi)$ directly as
         a ``virtual cycle" of the expected dimension in the
         (usually singular) orbifold local chart $s^{-1}(0)\subset V$
         of the moduli space, weighted by $1/|\Gamma|$.
        Equivalent Kuranishi neighborhoods give equivalent local
         virtual cycles.
        Transition data of a Kuranishi structure gives the patching data
         of these local cycles and defines a {\it virtual fundamental
         cycle} on the (usual singular) moduli orbifold space.
        Equivalent Kuranishi structures define the same virtual
         fundamental class on the moduli orbifold space.}
  via Kuranishi structures ${\cal K}_{\lambda}$.
Similarly, for
 $$
  \begin{array}{l}
   [\overline{\cal M}_{g,n}(W/B,[\beta])/B]^{\virt}\,,\hspace{2em}
   [\overline{\cal M}_{g,n}(Y,\underline{\beta})]^{\virt}\,,
                                                      \hspace{2em}
   [\overline{\cal M}(Y\,|\,\eta)]^{\virt}\,,\hspace{2em}
   [\overline{\cal M}(Y_1;D\,|\,\Gamma_1)]^{\virt}\,, \\[.6ex]
   [\overline{\cal M}(Y_2;D\,|\,\Gamma_2)]^{\virt}\,,\hspace{1em}
   [\overline{\cal M}(Y_1;D\,|\,\Gamma_1)
     \times_{D^l} \overline{\cal M}(Y_2;D\,|\,\Gamma_2)]^{\virt}\,,
                                                     \hspace{1em}
   [\overline{\cal M}((Y_1;D)\amalg (Y_2;D)\,|\,\eta)]^{\virt}
  \end{array}
 $$
 that are constructed from Kuranishi structures
  $$
   {\cal K}/B\,,\hspace{1em}
   {\cal K}_0\,,\hspace{1em}   {\cal K}_{0,\eta}\,,\hspace{1em}
   {\cal K}^{\prime}_1\,,\hspace{1em}
   {\cal K}^{\prime}_2\,,\hspace{1em}
   {\cal K}^{\prime}_1\times_{D^l}{\cal K}^{\prime}_2\,,\hspace{1em}
   {\cal K}_{\eta}
  $$
 respectively.
%%%%%%%%%%%%%%%%%%%%
% Fix a standard Kuranishi structure ${\cal K}/B$ on
%  $\overline{\cal M}_{(g,h),(n,\vec{m})}
%              (W/B,L\,|\,[\beta],\vec{\gamma},\mu)$.
% Then the Kuranishi neighborhoods from ${\cal K}_0^{(0)}$
%  ??????????????????.
%
% \bigskip
%
% \noindent $\bullet$
% The notion of a {\it maximal cluster} of
%  $\overline{\cal M}_{(g,h),(n,\vec{m})}
%            (Y,L\,|\,\underline{\beta},\vec{\gamma},\mu)$  and\\
%  $\overline{\cal M}_{(g,h),(n+l(\vec{s}),\vec{m})}
%    (Z,L;D\,|\,\beta^{\prime},\vec{\gamma},\mu^{\prime};\vec{s})$.
%
% \bigskip
%
% \noindent $\bullet$
% The multiplicity of a maximal cluster.
%
% \bigskip
%
% \noindent $\bullet$
% The reduced Kuranishi structure on a maximal cluster.
%
% \bigskip
%
% \noindent $\bullet$
% An  ${\cal K}_0$
%  as the gluing of Kuranishi structures on
%  the moduli space of relative stable maps
%  to $(Y_1,L_1;D)$ and to $(Y_2,L_2;D)$.
%
% \bigskip
%
% \noindent $\bullet$
%  ??????????????.
%
% \bigskip
%%%%%%%%%%%%%%%%%%%%%%%%%%%%%%
Since equivalent Kuranishi structures give identical
 virtual fundamental class,
 the above degeneration-gluing formula of Kuranishi structures
  can be reduced\footnote{In the Fukaya-Ono setting, the degeneration
           formulas of any form in Gromov-Witten theory should be regarded
            as the consequence of the more fundamental degeneration-gluing
            relations of Kuranishi structures and an assignment
            to each moduli space with a Kuranishi structure a virtual
            fundamental class or chain that is functorial,
            particularly with respect to restrictions to
            sub-moduli spaces, fibered product, and covering maps.
           Recall
            the layer-structure decompositions of the moduli spaces
             of stable or relative stable maps and
            the virtual co-dimension of each stratum.
           These notions extends to the fiber-products that occur
            in the problem.
           The functorial property of a virtual fundamental class
            $[{\cal M}]^{\it virt}$
            implies that $[{\cal M}]^{\it virt}$ is determined
            by its restriction to the depth-$0$
            (i.e.\ virtual codimension-$0$) stratum in the moduli space.
           As the depth-$0$ strata that occur in right-hand side of
             the decomposition
             $\overline{\cal M}_{g,n}(Y,\underline{\beta})\;
               =\; \cup_{ \eta \in \bar{\Omega}_{(g,n;\underline{\beta})} }
                          \overline{\cal M}(Y\,|\,\eta)$
             are disjoint from each other,
            the union becomes a disjoint union when restricted to
             depth-$0$ strata of the moduli spaces in the identity.
           This disjoint union is then turned into a summation
            of virtual fundamental classes on these strata
            when the degeneration-gluing relations of Kuranishi structures
            are applied.
           As recovering the whole moduli space by adding in strata
            of positive depth will extend the virtual fundamental
            class by only lower-dimensional classes in the strata
            of positive depth,
           the summation is not influenced. This gives thus the
           degeneration/gluing formula at the level of virtual
            fundamental classes.
          % As a comparison, to re-do [Li2] in this aspect, one would
          % first derive a degeneration-gluing formula of
          % perfect tangent-obstruction complex and construct
          % the virtual fundamental class functorially with respect to
          % fibered produce and flat morphism.
          % The construction of [Li2] merges the two.
           It is with this aspect that we state, as an example,
            the result of [L-R] and [Li2] as a corollary.}
  to the degeneration/gluing formulas of
  Li-Ruan [L-R] and Li [Li2]\footnote{Note
                              that the degeneration formulas of
                               Li-Ruan and Li are equivalent.
                              Here we use the expression in [Li2];
                               see [L-R] for the expression in terms
                               of integrals over virtual
                                neighborhoods [Ru] with Thom forms.
                              See also the Appendix of the current work
                               for a discussion on the equivalence
                               of the degeneration formulas of
                               Li-Ruan [L-R], Li [Li2], and the
                               formally different Ionel-Parker [I-P2].}:

\bigskip

\noindent
{\bf Corollary 7.1.15
     [degeneration-gluing$\,$: virtual fundamental class].}
{\it
 $[\overline{\cal M}_{g,n}(X,[\beta])]^{\virt}$  and
  $[\overline{\cal M}_{g,n}(Y,\underline{\beta})]^{\virt}$
  can be realized as the fibers of the flat class
  $[\overline{\cal M}_{g,n}(W_{\lambda},[\beta])]^{\virt}/B$
  over $B$.
 $$
  \begin{array}{l}
    [\overline{\cal M}_{g,n}(Y,\underline{\beta})]^{\virt}
     \hspace{1em} = \hspace{1em}
       \sum_{ \eta\in\bar{\Omega}_{(g,n;\underline{\beta})} }\,
           [\overline{\cal M}(Y\,|\,\eta)]^{\virt}     \\[1ex]
    \hspace{2em} =\;
    \sum_{ \eta\in\bar{\Omega}_{(g,n;\underline{\beta})} }\,
     m(\eta)\,
     [\overline{\cal M}((Y_1;D)\amalg(Y_2;D)\,|\,\eta)]^{\virt} \\[1ex]
    \hspace{2em} =\;
    \sum_{ \eta=(\Gamma_1,\Gamma_2,I)
                 \in\bar{\Omega}_{(g,n;\underline{\beta})} }\,
      \frac{m(\eta)}{|\smallPer_r(\eta)|}\,
       \Phi_{\eta\,\ast}\,
        [\overline{\cal M}(Y_1;D\,|\,\Gamma_1)  \times_{D^l}
          \overline{\cal M}(Y_2;D\,|\,\Gamma_2)]^{\virt}  \\[1ex]
    \hspace{2em} =\;
    \sum_{ \eta=(\Gamma_1,\Gamma_2,I)
                   \in\bar{\Omega}_{(g,n;\underline{\beta})} }\,
      \frac{m(\eta)}{|\smallPer_r(\eta)|}\,
       \Phi_{\eta\,\ast}\,
        \Delta_{\eta}^{!}(
          [\overline{\cal M}(Y_1;D\,|\,\Gamma_1)]^{\virt}  \times
          [\overline{\cal M}(Y_2;D\,|\,\Gamma_2)]^{\virt} )\,,
  \end{array}
 $$
 where, for $\eta$ with $l$-many roots,
  $\Delta_{\eta}:D^l \hookrightarrow D^l\times D^l$ is the diagonal map
    and
  $$
   \Delta_{\eta}^{!}\;:\;
    A_{\ast}(\overline{\cal M}(Y_1;D\,|\,\Gamma_1)
              \times \overline{\cal M}(Y_2;D\,|\,\Gamma_2))\;
   \longrightarrow\;
   A_{\ast}(\overline{\cal M}(Y_1;D\,|\,\Gamma_1)
             \times_{D^l} \overline{\cal M}(Y_2;D\,|\,\Gamma_2))
  $$
  is the Gysin homomorphism under
  $$
   \begin{array}{ccc}
    \overline{\cal M}(Y_1;D\,|\,\Gamma_1)
     \times_{D^l} \overline{\cal M}(Y_2;D\,|\,\Gamma_2)
     & \longrightarrow
     & \overline{\cal M}(Y_1;D\,|\,\Gamma_1)
        \times \overline{\cal M}(Y_2;D\,|\,\Gamma_2) \\[.6ex]
    \downarrow  & & \downarrow  \\
    D^l         & \stackrel{\Delta_{\eta}}{\longrightarrow}
                & D^l\times D^l\,.
   \end{array}
  $$
} % end-theorem

\bigskip

The Gromov-Witten invariants of $X$ associated to $(g,n;[\beta])$
 are defined by\footnote{All cohomologies in the definition of
                             Gromov-Witten and relative Gromov-Witten
                             invariants are over ${\footnotesizeBbb Q}$.}:
 $$
  \begin{array}{cccccl}
   \Psi^X_{(g,n;[\beta])} & :
    & H^{\ast}(X)^{\times n} \times H^{\ast}(\overline{\cal M}_{g,n})
    & \longrightarrow  & {\Bbb Q} \\[.6ex]
   && (\,\kappa\;,\;\varsigma\,)
    & \longmapsto
    &  \left[ \ev^{\ast}(\kappa)
               \cup \pi_{(g,n)}^{\ast}(\varsigma)\,
               [\overline{\cal M}_{g,n}(X,[\beta])]^{\virt}
       \right]_0\,,
  \end{array}
 $$
 where
  $\ev:\overline{\cal M}_{g,n}(X,[\beta])\rightarrow X^n$
   is the evaluation map associated to the ordered set of
   $n$ marked points,
  $\pi_{(g,n)}:\overline{\cal M}_{g,n}(X,[\beta])
                     \rightarrow \overline{\cal M}_{g,n}$
   is the domain-curve stabilization map\footnote{For
                   $X$ smooth, $\pi_{(g,n)}$ is a local complete
                   intersection morphism when extended to a map
                   on Kuranishi neighborhoods.},  and
  $[\,\cdot\,]_0$ means the degree-$0$ component of $\,\cdot\,$.

Given an admissible weighted graph $\Gamma$ with $n$ legs and $l$ roots,
let $\overline{\cal M}_{\Gamma}$ be the moduli space of
 stables curves
  with $|V(\Gamma)|$-many connected components
  in one-one correspondence with $V(\Gamma)$,
 $n$ ordinary marked points corresponding to legs and
 $l$ distinguished marked points.corresponding to roots accordingly.
The relative Gromov-Witten invariants of the pair $(Z,D)$
 associated to an admissible weighted graph $\Gamma$
 with $n$ legs and $l$ roots are defined by
 $$
  \begin{array}{cccccl}
   \Psi^{(Z,D)}_{\Gamma} & :
    & H^{\ast}(Z)^{\times n} \times H^{\ast}(\overline{\cal M}_{\Gamma})
    & \longrightarrow  & H_{\ast}(D^l) \\[.6ex]
   && (\,\kappa\;,\;\varsigma\,)
    & \longmapsto
    & {\mathbf q}_{\ast}
         \left(\, \rule{0em}{1.2em}
           \ev^{\ast}(\kappa)
            \cup \pi_{\Gamma}^{\ast}(\varsigma)\,
                  [\overline{\cal M}(Z;D\,|\,\Gamma)]^{\virt}\,
         \right)   &,
  \end{array}
 $$
 where
  $\ev:\overline{\cal M}_{g,n}(X,[\beta])\rightarrow X^n$
   is the evaluation map associated to the ordered set of
   ordinary $n$ marked points,
  $\pi_{\Gamma}:\overline{\cal M}(Z;D\,|\,\Gamma)
              \rightarrow \overline{\cal M}_{\Gamma}$
    is the domain-curve stabilization map,  and\\
  ${\bf q}:\overline{\cal M}(Z;D\,|\,\Gamma) \rightarrow D^l$
   is the evaluation map associated to the ordered set of $l$
   distinct marked points.

For an admissible triple
  $\eta=(\Gamma_1,\Gamma_2,I)$ with $(n_1,n_2)$-many legs and
  $l$-many roots,
gluing at the paired distinguished marked points defines
 an orbifold map
 $G_{\eta}:
  \overline{\cal M}_{\Gamma_1}\times \overline{\cal M}_{\Gamma_2}
                             \rightarrow \overline{\cal M}_{g,n}$,
  where $g=g(\eta)$ and $n=n_1+n_2$.
For $\varsigma\in H^{\ast}(\overline{\cal M}_{g,n};{\Bbb Q})$,
 we assume that the K\"{u}nneth decomposition
 $G_{\eta}^{\ast}(\varsigma)
    =\sum_{j\in N_{\eta}}
         \varsigma_{\eta,1,j}\bboxtimes \varsigma_{\eta,2,j}$
 exists.
Then the degeneration-gluing formula of Gromov-Witten invariants
 with respect to $X/B$ is given by: ([Li2])

\bigskip

\noindent
{\bf Corollary 7.1.16 [degeneration-gluing$\,$: invariant].} {\it
 Let
  $\kappa \in H_c^0(R^{\tinybullet}\pi_{\ast}{\Bbb Q}_W)^{\oplus n}$,
  $\varsigma \in H^{\ast}(\overline{\cal M}_{g,n})$,
  and
  $j_i:Y_i\hookrightarrow Y=W_0$, $i=1,\,2$.
 Then
 $$
  \begin{array}{l}
   \Psi_{(g,n;[\beta])}^{W_{\lambda}}(\kappa(\lambda),\varsigma)  \\[.6ex]
    \hspace{2em}
     =\; \sum_{\eta\in\bar{\Omega}_{(g,n;\underline{\beta})}}\;
          \frac{m(\eta)}{|\Per_r(\eta)|}\,
           \sum_{j\in N_{\eta}}\,
            \left[ \Psi_{\Gamma_1}^{(Y_1,D)}(j_1^{\ast}\kappa(0),
                     \varsigma_{\eta,1,j})\,
                   \bullet\,
                   \Psi_{\Gamma_2}^{(Y_2,D)}(j_2^{\ast}\kappa(0),
                     \varsigma_{\eta,2,j})
            \right]_0\,,
  \end{array}
 $$
 where
  $\kappa(\lambda)$ is the restriction of $\kappa$ to the fiber
    $W_{\lambda}$ of $W/B$,
  $\,\bullet\,$ is the intersection product on $H_{\ast}(D^l)$,
  $[\,\cdot\,]_0$ is the degree-$0$ component of $\,\cdot\,$.
} % end-theorem

\bigskip

\subsection{A degeneration axiom and a gluing axiom for
    open Gromov-Witten invariants under a symplectic cut.}

When $L$ is non-empty, the (real) codimension-$1$ boundary on
 the moduli space $\widetilde{\cal M}_{(g,h),(n,\vec{m})}$ of
 prestable labelled-bordered Riemann surfaces gives rise to
 the codimension-$1$ boundary $\partial {\cal K}_{(X,L)}$
 on the Kuranishi structure ${\cal K}_{(X,L)}$
 on the moduli space
 $\overline{\cal M}_{(g,h),(n,\vec{m})}
                          (X,L\,|\,\beta,\vec{\gamma},\mu)$
 of stable maps to $(X,L)$.
As a Gromov-Witten theory/invariant so far constructed is based
 on the intersection theory with the functorially constructed
 virtual fundamental class/chain on the moduli space,
the birth-'n-death of chain components along the codimension-$1$
 boundary $\partial{\cal K}_{(X,L)}$ of ${\cal K}_{(X,L)}$
 makes such construction not well-defined
 unless one has a way to fix the ambiguity.
Furthermore, it has been noticed ([K-L]) that
  to define meaningful open Gromov-Witten invariants and
  to match with the physicists' computation of
  open string instantons (e.g.\ [A-K-V]),
 a {\it decoration} $\alpha$ has to be imposed to the Lagrangian
 submanifold $L$, to which boundaries of
 Riemann surfaces/open string world-sheets are mapped.
Basic examples of decorations are
  a {\it group action} on $L$,
  a {\it framing} on $T_{\ast}L$, and
  an {\it involution} on $T_{\ast}X|_L$ that leave $T_{\ast}L$ fixed,
 if any of these structures on $L$ exists.
Denote a Lagrangian submanifold $L$ with a decoration $\alpha$
 by $L^{\alpha}$.
Thus:
 \begin{quote}
  \parbox[t]{6em}{\sc Problem$\,$:}\
  {\sl To define open Gromov-Witten invariants for $(X,L^{\alpha})$.}
 \end{quote}
Note that in general $\alpha$ on $L$ does not extend to
 a decoration on $X$.

With the above problem in mind,
 the degeneration and gluing of Kuranishi structures studied
  in this work  and
 the deformation-invariance requirement of open Gromov-Witten
  invariants
 propel us to impose the following two {\it axioms} on
 open Gromov-Witten invariants.

The {\it Gromov-Witten invariants} of $(X,L^{\alpha})$
 associated to $((g,h),(n,\vec{m})\,|\,\beta,\vec{\gamma},\mu)$
 are meant to be the evaluation of a map
 $$
  \Psi^{(X,L^{\alpha})}
      _{((g,h),(n,\vec{m})\,|\,\beta,\vec{\gamma},\mu)}\;
  :\;  H^{\ast}(X)^{\times n}
       \times H^{\ast}(L)^{\times |\vec{m}|}
       \times H^{\ast}(\overline{\cal M}_{(g,h),(n,\vec{m})})\;
      \longrightarrow\; {\Bbb Q}\,,
 $$
 where $|\vec{m}|=m_1+\,\cdots\,+m_h$, that satisfies
 a set of properties\footnote{Besides the interest in its own right,
                            open Gromov-Witten theory gives a mathematical
                             formulation for the problem of
                             {\it open string world-sheet instantons}
                             and their enumeration
                             in superstring theory;
                            it is closely related also to
                             {\it conformal field theory with boundary}
                             and {\it D-branes}.
                            Some of the properties
                             $\Psi^{\tinybullet}_{\tinybullet}
                                                      (\,\cdots\,)$
                             has to satisfy come from these subjects
                             in superstring theory.
                            The following
                             incomplete/intentionally-limited
                             additional stringy literatures only mean
                             to give unfamiliar readers a glimpse
                             of these diverse yet linked topics:
                             [B-C-O-V: Sec's 4, 5.5, 8.2],
                             [H-I-V], [K-K-L-MG], and reviews
                             [Dou], [Ga], [S-F-W], and [T-Z].},
 e.g.\ the list in [Ko-M].
The same holds with $X$ replaced by the singular $Y$.
Define also
 $$
  \Psi^{(X,L^{\alpha})}
      _{((g,h),(n,\vec{m})\,|\,[\beta],\vec{\gamma},\mu)}\;
  :=\; \sum_{\beta^{\prime\prime}\in[\beta]}\,
         \Psi^{(X,L^{\alpha})}
             _{((g,h),(n,\vec{m})\,|\,
                         \beta^{\prime\prime},\vec{\gamma},\mu)}\,.
 $$
Similarly,
given an admissible weighted graph $\Gamma$
 with $n$ legs, $m$ fingers, and $l$ roots,
the {\it relative Gromov-Witten invariants} of the relative pair
 $(Z,L^{\alpha};D)$ associated to $\Gamma$
 are meant to be the evaluation of a map
 $$
  \Psi^{(Z,L^{\alpha};D)}_{\Gamma}\; :\;
    H^{\ast}(Z)^{\times n}
    \times H^{\ast}(L)^{\times m}
    \times H^{\ast}(\overline{\cal M}_{\Gamma})\;
    \longrightarrow\;   H_{\ast}(D^l)\,,
 $$
 where
 $\overline{\cal M}_{\Gamma}$ is the moduli space of
  (not necessarily connected) labelled-bordered Riemann surfaces
  with marked points of combinatorial type specified by $\Gamma$.
For an admissible quadruple
  $\eta=(\Gamma_1,\Gamma_2,I_{\hand},I_{\leg})$
  with $(h_1,h_2)$-many hands, $(n_1,n_2)$-many legs, $l$-many roots,
  and  type $|\eta|=((g,h),(n,\vec{m})\,|\,
                                \underline{\beta},\vec{\gamma},\mu)$
gluing at the paired distinguished marked points defines
 an orbifold map
 $G_{\eta}:
  \overline{\cal M}_{\Gamma_1}\times \overline{\cal M}_{\Gamma_2}
                 \rightarrow \overline{\cal M}_{(g,h),(n,\vec{m})}$.
For $\varsigma \in
      H^{\ast}(\overline{\cal M}_{(g,h),(n,\vec{m})};{\Bbb Q})$,
 we assume that the K\"{u}nneth decomposition
 $G_{\eta}^{\ast}(\varsigma)
    =\sum_{j\in N_{\eta}}
         \varsigma_{\eta,1,j}\bboxtimes \varsigma_{\eta,2,j}$
 exists.

\bigskip

\noindent
{\bf Axiom OGW-degeneration.} {\it
 Let
  $W/B$ be a degeneration of $X$ associated to
   a symplectic cut $\xi:X \rightarrow Y=Y_1\cup _D Y_2 =W_0$ and
  $L^{\alpha}$ be a decorated Lagrangian submanifold of $X$
    disjoint from the cutting locus.
 The submanifold in $W_{\lambda}$ associated to $L$ is denoted
  also by $L$.
 Let
  $\kappa \in H_c^0(R^{\tinybullet}\pi_{\ast}{\Bbb Q}_W)^{\oplus n}$,
  $\upsilon \in H^{\ast}(L)^{|\vec{m}|}$, and
  $\varsigma \in H^{\ast}(\overline{\cal M}_{g,n})$.
 Then
  $$
   \Psi^{(W_{\lambda},L^{\alpha})}
     _{((g,h),(n,\vec{m})\,|\,[\beta],\vec{\gamma},\mu)}
      (\kappa(\lambda),\upsilon,\varsigma)\;
    =\; \Psi^{(Y,L^{\alpha})}
         _{((g,h),(n,\vec{m})\,|\,\underline{\beta},\vec{\gamma},\mu)}
        (\kappa(0),\upsilon,\varsigma)\,,
  $$
  where
   $\underline{\beta}=\xi_{\ast}([\beta])$ and
   $\kappa(\lambda)$ is the restriction of $\kappa$ to $W_{\lambda}$.
} % end-axiom

\bigskip

\noindent
{\bf Axiom OGW-gluing.} {\it
 Gromov-Witten invariants of
   $(Y,L^{\alpha})=(Y_1,L_1^{\alpha})\cup_D (Y_2,L_2^{\alpha})$
   can be expressed in terms of
  relative Gromov-Witten invariants of
   $(Y_i,L_i^{\alpha};D)$, $i=1,\,2$, by the identity:

 \vspace{-1ex}
 {\small
 $$
  \begin{array}{l}
   \Psi^{(Y,L^{\alpha})}
       _{((g,h),(n,\vec{m})\,|\,\underline{\beta},\vec{\gamma},\mu)}
    (\kappa,\upsilon,\varsigma) \\[.6ex]
   \;=\; \sum_{
          \eta\in\bar{\Omega}_{ ((g,h),(n,\vec{m})\,|\,
                               \underline{\beta},\vec{\gamma},\mu) } }\;
            \frac{m(\eta)}{|\Per_r(\eta)|}\,
           \sum_{j\in N_{\eta}}\,
            \left[ \Psi_{\Gamma_1}^{(Y_1,L_1^{\alpha};D)}
                      (j_1^{\ast}\kappa, j_1^{\ast}\upsilon,
                                        \varsigma_{\eta,1,j})\,
                   \bullet\,
                   \Psi_{\Gamma_2}^{(Y_2,L_2^{\alpha};D)}
                     (j_2^{\ast}\kappa,  j_2^{\ast}\upsilon,
                                          \varsigma_{\eta,2,j})
            \right]_0\,,
  \end{array}
 $$
 {\normalsize where}} % end-small
  $\,\bullet\,$ is the intersection product on $H_{\ast}(D^l)$, and
  $[\,\cdot\,]_0$ is the degree-$0$ component of $\,\cdot\,$.
} % end-axiom

\bigskip

\noindent
{\it Remark 7.2.1
    $[$selection of fundamental chains adapted to $\alpha$
       -- specialization$]$.}
 Concerning the ambiguity mentioned in the beginning of
   this subsection on the choices of virtual fundamental chains,
  below is how these two axioms are applied to this issue.
 For simplicity of presentation, we assume that
  $\xi:(X,L^{\alpha})  \rightarrow
    (Y,L^{\alpha})=(Y_1,L^{\alpha})\cup_D (Y_2,\emptyset)$.
 Suppose that
  \begin{quote}
   [$\,${\it assumption}$\,$]\hspace{1ex}
   the decoration $\alpha$ is full enough to select
    in a standard way a class of virtual fundamental chains
    $[\overline{\cal M}(Y_1,L^{\alpha};D\,|\,\Gamma_1)]^{\virt}$
    in $\overline{\cal M}(Y_1,L^{\alpha};D\,|\,\Gamma_1)$
    associated to a standard Kuranishi structure
    ${\cal K}^{\prime}_{(Y_1,L^{\alpha};D\,|\,\Gamma_1)}$
   for all $\Gamma_1$ in an
    $\eta\in\bar{\Omega}_{((g,h),(n,\vec{m})\,|\,
                           \underline{\beta},\vec{\gamma},\mu)}$,
  \end{quote}
 then it induces a class of virtual fundamental chains on
  $\overline{\cal M}_{(g,h),(n,\vec{m})}
     (W_{\lambda},L^{\alpha}\,|\,[\beta],\vec{\gamma},\mu)$
  as follows:
  \begin{itemize}
   \item[$\cdot$]
    the push-forward of the fibered product of
     $[\overline{\cal M}(Y_1,L^{\alpha};D\,|\,\Gamma_1)]^{\virt}$
      with
     $[\overline{\cal M}(Y_2;D\,|\,\Gamma_2)]^{\virt}$ over $D^l$
     weighted by $m(\eta)/|\Per_r(\eta)|$
    gives rise to a class of virtual fundamental subchains
    $[\overline{\cal M}(Y,L^{\alpha}\,|\,\eta)]^{\virt}$
    in  $\overline{\cal M}_{(g,h),(n,\vec{m})}(Y,L^{\alpha}\,|\,
                             \underline{\beta},\vec{\gamma},\mu)$;

   \item[$\cdot$]
    their summation over
    $\eta\in \bar{\Omega}
              _{((g,h),(n,\vec{m})\,|\,\underline{\beta},\vec{m},\mu)}$
    gives a class of virtual fundamental chains\\
    $[\overline{\cal M}_{(g,h),(n,\vec{m})}
      (Y,L^{\alpha}\,|\,\underline{\beta},\vec{\gamma},\mu)]^{\virt}$
    in  $\overline{\cal M}_{(g,h),(n,\vec{m})}
               (Y,L^{\alpha}\,|\,\underline{\beta},\vec{\gamma},\mu)$;

   \item[$\cdot$]
    deform the chains
    $[\overline{\cal M}_{(g,h),(n,\vec{m})}
      (Y,L^{\alpha}\,|\,\underline{\beta},\vec{\gamma},\mu)]^{\virt}$
    to over $\lambda\ne 0$  by a $2$-dimension-higher chain $c$ in
    $\overline{\cal M}_{(g,h),(n,\vec{m})}
                    (W/B,L^{\alpha}\,|\,[\beta],\vec{\gamma},\mu)$
    such that both $c$ and its restriction to\\
     $\partial\overline{\cal M}_{(g,h),(n,\vec{m})}
                    (W/B,L^{\alpha}\,|\,[\beta],\vec{\gamma},\mu)/B$
     are flat over $B$;
    this then defines a class of virtual fundamental chains
     $[\overline{\cal M}_{(g,h),(n,\vec{m})}
            (W/B,L^{\alpha}\,|\,[\beta],\vec{\gamma},\mu)]^{\virt}$
     in  $\overline{\cal M}_{(g,h),(n,\vec{m})}
             (W_{\lambda},L^{\alpha}\,|\,[\beta],\vec{\gamma},\mu)$.
 \end{itemize}
 In this prescription,
  $\partial\overline{\cal M}_{(g,h),(n,\vec{m})}
         (W/B,L^{\alpha}\,|\,[\beta],\vec{\gamma},\mu)/B$
  consists of stable maps to the fibers of
   $(\widehat{W},\widehat{L})/\widehat{B}$ of the given type
   such that the domain $\Sigma$ has either boundary nodes or
    free marked points landing on $\partial\Sigma$.
 The requirement of the flatness of the deformation of chains
  also on the restriction to
  $\partial\overline{\cal M}_{(g,h),(n,\vec{m})}
           (W/B,L^{\alpha}\,|\,[\beta],\vec{\gamma},\mu)/B$
  suppresses the birth-'n-death of chains from the
  codimension-$1$ boundary Kuranishi structure of
  the Kuranishi structure on
  $\overline{\cal M}_{(g,h),(n,\vec{m})}
             (X,L^{\alpha}\,|\,[\beta],\vec{\gamma},\mu)/B$.
 Here, we identify $X$ as some fiber $W_{\lambda_0}$ of $W/B$
  with $\lambda_0\ne 0$.
 This process is similar to the {\it specialization} technique
  in algebraic geometry.  %:
  % deform and specialize to a degenerate situation on which the
  % ``counting of points" can be done and then flatly deform back.

\bigskip

\noindent
{\bf Definition 7.2.2 [$L$-isolatable].} {\rm
 We call $(X,L^{\alpha})$ {\it $L$-isolatable}
  if there exists a symplectic cut
  $$
   X\; \longrightarrow\;
   (Y_0\,;\, \amalg_i\,D_i)
    \cup_{\,\cup_i\,D_i} \amalg_i (Z_i,L_i;D_i)\,,
  $$
  where $L_i$'s are the finitely many connected components of $L$,
  such that $Z_i$ is the symplectic manifold determined by $L_i$
  with the property that $Z_i-D_i$ is symplecto-isomorphic to
  a tubular neighborhood of the $0$-section of $T^{\ast}L_i$.
 Here, $T^{\ast}L$ is equipped with the canonical symplectic structure.
} % end-definition

\bigskip

Under Axiom OGW-degeneration and Axiom OGW-gluing,
the problem of the construction of open
 Gromov-Witten invariants of $L$-isolatable $(X,L^{\alpha})$
 is reduced to
 \begin{itemize}
  \item[]
   Step (2)$\,$:
    \parbox[t]{32em}{the
     construction of relative open Gromov-Witten invariants
     of $(Z,L^{\alpha};D)$\\  determined by $L^{\alpha}$.}
 \end{itemize}
Such class of $(X,L^{\alpha})$'s includes those that have occurred
 in the open/closed string duality.

%\bigskip
%\bigskip

\vspace{6em}
\begin{flushleft}
{\large\bf \parbox[t]{5em}{Appendix.}
  The equivalence of Li-Ruan/Li's degeneration formula and\\[.6ex]
  \hspace{5.4em}Ionel-Parker's degeneration formula.}
\end{flushleft}
The details of [L-R] and [Li1], [Li2],
 together with Comparison 3.2.4 in Sec.~3.2,
   % Comparison [refinement of [Li1] and [Li2]]
 imply that the degeneration formula of the (closed)
 Gromov-Witten invariants derived by
  A.-M.~Li and Y.~Ruan in [L-R] and J.~Li in [Li1], [Li2]
 are the same.
The Degeneration Axiom and the Gluing Axiom of open Gromov-Witten
 invariants we propose in Sec.~7.2 are of the Li-Ruan/Li form.
This form are formally different\footnote{
   See [L-R: p.~159], [I-P1: p.~48], and [Li1: Sec.~0]
     for a light comparison by these authors themselves.}
 from that derived by
 E.-N.\ Ionel and T.H.\ Parker in [I-P1], [I-P2].
Indeed, we can also adopt the discussion of [I-P2: Sec.~12] to
 give degeneration-gluing axioms of open Gromov-Witten invariants
 in the Ionel-Parker form,
though algebro-geometrically (cf.\ [Fu: Chap.~10])
 it is the Li-Ruan/Li form that we would choose,
 as it comes from a flat family construction.
This leads to the following question:
 \begin{itemize}
  \item[{\bf Q.}]
  {\sl Do Li-Ruan/Li and Ionel-Parker give different/independent sets
       of gluing/degeneration\\
       axioms for open Gromov-Witten invariants for a symplectic cut?}
 \end{itemize}

In this appendix,
 as a not-completely-irrelevant issue to our project,
 we explain the following conjecture, whose justification
 will answer the above question {\it negatively}\footnote{Indeed,
           from the algebro-geometric point of view,
           {\it any degeneration/gluing formula for intersection-theoretic
            type invariants that are constant under flat deformations
            must be re-derivable from a flat family construction}  and
           {\it any gluing/degeneration formula of Gromov-Witten invariants
            different from the one derived by a flat family construction
           {\rm (}i.e.\
            the Li-Ruan/Li formula in the case of symplectic cut{\rm )}
            must be convertible to the latter}
           unless they are indeed dealing with different invariants or
           different kinds of degenerations.}:

\bigskip

\noindent
{\bf Conjecture A.1 [Li-Ruan/Li = Ionel-Parker].} {\it
 Li-Ruan/Li's degeneration formula and Ionel-Parker's degeneration
   formula for closed Gromov-Witten invariants
  are equivalent/convert-ible to each other.
 Furthermore, the conversion is induced by
   $$
    \mbox{\sl Li-Ruan\,/\,Li formula}
     \hspace{2em}
    \begin{array}{c}
      \mbox{\footnotesize un-rigidifying $Y_{[k]}$'s} \\[-.6ex]
      \mbox{\rm ---------------------}\!\!\longrightarrow   \\[-1.6ex]
      \longleftarrow\!\!\mbox{\rm ---------------------}    \\[-.6ex]
      \mbox{\footnotesize rigidifying $Y_{[k]}$'s}
    \end{array}
     \hspace{2em}
    \mbox{\sl Ionel-Parker formula}\,.
   $$
} % end-claim

%\bigskip
%
\noindent
{\it Explanation.}
 Though, in format,
  \begin{itemize}
   \item[$\cdot$]
    [L-R] uses
     symplectic stretching similar to that in Floer homology theory and
     the notion of virtual neighborhoods construction in [Ru],

   \item[$\cdot$]
    [I-P2] uses the moduli space of $(J,\nu)$-holomorphic maps
     from the beginning and are thus dealing with a different
     moduli space from both [I-R] and [Li2],

   \item[$\cdot$]
    [Li2] uses the construction of virtual fundamental class from
     a perfect obstruction theory associated to the moduli problem
      of maps to fibers of a degeneration and
     is in the pure algebro-geometric setting in terms of
     Artin stacks and Deligne-Mumford stacks,
  \end{itemize}
  these differences should be only superficial as long as  the
  explicit form of the degeneration/gluing formula is concerned.
 The latter depends more on how objects in the moduli problem
  degenerate, i.e.\ on how maps in question break and how the target
  degenerates accordingly to keep the maps remain what we want.
 For this, what happens in the three are {\it the same},
  subject to the superficial difference of symplectic stretching
   in [L-R] versus the expansion of targets by a ruled
   manifold/variety in [I-P] and [Li2].

 The true cause of the difference of the formula of
  [L-R] and [Li2] versus [I-P2] is at [I-P2: Sec.~12].
 There, maps about to degenerate are pre-grouped by how many expansions
  it is going to take to remove degeneracy of maps in the limit.
 This gives rise to a covering of the moduli space of maps
   in question ([I-P2: Lemma 12.2]) and
  it is shown that the {\it inclusion-exclusion principle} way
  of counting does no harm (Identity (12.4) in [I-P2: Lemma 12.2]).
 It is this inclusion-exclusion identity of moduli spaces that leads
  to the form of the degeneration/gluing formula of [I-P2].
 Thus, to relate [I-P2] to [L-R] and [Li2], we should ask:
  \begin{itemize}
   \item[] {\bf Q.}
   {\it Is there an inclusion-exclusion identity in the setting
        of {\rm [L-R]} and {\rm [Li2]} as well?}
  \end{itemize}

 To investigate this, recall the layer-structure stratification
  $\{ {\cal M}_{(g,h),(n,\vec{m})}^{(i)}
           (Y,L\,|\,\underline{\beta},\vec{\gamma},\mu) \}
                         _{i\in {\scriptsizeBbb Z}_{\ge 0}}$  of
  $\overline{\cal M}_{(g,h),(n,\vec{m})}
           (Y,L\,|\,\underline{\beta},\vec{\gamma},\mu)$
  from Sec.~7.1.
 Stable maps in the depth-$i$ stratum
  ${\cal M}_{(g,h),(n,\vec{m})}^{(i)}
            (Y,L\,|\\[0ex]  \underline{\beta},\vec{\gamma},\mu)$
  is characterized by that their targets are all $Y_{[i]}$.
 The virtual codimension of
  ${\cal M}_{(g,h),(n,\vec{m})}^{(i)}
           (\\  Y,L\,|\,\underline{\beta},\vec{\gamma},\mu)$
  is $2i$.
 This is precisely $\dim {\Bbb G}_m[i]$.
 Indeed the occurrence of this virtual codimension comes exactly
  from the rigidification of the ${\Bbb G}_m[i]$-action
  on $W[i]/B[i]$ when we construct a standard Kuranishi
  neighborhood of an
  $f\in\overline{\cal M}_{(g,h),(n,\vec{m})}
                       (Y,L\,|\,\underline{\beta},\vec{\gamma},\mu)$
  that lies in
  ${\cal M}_{(g,h),(n,\vec{m})}^{(i)}
             (Y,L\,|\,\underline{\beta},\vec{\gamma},\mu)$.
 The situation is indeed analogous to what happens in
  [MD-S1: Remark A.5.3].

 In particular, if we
  take the depth-$i$ descendant Kuranishi structure
   ${\cal K}_0^{(i)}$ on
   ${\cal M}_{(g,h),(n,\vec{m})}^{(i)}
              (Y,L\,|\\  \underline{\beta},\vec{\gamma},\mu)$
   from the Kuranishi structure ${\cal K}_0$ on
   $\overline{\cal M}_{(g,h),(n,\vec{m})}
              (Y,L\,|\,\underline{\beta},\vec{\gamma},\mu)$   and
  consider the corresponding Kuranishi structure
   $\widetilde{\cal K}_0^{(i)}$ before rigidification, i.e.\
   Kuranishi structure for maps to the {\it rigid} $Y_{[i]}$,
 then we expect to have an open pseudo-embedding
  $$
   \tilde{\iota}^{(i)}:\widetilde{\cal K}_0^{(i)}\;
    \longrightarrow\; {\cal K}_0^{}\,,
  $$
  defined around $s^{-1}(0)$ on each Kuranishi neighborhood
  from $\widetilde{\cal K}_0^{(i)}$.
 (Recall a Kuranishi neighborhood data $(V,\Gamma_V,E_V;s,\psi)$.)
 We expect also that
  a resemble to [I-P2: Identity (12.4) in Lemma 12.1]
  $$
   {\cal K}_0\;
    =\; \tilde{\iota}^{(1)}({\cal K}_0^{(1)})\,
        -\, \tilde{\iota}^{(2)}({\cal K}_0^{(2)})\,
        +\, \tilde{\iota}^{(3)}({\cal K}_0^{(3)})\,
        -\, \cdots
  $$
  holds around $s^{-1}(0)$ on each Kuranishi neighborhood from
  ${\cal K}_0$.
 This should then reproduce the degeneration/gluing formula
  in the form of [I-P2].

 Since all the difference in [L-R], [Li2] versus [I-P2]
  that are related to the expression of the degeneration/gluing formula
  is whether or not and when and where to apply rigidification
  of targets of maps, we thus make the conjecture.

\noindent\hspace{15cm}$\Box$

\vspace{12em}
%\newpage
%references
%endfootnotesize

\end{document}